\documentclass[a4paper,11pt]{article}

\usepackage[
    a4paper,
    top=1.05in,
    bottom=1.05in,
    left=1.00in,
    right=1.00in
]{geometry}
\usepackage[T1]{fontenc}

\usepackage{amsmath,amsthm,mathtools}
\usepackage{amssymb,mathrsfs}
\usepackage{txfonts}
\usepackage{microtype}
\microtypesetup{protrusion=true,expansion=true}
\allowdisplaybreaks[4]

\usepackage{enumitem}
\setlist[enumerate]{leftmargin=2.5em,itemsep=0.25em,topsep=0.35em}
\setlist[itemize]{leftmargin=2.2em,itemsep=0.25em,topsep=0.35em}
\usepackage{booktabs,array,longtable,multicol,graphicx,float}

\usepackage{titlesec}
\titleformat{\section}{\normalfont\Large\bfseries}{\thesection}{0.75em}{}
\titleformat{\subsection}{\normalfont\large\bfseries}{\thesubsection}{0.75em}{}
\titlespacing*{\section}{0pt}{2.6ex plus 0.8ex minus 0.2ex}{1.2ex plus 0.2ex}
\titlespacing*{\subsection}{0pt}{2.0ex plus 0.6ex minus 0.2ex}{0.9ex plus 0.2ex}

\usepackage{fancyhdr}
\fancypagestyle{plain}{%
  \fancyhf{}%
  \fancyfoot[C]{\thepage}%
}

\usepackage[dvipsnames]{xcolor}
\definecolor{linkblue}{RGB}{20,67,128}
\definecolor{citegreen}{RGB}{25,105,80}
\definecolor{urlpurple}{RGB}{100,55,125}
\usepackage[colorlinks=true,
            linkcolor=linkblue,
            citecolor=citegreen,
            urlcolor=urlpurple,
            plainpages=false]{hyperref}
\usepackage{bookmark}
\numberwithin{equation}{section}
\newtheorem{theorem}{Theorem}[section]
\newtheorem{proposition}[theorem]{Proposition}
\newtheorem{lemma}[theorem]{Lemma}

\theoremstyle{definition}\newtheorem{definition}[theorem]{Definition}
\theoremstyle{remark}
\newcommand{\Ric}{\operatorname{Ric}}
\newcommand{\Rm}{\operatorname{Rm}}
\newcommand{\Tr}{\mathrm{Tr}}

\newcommand{\R}{\mathbb R}

\newcommand{\Vol}{\operatorname{Vol}}
\newcommand{\Id}{\operatorname{Id}}
\let\originalthebibliography\thebibliography
\renewcommand{\thebibliography}[1]{\originalthebibliography{#1}%
 \setlength{\itemsep}{0.4pt}\setlength{\parsep}{0pt}}

\title{Strong uniqueness of generalized cylindrical tangent flows}
\author{Yu Li}
\date{}
\hypersetup{pdftitle={Strong uniqueness of generalized cylindrical tangent flows (V81)},
 pdfauthor={Yu Li},
 pdfsubject={Strong uniqueness of generalized cylindrical tangent flows}}
\begin{document}
\maketitle
\begin{abstract}
We prove strong uniqueness of tangent flows at the first singular time of a closed Ricci flow
when one tangent flow is modeled on the product of an arbitrary compact positive Einstein
manifold and a Euclidean space. More precisely, the modified Ricci flow converges smoothly
in one time-independent embedding, with quantitative convergence on expanding regions.
\end{abstract}
\begingroup\small\tableofcontents\endgroup

\section{Introduction}\label{sec:introduction}
Tangent flows describe the geometry near a singularity of a geometric flow.
Rescaling and compactness produce subsequential limits. Uniqueness asks whether the limit is
independent of the rescaling sequence. For mean curvature flow, Colding--Minicozzi proved
uniqueness of multiplicity-one cylindrical tangent flows, including their axes~\cite{CM15}.
Strong uniqueness asks for convergence of the modified flow through one fixed identification,
together with quantitative estimates.

For Ricci flow, tangent flows in the noncollapsed measured framework are Ricci shrinkers,
possibly with singularities~\cite{BamStr,FLSt}. Chan--Ma--Zhang~\cite{CMZ24} proved uniqueness
when one tangent flow is smooth and closed. Noncompact models occur even for flows on closed
manifolds. Among smooth tangent shrinkers, geometry on the end plays an essential role in
quantitative comparison. Li--Yan--Zhang~\cite{LYZII} established strong uniqueness for smooth
asymptotically conical tangent shrinkers. Generalized cylinders, formed from a compact positive
Einstein factor and a Euclidean factor, constitute another fundamental class.

Li--Wang~\cite[Theorem~1.1]{LW24} proved pointed Gromov--Hausdorff rigidity of
$S^{n-1}\times\R$ among Ricci shrinkers. Colding--Minicozzi~\cite{CM25} proved rigidity of
round cylinders with arbitrary Euclidean dimension, and Fang--Li~\cite{FLC} proved strong
uniqueness of the corresponding tangent flows.
More recently, Bamler--Lai~\cite{BL} gave another proof of the uniqueness of cylindrical
tangent flows using their PDE--ODI principle.
Li--Zhang~\cite{LZ26} proved rigidity and tangent-flow uniqueness for generalized cylinders
under stability and obstruction assumptions. Fang--Li~\cite{FLQ} obtained strong uniqueness under
spectral and obstruction hypotheses. They also proved cylindrical rigidity whenever the
compact Einstein factor is rigid among normalized compact Ricci shrinkers~\cite{FLR}, removing
the earlier stability, nonresonance, and prescribed obstruction-order restrictions.
Here we prove strong uniqueness for an arbitrary smooth closed positive Einstein factor,
with no additional rigidity, stability, spectral, or integrability assumption.

Let $(N^m,g_N)$
be smooth, closed, and connected, with
\begin{equation}\label{eq:einstein}
 \Ric(g_N)=\frac12g_N,
\end{equation}
and let $d\ge1$ and $n:=m+d$. Put
\begin{equation}\label{eq:model}
 \bar M:=N\times\R^d,\qquad \bar g:=g_N+g_E,\qquad
 \Theta_N:=\log\frac{\Vol(N,g_N)}{(4\pi)^{m/2}}-\frac m2.
\end{equation}
\phantomsection\label{idx:euclidean-coordinates}
Write $\vec x=(x_1,\ldots,x_d)\in\R^d$ for the Euclidean coordinate vector.
The probability-normalized potential and measure are
\begin{equation}\label{eq:model-potential}
 \bar f:=\frac m2+\Theta_N+\frac{|\vec x|^2}{4},\qquad
 \mathrm{d}\bar\mu:=(4\pi)^{-n/2}e^{-\bar f}\,\mathrm{d}V_{\bar g}
 =\mathrm{d}\mu_N\,\mathrm{d}\gamma_d,
\end{equation}
where $\mu_N:=\Vol(N,g_N)^{-1}\,\mathrm{d}V_{g_N}$ and
$\mathrm{d}\gamma_d:=(4\pi)^{-d/2}e^{-|\vec x|^2/4}\,\mathrm{d}\vec x$.
\phantomsection\label{idx:product-flow}
Let $\mathcal C_N^d$ denote the noncollapsed Ricci flow limit space obtained by taking
the spacetime completion of the standard ancient solution associated with
$(\bar M,\bar g,\bar f)$, as in \cite[Section~2, before Definition~2.10]{FLQ}.
Its distinguished time-zero point is the limit of $(y,0,t)$ as $t\uparrow0$, for any $y\in N$.
On its smooth negative-time part, the metric at time $-a<0$ is $ag_N+g_E$ and the
terminal probability is
$\mu_N\otimes(4\pi a)^{-d/2}e^{-|\vec x|^2/(4a)}\,\mathrm{d}\vec x$.

Let $\mathcal X=\{(M^n,g(t)):-T\le t<0\}$ be a smooth closed Ricci flow, with first singular time $0$.
Assume that, for some $Y>0$,
\begin{equation}\label{eq:nash-hypothesis}
 \inf_{0<\tau\le T}\boldsymbol{\mu}(g(-T),\tau)\ge -Y,
\end{equation}
where $\boldsymbol{\mu}$ is Perelman's entropy. By monotonicity, this gives a common
lower bound $-Y$ for the pointed Nash entropies at all regular base points and admissible
backward scales; see \cite[Section~1, after Definition~1.1]{FLSt}.
We use the time-zero completion and measured tangent-flow framework
of~\cite[Section~2, especially Definitions~2.16 and~2.20 and Theorem~2.19]{FLC}. For a point $z_*$ in that completion, its terminal kernel is smooth and
positive on every negative-time slice. Denote its potential by $f(t)$, so that
\begin{equation}\label{eq:terminal-probability}
 \mathrm{d}\mu_t=(4\pi|t|)^{-n/2}e^{-f(t)}\,\mathrm{d}V_{g(t)},\qquad t<0.
\end{equation}
The measured tangent hypothesis below includes the terminal kernel and the time vector. On the smooth
product limit it supplies smooth convergence of the metrics and kernel densities on compact subsets.

\phantomsection\label{idx:modified-time}
Write $t=-e^{-s}$ and set $\tilde g(s):=|t|^{-1}g(t)$ and $\tilde f(s):=f(t)$.
Fix a sufficiently large $s_0$ and solve
\begin{equation}\label{eq:modified-diffeomorphism}
 \partial_s\Xi_s=-\nabla^{\tilde g(s)}\tilde f(s)\circ\Xi_s,\qquad \Xi_{s_0}=\Id.
\end{equation}
The closedness of $M$ gives a diffeomorphism on each finite interval. The modified flow is
\begin{equation}\label{eq:modified-gauge}
 g^{z_*}(s):=\Xi_s^*\tilde g(s),\qquad f^{z_*}(s):=\Xi_s^*\tilde f(s),\qquad
 \mu^{z_*}:=\Xi_s^*\mu_{-e^{-s}}.
\end{equation}
Its probability measure is independent of $s$.

\begin{theorem}[Strong uniqueness]\label{thm:main}
Let $\mathcal X$ be a smooth closed Ricci flow satisfying~\eqref{eq:nash-hypothesis}, and suppose that
$0$ is its first singular time. Let $z_*$ belong to its time-zero completion. If one measured tangent flow
at $z_*$ is $\mathcal C_N^d$, where $(N,g_N)$ satisfies~\eqref{eq:einstein}, then there is a smooth,
time-independent open embedding $\Psi:\bar M\longrightarrow M$ such that
\begin{equation}\label{eq:main-convergence}
 \Psi^*g^{z_*}(s)\longrightarrow\bar g,\qquad
 \Psi^*f^{z_*}(s)\longrightarrow\bar f\quad\text{in }C^\infty_{\mathrm{loc}}(\bar M),
 \qquad \Psi^*\mu^{z_*}=\bar\mu.
\end{equation}
Let $p\in(1,3/2]$ be the stationary exponent of
Theorem~\ref{thm:logarithmic}. The same embedding $\Psi$ has the following
property for every integer $k\ge0$ and every
\[
 0<\gamma<\frac{p-1}{2-p},\qquad 0<c<2\gamma.
\]
There are $C,\bar s<\infty$ such that, for every integer $j\ge\bar s$
and every $s\ge j$,
\begin{equation}\label{eq:main-quantitative}
 [\Psi^*g^{z_*}(s)-\bar g]_k
 +[\Psi^*f^{z_*}(s)-\bar f]_k
 \le Ce^{\bar f/2}(1+s)^{-\gamma}
 \quad\text{on }\{\bar f\le c\log j\},
\end{equation}
where $[T]_k:=\sum_{\ell=0}^k|\bar\nabla^\ell T|_{\bar g}$.
The constant $C$ depends on $(N,g_N)$, $d$, $Y$, the fixed norm
conventions, and $p,k,\gamma,c$; the waiting time $\bar s$ may also
depend on the particular flow, $z_*$, and the initial normalization of
the modified flow.
\end{theorem}

The convergence takes place in one time-independent embedding, as in~\cite{FLC},
with quantitative estimates on expanding potential sublevels; compare also
\cite[Theorem~1.2]{LYZII}. In particular, the tangent flow at $z_*$ is unique; see
Proposition~\ref{prop:time-preserving-reconstruction}. When $p=3/2$, every $\gamma<1$ is
permitted. Subsection~\ref{sec:finite-quotients} proves the corresponding statements for
smooth free finite quotients and for ancient flows with one prescribed tangent flow at infinity.

The main difficulty is that a small shrinker defect need not control displacement from a prescribed
compact metric. If the compact factor is nonrigid, nearby compact shrinkers give nearby cylinders
with zero defect. Thus the reduced entropy gradient controls the defect, while the compact parameter
itself must be controlled along the flow.

The extension--improvement strategy goes back to Colding--Ilmanen--Minicozzi's rigidity argument
for cylindrical self-shrinkers in mean curvature flow~\cite[Section~0.2 and Propositions~2.1--2.2]{CIM15}.
One first enlarges the region of control with some loss, then improves the estimates so that the
construction can be repeated. This strategy was developed in Colding--Minicozzi's cylindrical
mean curvature flow uniqueness theorem~\cite{CM15} and Chodosh--Schulze's treatment of
multiplicity-one asymptotically conical tangent flows~\cite{CS21}. Its Ricci-flow applications
include Colding--Minicozzi's cylinder rigidity theorem~\cite{CM25} and Li--Zhang's
work~\cite{LZ26}. The Ricci-flow arguments in~\cite{FLC,LYZII} combine
stationary estimates with parabolic propagation. Fang--Li~\cite{FLC} is the immediate source of
several propagation and comparison arguments used below.

The stationary input is a \L{}ojasiewicz--Simon inequality obtained from the analytic reduction
and radial obstruction method of~\cite{FLR}. The radial obstruction controls all positive
Hermite-degree kernel components, including additional resonances. The compact directions are
represented by an analytic finite-dimensional variational family $X_b$, with reduced entropy
$e(b)$. Its defect is controlled by $|\mathrm De(b)|$, and the finite-dimensional gradient
inequality applies even when the critical set is nonisolated. Theorem~\ref{thm:logarithmic}
controls the entropy and the metric--potential error relative to this family, together with the
positive Hermite modes and the reduced gradient.

The extension step propagates an incoming comparison over a fixed time block.
Two-sided pseudolocality and Shi's estimates give uniform regularity; compactness and
uniqueness of Ricci flow then yield qualitative metric control.
A localized metric--connection--curvature energy argument, following Kotschwar~\cite{Kot},
then gives quantitative metric comparison in the original coordinates.
The potential is treated by auxiliary caloric comparison, terminal-kernel estimates, and entropy
control; probability normalization fixes its additive constant. These estimates hold uniformly
over the compact variational family; the additional errors caused by its nonzero defect are
controlled by $|\mathrm De(b)|$.
Propagation followed by parabolic rescaling enlarges the comparison region in model coordinates.
At the endpoint, we place the pair in the reference slice and complete it outside the comparison
region. Applying the stationary inequality gives the improvement, and smoothing recovers the
same incoming regularity for the next block.

To compare successive blocks, we place the metric in a weighted divergence-free gauge with
respect to the model generalized cylinder: the metric perturbation has zero weighted divergence
with respect to $\bar g$ and $\bar\mu$. A similar gauge construction for asymptotically conical
Ricci shrinkers is given in~\cite[Theorem~1.1]{LYZI}. We choose these gauges continuously
and normalize them by model isometries to remove the symmetry components of their motion.
The resulting endpoint chart determines the data for the following block. Weak estimates for the chart
motion and the transitions between blocks then control finitely many localized metric pairings
that detect compact-factor displacement. Entropy summation bounds their cumulative change on
arbitrary finite initial segments. A first-exit argument confines the compact parameters to the
construction neighborhood and establishes continuation. This controls motion along the
variational family even when its critical set is nonisolated.

After continuation, the entropy decay and growing comparison regions give finite total variation
of the modified flow on fixed subsets. Potential barriers and integration of its equations give a
complete limit, which the
selected tangent sequence identifies with the prescribed generalized cylinder. The resulting
time-independent embedding is kept fixed when deriving the quantitative estimates for every
derivative order and admissible rate.

The improvement estimate in~\cite[Proposition~5.2]{FLC} is proved before the spatial extension theorems. Taylor expansion, combined with quadratic rigidity, gauge control,
and localization, gives an entropy bound in terms of $r_{B,\sigma}$
\cite[Theorems~4.17 and~5.4]{FLC}.
Alternating extension and improvement on one time slice yields
$\min\{r_A,r_{B,\sigma}\}\ge(1-3\sigma)r_E$, from which the discrete
\L{}ojasiewicz inequality follows
\cite[Theorems~6.18 and~6.24--6.27]{FLC}.
The propagation, terminal-kernel, caloric comparison, and normalization arguments in the present
extension step draw on~\cite[Section~6]{FLC}.

Here the dynamical construction propagates a chart over a parabolic block before applying the
independently proved stationary inequality at its endpoint. Improvement is relative to a new
member $X_{b_+}$ of the variational family; uniform family estimates and parameter confinement
are therefore needed even when the reduced critical set is nonisolated. The retained radii
depend on both current entropy dissipation and restrictions inherited from preceding blocks.
Their delayed recurrence supplies the finite-prefix summability needed for confinement.

In the asymptotically conical case, curvature decay permits control of retained exterior regions
under the flow~\cite[Lemma~6.4 and Proposition~6.5]{LYZII}. Weighted norms adapted to conical
scaling govern both geometric comparison and the analytic inequality; the geometric comparison
is converted into smallness in the stronger analytic admissibility space~\cite[Proposition~6.10 and Lemma~6.12]{LYZII}.
This gives a pointed entropy inequality on a single time slice~\cite[Proposition~6.14]{LYZII}.
Theorem~\ref{thm:logarithmic} is likewise stationary, but curvature of the compact factor
persists along the Euclidean directions. The asymptotically conical exterior argument therefore
does not transfer directly; the block construction and parameter confinement supply the
required control here.

The paper is organized as follows. Section~\ref{sec:stationary} proves the stationary inequality and comparison estimates.
Section~\ref{sec:analytic} sets up the extension--improvement construction, and
Sections~\ref{sec:physical}--\ref{sec:radii} establish propagation, endpoint improvement, and
the iteration estimates. Section~\ref{sec:motion} proves parameter confinement and continuation;
Section~\ref{sec:uniqueness} proves convergence in a fixed embedding and the quotient and ancient
extensions. Appendix~\ref{sec:inverse} supplies the gauge inverse.

\subsection*{Acknowledgements}
Yu Li is supported by the National Key Research and Development Program of China
(Grant No.~2025YFA1018200), NSFC-12522105, and research funds from the University of Science and
Technology of China and the Chinese Academy of Sciences.

\subsection*{AI declaration}
ChatGPT (GPT-6 Astra) was used to assist with drafting,
mathematical exposition, revision, and reference checking. The author takes full responsibility
for the mathematical content and correctness.

\section{A {\L}ojasiewicz--Simon inequality for generalized cylinders}
\label{sec:stationary}\label{sec:logarithmic}\label{sec:log}

\subsection{Reference spaces and variational formulas}
We use the weighted conventions of \cite[Sections~2.1--2.3]{FLR}.
\begin{definition}\label{def:weighted-quantities}
For a metric--potential pair $X=(g,f)$ on an $n$-manifold, define
\begin{enumerate}[label=\textup{(\roman*)},leftmargin=*]
\item (Weighted measure and mass) Set
\begin{equation}\label{eq:weighted-conventions}
 \mathrm d\mu_X:=(4\pi)^{-n/2}e^{-f}\;\mathrm dV_g.
\end{equation}
The pair is normalized when $\int\mathrm d\mu_X=1$.
\begin{subequations}\label{eq:physical-defect}
\item (Soliton defect) At scale one, set
\begin{equation}
 \Phi(X):=\frac12g-\Ric_g-\nabla_g^2f,
 \qquad \mathcal D(X):=\left(\int|\Phi(X)|_g^2\;\mathrm d\mu_X\right)^{1/2}.
\end{equation}
The tensor norm in $\mathcal D$ uses the current metric $g$.
\item (Entropy) The scalar Euler--Lagrange expression and weighted
entropy are
\begin{equation}
 \begin{split}
 w(X)&:=2\Delta_g f-|\nabla f|_g^2+\mathrm R_g+f-n,\\
 \mathcal W(X)&:=\int(\mathrm R_g+|\nabla f|_g^2+f-n)\;\mathrm d\mu_X.
 \end{split}
\end{equation}
For the pairs considered below, $\mathcal W(X)=\int w(X)\,\mathrm d\mu_X$.
\end{subequations}
\end{enumerate}
\end{definition}

\phantomsection\label{idx:model-potentials}
On the reference product,
$f_N:=m/2+\Theta_N$, $f_E:=|\vec x|^2/4$, and
$\bar X=(\bar g,\bar f)$ has $\bar f=f_N+f_E$ and
$\mathrm d\bar\mu=\mathrm d\mu_N\mathrm d\gamma_d$.

\phantomsection\label{idx:weighted-operators}
\emph{Weighted operators.} We use
\[
 \Delta_f:=\Delta_g-\nabla_{\nabla f},\qquad
 \mathcal L_f:=\Delta_f+2\Rm_g,
\]
where the weighted Laplacian acts on functions or tensors and
$(\Rm_g h)_{ij}=R_{ikjl}h^{kl}$.
\phantomsection\label{idx:weighted-divergence}
The divergence conventions are
\[
 \begin{gathered}
 (\delta_gh)_j=\nabla^ih_{ij},\qquad
 \delta_X:=\delta_{g,f},\qquad
 (\delta_Xh)_j=\nabla^ih_{ij}-h_{ij}\nabla^if,\\
 \delta_X^*\omega=\delta_g^*\omega=-\frac12L_{\omega^\sharp}g.
 \end{gathered}
\]
The last operator is the adjoint of $\delta_X$ in the $L^2(g,\mu_X)$
pairing, and $L_V$ denotes Lie differentiation. Reference operators use
$\bar X$.

\phantomsection\label{idx:pair-components}
\emph{Variations.} For a tensor--scalar pair $H=(h,k)$, write
$H_{\mathrm t}:=h$ and $H_{\mathrm s}:=k$.
In particular, $X_{\mathrm t}=g$ and $X_{\mathrm s}=f$;
for a variation of $(g,f)$, the scalar component varies the potential.
Writing $q=\frac12\Tr_g h-k$, we have
$\mathrm D(\mathrm d\mu_X)[H]=q\;\mathrm d\mu_X$.
The entropy first variation \cite[Section~2.3]{FLR} is
\begin{equation}\label{eq:W-first}
 \mathrm D\mathcal W_X[H]
 =\int\bigl(\langle\Phi(X),H_{\mathrm t}\rangle_g
             +(w(X)-1)q\bigr)\;\mathrm d\mu_X.
\end{equation}
\phantomsection\label{idx:product-regions}
On $\bar M=N\times\R^d$, reference derivatives, norms, projections,
and orthogonality use $(\bar g,\bar\mu)$; $\bar\nabla$ is its product
connection. We use $\mathrm D$ for differentials of maps and functionals
and $\mathcal D$ for the intrinsic defect norm.
Set $N_R:=N\times B_R(0)$. An embedding on $N_R$ is always
defined on an open neighborhood of $\overline N_R$, so that its boundary has a collar.
Unless additional dependencies are stated, constants depend only on
$(N,g_N)$, $d$, $Y$ in~\eqref{eq:nash-hypothesis}, and the fixed reference charts and norm conventions.
After Lemma~\ref{lem:compact-family}, dependence on its fixed
preliminary stationary family is also understood. Dependence on
prescribed derivative orders, retained fractions, Gaussian losses,
time horizons and comparison tolerances is specified in the relevant
statements. Constants are independent of varying radii unless stated
otherwise.

\begin{definition}\label{def:metric-potential-perturbations}
For a metric--potential pair on $\bar M$, use the following coordinates
and constraints.
\begin{enumerate}[label=\textup{(\roman*)},leftmargin=*]
\item (Perturbation and density) Write
\begin{equation}\label{eq:chart}
 X=(g,f)=(\bar g+h,\bar f+k).
\end{equation}
\phantomsection\label{idx:physical-density-variation}
For positive $\bar g+h$, the logarithmic density ratio is
\begin{equation}\label{eq:density-ratio}
 \ell(h,k):=\frac12\log\det(\Id+\bar g^{-1}h)-k,
 \qquad \mathrm d\mu_X=e^{\ell(h,k)}\,\mathrm d\bar\mu.
\end{equation}
\item (Normalized centered slice) For pairs with finite mass and
first moment, the slice is specified by
\begin{equation}\label{eq:actual-constraints}
 \delta_{\bar X}h=0,\qquad
 \int\mathrm d\mu_X=1,\qquad
 \int\vec x\,\mathrm d\mu_X=0.
\end{equation}
\end{enumerate}
\end{definition}

\begin{definition}\label{def:weighted-sobolev}\label{def:physical-sobolev}
Our norm conventions are as follows.
\begin{enumerate}[label=\textup{(\roman*)},leftmargin=*]
\item (Sobolev norms) For a positive smooth measure $\mu$ and an
integer $j\ge0$, set
\[
 \|u\|_{H^j(\mu)}^2
 :=\sum_{i=0}^j\int|\bar\nabla^iu|_{\bar g}^2\,\mathrm d\mu,
 \qquad
 \|U\|_{H^j(\mu)}^2
 :=\|U_{\mathrm t}\|_{H^j(\mu)}^2
                   +\|U_{\mathrm s}\|_{H^j(\mu)}^2.
\]
These spaces are the completions of smooth compactly supported fields.
Unqualified $H^j$ and $L^2$ norms use $\bar\mu$.
Derivatives and pointwise norms use $\bar g$, also when $\mu=\mu_X$;
the intrinsic defect $\mathcal D(X)$ instead uses $g$,
as in~\eqref{eq:physical-defect}.
\item (Uniform norms) Unqualified $C^j$ denotes the sum of the
supremum norms of reference covariant derivatives through order $j$.
The analytic proofs also use the usual $C^{j,1/2}$ norms on the fixed
uniform product atlas. For a tensor--scalar pair, either norm is the sum
of its component norms. Norms on $N$ use $g_N$ and $\mu_N$.
\end{enumerate}
\end{definition}
\begin{proposition}\label{prop:kernel}
\phantomsection\label{idx:hermite-spaces}
Let $\mathcal H_j:=\ker_{L^2(\gamma_d)}(\Delta_{f_E}+j/2)$ be the
space of degree-$j$ Hermite polynomials.
Put $\mathcal L_N=\Delta_N+2\Rm_N$ and
$\mathscr E_\lambda(N):=\{\eta:\delta_{g_N}\eta=0,
\ \mathcal L_N\eta=\lambda\eta\}$.
Then the kernel has the finite orthogonal decomposition
\[
 K:=\ker\mathcal L_{\bar f}\cap\ker\delta_{\bar X}
   =\bigoplus_{j\ge0}K_j,
 \qquad K_j:=\mathcal H_j\otimes\mathscr E_{j/2}(N).
\]
Its tensors have both indices tangent to $N$. Write
$K_N:=K_0$ and $K_+:=\bigoplus_{j\ge1}K_j$.
\end{proposition}
\begin{proof}
This is \cite[Proposition~2.4]{FLR}.
\end{proof}

\begin{lemma}\label{lem:gaussian-sobolev}\label{lem:physical-sobolev}
Let $\mathrm d\mu=e^\ell\,\mathrm d\bar\mu$ be a smooth positive measure
with $|\bar\nabla\ell|\le\epsilon\langle\vec x\rangle$, where
$\langle\vec x\rangle:=1+|\vec x|$, for a sufficiently small fixed
$\epsilon>0$. For every function or tensor $u\in H^2(\mu)$,
\begin{equation}\label{eq:graphnorm}
 \|\langle\vec x\rangle\bar\nabla u\|_{L^2(\mu)}
       +\|\langle\vec x\rangle^2u\|_{L^2(\mu)}
       \le C\|u\|_{H^2(\mu)}.
\end{equation}
The constant is uniform over these measures.
For the reference measure $\bar\mu$, the self-adjoint domains of
$\Delta_{\bar f}$ and $\mathcal L_{\bar f}$ are $H^2$ of the respective
tensor types, with
\[
 \|u\|_{H^2}\asymp\|u\|_2+\|\Delta_{\bar f}u\|_2,
 \qquad
 \|h\|_{H^2}\asymp\|h\|_2+\|\mathcal L_{\bar f}h\|_2.
\]
\end{lemma}
\begin{proof}
For a compactly supported field, Gaussian integration by parts gives
\[
 \int x_i^2|u|^2\;\mathrm d\mu
 =2\int|u|^2\;\mathrm d\mu
  +4\int x_i\langle u,\bar\nabla_i u\rangle\;\mathrm d\mu
  +2\int x_i\partial_i\ell\,|u|^2\;\mathrm d\mu.
\]
After summing over $i$, the last term is bounded by
$C\epsilon\|\langle\vec x\rangle u\|_{L^2(\mu)}^2$ and can be
absorbed. Cauchy--Schwarz then gives
\[
 \|\langle\vec x\rangle u\|_{L^2(\mu)}
 \le C\bigl(\|u\|_{L^2(\mu)}+\|\bar\nabla u\|_{L^2(\mu)}\bigr).
\]
Apply this estimate to $\bar\nabla u$ and $x_j u$, and then use
approximation in $H^2(\mu)$, to prove~\eqref{eq:graphnorm}.
The reference domain and graph-norm statements are
\cite[Lemma~3.6, Propositions~3.7 and~3.14, Remark~3.4]{LZ26}.
\end{proof}

\begin{definition}\label{def:internal-parameters}
For metric--density coordinate pairs $H=(h,\sigma)$, we use the following
fixed linear structure and chart.
\begin{enumerate}[label=\textup{(\roman*)},leftmargin=*]
\item (Pairing) Set
\begin{equation}\label{eq:pair-pairing}
 \langle U,V\rangle:=\int
  \big(\langle U_{\mathrm t},V_{\mathrm t}\rangle_{\bar g}
                  +U_{\mathrm s}V_{\mathrm s}\big)\,\mathrm d\bar\mu.
\end{equation}
\item (Slice and normal space) Set
\begin{equation}\label{eq:constraints}
 \mathcal S:=\ker\delta_{\bar X}\times S_0,
 \qquad S_0:=\left\{\sigma:\int\sigma\,\mathrm d\bar\mu=0,
                \ \int\vec x\,\sigma\,\mathrm d\bar\mu=0\right\}.
\end{equation}
Set $Z:=\mathcal S\cap(K\times\{0\})^\perp$.
The component in this kernel-orthogonal space is called the normal component.
\item (Density chart) For positive metric $\bar g+h$ and positive
density $1+\sigma$, define
\begin{equation}\label{eq:density-lift}
 \begin{split}
 \mathcal C(H)&:=\left(\bar g+h,
  \bar f+\frac12\log\det(\Id+\bar g^{-1}h)-\log(1+\sigma)\right),
 \\
 \mathrm d\mu_{\mathcal C(H)}&=(1+\sigma)\,\mathrm d\bar\mu.
 \end{split}
\end{equation}
For $X=\mathcal C(H)$, one has $h=g-\bar g$ and
$\sigma=\mathrm d\mu_X/\mathrm d\bar\mu-1$;
the potential difference $k=f-\bar f$ satisfies
$1+\sigma=\exp(\ell(h,k))$.
\end{enumerate}
\end{definition}
The chart is defined also off $\mathcal S$. On $\mathcal S$ it imposes
the exact conditions~\eqref{eq:actual-constraints}. The relative measure variation in the direction $U$ is
$U_{\mathrm s}/(1+\sigma)$.

\phantomsection\label{idx:slice-projections}
The fixed orthogonal projections are
\[
 \begin{split}
 \Pi_{S_0}u&=u-\int u\,\mathrm d\bar\mu
           -\frac12\sum_i x_i\int x_i u\,\mathrm d\bar\mu,\\
 \Pi_ZU&=\bigl((\Pi_{\mathrm{sl}}-\Pi_K)U_{\mathrm t},
                         \Pi_{S_0}U_{\mathrm s}\bigr),
 \end{split}
\]
where $\Pi_{\mathrm{sl}}$ projects onto the full metric slice.
Compact-factor isometries and Euclidean rotations preserve the reference
pairing, slice, kernel and projections.
All these projections are bounded on each $H^j$
by \cite[Lemma~2.9]{FLR}. The parameters
$(b,a)=\Pi_Kh\in K_N\oplus K_+$ use the fixed metric $L^2$ norm,
and reduced derivatives use its dual norm.

The constrained gradient of $\mathcal W\circ\mathcal C$ for the fixed pairing is
\begin{equation}\label{eq:density-lift-gradient}
 \mathcal G(H):=
 \left(\Pi_{\mathrm{sl}}[(1+\sigma)\mathrm A\Phi(\mathcal C(H))\mathrm A],
                         \Pi_{S_0}[w(\mathcal C(H))-1]\right),
 \qquad \mathrm A:=g^{-1}\bar g.
\end{equation}
Here $(\mathrm A\Phi\mathrm A)_{ij}
=\mathrm A^p{}_i\mathrm A^q{}_j\Phi_{pq}$.
Formula~\eqref{eq:density-lift-gradient} also defines the projected
ambient gradient off $\mathcal S$; it follows from the first variation
and \cite[Lemma~2.8]{FLR}. On $\mathcal S$ the reference Hessian is
\begin{equation}\label{eq:hessbase}
 \mathfrak B_0U:=\mathrm D\mathcal G(0)[U]
 =\left(\frac12\mathcal L_{\bar f}U_{\mathrm t},
                       -(2\Delta_{\bar f}+1)U_{\mathrm s}\right).
\end{equation}
The compact scalar gap exceeds $1/2$, and centering removes the
Euclidean scalar $1/2$-eigenfunctions. Hence
$\ker\mathfrak B_0=K\times\{0\}$. By \cite[Lemma~2.9]{FLR}, the normal
restriction is an isomorphism $Z\cap H^{j+2}\to Z\cap H^j$ and
\begin{equation}\label{eq:finiteinverse}
 \|U\|_{H^{j+2}}\le C_j\|\mathfrak B_0U\|_{H^j},
 \qquad U\in Z\cap H^{j+2}.
\end{equation}

For a variation $J=(\eta,\kappa)$ at an arbitrary pair $X=(g,f)$,
put $q:=\frac12\Tr_g\eta-\kappa$. By \cite[Lemma~2.2]{FLR}, the linearizations in
metric--potential coordinates are
\begin{equation}\label{eq:physical-linearizations}
 \begin{split}
 \mathrm D\Phi_X[J]&=\frac12\mathcal L_f\eta
       +\delta_X^*\delta_X\eta+\nabla_g^2q
       +\frac12(\Phi\circ_g\eta+\eta\circ_g\Phi),\\
 \mathrm Dw_X[J]&=-(2\Delta_f+1)q+\delta_X\delta_X\eta
                                      +\langle\Phi,\eta\rangle_g,
 \end{split}
\end{equation}
where $\circ_g$ denotes contraction of one index with $g$.

\subsection{Finite-dimensional reduction and positive cutoff family}
\label{subsec:finite-jets}
\begin{lemma}[Compact-factor variational family]\label{lem:compact-family}
For $b$ in a neighborhood of zero in $K_N$, there is a canonical
analytic family of compact metric--potential pairs $(g_b,f_b)$ such that
\begin{equation}\label{eq:compactfamily}
 X_b:=(g_b+g_E,f_b+f_E),\qquad
 e(b):=\mathcal W(X_b)=w(X_b),\qquad e(0)=\Theta_N.
\end{equation}
Here $f_b$ is the canonical normalized minimizing potential and,
at every prescribed finite order,
\[
 \|g_b-g_N-b\|_{C^j(N)}\le C_j|b|^2,
 \qquad \|\Phi(X_b)\|_{C^j(N)}\le C_j|\mathrm De(b)|.
\]
Its metric perturbation lies in the reference slice and has kernel projection $b$.
The family has mass one and center zero, is equivariant under
center-preserving reference isometries, and agrees on overlapping
parameter neighborhoods at different regularities.
\end{lemma}
\begin{proof}
Apply \cite[Lemma~3.1 and Proposition~3.2]{FLR} to the compact
metric--density slice. Its unique small analytic normal graph has
kernel projection $b$, correction $O(|b|^2)$ in every prescribed
$C^j$ norm, and the canonical minimizing potential. Extend its
coordinate pair constantly in the Euclidean factor and write it as
$H_b:=\mathcal C^{-1}(X_b)$. The normal equation gives
\[
 H_b-(b,0)\in Z,\qquad
 \mathcal G(H_b)=(\nabla e(b),0).
\]
The cited compact estimates give the displayed metric and defect bounds.
Choose a parameter ball invariant under compact-factor isometries and
Euclidean rotations. The pairings, projections, density chart and
normal equation commute with this action, so uniqueness gives
$R^*X_b=X_{R\cdot b}$. Compact normalization and the Euclidean
Gaussian give mass one and center zero for the product family.
\end{proof}
Fix a relatively compact preliminary parameter neighborhood for this
stationary family and its reference norms.

\begin{proposition}[Second-order reduction]
\label{thm:jets}\label{lem:closing}
For every prescribed finite spatial and parameter regularity, there
is an analytic reduced energy $F_2(b,a)$ on a sufficiently small
parameter neighborhood, with
\begin{equation}\label{eq:reduced-energy}
 F_2(b,0)=e(b),\qquad \mathrm D_aF_2(b,0)=0.
\end{equation}
On a smaller neighborhood,
\begin{equation}\label{eq:closing}
 |a|^2+|\mathrm De(b)|\le C|\mathrm DF_2(b,a)|,
\end{equation}
and there is $p_N\in(1,2]$ such that
\begin{equation}\label{eq:finiteLS}
 |F_2(b,a)-\Theta_N|\le C|\mathrm DF_2(b,a)|^p,
 \qquad p:=\min\{3/2,p_N\}>1.
\end{equation}
The neighborhood and constants depend on the prescribed orders.
\end{proposition}
\begin{proof}
By \cite[Proposition~3.11]{FLR}, there is a coordinate pair
$H^{[2]}(b,a)$, quadratic in $a$ and analytic in $b$, with smooth
coefficients on $N$ and a fixed finite Euclidean polynomial degree,
such that
\begin{equation}\label{eq:formalnormal}
 H^{[2]}-H_b-(a,0)\in Z,\qquad H^{[2]}\in\mathcal S,
 \qquad \Pi_K(H^{[2]})_{\mathrm t}=b+a.
\end{equation}
For fixed $b$, the terms of degrees zero, one and two in $a$ in the
formal Taylor expansion of $\Pi_Z\mathcal G(H^{[2]}(b,a))$ vanish.
The correction coefficients are bounded at order $j$ by
$C_j(|b||a|+|a|^2)$. In particular, for a fixed polynomial degree $Q$,
\[
 |\bar\nabla^j(H^{[2]}-H_b)|\le C_j|a|\langle\vec x\rangle^Q.
\]
Define $F_2$ to be the cubic Taylor polynomial in $a$ of
$\mathcal W(\mathcal C(H^{[2]}))$. Its coefficients are defined by
cutoff approximation and Gaussian integration. The cited construction
gives~\eqref{eq:reduced-energy}, and~\eqref{eq:closing} is
\cite[(4.07)]{FLR}.

The finite-dimensional \L{}ojasiewicz gradient inequality
\cite[Section~18, Proposition~1]{Loj65}
for the analytic function $e$ at its critical point $0$ gives, for
some $p_N\in(1,2]$,
\begin{equation}\label{eq:compactLS}
 |e(b)-e(0)|\le C|\mathrm De(b)|^{p_N}.
\end{equation}
If $e$ is locally constant, take $p_N=2$. Taylor expansion in $a$ gives
\[
 F_2(b,a)-e(b)=\frac12\mathrm D_aF_2(b,a)[a]+O(|a|^3).
\]
By \eqref{eq:closing}, its absolute value is at most
$C|\mathrm DF_2|^{3/2}$, while \eqref{eq:compactLS} bounds
$|e(b)-e(0)|$ by $C|\mathrm DF_2|^{p_N}$.
Adding proves \eqref{eq:finiteLS}.
\end{proof}

\begin{proposition}\label{prop:cutoffjets}
Prescribe an integer $j\ge0$ and finitely many parameter derivatives.
For sufficiently small $(b,a)$ there are coordinate pairs $\hat H_{b,a}$
defining positive metrics and densities, and complete normalized,
centered pairs $\hat X_{b,a}=\mathcal C(\hat H_{b,a})$,
such that $\hat H_{b,0}=H_b$ and $\hat X_{b,0}=X_b$.
The uniform spatial bound is
\[
 \|\hat H_{b,a}-H_b\|_{C^{j+2}}\le C_j|a|^{1/2}.
\]
On every fixed $N_r$ the family admits the prescribed parameter
derivatives in $C^{j+2}(N_r)$. Writing $\hat h=(\hat H_{b,a})_{\mathrm t}$,
for every fixed $A>0$ the following estimate holds, also when the
prescribed derivatives are applied before taking the norms:
\begin{equation}\label{eq:representative-flat-errors}
 \|(\Id-\Pi_{\mathrm{sl}})\hat h\|_{H^{j+2}}
 +|\Pi_K\hat h-(b+a)|
 +\|\delta_{\bar X}\hat h\|_{H^{j+1}}
       \le C_{j,A}|a|^A.
\end{equation}
The variational and entropy estimates are
\begin{gather}\label{eq:cutresidual}
 \|\Pi_Z\mathcal G(\hat H_{b,a})\|_{H^j}
 +|\Pi_K\mathcal G(\hat H_{b,a})_{\mathrm t}-\mathrm DF_2(b,a)|
       \le C_j|a|^3,\\*
 |\mathcal W(\hat X_{b,a})-F_2(b,a)|\le C|a|^4.
 \label{eq:cutenergy}
\end{gather}
For every prescribed $j$ there is a fixed polynomial degree $Q_j$ such that
\begin{equation}\label{eq:representative-cutoff-core}
 |\bar\nabla^j(\hat H_{b,a}-H_b)|
 +|\bar\nabla^j(\hat X_{b,a}-X_b)|
       \le C_j|a|\langle\vec x\rangle^{Q_j}.
\end{equation}
All constants are uniform for $b$ in a smaller fixed neighborhood.
The constructions agree on overlaps at different regularities and
are equivariant under center-preserving reference isometries.
\end{proposition}
\begin{proof}
Use the polynomial pair $H^{[2]}$ from the proof of
Proposition~\ref{thm:jets}. The positive cutoff construction and
\eqref{eq:cutresidual} are \cite[Lemmas~3.12--3.13]{FLR}.
The cutoff errors are flat in $a$, meaning $O(|a|^A)$ for every $A>0$
in the indicated norm, uniformly for $b$ in a fixed smaller neighborhood.
Fix a radial function $\vartheta$, equal to one on
$B_1$ and zero outside $B_2$, and a polynomial degree $Q$ for the
preceding approximation. For $a\ne0$, take
$R=|a|^{-1/(4(Q+1))}$ and use
\[
 H_b+\vartheta(\vec x/R)(H^{[2]}-H_b)
\]
with a scalar correction supported in a fixed compact set to restore
mass and first moment. Its coefficients are Gaussian tail integrals
outside $R$, hence are flat with every prescribed derivative.
The cutoff perturbation is uniformly $C^{j+2}$-small, with bound
$C_j|a|^{1/2}$, so its metric and density remain positive.
Outside a compact set its metric is $g_b+g_E$, which proves completeness.
The polynomial coefficient bounds give the first pointwise estimate
in~\eqref{eq:representative-cutoff-core}; the density chart and uniform
positivity give the metric--potential estimate.

On a fixed $N_r$ the cutoff equals one for all sufficiently small $|a|$.
After subtracting the scalar correction, $\hat H_{b,a}-H^{[2]}$ is
supported where $|\vec x|\ge R$. Each parameter derivative of the cutoff
costs at most a power of $|a|^{-1}$ there, which is absorbed by the
Gaussian tail. The difference is thus flat in every prescribed
polynomially weighted Gaussian Sobolev norm, also after the prescribed
parameter derivatives. Since $H^{[2]}$ lies in the exact slice and has metric kernel
projection $b+a$, boundedness of the fixed projections and weighted
divergence gives~\eqref{eq:representative-flat-errors}.
The radial cutoff and the canonical compact construction commute
with reference isometries. For the scalar correction use a radial
Euclidean function and its products with the coordinate functions.
Their coefficients transform as the mass and first moment, so the
cutoff family is equivariant.
The fixed choices agree at different regularities.

For \eqref{eq:cutenergy}, put $V:=\hat H_{b,a}-H_b$ and Taylor-expand
the ambient entropy along $H_b+tV$ through degree three. Metric and
density are uniformly positive on this segment. Each term in the
fourth-order integral remainder contains four factors of $V$ or its
spatial derivatives; the denominators are uniformly bounded and
the background has only fixed polynomial growth. The pointwise
bounds above therefore make its Gaussian integral $O(|a|^4)$.
Replacing $V$ in the Taylor coefficients by $H^{[2]}-H_b$ changes
them only by flat errors. Their terms of $a$-degree at most three
are the defining cubic Taylor polynomial $F_2$; every remaining term
is $O(|a|^4)$. This proves~\eqref{eq:cutenergy}.
\end{proof}

\subsection{Logarithmic bounds and Gaussian estimates}\label{part5-diagnostic-and-choice-of-topology}

In this subsection write $\mu=\mu_X=e^{\ell(h,k)}\bar\mu$.
In the variables of \eqref{eq:chart}--\eqref{eq:density-ratio}, we allow
quadratic growth of the potential difference. The logarithmic smallness condition is
\begin{equation}\label{eq:L-1-2}
 \|h\|_{C^2}+
 \sup_{\bar M}\left\{\frac{|k|}{1+|\vec x|^2}
       +\frac{|\bar\nabla k|}{1+|\vec x|}+|\bar\nabla^2 k|\right\}
 \le\epsilon.
\end{equation}
All pairs are smooth. For sufficiently small $\epsilon$, the determinant formula gives
$|\bar\nabla\ell|\le C\epsilon\langle\vec x\rangle$, so
Lemma~\ref{lem:gaussian-sobolev} applies with measure $\mu$.

\subsubsection{Scalar estimate}\label{part5-uniform-physical-scalar-inverse-without-bounded-density-ratio}

\begin{lemma}\label{lem:log-scalar}
Suppose that $\mu(\bar M)=1$ and \eqref{eq:L-1-2} holds with sufficiently
small tolerance. Then
\begin{equation}\label{eq:L-2-1}
 \|w-\mathcal W(X)\|_{L^2(\mu)}\le C\mathcal D.
\end{equation}
The constants are uniform over this class of pairs.
\end{lemma}
\begin{proof}
For $F\in L^2(\mu)$, consider the mean-corrected equation
\[
 \Delta_f v-\int_{\bar M}v\;\mathrm d\mu=F,
 \qquad \|v\|_{H^2(\mu)}\le C\|F\|_{L^2(\mu)}.
\]
Set $\varphi:=\bar f-\ell$ and
$L_\ell:=\Delta_{\bar g}-\langle\bar\nabla\varphi,\bar\nabla\cdot\rangle$,
which is self-adjoint for $\mu$. For small $\epsilon$,
\[
 \Ric_{\bar g}+\bar\nabla^2\varphi
 =\frac12\bar g-\bar\nabla^2\ell
 \ge\frac14\bar g.
\]
The Bakry--\'Emery Poincar\'e inequality
\cite[Theorem, part~1, p.~201]{BE85} and the integrated Bochner identity give
\[
 \left\|v-\int_{\bar M}v\;\mathrm d\mu\right\|_{L^2(\mu)}^2
 \le4\|\bar\nabla v\|_{L^2(\mu)}^2,
 \qquad
 \|\bar\nabla^2v\|_{L^2(\mu)}^2
 +\frac14\|\bar\nabla v\|_{L^2(\mu)}^2
 \le\|L_\ell v\|_{L^2(\mu)}^2.
\]
The Dirichlet form therefore gives an inverse of $L_\ell$ on the
mean-zero subspace. Local elliptic regularity and the Bochner estimate
place its solutions in $H^2(\mu)$. These estimates follow first for
compactly supported smooth functions and extend to the operator domain
by regularization and cutoffs: $|\bar\nabla\varphi|\le C\langle\vec x\rangle$,
and cutoffs on $B_{2R}\setminus B_R$ have derivatives of orders one and
two bounded by $C/R$ and $C/R^2$. The resulting boundary terms tend to
zero in the form and graph norms. Consequently
\[
 v\longmapsto L_\ell v-\int_{\bar M}v\;\mathrm d\mu
\]
is an isomorphism from $H^2(\mu)$ to $L^2(\mu)$ with uniform inverse
bound; on constants it acts as minus the identity.

We compare $\Delta_f$ with $L_\ell$ on these spaces. Its divergence
expression in reference indices gives
\[
  (\Delta_f-L_\ell)v
 =(g^{ij}-\bar g^{ij})\bar\nabla_i\bar\nabla_jv
 +\left[\bar\nabla_i g^{ij}
 +(g^{ij}-\bar g^{ij})(\partial_i\ell-\partial_i\bar f)\right]
       \bar\nabla_jv.
 \]
The second-order coefficient is bounded by $C\epsilon$, and the
first-order coefficient by $C\epsilon\langle\vec x\rangle$.
Lemma~\ref{lem:physical-sobolev} yields
\[
 \|(\Delta_f-L_\ell)v\|_{L^2(\mu)}
 \le C\epsilon\|v\|_{H^2(\mu)}.
\]
A Neumann series proves the mean-corrected inverse.

Put $u:=w-\mathcal W(X)$. The inverse supplies the
mean-zero solution of $\Delta_f v=u$. Indeed, integration of the
inverse equation shows that $\int_{\bar M}v\;\mathrm d\mu=0$.
The weighted contracted Bianchi identity used in the proof of \cite[Lemma~2.10]{FLR} is
$\mathrm dw=-2\delta_X\Phi$. Hence
\[
 \begin{split}
 \|u\|_{L^2(\mu)}^2
 &=\int_{\bar M}u\Delta_f v\;\mathrm d\mu
   =-\int_{\bar M}\langle\mathrm dw,\mathrm dv\rangle_g\;\mathrm d\mu\\
 &=-2\int_{\bar M}\langle\Phi,\nabla_g^2v\rangle_g\;\mathrm d\mu
 \le C\mathcal D\,\|v\|_{H^2(\mu)}
 \le C\mathcal D\,\|u\|_{L^2(\mu)}.
 \end{split}
\]
The weak identity is first tested against compactly supported smooth
functions and then extended by $H^2(\mu)$ approximation. Both $u$ and
$\Phi$ belong to $L^2(\mu)$ by their polynomial bounds, and
$\|\nabla_g^2v\|_{L^2(\mu;g)}\le C\|v\|_{H^2(\mu)}$.
This proves \eqref{eq:L-2-1}.
\end{proof}

\subsubsection{Comparison of Gaussian weights}\label{part5-two-power-comparisons-for-gaussian-weights}

Suppose $|\ell|\le\delta(1+|\vec x|^2)$, where $0\le\delta<1/4$.
Let $F$ satisfy $|F|\le C_0\langle\vec x\rangle^Q$ for fixed $C_0,Q$.
Constants in the comparisons depend on these bounds, on $\delta$,
and on the prescribed exponent within its displayed strict range. Then
\begin{equation}\label{eq:L-3-1}
 0<p_2<\frac1{1+4\delta},\qquad
 \|F\|_{L^2(\bar\mu)}\le C_{p_2}\|F\|_{L^2(\mu)}^{p_2}
 \quad\text{when }\|F\|_{L^2(\mu)}\le1.
\end{equation}
Indeed, H\"older's inequality gives
\[
 \int|F|^2\;\mathrm d\bar\mu
 \le\left(\int|F|^2e^\ell\;\mathrm d\bar\mu\right)^{p_2}
 \left(\int|F|^2e^{-p_2\ell/(1-p_2)}\;\mathrm d\bar\mu\right)^{1-p_2}.
\]
The second factor is uniformly bounded because
$p_2\delta/(1-p_2)<1/4$ and the reference Gaussian density decays as
$e^{-|\vec x|^2/4}$.

If $u$ and its first two derivatives have uniform polynomial bounds,
then
\begin{equation}\label{eq:L-3-2}
 0<b_0<1-4\delta,\qquad
 \|u\|_{H^2(\mu)}\le C_{b_0}\|u\|_{H^2}^{b_0}
 \quad\text{when }\|u\|_{H^2}\le1.
\end{equation}
For each $V=\bar\nabla^ju$, $0\le j\le2$, apply
\[
 \int|V|^2e^\ell\;\mathrm d\bar\mu
 \le\left(\int|V|^2\;\mathrm d\bar\mu\right)^{b_0}
 \left(\int|V|^2e^{\ell/(1-b_0)}\;\mathrm d\bar\mu\right)^{1-b_0}.
\]
The second factor is uniformly bounded since
$\delta/(1-b_0)<1/4$. Taking square roots and summing over the three
derivative orders proves~\eqref{eq:L-3-2}.

\subsection{Endpoint moments and residual comparison}
\label{subsec:physical-comparison-preparation}
Fix a small tolerance $\epsilon_0>0$. Let $X=(\bar g+h,\bar f+k)$ be
a smooth complete mass-one pair satisfying \eqref{eq:L-1-2} with this
tolerance.
\begin{enumerate}[label=\textup{(\roman*)},leftmargin=*]
\item \emph{Endpoints.} Set $(b,a)=\Pi_Kh$ and take
$\hat X=\hat X_{b,a}$ from Proposition~\ref{prop:cutoffjets}. Write
\[
 \hat X=(\hat g,\hat f)=(\bar g+\hat h,\bar f+\hat k)
                    =\mathcal C(\hat H_{b,a})
\]
and assume $\|\hat H_{b,a}\|_{C^2}\le\epsilon_0$.
\item \emph{Differences and errors.} Set
\[
 \begin{gathered}
 a_{\mathrm{sl}}:=\|(\Id-\Pi_{\mathrm{sl}})h\|_{H^2},\qquad
 \vec m:=\int\vec x\,\mathrm d\mu_X,\qquad
 J:=X-\hat X=(\eta,\kappa),\\
 U:=(\eta,\chi),\qquad
 \chi:=\frac12\Tr_{\bar g}\eta-\kappa,\qquad
 \tau:=\|J\|_{H^2}.
 \end{gathered}
\]
\item \emph{Comparison segment.} Define
\begin{equation}\label{eq:physical-comparison-segment}
 X_t:=\hat X+t J=(g_t,f_t),\qquad
 \mu_t:=\mu_{X_t},\qquad
 \ell_t:=\log\frac{\mathrm d\mu_t}{\mathrm d\bar\mu},
 \qquad 0\le t\le1.
\end{equation}
\end{enumerate}
The representative has mass one and center zero. The density chart
gives $\|(\hat h,\hat k)\|_{C^2}\le C\epsilon_0$. Since the reference trace is parallel,
$\tau\asymp\|U\|_{H^2}$.
Since the metric differences vary affinely and are uniformly $C^2$-small, the
metrics $g_t$ are uniformly equivalent to $\bar g$ and hence complete. The potential differences also vary affinely and obey
\eqref{eq:L-1-2} with tolerance $C\epsilon_0$. The determinant formula
therefore gives, uniformly in $t$,
\begin{equation}\label{eq:comparison-log-bounds}
 |\ell_t|\le\delta(1+|\vec x|^2),\qquad
 |\bar\nabla\ell_t|\le C\epsilon_0\langle\vec x\rangle,
 \qquad |\bar\nabla^2\ell_t|\le C\epsilon_0,
 \qquad \delta=C\epsilon_0.
\end{equation}
The endpoint hypotheses also imply the polynomial bounds
\begin{equation}\label{eq:comparison-polynomial-bounds}
 |\bar\nabla^j\eta|\le C\epsilon_0,\qquad
 |\bar\nabla^j\kappa|\le C\epsilon_0
                       \langle\vec x\rangle^{2-j},
 \qquad 0\le j\le2.
\end{equation}
All the measures $\mu_t$ have uniform Gaussian moments.
For each fixed integer $j\ge0$,
\begin{equation}\label{eq:endpoint-density-moments}
 \|\langle\vec x\rangle^j(e^{\ell_t}-1)\|_2
       \le C_j\epsilon_0,
 \qquad
 \|\langle\vec x\rangle^je^{\ell_t}\|_2\le C_j.
\end{equation}
Equation~\eqref{eq:comparison-log-bounds} and Gaussian integrability
give~\eqref{eq:endpoint-density-moments}.

Differentiating the measure along this segment gives
\begin{equation}\label{eq:affine-density-velocity}
 q_t=\dot\ell_t
       =\frac12\Tr_{g_t}\eta-\kappa,
 \qquad
 \dot q_t=-\frac12\Tr\bigl((g_t^{-1}\eta)^2\bigr).
\end{equation}
In particular, the moving-trace correction is
\[
 q_t-\chi
       =\frac12(g_t^{ij}-\bar g^{ij})\eta_{ij},\qquad
 |q_t-\chi|\le C\epsilon_0|\eta|.
\]
The endpoint masses are one and the representative's center is zero;
intermediate masses and centers remain unconstrained. The uniform
Gaussian and polynomial bounds justify differentiation of their
integrals by dominated convergence.

\begin{lemma}
\label{lem:endpoint-moments}
Under the preceding assumptions, for sufficiently small $\epsilon_0$,
\begin{equation}\label{eq:endpoint-moments}
 \left|\int\chi\,\mathrm d\bar\mu\right|
 +\left|\int\vec x\,\chi\,\mathrm d\bar\mu-\vec m\right|
 \le C\epsilon_0\tau.
\end{equation}
Consequently, for every prescribed fixed $A>0$,
\begin{equation}\label{eq:endpoint-normal-decomposition}
 U=v+\zeta,\qquad v=\Pi_ZU\in Z,\qquad
 \|\zeta\|_{H^2}
 \le C(a_{\mathrm{sl}}+|\vec m|)
              +C\epsilon_0\tau+C_A|a|^A.
\end{equation}
\end{lemma}
\begin{proof}
The endpoint mass and center conditions, together with
\eqref{eq:affine-density-velocity}, give exactly
\[
 \int_0^1\int q_t e^{\ell_t}
                    \,\mathrm d\bar\mu\,\mathrm dt=0,
 \qquad
 \int_0^1\int\vec x\,q_t e^{\ell_t}
                    \,\mathrm d\bar\mu\,\mathrm dt=\vec m.
\]
For either the constant or a coordinate function $P\in\{1,x_1,\ldots,x_d\}$,
Cauchy--Schwarz and \eqref{eq:endpoint-density-moments} give
\[
 \left|\int P(q_t e^{\ell_t}-\chi)\,\mathrm d\bar\mu\right|
 \le\|\chi\|_2\|P(e^{\ell_t}-1)\|_2
       +C\epsilon_0\|\eta\|_2\|Pe^{\ell_t}\|_2
 \le C\epsilon_0\tau.
\]
Integration in $t$ proves \eqref{eq:endpoint-moments}.

The fixed orthogonal splitting gives
\[
 (\Id-\Pi_Z)U
 =\big((\Id-\Pi_{\mathrm{sl}})\eta,0\big)
       +(\Pi_K\eta,0)
       +\big(0,(\Id-\Pi_{S_0})\chi\big).
\]
By~\eqref{eq:representative-flat-errors} and the fixed projection bounds,
\[
 \|(\Id-\Pi_{\mathrm{sl}})\hat h\|_{H^2}
       +\|\Pi_K\hat h-(b+a)\|_{H^2}
       \le C_A|a|^A.
\]
It follows that the first two terms in the splitting are bounded
by $a_{\mathrm{sl}}+C_A|a|^A$.
Finally,
\[
 (\Id-\Pi_{S_0})\chi
 =\int\chi\,\mathrm d\bar\mu
       +\frac12\sum_i x_i\int x_i\chi\,\mathrm d\bar\mu.
\]
The $H^2$ norms of these fixed Hermite functions, together with
\eqref{eq:endpoint-moments}, bound this last term by
$C|\vec m|+C\epsilon_0\tau$. This proves
\eqref{eq:endpoint-normal-decomposition}.
\end{proof}

For the residual comparison freeze the density and inverse metric
of the cutoff comparison pair. Put
$\hat{\mathrm A}:=\hat g^{-1}\bar g$ and
$\hat\ell:=\log(\mathrm d\mu_{\hat X}/\mathrm d\bar\mu)$, and define
\begin{equation}\label{eq:fixed-representative-residual}
 \mathcal R_{\hat X}(Y):=
 \left(
 \Pi_{\mathrm{sl}}[e^{\hat\ell}\hat{\mathrm A}
                          \Phi(Y)\hat{\mathrm A}],
 \ \Pi_{S_0}[w(Y)-1]\right).
\end{equation}

\begin{lemma}
\label{lem:fixed-residual-comparison}
Under the preceding assumptions,
\begin{equation}\label{eq:fixed-residual-comparison}
 \mathcal R_{\hat X}(X)-\mathcal G(\hat H_{b,a})
       =\mathrm D\mathcal G(0)[U]+\mathcal E,
 \qquad \|\mathcal E\|_2\le C\epsilon_0\tau.
\end{equation}
The linearization of $\mathcal G$ at zero on smooth compactly supported
pairs extends to a bounded operator $H^2\to L^2$. For
$U=(U_{\mathrm t},U_{\mathrm s})\in H^2$ it is given by
\begin{equation}\label{eq:ambient-reference-derivative}
 \mathrm D\mathcal G(0)[U]=
 \left(
 \frac12\Pi_{\mathrm{sl}}\mathcal L_{\bar f}U_{\mathrm t},
 \ \Pi_{S_0}\big[-(2\Delta_{\bar f}+1)U_{\mathrm s}
                   +\delta_{\bar X}\delta_{\bar X}U_{\mathrm t}\big]
 \right).
\end{equation}
Its tensor component is orthogonal to $K$.
\end{lemma}
\begin{proof}
Since $\hat X=\mathcal C(\hat H_{b,a})$,
$e^{\hat\ell}=1+(\hat H_{b,a})_{\mathrm s}$.
The bounded $C^2$ smallness of this input gives
\[
 \|e^{\hat\ell}-1\|_{C^2}
       +\|\hat{\mathrm A}-\Id\|_{C^2}\le C\epsilon_0.
\]
Thus $Q\mapsto e^{\hat\ell}\hat{\mathrm A}Q\hat{\mathrm A}$
differs from the identity on $L^2$ by at most $C\epsilon_0$. Formula~\eqref{eq:density-lift-gradient} gives the
exact endpoint identity
\[
 \mathcal R_{\hat X}(\hat X)=\mathcal G(\hat H_{b,a}).
\]

Along~\eqref{eq:physical-comparison-segment}, express tensor
components using the fixed reference connection.
Let $C_t=\nabla^{g_t}-\bar\nabla$. The metric hypotheses give
\[
 \|g_t^{-1}-\bar g^{-1}\|_{C^2}
       +\|C_t\|_{C^1}\le C\epsilon_0,
 \qquad |\Rm_{g_t}-\Rm_{\bar g}|\le C\epsilon_0.
\]
Write $f_t=\bar f+k_t$. The logarithmic bounds in the
common setup give $|\mathrm df_t|\le C\langle\vec x\rangle$,
and the difference of its reference Hessian from that of $\bar f$
is $O(\epsilon_0)$. When expressed using $g_t$, its Hessian
also contains $-C_t\,\mathrm df_t$, bounded by
$C\epsilon_0\langle\vec x\rangle$.
In particular
$|\Phi(X_t)|\le C\epsilon_0\langle\vec x\rangle$.

The moving-trace correction from the common setup obeys, for
$0\le r\le2$,
\[
 |\bar\nabla^r(q_t-\chi)|
       \le C\epsilon_0\sum_{s=0}^r|\bar\nabla^s\eta|.
\]
Apply \eqref{eq:physical-linearizations} to $J$ at $X_t$.
The second-order terms in
$\mathcal L_{f_t}\eta$, $\delta_{X_t}^*\delta_{X_t}\eta$,
$\nabla_{g_t}^2q_t$, $\Delta_{f_t}q_t$,
and $\delta_{X_t}\delta_{X_t}\eta$ have coefficients
formed from inverse metrics. Their differences from the reference
coefficients are $O(\epsilon_0)$. The first-order terms additionally
contain one connection coefficient or one factor of $\mathrm df_t$;
their differences are bounded by
$C\epsilon_0\langle\vec x\rangle$.
The zeroth-order terms contain curvature, a differentiated connection,
a Hessian of $f_t$, products of connection and drift coefficients,
or, in the double weighted divergence, two drift factors.
Their differences are bounded by
$C\epsilon_0\langle\vec x\rangle^2$.
For example the drift-square term has difference
\[
 \left|g_t^{-1}\mathrm df_t
            \otimes g_t^{-1}\mathrm df_t
       -\bar g^{-1}\mathrm d\bar f\otimes
                          \bar g^{-1}\mathrm d\bar f\right|
       \le C\epsilon_0\langle\vec x\rangle^2.
\]
Terms containing $\Phi(X_t)$ obey the same zeroth-order
bound. These coefficients involve at most two spatial derivatives of
$g_t$ and $f_t$; the moving-trace terms are controlled
by the preceding bounds for $q_t-\chi$.

Consequently the differences of the two differentials are
sums of contractions of the form
\begin{equation}\label{eq:physical-differential-coefficients}
 a_2\bar\nabla^2J+a_1\bar\nabla J+a_0J,
 \qquad |a_2|\le C\epsilon_0,\quad
 |a_1|\le C\epsilon_0\langle\vec x\rangle,\quad
 |a_0|\le C\epsilon_0\langle\vec x\rangle^2.
\end{equation}
Lemma~\ref{lem:gaussian-sobolev} therefore proves
\begin{equation}\label{eq:physical-differential-comparison}
 \|\mathrm D\Phi_{X_t}[J]-\mathrm D\Phi_{\bar X}[J]\|_2
 +\|\mathrm Dw_{X_t}[J]-\mathrm Dw_{\bar X}[J]\|_2
       \le C\epsilon_0\|J\|_{H^2}.
\end{equation}
The reference differentials satisfy the same coefficient bounds
without the factor $\epsilon_0$. By Lemma~\ref{lem:gaussian-sobolev}
they extend by density from smooth compactly supported pairs to
bounded operators $H^2\to L^2$.

At the reference pair, $\frac12\Tr_{\bar g}\eta-\kappa=\chi$
and $\Phi(\bar X)=0$.
The tensor projection annihilates
$\delta_{\bar X}^*\delta_{\bar X}\eta$ and
$\bar\nabla^2\chi$: both are in the reference gauge range, in the
weak sense by integration by parts. The scalar row retains its double
divergence. Thus the fixed projections of
$(\mathrm D\Phi_{\bar X}[J],\mathrm Dw_{\bar X}[J])$ give exactly
\eqref{eq:ambient-reference-derivative}, also as stated in
\cite[Lemma~2.8]{FLR}.
For any $a_0\in K$, self-adjointness on the reference domain gives
\[
 (\Pi_{\mathrm{sl}}\mathcal L_{\bar f}U_{\mathrm t},a_0)_2
   =(U_{\mathrm t},\mathcal L_{\bar f}a_0)_2=0.
\]
Hence the tensor kernel projection vanishes for every ambient $U$.

Finally integrate \eqref{eq:physical-differential-comparison} in
$t$, multiply its tensor row by the frozen representative
factors, and apply the fixed $L^2$ projections. The multiplier's
$O(\epsilon_0)$ difference from the identity, applied to the reference
tensor differential, contributes at most $C\epsilon_0\tau$.
This proves \eqref{eq:fixed-residual-comparison}.
\end{proof}

\subsection{Affine entropy comparison}
\label{part5-affine-entropy-comparison}

\begin{lemma}\label{lem:affine-entropy}
Under the assumptions of
Subsection~\ref{subsec:physical-comparison-preparation}, for every
fixed $0<b_0<1$, sufficiently small logarithmic and representative
neighborhoods give
\begin{equation}\label{eq:L-6-4}
 |\mathcal W(X)-\mathcal W(\hat X)|
 \le C\mathcal D \tau^{b_0}+C\tau^{2b_0},
 \qquad \mathcal D=\mathcal D(X).
\end{equation}
The neighborhood thresholds and the constant depend only on $b_0$,
the fixed background, and the prescribed polynomial bounds.
\end{lemma}
\begin{proof}
Choose the neighborhoods so that $\tau\le1$ and $b_0<1-4\delta$.
The polynomial bounds \eqref{eq:comparison-polynomial-bounds} and
the logarithmic bounds \eqref{eq:comparison-log-bounds} allow us to apply
\eqref{eq:L-3-2} uniformly along
\eqref{eq:physical-comparison-segment}. Thus
\[
 \|J\|_{H^2(\mu_t)}\le C\tau^{b_0},
 \qquad 0\le t\le1.
\]
The power comparison and Lemma~\ref{lem:physical-sobolev} use only
the logarithmic coefficient bounds, so apply also to the
intermediate measures, whose masses need not equal one.

Put
\[
 F(t)=\mathcal W(X_t),\qquad
 V(t)=\int\mathrm d\mu_t.
\]
The formulas for $q_t$ and $\dot q_t$ give the exact
mass Hessian
\[
 V''(t)
 =\int\left(q_t^2
       -\frac12\Tr\bigl((g_t^{-1}\eta)^2\bigr)\right)
       \;\mathrm d\mu_t.
\]
Uniform metric equivalence and the power comparison imply
\[
 |V''(t)|\le C\|J\|_{L^2(\mu_t)}^2
                   \le C\tau^{2b_0}.
\]
Since $V(1)=V(0)=1$, the endpoint Taylor identity yields
\[
 V'(1)=\int_0^1t V''(t)\;\mathrm dt,
 \qquad |V'(1)|\le C\tau^{2b_0}.
\]

The common metric and potential bounds give
$|w(X_t)|\le C\langle\vec x\rangle^2$: the curvature is
bounded, $|\mathrm df_t|\le C\langle\vec x\rangle$, and
\[
 \nabla_{g_t}^2f_t
 =\bar\nabla^2f_t
   -(\Gamma_{g_t}-\Gamma_{\bar g})\bar\nabla f_t
\]
has at most linear growth. Uniform Gaussian moments therefore give
$|\mathcal W(X)|\le C$. At $t=1$, \eqref{eq:W-first} reads
\[
 F'(1)=\int\bigl(\langle\Phi(X),\eta\rangle_g
       +(w(X)-\mathcal W(X))q_1\bigr)\;\mathrm d\mu_X
       +(\mathcal W(X)-1)V'(1).
\]
Lemma~\ref{lem:log-scalar}, the estimate
$\|q_1\|_{L^2(\mu_X)}\le C\|J\|_{L^2(\mu_X)}$, and the mass
derivative bound consequently give
\[
 |F'(1)|\le C\mathcal D\|J\|_{H^2(\mu_X)}+C\tau^{2b_0}
          \le C\mathcal D \tau^{b_0}+C\tau^{2b_0}.
\]

Differentiate~\eqref{eq:W-first} along the affine path. Since
$\partial_tg_t=\eta$ and $\partial_t\mathrm d\mu_t=q_t\mathrm d\mu_t$,
\[
 \begin{split}
 F''(t)=\int\bigl(&\langle\mathrm D\Phi_{X_t}[J],\eta\rangle_{g_t}
       +\mathrm Dw_{X_t}[J]q_t
       -2\langle\Phi(X_t),\eta\circ_{g_t}\eta\rangle_{g_t}\\
       &+q_t\langle\Phi(X_t),\eta\rangle_{g_t}
       +(w(X_t)-1)(\dot q_t+q_t^2)\bigr)\;\mathrm d\mu_t.
 \end{split}
\]
The coefficient calculation in Lemma~\ref{lem:fixed-residual-comparison}
gives the pointwise bound
\[
 |\mathrm D\Phi_{X_t}[J]|+|\mathrm Dw_{X_t}[J]|
 \le C\bigl(|\bar\nabla^2J|
       +\langle\vec x\rangle|\bar\nabla J|
       +\langle\vec x\rangle^2|J|\bigr).
\]
Apply Lemma~\ref{lem:gaussian-sobolev} with the moving measure
$\mu_t$, whose logarithmic coefficients satisfy the same smallness
bounds. Thus both differentials have $L^2(\mu_t)$ norm at most
$C\|J\|_{H^2(\mu_t)}$. Moreover,
$|\Phi(X_t)|\le C\langle\vec x\rangle$,
$|w(X_t)-1|\le C\langle\vec x\rangle^2$,
$|q_t|\le C|J|$ and $|\dot q_t|\le C|\eta|^2$.
Cauchy--Schwarz and the same Gaussian multiplication estimate give
\[
 |F''(t)|\le C\|J\|_{H^2(\mu_t)}^2
                 \le C\tau^{2b_0}.
\]
The common polynomial bounds and uniform Gaussian tails justify
differentiation under the integrals and passage to the limit in the
cutoff first-variation formula, including its integrations by parts.

Finally, the endpoint identity
\[
 F(1)-F(0)=F'(1)-\int_0^1t F''(t)\;\mathrm dt
\]
proves \eqref{eq:L-6-4}.
\end{proof}

\subsection{Proof of the {\L}ojasiewicz--Simon inequality}
\label{part5-theorem-s-logarithmic-entropy-and-normal-comparison}

\begin{theorem}[Entropy and comparison estimates]\label{thm:logarithmic}
Put $p=\min\{3/2,p_N\}$. For every fixed $1<p_1<p$, there are
$\epsilon,C>0$ such that every smooth complete mass-one pair $X$ satisfying
\eqref{eq:L-1-2} obeys the following estimates, with the metric kernel
projection, positive cutoff pair, slice error, and center defined
in Subsection~\ref{subsec:physical-comparison-preparation}:
\begin{equation}\label{eq:L-6-1}
 |\mathcal W(X)-\Theta_N|
       \le C(\mathcal D+a_{\mathrm{sl}}+|\vec m|)^{p_1}.
\end{equation}
With $p_2=p_1/p$, the comparison in metric--potential variables is
\begin{equation}\label{eq:L-6-2}
 \begin{split}
 |a|^2+|\mathrm De(b)|+|\mathrm DF_2(b,a)|
       +\|(h,k)-(\hat h,\hat k)\|_{H^2}
 \le C(\mathcal D^{p_2}+a_{\mathrm{sl}}+|\vec m|).
 \end{split}
\end{equation}
The thresholds and estimate constants depend only on the prescribed
$p_1$ and the fixed background data.
\end{theorem}
\begin{proof}
Fix $b_0$ with $p/2<b_0<1$ and put $p_2:=p_1/p$.
Choose a coefficient tolerance $\epsilon_0>0$ small enough for the
normal absorption below and for all the preceding lemmas. Also require
that the resulting bound $\delta=C\epsilon_0$ satisfies
\begin{equation}\label{eq:unified-exponent-choices}
 p_2<\frac1{1+4\delta},\qquad b_0<1-4\delta.
\end{equation}
Since $p_2,b_0<1$, these conditions hold for small $\epsilon_0$.
Fix this tolerance and then choose the input tolerance $\epsilon$. We have
$|b|+|a|\le C\|h\|_2\le C\epsilon$, whereas
Proposition~\ref{prop:cutoffjets} gives
\[
 \|\hat H_{b,a}\|_{C^2}\le C(|b|+|a|^{1/2}).
\]
For sufficiently small $\epsilon$, both $\hat H_{b,a}$ and
$(\hat h,\hat k)$ have bounded $C^2$ norm at most $\epsilon_0$.
Thus Subsection~\ref{subsec:physical-comparison-preparation} applies,
with~\eqref{eq:comparison-log-bounds} at the chosen $\delta$.
The logarithmic coefficient bounds and Gaussian moments also make
$\mathcal D$, $a_{\mathrm{sl}}$, $|\vec m|$, and $\tau$ small.

Use the fixed representative residual
$\mathcal R_{\hat X}$ from \eqref{eq:fixed-representative-residual}.
Its tensor multiplier is bounded uniformly, since
$e^{\hat\ell}=1+(\hat H_{b,a})_{\mathrm s}$ and
$\hat g^{-1}\bar g$ are bounded and close to the identity.
Lemma~\ref{lem:log-scalar} controls
$w(X)-\mathcal W(X)$ in
$L^2(\mu_X)$ by $C\mathcal D$. This scalar field and $\Phi(X)$
have fixed polynomial bounds. Apply \eqref{eq:L-3-1} to each,
then apply the fixed projections. Since $\Pi_{S_0}$ annihilates
constants, this gives
\begin{equation}\label{eq:actual-fixed-residual-bound}
 \|\mathcal R_{\hat X}(X)\|_2\le C\mathcal D^{p_2}.
\end{equation}

Set $v:=\Pi_ZU$ and $\zeta:=(\Id-\Pi_Z)U$.
Lemma~\ref{lem:endpoint-moments}, with a fixed choice $A=4$, yields
\begin{equation}\label{eq:unified-off-normal-bound}
 \|\zeta\|_{H^2}
 \le C(a_{\mathrm{sl}}+|\vec m|)+C\epsilon_0\tau+C|a|^4.
\end{equation}
Lemma~\ref{lem:fixed-residual-comparison} gives
\[
 \mathcal R_{\hat X}(X)-\mathcal G(\hat H_{b,a})
   =\mathfrak B_0v+\mathrm D\mathcal G(0)[\zeta]+\mathcal E,
 \qquad \|\mathcal E\|_2\le C\epsilon_0\tau.
\]
Here $\mathrm D\mathcal G(0)$ is the ambient derivative in
\eqref{eq:ambient-reference-derivative}; it is bounded from
$H^2$ to $L^2$. Its restriction to $Z$ is $\mathfrak B_0$.
Project onto $Z$ and apply~\eqref{eq:finiteinverse} to $v\in Z$. Equations
\eqref{eq:cutresidual}, \eqref{eq:actual-fixed-residual-bound},
and \eqref{eq:unified-off-normal-bound}, together with
$\|U\|_{H^2}\asymp\tau$, imply
\[
 \tau\le C(\mathcal D^{p_2}+a_{\mathrm{sl}}+|\vec m|+|a|^3)
             +C\epsilon_0\tau.
\]
The tensor kernel projection of $\mathrm D\mathcal G(0)[U]$
vanishes for every ambient $U$, as proved in
Lemma~\ref{lem:fixed-residual-comparison}. Thus the kernel
projection of the same identity gives
\[
 |\mathrm DF_2|
 \le C(\mathcal D^{p_2}+|a|^3)+C\epsilon_0\tau.
\]
Combining these inequalities gives
\begin{equation}\label{eq:unified-preabsorption}
 \tau+|\mathrm DF_2|
 \le C_1(\mathcal D^{p_2}+a_{\mathrm{sl}}+|\vec m|+|a|^3)
                  +C_2\epsilon_0\tau.
\end{equation}
The initial choice $C_2\epsilon_0\le1/2$ absorbs the final term.
Put $\Delta:=\mathcal D^{p_2}+a_{\mathrm{sl}}+|\vec m|$.
By \eqref{eq:closing},
$|a|^3\le C|a|\,|\mathrm DF_2|$. Shrink the parameter ball so
that the resulting coefficient of $|\mathrm DF_2|$ is at most
$1/2$, and absorb it as well. We obtain
\begin{equation}\label{eq:unified-normal-bound}
 \tau+|\mathrm DF_2|+|a|^2+|\mathrm De(b)|\le C\Delta.
\end{equation}
This is \eqref{eq:L-6-2}.

Since $X$ and $\hat X$ have mass one,
Lemma~\ref{lem:affine-entropy} and~\eqref{eq:unified-exponent-choices} give
\[
 |\mathcal W(X)-\mathcal W(\hat X)|
       \le C\mathcal D \tau^{b_0}+C\tau^{2b_0}.
\]
By \eqref{eq:unified-normal-bound}, $\tau\le C\Delta$ and
$\mathcal D\le\Delta^{1/p_2}$. Since
$1/p_2>1>b_0$ and $2b_0>p$, for $0\le\Delta\le1$,
\[
 \mathcal D \tau^{b_0}+\tau^{2b_0}
 \le C\bigl(\Delta^{1/p_2+b_0}+\Delta^{2b_0}\bigr)
 \le C\Delta^p.
\]
Proposition~\ref{prop:cutoffjets} and
Proposition~\ref{lem:closing} supply the remaining two comparisons:
\[
 |\mathcal W(\hat X)-F_2|\le C|a|^4\le C\Delta^2,
 \qquad |F_2-\Theta_N|\le C|\mathrm DF_2|^p\le C\Delta^p.
\]
Because $p\le3/2<2$, their sum is bounded by $C\Delta^p$.
Finally, for $\mathcal D+a_{\mathrm{sl}}+|\vec m|\le1$,
\[
 \Delta\le2(\mathcal D+a_{\mathrm{sl}}+|\vec m|)^{p_2}.
\]
Raising to the power $p$ proves \eqref{eq:L-6-1}.
\end{proof}

\section{Dynamical framework and extension--improvement}\label{sec:analytic}
We work on the fixed manifold $M$ with the modified pair defined by
\eqref{eq:modified-diffeomorphism}--\eqref{eq:modified-gauge}. Write
\begin{equation}\label{eq:modified-flow}
 \begin{gathered}
 \mathcal W(s):=\mathcal W(g^{z_*}(s),f^{z_*}(s)),\qquad
 E(s):=\Theta_N-\mathcal W(s),\qquad
 \mathcal D(s):=\|\Phi(s)\|_{L^2(\mu^{z_*})},\\
 \partial_sg^{z_*}=2\Phi,\qquad
 \partial_sf^{z_*}=\Tr_{g^{z_*}}\Phi,\qquad
 \partial_s\mathrm{d}\mu^{z_*}=0,\qquad
 E'(s)=-2\mathcal D(s)^2.
 \end{gathered}
\end{equation}
Moving comparison charts
will establish continuation; final convergence will be obtained by
integration on fixed subsets of $M$.

\phantomsection\label{subsec:terminal-tails}
For $0<\tau\le T/2$, put
$g_\tau=\tau^{-1}g(-\tau)$, $f_\tau=f(-\tau)$ and
$\mathrm d\mu_\tau=(4\pi)^{-n/2}e^{-f_\tau}\,\mathrm dV_{g_\tau}$.
Recall the definition of an $H_n$-center from
\cite[Definition~3.10]{Bam21}.

\begin{lemma}\label{lem:terminal-tails}
For every $\eta>0$, there is $C=C(n,Y,\eta)$ such that, for
$0<\tau\le T/2$, every $H_n$-center $(z_H,-\tau)$ of $z_*$, and every
$L>0$,
\begin{equation}\label{eq:terminal-tail}
 \int_{M\setminus B_{g_\tau}(z_H,L)}
 (1+|f_\tau|+|\mathrm R_{g_\tau}|+|\nabla f_\tau|^2)
 \,\mathrm d\mu_\tau
 \le Ce^{-L^2/(4+\eta)}.
\end{equation}
Moreover,
$f_\tau(x)\ge d_{g_\tau}(z_H,x)^2/(4+\eta)-C$.
All derivatives and norms here are taken with respect to $g_\tau$.
\end{lemma}
\begin{proof}
Rescale by $\tau$. The scalar-curvature maximum principle gives
\[
 \mathrm R_{g_\tau}\ge-\frac{n\tau}{2(T-\tau)}\ge-\frac n2.
\]
Thus the lower potential bound in \cite[Theorem~5.30]{FLSt} and
the integral and mass tails in
\cite[Proposition~3.1 and Theorem~2.13(i)]{FLC} have uniform constants.
For regular terminal points these estimates apply directly. Approximate
$z_*$ by regular terminal points. For each fixed $\tau$, their Nash
entropies are at least $-Y-1$ for sufficiently close approximations
by \cite[Lemma~7.2]{FLSt}. Smooth kernel convergence on negative-time
slices from \cite[Theorem~5.20]{FLSt} and Fatou's lemma pass the
integral estimates to the terminal kernel. The approximating $H_n$-centers have bounded rescaled distance from
$z_H$; starting with a strictly smaller Gaussian loss absorbs this
displacement.
\end{proof}

\phantomsection\label{lem:entropy-limit}
The selected tangent and \cite[Lemma~7.2 and proof of Proposition~7.9]{FLSt}
identify the limiting Nash entropy as $\Theta_N$; the entropy relation
and monotonicity in \cite[Proposition~2.8]{FLC} give the same limit for
the $\mathcal W$-entropy:
\begin{equation}\label{eq:entropy-limit}
 \mathcal W(s)\longrightarrow\Theta_N\qquad(s\longrightarrow\infty).
\end{equation}
Consequently, $E(s)\ge0$ decreases to zero by \eqref{eq:modified-flow}.

\phantomsection\label{idx:product-background}
In this section and below, $X_b$ denotes the product of the compact-factor
variational family with the Gaussian factor, with
$G_b:=(X_b)_{\mathrm{t}}=g_b+g_E$ and
$F_b:=(X_b)_{\mathrm{s}}=f_b+f_E$. Its associated probability is $\mu_b:=\mu_b^N\otimes\gamma_d$, where
$\mathrm d\mu_b^N:=(4\pi)^{-m/2}e^{-f_b}\;\mathrm dV_{g_b}$.
Its soliton defect $\Phi_b:=G_b/2-\Ric-\nabla^2F_b$ obeys $\|\Phi_b\|_{C^j}\le C_j|\mathrm{D}e(b)|$.
It is independent of the Euclidean variable.

\subsection{Overview of the iteration process}\label{subsec:iteration-overview}
Normalize a sufficiently late selected tangent chart so that the Ricci
flow is defined on $[-2,0)$. A block uses $[-1,-1/2]$, or
$[0,\log2]$ in modified time, and begins with an order-$100n$
comparison chart near $X_b$ on a sufficiently large region. We follow
the passage from the variational entropy estimate
\cite[Theorem~4.17]{FLC} to entropy drop in \cite[Sections~5--6]{FLC}.

Section~\ref{sec:physical} compares the given flow with
$\gamma_b(t)+g_E$, where $\gamma_b(-1)=g_b$, and compares the terminal
potential with an auxiliary forward scalar. Probability normalization
fixes the additive constant. Parabolic rescaling enlarges the available
region, while the entropy drop up to the current time determines the
radius on which the spatial gauge and stationary estimates apply.
Section~\ref{sec:restart} puts the endpoint in the reference divergence
slice and applies Theorem~\ref{thm:logarithmic} to its complete cutoff
extension. This improves comparison with $X_{b_+}$ and controls
$|\mathrm De(b_+)|$ and the positive Hermite component separately.
The endpoint estimates return order-$100n$ comparison data for
$g_+(t')=2g(t'/2)$ and $f_+(t')=f(t'/2)$.

The compact parameter may move along a nonisolated critical family,
so the endpoint estimate alone does not confine it to the construction
neighborhood. Proposition~\ref{prop:coherent-maps} normalizes the gauged
charts by a path of reference isometries starting at the identity.
The kernel projection of the cutoff metric perturbation in this chart defines the next
candidate parameter. Restricting the normalized chart gives the next
incoming chart, as in
Subsection~\ref{subsec:inherited-endpoint-coordinates}. If the candidate
belongs to the construction neighborhood, the next spatial gauge is
constructed with that parameter fixed. Proposition~\ref{prop:gauge-junction}
controls the transition on the common collar.

Proposition~\ref{prop:causal-radii} verifies the radii and endpoint
estimates. The resulting recurrence retains the spatial restriction
inherited from earlier blocks and gives entropy summability on every
finite initial segment, including a final partial block. The weak
motion estimates of Section~\ref{sec:motion} then control localized
metric pairings. A first-exit argument confines the compact parameter
and proves continuation. Entropy summation and radius growth allow
Section~\ref{sec:uniqueness} to integrate the modified-flow equations
on fixed subsets, yielding one embedding and quantitative convergence
on expanding potential sublevels.

\subsection{Comparison charts and one time block}\label{subsec:iteration-data}
The background product, entropy bound, reference charts and preliminary
stationary family are fixed. Dependence on these data is suppressed.
Each block starts with a $C^{100n}$ chart. Intrinsic smoothing and the
elliptic estimate in the endpoint slice (Lemma~\ref{lem:restart-coordinate-regularity})
give the higher coordinate bounds. The interpolation powers for metric
propagation, actual-pair comparison and restart are $1/2$.

The estimate constants are uniform in the compact parameter within
its stated neighborhood, the chart radius, the block number and the
duration of a partial block. A constant attached to a fixed core $N_r$
may depend on $r$.

Sections~\ref{sec:physical}--\ref{sec:restart} retain their local
parameters. Subsection~\ref{subsec:parameter-selection} introduces
and fixes the compatible iteration data in the order summarized in
Table~\ref{tab:iteration-data}.

\noindent\emph{Comparison data.}
The stationary and construction neighborhoods are defined in Stages~\ref{stage:parameter-pairings}--\ref{stage:iteration-neighborhoods}
of Subsection~\ref{subsec:parameter-selection}.
\begin{enumerate}[label=\textup{(\roman*)},leftmargin=*]
\item An \emph{initial comparison chart of order $100n$} at normalized
Ricci-flow time $t=-1$ is an embedding $\iota$ of a neighborhood of
$\overline N_S$, with parameter $b$ and error
\begin{equation}\label{eq:seed-error}
 \boldsymbol{\epsilon}(S):=\|\iota^*g(-1)-G_b\|_{C^{100n}(N_S)}
       +\|\iota^*f(-1)-F_b\|_{C^{100n}(N_S)}.
\end{equation}
Small error gives a positive lower bound for the actual terminal
probability of a fixed core.
\item A \emph{raw chart} is a comparison chart before the spatial
gauge change for the current block. This low-order construction uses
raw $C^8$ data
with $C^5$ displacement, and imposes the exact reference
metric slice and a centered auxiliary probability on an inner domain.
A \emph{reset} is this gauge change at a block boundary; Ricci-flow
propagation uses the initial order-$100n$ chart before the reset.
\item A finite initial segment of the iteration may end during a
block. For a block starting at $s_j$, its interval through the current
time is $[s_j,s_j+u]$, $0\le u\le\log2$. The comparison radius
in~\eqref{eq:E-10} uses only the incoming radius and entropy drop
on this interval.
\end{enumerate}

\phantomsection\label{idx:vanishing-modulus}
We write $\omega$ for a nondecreasing function on $[0,\epsilon_*]$
with $\omega(\epsilon)\to0$ as $\epsilon\to0$; it may be enlarged
when combining estimates. Its dependence is on the background, the
parameters in its estimate and, when indicated, a fixed support region $N_r$.
It is uniform in the number of blocks, the duration of a final
partial block and any prospective exit time, and does not depend on
the growing chart radius.

Initialize with one sufficiently late tangent chart of radius $S_0$
from the selected sequence, with $b_0=0$. Include the first gauge
change in the initial pairing error.

\section{Propagation along the Ricci flow and terminal-density comparison}
\label{sec:physical}
Propagation along the Ricci flow, terminal-kernel bounds and caloric comparison
extend one product chart over a time block. We follow
\cite[Propositions~6.2--6.3, Corollary~6.6 and Proposition~6.7]{FLC},
uniformly over the compact-factor variational family. The metric is first
compared with a product Ricci flow, and the actual terminal
potential is compared with an auxiliary forward scalar. Probability
normalization then determines the additive scalar constant.

\subsection{Common setup}\label{subsec:physical-common-setup}
Normalize the incoming Ricci-flow time to $t=-1$, and suppose that
the actual closed flow exists on $[-2,0)$ with the standing entropy
and terminal-point hypotheses. Put
\[
 s=-\log|t|,\qquad F(t):=|t|f(t).
\]
The main block is $[-1,-1/2]$, or $[0,\log2]$ in normalized time.
Use raw comparison order $8$. Prescribe a retention
$0<\rho<1$, a growth rate $1/8<c_1<1/4$, and an even integer
\[
 Q_*\ge1000n.
\]
Intrinsic smoothing supplies order $Q_*$ from incoming order $100n$.

The intermediate retentions and the common Gaussian loss are determined by
\[
 \rho_1:=\frac{3\rho+1}{4},\qquad
 \rho_2:=\frac{\rho+1}{2},\qquad
 \rho_3:=\frac{\rho+3}{4},\qquad
 c_3:=\frac{c_1-1/8}{32}.
\]
In particular $\rho<\rho_1<\rho_2<\rho_3<1$, and
\begin{equation}\label{eq:C-14}
 \frac{1+c_3}{16}+\frac{c_3}{2}<\frac{c_1}{2},
 \qquad \frac{c_1}{2}+c_3<\frac14.
\end{equation}
Proposition~\ref{prop:terminal-comparison} first chooses
$\eta=\eta(\rho,c_1,Q_*)\in(0,1/3)$ and then a metric
accuracy $\epsilon_{\rm met}>0$. If the additional intrinsic
estimate~\eqref{eq:terminal-arbitrary-rough} is required, its finite
order and Gaussian loss are prescribed before these choices.
The time and metric accuracy are then fixed throughout the section.
The comparison neighborhood is chosen afterward, so all its reference
flows exist on $[-1,-\eta]$.

For an embedding of a neighborhood of $\overline N_S$, the hypotheses
below are $b$ in the indicated common neighborhood,
$\boldsymbol{\epsilon}(S)\le\bar\epsilon$ and $S\ge S_*$, where the threshold
and minimum radius may be reduced or enlarged finitely many times.
For the current endpoint $t=-e^{-s}$, set
\[
 A(s)^2:=\mathcal W(|t|^{-1}g(t),f(t))
                     -\mathcal W(g(-1),f(-1)).
\]
Dependence on the comparison inputs means dependence on
$(\rho,c_1,Q_*)$ and the prescribed finite derivative orders
and Gaussian losses, with the time and metric accuracy chosen above.
All comparison neighborhoods and thresholds are independent of $S$;
constants are uniform in $b,S,\boldsymbol{\epsilon}(S)$ and the duration of a
partial block. The proofs supply positive Gaussian boundary-decay
rates; at the end of the section they are decreased to one common
output rate $c_4>0$.

\subsection{Metric propagation}\label{subsec:metric-propagation}
We first establish qualitative propagation in the given identification
following the argument of \cite[Proposition~6.2]{FLC}, then obtain uniform regularity and
quantitative comparison.

\medskip\noindent\textit{Qualitative propagation.}
For each fixed $\eta\in(0,1/3)$, continuous dependence on the
compact Einstein flow $(-t)g_N$ gives a parameter neighborhood on
which the compact Ricci flows $\gamma_b$, starting at $g_b$ at
$t=-1$, exist up to $-\eta$ with uniform finite-order bounds.
Write $\hat g_b(t):=\gamma_b(t)+g_E$.

\begin{lemma}\label{lem:qualitative-propagation}
For each $\eta\in(0,1/3)$ and each prescribed sufficiently small
$\epsilon_{\rm met}>0$, there are a parameter neighborhood
$B_{\rm met}=B_{\rm met}(\eta,\epsilon_{\rm met})$, an error
threshold $\bar\epsilon=\bar\epsilon(\eta,\epsilon_{\rm met})>0$,
a collar $C_0$ and a radius $S_*>C_0$, also depending only on
$(\eta,\epsilon_{\rm met})$, with the following property.
Every $b\in B_{\rm met}$, $S\ge S_*$ and initial comparison chart
with $\boldsymbol{\epsilon}(S)\le\bar\epsilon$ satisfy
\begin{equation}\label{eq:A-2}
 \sup_{-1\le t\le-\eta}
 \|\iota^*g(t)-\hat g_b(t)\|_{C^2(N_{S-C_0})}\le\epsilon_{\rm met},
\end{equation}
and, on the same cylinder,
\begin{equation}\label{eq:B-4}
 (1-\epsilon_{\rm met})\hat g_b(t)\le\iota^*g(t)
                         \le(1+\epsilon_{\rm met})\hat g_b(t).
\end{equation}
\end{lemma}
\begin{proof}
Apply the limiting argument of \cite[Proposition~6.2]{FLC} to the
compact family of reference flows. These have uniform metric,
curvature and volume bounds. At a time of sufficiently small $C^2$
comparison error, a first-exit argument places fixed-size ambient
balls inside the chart, with uniform volume lower bounds and curvature
upper bounds. The two-sided pseudolocality theorem of
Perelman~\cite[Theorem~10.1]{Per02} and
Bamler~\cite[Theorem~2.47]{BamStr}, in the combined form
\cite[Theorem~2.14]{FLC}, and Shi's local estimates~\cite{Shi89}
give intrinsic curvature derivative bounds on smaller balls over
uniform short windows inside $[-2,0)$.
To obtain convergence in the given coordinates, put
$g_0=\iota^*g$ and $\mathcal A=\nabla^{g_0}-\bar\nabla$. The identities
\[
 \partial_t\mathcal A^k{}_{ij}
 =-g_0^{kq}(\nabla_i\Ric_{jq}+\nabla_j\Ric_{iq}
                         -\nabla_q\Ric_{ij}),
 \qquad \bar\nabla g_0=g_0*\mathcal A
\]
propagate the finite coordinate bounds by induction; here
$\nabla=\nabla^{g_0}$ and $*$ denotes fixed contractions.
The incoming $C^{100n}$ bound gives spatial $C^6$ bounds and bounds for
$\partial_tg_0$ in $C^4$ on a smaller collar. Thus compactness holds
in $C^4$ on compact product regions, uniformly in time and in the
original identification.

If comparison failed on a forward window while the initial errors
tended to zero and the collars tended to infinity, recenter at a
failure point by Euclidean translation. After passing to a subsequence
$b_i\to b_\infty$, the limit is a complete bounded-curvature Ricci
flow starting from $\hat g_{b_\infty}(t_0)$ at the comparison time.
The complete bounded-curvature uniqueness theorem of
Chen--Zhu~\cite[Theorem~1.1]{CZ06} identifies the limit with
$\gamma_{b_\infty}(t_0+\cdot)+g_E$, a contradiction.
Cover $[-1,-\eta]$ by finitely many forward windows. Choose the
comparison tolerances backwards from the prescribed final accuracy
and sum their fixed collars to obtain $C_0$. Adjusting the final
accuracy gives both~\eqref{eq:A-2} and~\eqref{eq:B-4}.
The uniform reference and coordinate bounds make these choices
uniform on the closure of a sufficiently small parameter neighborhood.
\end{proof}

\medskip\noindent\textit{Intrinsic and coordinate regularity.}
Apply Lemma~\ref{lem:qualitative-propagation} with the prescribed
$\eta$ and with accuracy sufficient for $\epsilon_{\rm met}$.
Pseudolocality and Shi's estimates on the resulting cylinder, with
fixed spatial collars and larger time intervals inside $[-2,0)$, give
\begin{equation}\label{eq:A-4}
 |(\nabla^{g(t)})^j\Rm(g(t))|\le C_j
 \qquad\text{for every fixed }j.
\end{equation}
These bounds include $t=-1$; the earlier Ricci-flow interval
provides the interior estimates at that time.
For the coordinate bounds, differentiate the evolution of
$\mathcal A=\nabla^{\iota^*g}-\bar\nabla$ with $\bar\nabla$.
At order $r\le100n-3$, conversion to intrinsic derivatives involves
curvature derivatives through order $r+1$ and reference derivatives
of $\mathcal A$ of strictly lower order. The incoming $C^{100n}$
bound controls $\mathcal A(-1)$ through order $100n-3$.
Integrating successively and using
$\bar\nabla(\iota^*g)=(\iota^*g)*\mathcal A$ gives, on a fixed further collar,
\[
 \sup_{-1\le t\le-\eta}
 \|\iota^*g(t)\|_{C^{100n-2}(N_{S-C})}\le C_{100n-2}.
\]
Here $C>C_0$ is a fixed total collar, and $C_{100n-2}$ depends on the
background, compact family and fixed time interval, but not on
$S,b$ or the initial comparison chart. The reference flow has the
same bounds.

\begin{lemma}[Local metric comparison]\label{lem:local-metric-comparison}
Fix $0<\eta<1$. For $D\ge4$ and $\vec x_0\in\R^d$, put
$\Omega_D=N\times B_D(\vec x_0)$. Let $g,\tilde g$ be Ricci flows
on a neighborhood of $\overline\Omega_D\times[-1,-\eta]$, with
fixed uniform metric equivalence to $\bar g$, $C^{100n-2}$ coordinate
bounds, and intrinsic curvature derivative bounds through order two.
Suppose that, in reference coordinate norms,
\[
 e_{\rm in}=\|g(-1)-\tilde g(-1)\|_{C^2(\Omega_D)}\le1.
\]
Then
\begin{equation}\label{eq:A-7}
 \sup_{-1\le t\le-\eta}
 \|g(t)-\tilde g(t)\|_{C^{12}(\Omega_1)}
 \le C_*(e_{\rm in}^{1/2}+D^{d/2}e^{-c_*D^2}).
\end{equation}
The constants $C_*,c_*>0$ depend only on $\eta,n$ and the stated
metric-equivalence, coordinate and curvature bounds. They are
independent of $D$ and $\vec x_0$.
\end{lemma}
\begin{proof}
For the metric--connection--curvature energy of~\cite[Sections~1.1 and~2]{Kot}, set
\[
 h:=g-\tilde g,\qquad A:=\nabla-\tilde\nabla,
 \qquad \mathcal R:=\Rm(g)-\widetilde{\Rm},
 \qquad \widetilde{\Rm}:=\Rm(\tilde g),
\]
where $\nabla=\nabla^g$ and $\tilde\nabla=\nabla^{\tilde g}$.
Use curvature of type $(1,3)$ and measure tensors with $g(t)$.
Then Ricci curvature is a fixed trace, so
$h_t=-2\operatorname{tr}\mathcal R$ and
$|h_t|\le C|\mathcal R|$. Subtracting the two connection evolution
identities gives, schematically,
\[
 A_t=g^{-1}*\nabla\mathcal R
      +(g^{-1}-\tilde g^{-1})*\tilde\nabla\widetilde{\Rm}
      +g^{-1}*A*\widetilde{\Rm}.
\]
Since $g^{-1}-\tilde g^{-1}=-g^{-1}h\tilde g^{-1}$, this proves
$|A_t|\le C(|h|+|A|+|\nabla\mathcal R|)$.

To treat the difference of curvature Laplacians, define
\[
 V^p:=g^{pq}\nabla_q\widetilde{\Rm}
             -\tilde g^{pq}\tilde\nabla_q\widetilde{\Rm}.
\]
Metric compatibility gives the exact tensor identity
\[
 \begin{split}
 (\Delta_g-\Delta_{\tilde g})\widetilde{\Rm}
 &=\nabla_pV^p
   +(\nabla_p-\tilde\nabla_p)
           (\tilde g^{pq}\tilde\nabla_q\widetilde{\Rm}).
 \end{split}
\]
The second term is the algebraic action of $A$ on all tensor
indices of $\tilde g^{pq}\tilde\nabla_q\widetilde{\Rm}$, including
the index $p$. Also
\[
 V=(g^{-1}-\tilde g^{-1})*\tilde\nabla\widetilde{\Rm}
                           +g^{-1}*A*\widetilde{\Rm},
 \qquad |V|\le C(|h|+|A|).
\]
The curvature evolution is $\partial_t\Rm=\Delta\Rm$
plus algebraic quadratic curvature contractions. Subtracting
these contractions costs at most $C(|h|+|\mathcal R|)$.
Consequently
\[
 |\partial_t\mathcal R-\Delta_g\mathcal R
                         -\operatorname{div}_gV|
 \le C(|h|+|A|+|\mathcal R|).
\]
Choose a smooth radial cutoff $0\le\psi_D\le1$, equal to one on
$\Omega_{D-2}$ and supported in $\Omega_{D-1}$, with
$|\bar\nabla\psi_D|\le C$. Its derivative is supported where
$D-2\le|\vec x-\vec x_0|\le D-1$; these bounds are uniform for
$D\ge4$. For a phase $\varphi$ to be chosen below, set
\[
 Q:=|h|^2+|A|^2+|\mathcal R|^2,\qquad
 w:=\psi_D^2e^{-\varphi},\qquad
 \mathcal E(t):=\int_{\Omega_D}wQ\,\mathrm dV_g.
\]
The derivatives of the tensor norms contribute at most $CQ$,
because $\partial_tg=-2\Ric$ and $|\Ric|\le C$; the volume
derivative is $\partial_t\mathrm dV_g=-\mathrm R_g\mathrm dV_g$
and has the same bound. The $h$ and $A$ equations therefore give
an integrand bounded by $CQ+C|A||\nabla\mathcal R|$ in addition
to the curvature equation and $-\varphi_tQ$.
For the principal curvature terms, integration by parts yields
\[
 \begin{split}
 2\int w\langle\mathcal R,\Delta_g\mathcal R\rangle
 &=-2\int w|\nabla\mathcal R|^2
       -2\int\langle\nabla w\otimes\mathcal R,
                                      \nabla\mathcal R\rangle,\\
 2\int w\langle\mathcal R,\operatorname{div}_gV\rangle
 &=-2\int w\langle\nabla\mathcal R,V\rangle
       -2\int\langle\nabla w\otimes\mathcal R,V\rangle.
 \end{split}
\]
All integrals use $\mathrm dV_g$; the cutoff removes the boundary
terms. Since
$\nabla w=e^{-\varphi}(2\psi_D\nabla\psi_D-\psi_D^2\nabla\varphi)$,
Young's inequality bounds the terms containing
$\nabla\mathcal R$ by an arbitrarily small fixed multiple of
$\int w|\nabla\mathcal R|^2$, plus
\[
 C\mathcal E+C\int w|\nabla\varphi|^2Q
                  +C\int e^{-\varphi}|\nabla\psi_D|^2Q.
\]
The bound $|V|\le C(|h|+|A|)$ controls the last divergence term
by the same expression. Choosing the small multiples with sum at
most one gives
\[
  \mathcal E'\le{}C\mathcal E-\int w|\nabla\mathcal R|^2\,\mathrm dV_g
 +\int w(C_1|\nabla\varphi|^2-\varphi_t)Q\,\mathrm dV_g
    +C\int e^{-\varphi}|\nabla\psi_D|^2Q\,\mathrm dV_g.
 \]
Here $C_1$ is a fixed local energy constant.

Let $\varphi:=a(t)|\vec x-\vec x_0|^2$. Write $C_2$ for a
fixed upper bound of the Euclidean covector norm squared with
respect to $g^{-1}$. Then
$|\nabla\varphi|^2\le4C_2a(t)^2|\vec x-\vec x_0|^2$.
Choose a fixed $a_0>0$ sufficiently small that
$4C_1C_2 a_0(1-\eta)\le1/2$, and solve
\[
 a'=4C_1C_2 a^2,\qquad a(-1)=a_0,
 \qquad a(t)=\frac{a_0}{1-4C_1C_2 a_0(t+1)}.
\]
Thus $a_0\le a(t)\le2a_0$ on the entire Ricci-flow interval,
and $C_1|\nabla\varphi|^2-\varphi_t\le0$.
The assumed coordinate and curvature bounds give
$Q\le C$. The cutoff annulus has volume at most $CD^d$, so
\[
 \int e^{-\varphi}|\nabla\psi_D|^2Q\,\mathrm dV_g
       \le CD^d e^{-a_0(D-2)^2}.
\]
At the original initial time, the formulas for the difference
of inverse metrics, connections and curvatures imply
$|h|+|A|+|\mathcal R|\le Ce_{\rm in}$. Uniform metric
equivalence and compactness of $N$ give
\[
 \mathcal E(-1)\le Ce_{\rm in}^2
 \int_{N\times\R^d}e^{-a_0|\vec x-\vec x_0|^2}\,\mathrm dV_{\bar g}
 =Ce_{\rm in}^2\Vol(N,g_N)(\pi/a_0)^{d/2}\le Ce_{\rm in}^2.
\]
Gronwall on $[-1,-\eta]$ yields
\[
 \sup_t\mathcal E(t)
 \le C\{e_{\rm in}^2+D^d e^{-a_0(D-2)^2}\}.
\]
On $\Omega_2$ the cutoff is one and the weight is bounded below
by $e^{-8a_0}$. Metric equivalence therefore gives the ordinary,
unweighted estimate
\[
 \sup_t\|h(t)\|_{L^2(\Omega_2;\mathrm dV_{\bar g})}
 \le C\{e_{\rm in}+D^{d/2}e^{-a_0(D-2)^2/2}\}.
\]

Interpolation on fixed product regions, using finitely many
fixed coordinate charts on $N$, gives
\[
 \|h\|_{C^{12}(\Omega_1)}
 \le C\|h\|_{L^2(\Omega_2)}^{3/4}
             \|h\|_{C^{100n-2}(\Omega_2)}^{1/4}
       +C\|h\|_{L^2(\Omega_2)}.
\]
Indeed, $100n-2>48+3n/2$, so $3/4$ is strictly below the
interpolation endpoint. Mollification in the fixed charts bounds
the low-frequency term by $r^{-12-n/2}\|h\|_2$ and the
remainder by $r^{100n-2-12}\|h\|_{C^{100n-2}}$;
optimizing $r$ proves the estimate. The same charts and constants
work after every Euclidean translation. Since the stronger norm
is uniformly bounded and $e_{\rm in}\le1$, the initial-error
contribution is at most $Ce_{\rm in}^{1/2}$. The interpolated tail
is bounded by $CD^{3d/8}e^{-3a_0(D-2)^2/8}$, and the remaining
$L^2$ term has factor $D^{d/2}$ and a faster Gaussian rate.
Since $D\ge4$ and $(D-2)^2\ge D^2/4$, both are bounded by
$CD^{d/2}e^{-c_*D^2}$ with $c_*=3a_0/32>0$.
This proves~\eqref{eq:A-7}.

\end{proof}

\begin{proposition}[Quantitative metric comparison]
\label{prop:metric-comparison}
Fix $\rho\in(0,1)$ and $\eta,\epsilon_{\rm met}$ as in
Lemma~\ref{lem:qualitative-propagation}. There are a parameter
neighborhood, $\bar\epsilon>0$, $S_*$, fixed collars and constants
$C,c_4>0$, all depending only on $(\eta,\epsilon_{\rm met},\rho)$,
such that for every parameter in this neighborhood and every initial
comparison chart with $S\ge S_*$ and $\boldsymbol{\epsilon}(S)\le\bar\epsilon$,
the qualitative estimates~\eqref{eq:A-2}--\eqref{eq:B-4}, the
regularity bounds above, and the following estimate hold:
\begin{equation}\label{eq:A-3}
 \sup_{-1\le t\le-\eta}
 \|\iota^*g(t)-\hat g_b(t)\|_{C^{12}(N_{\rho_3S})}
       \le C\{\boldsymbol{\epsilon}(S)^{1/2}+S^{d/2}e^{-c_4S^2}\}.
\end{equation}
Here $c_4>0$ may be any sufficiently small Gaussian boundary-decay rate
furnished by the proof, and may subsequently be decreased.
\end{proposition}
\begin{proof}
Apply Lemma~\ref{lem:local-metric-comparison} to
$\iota^*g$ and $\hat g_b$ over the whole interval $[-1,-\eta]$,
using the regularity just established. For each center
$|\vec x_0|\le\rho_3S$, choose
\[
 D=\frac{1-\rho_3}{2}S.
\]
If $C$ is the fixed regularity collar above, increase $S_*$ so that
$D\ge4$ and $D\ge C+5$. Then
\[
 N\times B_{D+4}(\vec x_0)\subset N_{S-C},
 \qquad |\vec x_0|\le\rho_3S.
\]
Thus the comparison region, cutoff annulus and regularity collars
lie in the controlled cylinder. On each region
$e_{\rm in}\le\boldsymbol{\epsilon}(S)$, by the incoming $C^{100n}(N_S)$ bound.
The initial energy in the lemma is at most $C\boldsymbol{\epsilon}(S)^2$,
independently of $S$, because its Gaussian is centered at $\vec x_0$.
The lemma therefore gives, uniformly in the center,
\[
 \sup_t\|\iota^*g(t)-\hat g_b(t)\|_{C^{12}
                   (N\times B_1(\vec x_0))}
 \le C\{\boldsymbol{\epsilon}(S)^{1/2}+S^{d/2}e^{-c_*'S^2}\},
 \qquad c_*':=\frac{c_*(1-\rho_3)^2}{8}>0.
\]
The fixed compact atlas and Euclidean translations give uniform
constants. Taking the supremum over the centers and choosing
$0<c_4\le c_*'$ proves~\eqref{eq:A-3}.
\end{proof}

\subsection{Terminal probability estimates}
\label{subsec:terminal-comparison}
\label{subsec:coverage}
The following first-exit argument converts comparison domains into
ambient balls, so the terminal probability estimates apply to their complements.

\begin{lemma}\label{lem:coverage}
Suppose $\iota$ embeds a neighborhood of $\overline N_R$, $R>0$,
into a complete $n$-manifold and
$(1-\delta)\bar g\le\iota^*g\le(1+\delta)\bar g$ with $0\le\delta\le1/2$.
Fix $y_0\in N$, and set $p_0:=\iota(y_0,0)$ and $\Omega_R:=\iota(N_R)$.
Then
\begin{equation}\label{eq:coverage}
 B_g(p_0,\sqrt{1-\delta}\,R)\subset\Omega_R.
\end{equation}
If $R\ge1$, $\mu$ is a probability with
$\mu(\iota(N_1))\ge q>0$, and
$\int_M d_g(z_H,x)^2\,\mathrm d\mu(x)\le H_n$, then
$d_g(z_H,p_0)\le C_1+\sqrt{H_n/q}=:C_0$ and
\begin{equation}\label{eq:center-coverage}
 B_g(z_H,\sqrt{1-\delta}\,R-C_0)\subset\Omega_R.
\end{equation}
Here $C_1$ depends only on $(N,g_N)$; the last inclusion is used
when its radius is positive.
\end{lemma}
\begin{proof}
The image is open and its boundary is the image of
$N\times\partial B_R$. Every path from $p_0$ to the complement has a
first boundary intersection. Before that intersection its Euclidean
radial displacement is $R$, so its length is at least
$\sqrt{1-\delta}\,R$. This proves~\eqref{eq:coverage}.
The upper metric bound gives core diameter at most $C_1$. Thus
\[
 H_n\ge\int_{\iota(N_1)}d_g(z_H,x)^2\,\mathrm d\mu(x)
       \ge q\,(d_g(z_H,p_0)-C_1)_+^2,
\]
which proves the center bound and~\eqref{eq:center-coverage}.
\end{proof}

For every prescribed Gaussian loss $\eta>0$ and core-mass bound
$q_0>0$, \eqref{eq:terminal-tail} and Lemma~\ref{lem:coverage}
give a metric threshold $\bar\delta(\eta)>0$ such that, for an actual
normalized terminal probability with $\delta\le\bar\delta(\eta)$
and $\mu(\iota(N_1))\ge q_0$, the mass and entropy outside
$\Omega_R$ are at most
$C_{\eta,q_0}e^{-(1-\eta)R^2/4}$. To absorb the center displacement use
$(aR-C_0)^2\ge(1-\epsilon)a^2R^2-C_0^2/\epsilon$ and begin with a
smaller Gaussian loss. The reference product has the same type of tails.

The next proposition gives the family-uniform forms of
\cite[Proposition~6.3 and Corollary~6.6]{FLC}. It uses the qualitatively
propagated cylinder and its intrinsic regularity.

\begin{proposition}[Terminal probability comparison]
\label{prop:terminal-comparison}\label{thm:terminal-kernel}
Prescribe $(\rho,c_1,Q_*)$ as in
Subsection~\ref{subsec:physical-common-setup}. There is
$\eta=\eta(\rho,c_1,Q_*)\in(0,1/3)$, and then a permitted
metric accuracy
$\epsilon_{\rm met}=\epsilon_{\rm met}(\rho,c_1,Q_*,\eta)>0$,
for which the following holds. There are a parameter neighborhood,
$\bar\epsilon>0$ and $S_*$, depending only on these comparison inputs,
such that the conclusions hold whenever $b$ belongs to this
neighborhood, $\boldsymbol{\epsilon}(S)\le\bar\epsilon$ and $S\ge S_*$.
The constants below have the same dependence, together with any
order or fixed core radius displayed in their subscripts; they are
uniform in $b,S$ and the duration of a partial block.

There are $R_c,q_c>0$ such that
\begin{equation}\label{eq:B-1}
 \mu_t(\iota(N_{R_c}))\ge q_c
       \qquad(-1\le t\le-\eta).
\end{equation}
For $R_c\le R\le\rho_2S$ and $-1\le t\le-\eta$,
\begin{equation}\label{eq:B-6}
 \int_{M\setminus\iota(N_R)}
  (1+|f(t)|+|t||\mathrm R|+|t||\nabla f(t)|^2)\,\mathrm d\mu_t
 \le C e^{-7R^2/(32|t|)}.
\end{equation}
Writing $r=|\vec x|$, the actual terminal potential satisfies on
$\iota(N_{\rho_2S})$
\begin{equation}\label{eq:B-7}
 (1-c_3)\frac{r^2}{4}-C
       \le F(\iota(y,\vec x),t)
       \le(1+c_3)\frac{r^2}{4}+C
       \qquad(-1\le t\le-1/3),
\end{equation}
and
\begin{equation}\label{eq:B-8}
 |(\nabla^g)^jF|+|(\nabla^g)^jf(t)|
       \le C_j(1+r)^{Q_j}e^{c_3 r^2}
       \quad(0\le j\le Q_*+4,\ -1\le t\le-1/2).
\end{equation}
The reference-coordinate version of \eqref{eq:B-8} holds through
$j\le 100n-2$, the fixed coordinate order in
Proposition~\ref{prop:metric-comparison}.
If a finite integer $l\ge0$ and a loss $\delta>0$ are also prescribed
before choosing $\eta$, the same conclusions hold together with
\begin{equation}\label{eq:terminal-arbitrary-rough}
 \begin{aligned}
 \sup_{-1\le t\le-1/2}\sum_{j=0}^{l}
       |(\nabla^{g(t)})^j\Rm(g(t))|&\le C_l,\\
 \sup_{-1\le t\le-1/2}\sum_{j=0}^{l}
       |(\nabla^{g(t)})^jf(t)|
       &\le C_{l,\delta}(1+r)^{Q_l}e^{\delta r^2}
       \quad\text{on }\iota(N_{\rho_2S}).
 \end{aligned}
\end{equation}
For this additional estimate, the choices of $\eta,\epsilon_{\rm met}$,
the parameter neighborhood and thresholds also depend on $(l,\delta)$.
After the comparison data and derivative order are fixed, the constants
are independent of $r,S,b$.
For each fixed $R_0$ and every finite $k$, on
$\iota(N_{R_0})\subset\iota(N_{\rho_2S})$ one also has
\begin{equation}\label{eq:fixed-core-rough}
 \sup_{-1\le t\le-1/2}\sum_{j=0}^k
  (|(\nabla^{g(t)})^j\Rm|+|(\nabla^{g(t)})^jf(t)|)
       \le C_{R_0,k}.
\end{equation}
Here $S_*$ is enlarged to include the fixed core if necessary.
All local estimates above use the outer
collar $N_{\rho_3S}$: one fixed $C_T$ suffices, and $S_*$ is
chosen so that $(\rho_3-\rho_2)S\ge C_T$ and $R_c<\rho_2S$.
\end{proposition}
\begin{proof}
For \eqref{eq:B-8}, choose an internal loss $\delta'\le c_3$ with
$C(Q_*+4)\delta'\le c_3$.
If \eqref{eq:terminal-arbitrary-rough} is also required, first fix
$l,\delta$ and make the stronger choice
\[
 C(\max\{l,Q_*+4\})\delta'
                         \le\min\{c_3,\delta\}.
\]
The constant bounds the products arising in logarithmic
differentiation through the prescribed order.
Choose $0<\eta<1/6$ so that
$(1-3\eta)^{-1}<1+\delta'/8$, and then choose $\epsilon_{\rm met}>0$
so small that
\[
 \frac{1+\epsilon_{\rm met}}{1-3\eta}<1+\frac{\delta'}{4},
 \qquad \epsilon_{\rm met}<\frac18\min\{\delta',1/8\}.
\]
Choose the kernel and bounded-center losses smaller still. Apply
Lemma~\ref{lem:qualitative-propagation},
with output accuracy sufficient for \eqref{eq:B-4} with this
$\epsilon_{\rm met}$, and use its intrinsic bounds \eqref{eq:A-4}.
All reference-flow bounds, collars and constants are allowed to
depend on this fixed $\eta$.

The initial pair condition gives
$\mu_{-1}(\iota(N_{1/2}))\ge q>0$ uniformly in $b$. Let $v$ solve
the forward heat equation on the closed flow from a cutoff equal
to one on this core and supported in $\iota(N_1)$. Then
$0\le v\le1$ and $\int v(t)\,\mathrm d\mu_t\ge q$ by duality.
To localize this mass, choose a radial function $a_R$ with uniformly
bounded first two derivatives, zero on $N_1$ and equal to $R$
outside $N_{R+3}$. Extend $e^{-a_R}$ as the constant $e^{-R}$
beyond the chart. The controlled metric and connection give
$|\Delta_g e^{-a_R}|\le C e^{-a_R}$, so the maximum principle gives
$v\le e^{C(t+1)}e^{-a_R}$.
Choose a fixed sufficiently large $R$, and then $R_c=R+3$.
Since $\mu_t$ has mass one, duality yields \eqref{eq:B-1}, after
increasing $S_*$ to contain the cutoff collar. The second-moment
bound for any terminal $H_n$-center $z_H(t)$ then gives
$d_{g(t)}(z_H(t),p)\le C_H$, where $p=\iota(y_0,0)$.
Indeed, the second moment is at most $H_n|t|$, and
the fixed core has uniformly bounded diameter.

The Euclidean factor of $\hat g_b$ is exactly $g_E$. By
\eqref{eq:B-4}, a path leaving $\iota(N_R)$ first crosses
Euclidean distance $R$. Consequently, for the radii in the
statement and points of the retained cylinder,
\begin{equation}\label{eq:B-5}
 d_{g(t)}(p,\iota(y,\vec x))\ge\sqrt{1-\epsilon_{\rm met}}\,r,
 \qquad B_{g(t)}(p,\sqrt{1-\epsilon_{\rm met}}R)
                         \subset\iota(N_R).
\end{equation}
The upper distance bound is $C+\sqrt{1+\epsilon_{\rm met}}\,r$.
Apply the terminal estimates stated with~\eqref{eq:terminal-tail}
about $z_H(t)$, scaled to Ricci-flow time, and use~\eqref{eq:B-5}. Absorbing $C_H$ with an arbitrarily small
additional Gaussian loss proves \eqref{eq:B-6} and the lower
bound in \eqref{eq:B-7}, first with the finer loss $\delta'$.
The constants in this absorption may depend on $\eta$, because
$|t|\ge\eta$; they are independent of the chart radius.

For the opposite bound, fix $-1\le t\le-1/3$, reproduce the terminal
probability from $t_a=-\eta$ and restrict to the core in \eqref{eq:B-1}. If
$u_t=(4\pi|t|)^{-n/2}e^{-f(t)}$, then
\[
 u_t(x)=\int_M K(w,t_a;x,t)\,\mathrm d\mu_{t_a}(w),
 \qquad
 K(w,t_a;x,t)\ge(4\pi\Delta)^{-n/2}
                   e^{-\ell_{w,t_a}(x,t)},
 \qquad \Delta=|t|-\eta.
\]
The kernel lower bound is \cite[Corollary~9.5]{Per02}.
Write $w=\iota(y_w,\vec x_w)$ and $x=\iota(y,\vec x)$, so
$\vec x_w$ is the Euclidean coordinate of $\iota^{-1}(w)$ and
$|\vec x_w|\le R_c$. Let $\sigma:[0,1]\to N$ be a constant-speed
minimizing $g_N$-geodesic from $y_w$ to $y$. In backward time
$0\le v\le\Delta$, use the curve
\[
 \zeta(v)=\iota\bigl(\sigma(v/\Delta),\vec z(v)\bigr),\qquad
 \vec z(v)=\vec x_w+\sqrt{v/\Delta}\,(\vec x-\vec x_w),
\]
at time $-\eta-v$. Its Euclidean component lies on the segment
joining $\vec x_w$ to $\vec x$, so the curve stays in the retained
cylinder. The Euclidean contribution to its reduced action is
\[
 \frac1{2\sqrt\Delta}\int_0^\Delta
       \sqrt v\,|\dot{\vec z}(v)|^2\,\mathrm dv
       =\frac{|\vec x-\vec x_w|^2}{4\Delta}.
\]
The compact contribution is bounded by $C/\Delta$, using the
uniform bound for $\gamma_b(t)$ and
$|\dot\sigma|\le\operatorname{diam}(N,g_N)$; the curvature
contribution is bounded by $C\Delta$. Both are bounded uniformly
because $1/6\le\Delta\le1$. Smooth approximation at $v=0$
justifies the square-root parametrization in the reduced-action
infimum. The metric comparison~\eqref{eq:B-4} therefore gives
\[
 \ell_{w,-\eta}(x,t)
 \le\frac{(1+\epsilon_{\rm met})|\vec x-\vec x_w|^2}{4\Delta}+C.
\]
Thus
\[
 u_t(x)\ge c\exp\left(-\frac{(1+\epsilon_{\rm met})(r+R_c)^2}
                                      {4(|t|-\eta)}\right).
\]
For $|t|\ge1/3$, the choices of $\eta$ and $\epsilon_{\rm met}$ control
$|t|/(|t|-\eta)$; a further small Gaussian loss absorbs the
fixed cross term. This gives the upper bound with loss
$\delta'$ and hence proves \eqref{eq:B-7} as stated.
The same argument gives these bounds on the fixed neighborhood
$N_{\rho_2S+C_T/2}$, with $C_T$ large enough for the local balls
used next and $(\rho_3-\rho_2)S\ge C_T$.

For derivatives at a main-block point $(x,t)$, put
$\varrho=c/(1+r)$ and use the later-time cylinder
\[
 Q_{2\varrho}=B_{g(t)}(x,2\varrho)\times[t,t+4\varrho^2],
\]
with the fixed spatial identification. Choose $c>0$ small enough
that this cylinder and its fixed larger collars lie in the controlled
region with time interval $[-1,-1/3]$; increasing $S_*$ puts their
spatial collars inside $N_{\rho_3S}$. Reading time backward, the
terminal density solves a uniformly parabolic equation. The local estimates used in
\cite[Corollary~6.6]{FLC}, with the intrinsic bounds
\eqref{eq:A-4} and a uniform local volume lower bound, give
\begin{equation}\label{eq:B-13}
 |(\nabla^g)^ju_t(x,t)|
       \le C_j\varrho^{-j}\sup_{Q_{2\varrho}}u.
\end{equation}
Here one fixed geometric cylinder supports every prescribed
intrinsic order; increasing the order changes only the constants.
Radius and time vary on this cylinder by $O(\varrho)$ and
$O(\varrho^2)$, respectively. Since $r\varrho$ and
$r^2\varrho^2$ are bounded and $|t|\ge1/3$, the quadratic bounds
with loss $\delta'$ imply
\[
 \frac{\sup_{Q_{2\varrho}}u}{u_t(x)}
                \le C e^{C\delta' r^2}.
\]
Each derivative of $\log u$ is a sum of products of at most $j$
factors $(\nabla^g)^a u/u$, of total derivative order $j$.
Thus \eqref{eq:B-13} gives a polynomial in $1+r$ times
$e^{C(Q_*+4)\delta' r^2}$ through order $Q_*+4$.
The internal choice of $\delta'$ proves \eqref{eq:B-8}; the
order-zero estimate follows directly from \eqref{eq:B-7}.
For~\eqref{eq:terminal-arbitrary-rough}, the curvature bound follows from
\eqref{eq:A-4} on the reserved collars, using the earlier interval
inside $[-2,0)$ to include $t=-1$. Logarithmic differentiation through
order $l$ gives the potential bound with the prescribed Gaussian loss. Conversion to reference derivatives
uses the metric and connection bounds through order $100n-2$
in the fixed comparison chart.

On each fixed core the kernel quotient is bounded, giving every
finite intrinsic derivative of $f(t)$. Curvature uses the earlier
Ricci-flow interval in~\eqref{eq:A-4}, while the terminal potential
uses the later interval constructed above. These intervals support
all prescribed orders and prove~\eqref{eq:fixed-core-rough}.
Metric equivalence, bounded curvature and the local volume lower
bound give uniform injectivity and interpolation constants there.
\end{proof}

\subsection{Auxiliary forward comparisons}\label{subsec:auxiliary-forward}
The auxiliary product comparison is constructed in three stages.
\begin{enumerate}[label=\textup{(\roman*)},leftmargin=*]
\item \emph{Compact Ricci flow and forward scalar.}
On the fixed auxiliary interval $[-1,-\eta]$, define
\begin{equation}\label{eq:compact-physical-forward}
 \begin{gathered}
 \partial_t\gamma_b=-2\Ric(\gamma_b),\qquad \gamma_b(-1)=g_b,\\
 (\partial_t-\Delta_{\gamma_b})\hat f_b^N=-m/2,\qquad \hat f_b^N(-1)=f_b,
 \qquad -1\le t\le-\eta,\\
 \hat g_b(t)=\gamma_b(t)+g_E,\qquad
 \hat f_b(y,\vec x,t):=\hat f_b^N(y,t)+|\vec x|^2/4.
 \end{gathered}
\end{equation}
The product scalar satisfies $\hat f_b(-1)=F_b$ and
$(\partial_t-\Delta_{\hat g_b})\hat f_b=-n/2$. We compare it below
with the forward solution $\tilde F$ and then the scaled terminal
potential $F(t)=|t| f(t)$.
\item \emph{Parabolic rescaling and compact pullback.}
With $t=-e^{-s}$, set
\begin{equation}\label{eq:compact-rescaling}
 \check g_b(s):=|t|^{-1}\gamma_b(t),\qquad
 \check f_b(s):=|t|^{-1}\hat f_b^N(t),
 \qquad
 \dot\psi_b=-\nabla_{\check g_b}\check f_b\circ\psi_b,\quad
 \psi_b(0)=\Id .
\end{equation}
The map $\psi_b$ removes the compact gradient drift from the
rescaled comparison.
\item \emph{Probability normalization.}
Define the compact mass and normalized product pair by
\begin{equation}\label{eq:compact-forward-definition}
 \begin{split}
 \mathcal M_b(s)&:=(4\pi)^{-m/2}
       \int_Ne^{-\check f_b(s)}\,\mathrm dV_{\check g_b(s)},\\
 Z_b(s)&:=\bigl(\psi_b^*\check g_b+g_E,\,
       \psi_b^*\check f_b+\log\mathcal M_b+|\vec x|^2/4\bigr).
 \end{split}
\end{equation}
\end{enumerate}
Since the compact input has mass one, $\mathcal M_b(0)=1$ and
$Z_b(0)=X_b$. The weighted measure of $Z_b(s)$ is a probability:
its compact factor is $\mathcal M_b(s)^{-1}$ times the pullback of
the measure defining $\mathcal M_b$, and its Euclidean factor is
the standard Gaussian.

\begin{proposition}\label{prop:compact-forward}
For each fixed $\eta\in(0,1/3)$ and finite integer $j\ge0$, there
are a parameter neighborhood $B=B(\eta,j)$ and a constant
$C_j=C_j(\eta)$ such that the constructions above exist for
$b\in B$ and $0\le s\le\log(1/\eta)$, and
\begin{equation}\label{eq:G-3-1}
 \sup_{0\le s\le\log(1/\eta)}\|Z_b(s)-X_b\|_{C^j(N)}
       \le C_j|\mathrm De(b)|.
\end{equation}
The norm denotes the two compact components; the Euclidean
components coincide. The maps $\psi_b,\psi_b^{-1}$ have uniform
$C^{j+1}$ bounds and $|\log\mathcal M_b|\le C_j$ on this
interval.
\end{proposition}
\begin{proof}
The reference compact flow $\gamma_0(t)=|t|g_N$ is smooth on
$[-1,-\eta]$. By compact Ricci-flow continuous dependence and the
linear scalar equation~\eqref{eq:compact-physical-forward}, the nearby
solutions $(\gamma_b,\hat f_b^N)$ exist on this interval with uniform
bounds at any prescribed finite order, after restricting $b$.
The rescaling, the ODE in~\eqref{eq:compact-rescaling}, and
\eqref{eq:compact-forward-definition} then give the asserted bounds
for $\psi_b$, $\psi_b^{-1}$, and $\log\mathcal M_b$.

On $N$ set
\[
 T=|t|\Ric_{\gamma_b}+\nabla_{\gamma_b}^2\hat f_b^N-\frac12\gamma_b.
\]
Since $(\partial_t-\Delta_{\gamma_b})\hat f_b^N=-m/2$, the Ricci
and Hessian evolution equations give $(\partial_t-\Delta_L)T=0$,
where $\Delta_L$ is the Lichnerowicz Laplacian of $\gamma_b(t)$
\cite[proof of Lemma~6.11]{FLC}. At $t=-1$ we have
$T=-\Phi(X_b)|_N$. By Lemma~\ref{lem:compact-family} its norm in every
prescribed finite order is at most $C|\mathrm De(b)|$. The compact
energy estimate for this linear equation at an integer order
greater than $j+1+m/2$, followed by Sobolev embedding, gives
\[
 \sup_{-1\le t\le-\eta}\|T(t)\|_{C^{j+1}(N)}
 \le C_j|\mathrm De(b)|.
\]

On $N$, write $g=\psi_b^*\check g_b$, $f^0=\psi_b^*\check f_b$,
$c=\log\mathcal M_b$, and $f=f^0+c$, using dimension $m$ in $w$. Then
$\Phi(g,f)=-|t|^{-1}\psi_b^*T$, so the same estimate holds for
$\Phi(g,f)$ in $C^{j+1}(N)$. Differentiating the rescaling and pullback gives
\[
 \partial_sg=2\Phi(g,f),\qquad
 \partial_sf^0=\Delta f^0-|\mathrm df^0|^2+f^0-m/2.
\]
Thus the logarithmic derivative of
$\mathrm d\mu^0=(4\pi)^{-m/2}e^{-f^0}\mathrm dV_g$ is $-w(g,f^0)$.
Since $\mathrm d\mu^N=e^{-c}\mathrm d\mu^0$ has mass one and
$w(g,f)=w(g,f^0)+c$, we obtain
\begin{equation}\label{eq:compact-mass-evolution}
 c'=c-\int_Nw(g,f)\,\mathrm d\mu^N .
\end{equation}
Consequently the two compact components of $Z_b(s)$ satisfy
\begin{equation}\label{eq:G-3-2}
 \partial_sg=2\Phi(g,f),\qquad
 \partial_sf=\Tr_g\Phi(g,f)+w(g,f)-\int_Nw(g,f)\,\mathrm d\mu^N .
\end{equation}
The contracted Bianchi identity gives $\mathrm dw=-2\delta_{g,f}\Phi$.
The uniform diameter bound therefore controls
$\|w-\int_N w\,\mathrm d\mu^N\|_{C^0}$ by $C\|\mathrm dw\|_{C^0}$,
and differentiating the identity bounds the higher derivatives through
order $j$ by $C_j\|\Phi\|_{C^{j+1}}$.
Both velocities in~\eqref{eq:G-3-2} are thus bounded in $C^j$ by
$C_j|\mathrm De(b)|$. Integrating from $Z_b(0)=X_b$ proves
\eqref{eq:G-3-1}.
\end{proof}

\phantomsection\label{idx:caloric-data}
On the actual closed flow solve the global forward equation
\begin{equation}\label{eq:C-4}
 (\partial_t-\Delta_g)\tilde F=-n/2,
 \qquad \tilde F(-1)=F(-1).
\end{equation}
The solution is smooth on every closed negative-time interval.
In comparisons on the product we suppress the pullback by $\iota$.

\begin{proposition}
\label{prop:scalar-comparison}
Fix $(\rho,c_1,Q_*)$ and use the time and metric accuracy
chosen in Proposition~\ref{prop:terminal-comparison}. There are a
common parameter neighborhood and constants
$\bar\epsilon,S_*,C,Q,c_4>0$, depending only on these comparison
inputs, such that every initial comparison chart with parameter in
that neighborhood, $\boldsymbol{\epsilon}(S)\le\bar\epsilon$ and $S\ge S_*$
satisfies
\begin{equation}\label{eq:C-9}
 \sup_{-1\le t\le-1/2}
 \|\tilde F(t)-\hat f_b(t)\|_{C^{9}(N_{\rho_1S})}
 \le CS^Q\bigl(\boldsymbol{\epsilon}(S)^{1/2}+e^{-c_4S^2}\bigr).
\end{equation}
The positive rate $c_4$ is furnished by this proof and may be
decreased to agree with the preceding boundary rates.
\end{proposition}
\begin{proof}
For a regular terminal point at time $t_i$ with kernel potential $f_i$,
Perelman's differential Harnack inequality~\cite[Corollary~9.3]{Per02} gives
$(\partial_t-\Delta_g)((t_i-t)f_i)\ge-n/2$. As the terminal points tend to
$z_*$, smooth kernel convergence on negative-time slices~\cite[Theorem~5.20]{FLSt}
yields $(\partial_t-\Delta_g)F\ge-n/2$. The global lower bound
$F(-1)\ge-C$, which follows from the entropy lower bound and the
terminal heat-kernel upper bound, and the maximum principle on the
closed flow imply
\[
 -C-n(t+1)/2\le\tilde F(t)\le F(t).
\]
The quadratic upper bound for $F$ in the preceding subsection
therefore proves, on $N_{\rho_2S}\times[-1,-1/2]$,
\begin{equation}\label{eq:C-5}
 |\tilde F|\le C(1+r^2).
\end{equation}
The compact forward solution and its derivatives are uniformly
bounded on the prescribed time interval, so $\hat f_b=\hat f_b^N+r^2/4$
has the same quadratic bound.

Set $v:=\tilde F-\hat f_b$ and write $\hat\nabla$ for the connection of
$\hat g_b(t)$. Both scalars have source $-n/2$, and hence
\begin{equation}\label{eq:C-6}
 \begin{split}
 (\partial_t-\Delta_g)v&=P,\\
 P&=(g^{ij}-\hat g_b^{ij})(\hat\nabla^2\hat f_b)_{ij}
       -g^{ij}(\Gamma_g-\Gamma_{\hat g_b})^a{}_{ij}
                                                   \partial_a\hat f_b.
 \end{split}
\end{equation}
Put
$E_{\rm met}:=\boldsymbol{\epsilon}(S)^{1/2}+S^Q e^{-c_*'S^2}$,
enlarging the fixed polynomial exponent when necessary.
Here $c_*'=c_*(1-\rho_3)^2/8$ is the Gaussian
boundary-decay rate obtained in the proof of
Proposition~\ref{prop:metric-comparison}. The comparison at order $12$ gives
\[
 \sup_t\|P(t)\|_{C^{10,1/2}(N_{\rho_2S})}
       \le C(1+S)^Q E_{\rm met}.
\]
Here one derivative is used by the connection difference and one
controls the H\"older seminorm on fixed reference charts. In
particular this bounds the source in $C^0$. The initial
$C^{9,1/2}$ norm of $v$ is at most $C\boldsymbol{\epsilon}(S)$:
the incoming $C^{100n}$ error controls it on unit charts, with
constants independent of $S$. By~\eqref{eq:C-5} and the quadratic bound for
$\hat f_b$, we have $|v|\le C(1+S)^2$ on $N_{\rho_2S}$.

Fix $\rho_1<\rho'<\rho_2$. For $|\vec x_0|\le\rho'S$, take
\[
 D=\frac{\rho_2-\rho'}2S.
\]
The region $N\times B_D(\vec x_0)$ and every required fixed collar
lie in $N_{\rho_2S}$ for sufficiently large $S$. Metric equivalence
and the first coordinate derivative bounds give, uniformly in the center,
\[
 |\mathrm d|\vec x-\vec x_0|^2|_g^2\le C_0|\vec x-\vec x_0|^2,
 \qquad
 \Delta_g|\vec x-\vec x_0|^2
       \le C_0(1+|\vec x-\vec x_0|^2).
\]
Choose $C\ge2C_0+1$ and then $a>0$ with $ae^{C/2}\le1$.
Direct differentiation shows that
\[
 \begin{split}
 &e^{aD^2-ae^{C(t+1)}|\vec x-\vec x_0|^2-C(t+1)}
 (\partial_t-\Delta_g)
 e^{-aD^2+ae^{C(t+1)}|\vec x-\vec x_0|^2+C(t+1)}\\
 &\quad=C+Ca e^{C(t+1)}|\vec x-\vec x_0|^2
       -ae^{C(t+1)}\Delta_g|\vec x-\vec x_0|^2
       -a^2e^{2C(t+1)}|\mathrm d|\vec x-\vec x_0|^2|_g^2\\
 &\quad\ge C-C_0+(C-2C_0)ae^{C(t+1)}|\vec x-\vec x_0|^2
       \ge0.
 \end{split}
\]
This exponential is at least one when $|\vec x-\vec x_0|=D$,
and at most $e^{1+C/2}e^{-aD^2}$ when
$|\vec x-\vec x_0|\le1$. Comparing both signs of $v$ with the
initial, source and boundary bounds gives
\[
 \begin{split}
 |v(y,\vec x,t)|\le{}&C\boldsymbol{\epsilon}(S)
       +C(t+1)(1+S)^QE_{\rm met}\\
       &+C(1+S)^2
       e^{-aD^2+ae^{C(t+1)}|\vec x-\vec x_0|^2+C(t+1)}.
 \end{split}
\]
The scalar maximum principle applies on
$N\times B_D(\vec x_0)\times[-1,t]$ for every
$-1\le t\le-1/2$. Choose the output rate $c_4$ so that
\[
 0<c_4<\min\left\{c_*',
                  \frac{a(\rho_2-\rho')^2}{4}\right\}.
\]
Evaluating on the inner unit regions yields
\[
 \sup_{-1\le t\le-1/2}\|v(t)\|_{C^0(N_{\rho'S})}
       \le CS^Q(\boldsymbol{\epsilon}(S)^{1/2}+e^{-c_4S^2}).
\]

For fixed nested coordinate balls $B_{r_*}(p)\Subset B_{2r_*}(p)$
and $Q_2=B_{2r_*}(p)\times[-1,-1/2]$, the spatial-interior,
initial-time Schauder estimate is
\[
 \sup_{-1\le t\le-1/2}\|v(t)\|_{C^{9,1/2}(B_{r_*}(p))}
 \le C\left(\|v\|_{C^0(Q_2)}
       +\|v(-1)\|_{C^{9,1/2}(B_{2r_*}(p))}
       +\|P\|_{C^{15/2,15/4}(Q_2)}\right).
\]
The spatial-interior Schauder estimate, including the initial trace,
gives this inequality with $C$ depending on the coefficient bounds
and the fixed time interval \cite[Chapter~3, Sections~2 and~5]{Fri64}.

For $h=g-\hat g_b$, the Ricci-flow equations give
$\partial_th=-2(\Ric(g)-\Ric(\hat g_b))$. In fixed coordinates,
the right side is a sum of uniformly bounded coefficients times
$h,\bar\nabla h,\bar\nabla^2h$. Repeated time differentiation
and the two Ricci-flow equations show that a mixed derivative
$\partial_t^a\bar\nabla^j h$ uses spatial derivatives through
$j+2a$, with every difference term containing a derivative of $h$.
The same statement for $P$ costs one further spatial derivative;
the compact factors in $\hat f_b$ have uniform mixed bounds.
Consequently the quantitative $C^{12}$ bound for $h$
controls every mixed derivative required in
$C^{15/2,15/4}$, as well as one additional spatial
and time derivative for the relevant H\"older seminorms, by
$C(1+S)^Q E_{\rm met}$. The finite coordinate bounds supply the coefficient factors,
and the initial comparison controls the trace norm. Applying the
estimate on fixed balls inside $N_{\rho'S}$ and taking derivative
suprema gives~\eqref{eq:C-9}. The fixed atlas and Euclidean
translations give uniform constants.

\end{proof}

\subsection{Entropy drop and local derivatives}
Put
\[
 u:=\tilde F-F,\qquad
 T:=|t|\Ric+\nabla^2F-g/2,\qquad
 \tilde T:=T+\nabla^2u.
\]
The following identity is the caloric comparison of
\cite[Proposition~6.7]{FLC}, applied to each partial block.

\begin{lemma}\label{thm:caloric}
Let $t=-e^{-s}$, $0\le s\le\log2$, and set
\[
 J(t):=\int|\mathrm du|^2\,\mathrm d\mu_t,\qquad
 I(t):=\int_{-1}^t\int|T|^2\,\mathrm d\mu_r\,\mathrm dr,\qquad
 \tilde I(t):=\int_{-1}^t\int|\tilde T|^2\,\mathrm d\mu_r\,\mathrm dr.
\]
Then
\begin{equation}\label{eq:G-2-1}
 \tilde I(t)+J(t)/2=I(t),\qquad
 J(t)\le A(s)^2,\qquad \tilde I(t)\le A(s)^2/2,
\end{equation}
and, with
$\Delta_Lh=\Delta h+2\Rm(h)-\Ric\circ h-h\circ\Ric$,
\begin{equation}\label{eq:G-2-3}
 (\partial_t-\Delta_L)\tilde T=0.
\end{equation}
\end{lemma}
\begin{proof}
The calculation in \cite[proof of Proposition~6.7]{FLC}, starting
at $t=-1$, gives
\[
 J'(t)=2\int_M\bigl(|T|^2-|\tilde T|^2\bigr)\,\mathrm d\mu_t.
\]
Since $J(-1)=0$, integration gives the identity
in~\eqref{eq:G-2-1}. Entropy monotonicity and $|t|\le1$ imply
$2I(t)\le A(s)^2$, hence the remaining inequalities.
The evolution equation~\eqref{eq:G-2-3} follows from
\cite[proof of Lemma~6.11]{FLC}.
\end{proof}

To interpolate the entropy bound for $\mathrm du$, we first establish a local high derivative bound for the auxiliary scalar.

\begin{lemma}
\label{lem:caloric-local-derivatives}
Under the hypotheses of Proposition~\ref{prop:scalar-comparison},
there are a fixed chart radius $r_*>0$ and constants $C,Q$ such
that, for every
$p=(y,\vec x)\in N_{\rho_1S}$, with $r=|\vec x|$,
\begin{equation}\label{eq:caloric-local-derivatives}
 \sup_{-1\le t\le-1/2}
 \|\mathrm d\tilde F(t)\|_{C^{100n-7}(B_{r_*}(p))}
 \le C(1+r)^{Q}.
\end{equation}
Balls and norms use the fixed reference charts; the enlarged balls
$B_{2r_*}(p)$ lie in $N_{\rho_2S}$ after increasing $S_*$.
The radius $r_*$ and constants $C,Q$ depend only on
$(\rho,c_1,Q_*)$, with the preceding time and metric accuracy
fixed.
At the endpoint $t=-1/2$, every prescribed finite order $j$ also satisfies
\[
 |(\nabla^{g(-1/2)})^j\tilde F(-1/2)|(\iota(p))
       \le C_j(1+r)^2\qquad(p\in N_{\rho_1S}),
\]
after enlarging the radius threshold for a fixed spatial collar.
\end{lemma}
\begin{proof}
Apply the initial-time Schauder estimate used in
Proposition~\ref{prop:scalar-comparison} to $\tilde F$ through
spatial order $100n-6$. Its source is the constant $-n/2$.
On the enlarged ball, \eqref{eq:C-5} bounds its supremum by
$C(1+r)^2$. The initial $C^{100n-6,1/2}$ norm has the same
polynomial bound, since $\tilde F(-1)$ differs from $F_b$ by the
initial comparison error. The coordinate bound at order $100n-2$
and the Ricci-flow equations provide the mixed coefficient bounds;
the incoming $C^{100n}$ bound controls the required lower-order
H\"older trace on fixed reference charts, uniformly in $S$.
Schauder estimates give a $C^{100n-7,1/2}$ bound for
$\mathrm d\tilde F$, and hence~\eqref{eq:caloric-local-derivatives}.

For the endpoint estimates, use fixed-radius intrinsic balls and a
fixed backward time window ending at $-1/2$, contained in
$(-1,-\eta)$. Metric comparison places their spatial collars inside
$\iota(N_{\rho_2S})$ after increasing $S_*$, uniformly in the center.
On each such cylinder,~\eqref{eq:C-5} gives
$|\tilde F|\le C(1+r)^2$: its radial coordinate differs from that of
the center by a bounded amount. The function $\tilde F+nt/2$ solves
the heat equation. The intrinsic curvature derivative estimates
\eqref{eq:A-4} and local noncollapsing give uniform interior derivative
estimates of every prescribed finite order. They yield the asserted
bound at $-1/2$, including order $Q_*+4$, by interior smoothing
after time $-1$.
\end{proof}

\begin{lemma}
\label{lem:entropy-local-derivatives}
Fix the comparison inputs $(\rho,c_1,Q_*)$ and assume the
hypotheses and conclusions of Propositions~\ref{prop:terminal-comparison}
and~\ref{prop:scalar-comparison}, with their time and metric accuracy.
Let $p=(y,\vec x)\in N_{\rho_1S}$, put
$r=|\vec x|$, and choose a fixed sufficiently small $c>0$ so that
the balls below and their prescribed collars lie in $N_{\rho_2S}$.
All local quantities are pulled back by $\iota$.
For $q_r=c/(1+r)$ and $-1\le t\le-1/2$,
\begin{gather}
 \|\mathrm du\|_{L^2(B_{q_r}(p);\mathrm dV_{\bar g})}
       \le C(1+r)^Q e^{(1+c_3)r^2/(8|t|)}A(s),
       \label{eq:C-12}\\*[2pt]
 \|\mathrm du\|_{C^{7}(B_{q_r/2}(p))}
 \le C(1+r)^Q A(s)^{1/2}
    \exp\!\left[\left(\frac{1+c_3}{16|t|}
                          +\frac{c_3}2\right)r^2\right].
       \label{eq:C-13}
\end{gather}
The constants $c,C,Q$ depend only on $(\rho,c_1,Q_*)$,
with those choices fixed, and are uniform in $b,S,t$ and the duration
of a partial block. In normalized radius $l=r/\sqrt{|t|}$, the
exponential in~\eqref{eq:C-13} is bounded by $e^{c_1l^2/2}$.
\end{lemma}
\begin{proof}
The metric equivalence and the upper bound for $F$ in~\eqref{eq:B-7}
imply, throughout $B_{q_r}(p)$,
\[
 \frac{\mathrm d\mu_t}{\mathrm dV_{\bar g}}
 \ge C^{-1}\exp\!\left(-\frac{(1+c_3)r^2}{4|t|}\right).
\]
Indeed, the Euclidean radius changes by $O((1+r)^{-1})$ on this
ball, so its square changes by a bounded amount; also
$1/2\le|t|\le1$. The norms of one-forms for $g(t)$ and $\bar g$
are uniformly equivalent. Applying~\eqref{eq:G-2-1} and taking
square roots proves~\eqref{eq:C-12}, with any polynomial factors
from the fixed local charts included in $(1+r)^Q$.

The local polynomial estimate~\eqref{eq:caloric-local-derivatives}
for $\mathrm d\tilde F$ has constants independent of the outer
radius $S$. The terminal estimates~\eqref{eq:B-8} through order
$Q_*+4>100n-6$, converted to reference derivatives by the
available coordinate connection bounds, consequently give
\[
 \|\mathrm du\|_{C^{100n-7}(B_{q_r}(p))}
       \le C(1+r)^Q e^{c_3 r^2}.
\]
Increasing the fixed polynomial power accommodates the lower
orders and the conversion of covariant derivatives.

Rescale the balls of radii $q_r/2$ and $q_r$ to $B_1\Subset B_2$.
For the rescaled components, interpolation gives
\[
 \|\mathrm du\|_{C^7(B_1)}
 \le C\|\mathrm du\|_{L^2(B_2)}^{1/2}
         \|\mathrm du\|_{C^{100n-7}(B_2)}^{1/2}.
\]
The inequality $100n-7>14+n/2$ leaves a strict interpolation
margin; the lower-order term is absorbed using
$\|\mathrm du\|_{L^2(B_2)}\le C\|\mathrm du\|_{C^{100n-7}(B_2)}$.
Norm conversion costs only powers of $q_r^{-1}$, hence powers of
$1+r$. Substituting the two bounds above and returning to the
original balls proves~\eqref{eq:C-13}.

Finally $r^2=|t| l^2$, and the strict margin~\eqref{eq:C-14}
gives
\[
 \left(\frac{1+c_3}{16|t|}
                 +\frac{c_3}2\right)r^2
 =\left(\frac{1+c_3}{16}
                 +\frac{|t| c_3}2\right)l^2
 <\frac{c_1l^2}{2}\qquad(l>0).
\]
\end{proof}

\subsection{Probability normalization and comparison of the terminal potential}
The total mass determines the additive constant in the actual potential.

\begin{proposition}[Actual metric--potential comparison]
\label{prop:actual-comparison}\label{thm:raw}
Fix $(\rho,c_1,Q_*)$ and make the time and metric-accuracy
choices of Proposition~\ref{prop:terminal-comparison}. There are a
compact-family neighborhood and constants
$\bar\epsilon,S_*,C_1,C,Q,c_4>0$, all depending only on these
comparison inputs, with the following property. For every initial
comparison chart with parameter in this neighborhood, $S\ge S_*$
and $\boldsymbol{\epsilon}(S)\le\bar\epsilon$, the embeddings
\begin{equation}\label{eq:C-15}
 \iota_s(y,\vec x):=\iota(\psi_b(s)y,e^{-s/2}\vec x),
       \qquad 0\le s\le\log2,
\end{equation}
are defined on neighborhoods of $\overline{N_{r_s}}$, where
\begin{equation}\label{eq:C-16}
 r_s:=e^{s/2}(\rho S-C_1).
\end{equation}
Here $\psi_b$ is the compact map defining $Z_b(s)$. Set
$g_s:=|t|^{-1}\iota_s^*g(t)$, $f_s:=\iota_s^*f(t)$ and
\[
 \mu_s:=\iota_s^*\mu_t=(4\pi)^{-n/2}e^{-f_s}\,\mathrm dV_{g_s},
 \qquad t=-e^{-s}.
\]
Write $\Phi_s:=\Phi(g_s,f_s)$ for the soliton defect in this raw chart.
For every partial block with $A(s)^2\le1$ and every
$0<l\le r_s$,
\begin{equation}\label{eq:C-21}
 \begin{split}
 &\|g_s-G_b\|_{C^{8}(N_l)}
        +\|f_s-F_b\|_{C^{8}(N_l)}\\
 &\qquad\le C|\mathrm De(b)|
       +CS^Q\{\boldsymbol{\epsilon}(S)^{1/2}+e^{-c_4S^2}\}
       +C(1+l)^Q e^{c_1l^2/2}A(s)^{1/2}.
 \end{split}
\end{equation}
The constants are uniform in $b,S,\boldsymbol{\epsilon}(S)$ and the duration of
the current partial block. The coefficient $c_4>0$ is the common
output boundary rate from the preceding comparisons and the
probability normalization in the proof below.
\end{proposition}
\begin{proof}
Choose $C_1$ larger than the fixed collars in the preceding results
and then $S_*$ sufficiently large. Proposition~\ref{prop:metric-comparison}
and Proposition~\ref{prop:scalar-comparison} then apply on a
neighborhood of $\overline{N_{\rho S-C_1}}$. The compact flow maps
$\psi_b$ and their inverses have uniform bounds through every
prescribed order on $0\le s\le\log2$. Formula~\eqref{eq:C-15}
therefore gives the stated embeddings and, exactly,
\begin{equation}\label{eq:actual-chart-image}
 \iota_s(N_{r_s})=\iota(N_{\rho S-C_1}).
\end{equation}
Composition with the inverse modified-flow diffeomorphism gives
comparison charts of the modified flow with the same pulled-back pairs.

Write $G_{Z_b(s)}$ and $f_{Z_b(s)}$ for the components of the
auxiliary product.
For the product map
$\Psi_{b,s}(y,\vec x):=(\psi_b(s)y,\sqrt{|t|}\,\vec x)$, we have
$\iota_s=\iota\circ\Psi_{b,s}$ and, directly from the definition of
$Z_b(s)$,
\[
 \begin{split}
 g_s-G_{Z_b(s)}&=|t|^{-1}\Psi_{b,s}^*(\iota^*g(t)-\hat g_b(t)),\\
 f_s-f_{Z_b(s)}&=|t|^{-1}\Psi_{b,s}^*(\iota^*F-\hat f_b)
                              -\log\mathcal M_b(s).
 \end{split}
\]
Identify scalars with their pullbacks by $\iota$. Since $\log\mathcal M_b(s)$
is spatially constant and $F-\hat f_b=(\tilde F-\hat f_b)-u$, we have
\[
 \mathrm d(f_s-f_{Z_b(s)})
 =|t|^{-1}\Psi_{b,s}^*\bigl(\mathrm d(\tilde F-\hat f_b)-\mathrm du\bigr).
\]
Estimate the two terms by~\eqref{eq:C-9} and
Lemma~\ref{lem:entropy-local-derivatives}, respectively, and the metric
difference by~\eqref{eq:A-3}. Put
\[
 \ell_s=\log\frac{\iota_s^*\mu_t}{\mu_{Z_b(s)}}.
\]
The exact determinant identity is
\begin{equation}\label{eq:aux-density-identity}
 \ell_s
 =-(f_s-f_{Z_b(s)})
       +\frac12\log\det\!\left(G_{Z_b(s)}^{-1}g_s\right).
\end{equation}
On a normalized ball centered at Euclidean radius $r$ and of
radius $c/(1+r)$, these estimates and the uniform product-map bounds give
\begin{equation}\label{eq:C-17}
 \|\mathrm d\ell_s\|_{C^{7}}
 \le C S^Q\{\boldsymbol{\epsilon}(S)^{1/2}+e^{-c_4S^2}\}
       +C(1+r)^Q e^{c_1r^2/2}A(s)^{1/2},
\end{equation}
where $c_4$ is no larger than the Gaussian boundary rates in the
metric and scalar comparisons. Fixed changes in the ball radius
only change $C,Q$. The metric difference has the first term
alone.

We next determine the additive constant on the full chart $N_{r_s}$.
The compact part of $f_{Z_b(s)}$ is uniformly bounded, as is the
metric determinant in~\eqref{eq:aux-density-identity}. The terminal
quadratic estimate~\eqref{eq:B-7}, under the dilation $\vec x\mapsto\sqrt{|t|}\,\vec x$,
therefore gives
\begin{equation}\label{eq:C-18}
 |\ell_s(y,\vec x)|\le C+c_3 r^2,
       \qquad |\ell_s(p)|\le C,
       \qquad p=(y_0,0).
\end{equation}
The quadratic coefficient is $c_3/4<c_3$;
$|t|^{-1}\le2$ bounds the remaining constants. Set
$c_{\rm mass}=\ell_s(p)$ and
$z=\ell_s-c_{\rm mass}$. Increasing the fixed constant
in~\eqref{eq:C-18} gives the required rough bound for $z$ as well.
The reference probability is the product of the normalized compact
measure in~\eqref{eq:compact-forward-definition} and $\gamma_d$;
its compact density and metric have uniform bounds.

Both exterior masses are exponentially small in the original
radius $S$. Indeed, \eqref{eq:actual-chart-image} and~\eqref{eq:B-6}
give
\[
 \mu_t(M\setminus \iota_s(N_{r_s}))
 \le C\exp\!\left(-\frac{7(\rho S-C_1)^2}{32|t|}\right),
\]
while the Gaussian factor of $\mu_{Z_b(s)}$ gives the same bound
for $\mu_{Z_b(s)}(\bar M\setminus N_{r_s})$, after changing
$C$. For $S_*\ge2C_1/\rho$, both are at most
$Ce^{-7\rho^2S^2/128}$. Decrease $c_4$ further so that
$c_4<7\rho^2/128$, preserving all preceding estimates.
Join $p=(y_0,0)$ to $(y,0)$ by a compact path of uniformly
bounded length and then follow the radial segment to $(y,\vec x)$.
Both paths lie in $N_{r_s}$. Integrating~\eqref{eq:C-17} gives
\[
 \begin{split}
 |z(y,\vec x)|\le{}&CS^Q(1+r)^{Q+1}
       \{\boldsymbol{\epsilon}(S)^{1/2}+e^{-c_4S^2}\}\\
       &+C(1+r)^{Q+1}e^{c_1r^2/2}A(s)^{1/2}.
 \end{split}
\]
Use this bound for $|z|$ and~\eqref{eq:C-18} for the exponential
in $|e^z-1|\le |z|e^{|z|}$. The Gaussian factor of $\mu_{Z_b(s)}$
and the strict inequality $c_1/2+c_3<1/4$ yield
\[
 \begin{split}
 \int_{N_{r_s}}|e^z-1|\,\mathrm d\mu_{Z_b(s)}
 &\le CS^Q\{\boldsymbol{\epsilon}(S)^{1/2}+e^{-c_4S^2}\}
                         +CA(s)^{1/2},\\
 \left|\int_{N_{r_s}}e^z\,\mathrm d\mu_{Z_b(s)}-1\right|
 &\le CS^Q\{\boldsymbol{\epsilon}(S)^{1/2}+e^{-c_4S^2}\}
                         +CA(s)^{1/2}.
 \end{split}
\]
The second line also uses the auxiliary exterior-mass bound.
The constants are independent of $r_s$: the integrals of
$(1+r)^{Q+1}e^{-(1/4-c_3)r^2}$ and
$(1+r)^{Q+1}e^{-(1/4-c_3-c_1/2)r^2}$ over $\mathbb R^d$ are finite.
Total probability gives the exact identity
\[
 e^{c_{\rm mass}}\int_{N_{r_s}}e^z\,\mathrm d\mu_{Z_b(s)}
       =1-\mu_t(M\setminus\iota_s(N_{r_s})).
\]
Increase $S_*$ so that the actual exterior mass is at most $1/2$.
Together with $|c_{\rm mass}|\le C$, this identity places the
integral in $[e^{-C}/2,e^C]$ and its right side in $[1/2,1]$.
The logarithm is Lipschitz on a fixed positive interval containing
both and $1$. Taking logarithms therefore gives
\[
 |c_{\rm mass}|=|\ell_s(p)|
 \le CS^Q\{\boldsymbol{\epsilon}(S)^{1/2}+e^{-c_4S^2}\}
                         +CA(s)^{1/2}.
\]

Since $\ell_s=z+c_{\rm mass}$, the path estimate gives its
zeroth-order bound on $N_l$. The higher derivatives
through order eight follow directly from~\eqref{eq:C-17}, by
taking the supremum over the local balls. Along the paths the
entropy exponential is at most
$e^{c_1l^2/2}$.
Since $l\le r_s\le CS$, every polynomial multiplying the initial
and Gaussian boundary errors can be included in $S^Q$; the entropy error
retains the factor $(1+l)^Q$.
Equation~\eqref{eq:aux-density-identity} and the metric estimate
now give
\[
 \|g_s-G_{Z_b(s)}\|_{C^{8}(N_l)}
 +\|f_s-f_{Z_b(s)}\|_{C^{8}(N_l)}
 \le CS^Q\{\boldsymbol{\epsilon}(S)^{1/2}+e^{-c_4S^2}\}
       +C(1+l)^Q e^{c_1l^2/2}A(s)^{1/2}.
\]
Since $Z_b(s)-X_b$ is independent of the Euclidean variable,
Proposition~\ref{prop:compact-forward} bounds its $C^{8}(N_l)$ norm
by $C|\mathrm De(b)|$. Adding this bound proves~\eqref{eq:C-21}.
\end{proof}

\section{Improvement toward the compact-factor family}\label{sec:restart}
Put the endpoint pair from Section~\ref{sec:physical} into the reference
slice and apply stationary comparison to a complete auxiliary extension.
The resulting comparison with $X_{b_+}$ has the same order $100n$ as
the incoming chart and satisfies
\[
 E_+^{1/p_1}\le Ce^{-c_2R^2},\qquad
 |\mathrm De(b_+)|\le Ce^{-c_2p_2R^2},\qquad
 \boldsymbol{\epsilon}_+(S')
 \le C(1+S')^Q e^{(c_1S'^2-c_2p_2R^2)/2}.
\]
Theorem~\ref{thm:restart} proves these estimates for one block.
Sections~\ref{sec:radii}--\ref{sec:motion} choose the parameters and
radii and define the recursive charts.

Throughout this section the Ricci-flow block is $[-1,-1/2]$ and
$A=A(\log2)$ is its completed entropy-drop square root. We use the
common incoming order $100n$ and the Ricci-flow hypotheses of
Subsection~\ref{subsec:physical-common-setup}. The four radii have
different roles:
\begin{center}
\begin{tabular}{cl}
 $S$ & radius of the incoming comparison chart,\\
 $R_+$ & endpoint gauge radius,\\
 $R$ & core radius of the stationary extension,\\
 $S'$ & radius of the endpoint comparison chart.
\end{tabular}
\end{center}
Their geometric requirements are
\[
 R_+\le\sqrt2(\rho S-C_1),\qquad
 \Lambda R+10\le R_+,\qquad S'+4\le R,
 \qquad \Lambda>1.
\]
The numerical collars below reserve room for cutoffs, local estimates,
and displacements smaller than $1/10$.

For the combined construction, prescribe $0<p_3<1$ and choose an even
order $Q_*\ge1000n$ such that
\[
 p_3<\frac{Q_*+1}{Q_*+2+n/2}.
\]
Prescribe the comparison estimates of
Subsection~\ref{subsec:physical-common-setup} through order $Q_*+4$
and the stationary representatives of Proposition~\ref{prop:cutoffjets}
through order $Q_*$ before choosing the auxiliary time and
metric accuracy in Proposition~\ref{prop:terminal-comparison}.
The remaining prescribed inputs include $\rho,c_1$, the stationary
exponent $1<p_1<p$, the retention $\Lambda>1$, and the gauge loss
$0<\beta<1$. Constants and thresholds depend on these inputs and
the resulting comparison bounds, uniformly in the compact parameter
and admissible radii. Section~\ref{sec:radii} makes compatible choices
for the iteration.

\subsection{Endpoint gauge and centering}
\label{part8-low-slice-first-with-an-explicit-numerical-hypothesis}

We construct the gauge at order $4$ using the $C^8$ estimate before
gauging from Proposition~\ref{prop:actual-comparison}.

\begin{definition}\label{def:gauge-reference-data}
\phantomsection\label{idx:killing-spaces}
Let $\mathcal K_0$ be the space of compact Killing fields and
Euclidean rotations, and let
$\mathcal K:=\mathcal K_0\oplus\mathbb R^d$ be the full reference
Killing space, with $\mathbb R^d$ denoting constant Euclidean
translations. There are no mixed Killing fields, since the positive
Einstein factor has no parallel one-forms. Projections use
$L^2(\bar\mu)$. Identifying vector fields and one-forms with
$\bar g$, set
$\bar{\mathcal P}:=\delta_{\bar X}\delta_{\bar X}^*$, where
$\delta_{\bar X}^*=-\operatorname{sym}\bar\nabla$.
A subscript $E$ on a vector field denotes its Euclidean component.

\phantomsection\label{idx:adjustable-weights}
For $0<\beta\le1$, define the Gaussian measure
\[
 \mathrm d\mu_\beta
 :=\mathrm d\mu_N(\beta/4\pi)^{d/2}
                    e^{-\beta|\vec x|^2/4}\,\mathrm d\vec x.
\]
Thus $\mu_1=\bar\mu$, and $\beta<1$ imposes stronger Gaussian
integrability. Sobolev norms use the reference connection and metric.
\end{definition}

\begin{proposition}
\label{prop:restart-gauge}\label{prop:endpoint-gauge}
Fix $0<\beta<1$. The exponent $Q_\beta$ depends only on $\beta$ and
the reference product.
Set
\begin{equation}\label{eq:D-1}
 \mathcal C_\beta(R):=(1+R)^{Q_\beta}e^{\beta(R+1)^2}\ge1.
\end{equation}
For $0<\epsilon_{\rm sl}\le\overline\epsilon_{\rm sl}$, there are a
parameter neighborhood $B_{\rm gauge}=B_{\rm gauge}(\beta,\epsilon_{\rm sl})$
and a minimum radius
$R_{+,0}=R_{+,0}(\beta,\epsilon_{\rm sl})>10$ with the following
property. The upper threshold $\overline\epsilon_{\rm sl}>0$ depends
only on the background, and the neighborhood is independent of radius.

Let $b\in B_{\rm gauge}$ and $R_+\ge R_{+,0}$. Let $(G,F)$ be a
smooth metric--potential pair on a neighborhood of
$\overline N_{R_+}$, with associated density $\mu$, and set
$\ell:=\log(\mathrm d\mu/\mathrm d\mu_b)$. Suppose
\[
 \|G-G_b\|_{C^8(N_{R_+})}
 +\|F-F_b\|_{C^8(N_{R_+})}\le \epsilon_+
\]
and
\begin{equation}\label{eq:D-3}
 \mathcal C_\beta(R_+)^2 \epsilon_+\le\frac1{100}.
\end{equation}
Here $\epsilon_+$ bounds the input $C^8$ error;
$\boldsymbol{\epsilon}_+(S')$ in Theorem~\ref{thm:restart}
denotes the outgoing order-$100n$ comparison error.
There is a smooth diffeomorphism $P$ of $\bar M$, equal to the
identity outside $N_{R_+-1}$, whose reference displacement is
less than $1/10$. The actual pulled-back pair
$X^+:=(G^+,F^+)=(P^*G,P^*F)$ and its pulled-back density
$\mu^+:=P^*\mu$ satisfy
$\delta_{\bar X}(G^+-G_b)=0$ on $N_{R_+-4}$. In particular
$\delta_{\bar X}(G^+-\bar g)=0$ there. With
$\ell^+:=\log(\mathrm d\mu^+/\mathrm d\bar\mu)$,
\begin{equation}\label{eq:D-4}
 \|G^+-\bar g\|_{C^3(N_{R_+-4})}
 +\|F^+-\bar f\|_{C^3(N_{R_+-4})}
 +\|\ell^+\|_{C^3(N_{R_+-4})}
 \le\epsilon_{\rm sl}.
\end{equation}
There is also an auxiliary probability $\mu'$ on $\bar M$,
with normalization factor $Z'>0$ satisfying $|\log Z'|\le C\epsilon_+$, such that
\begin{equation}\label{eq:gauge-centered-probability}
 \nu':=P^*\mu',\qquad
 \int\mathrm d\nu'=1,\qquad
 \int\vec x\,\mathrm d\nu'=0,\qquad
 \nu'=(Z')^{-1}\mu^+
       \quad\hbox{on }N_{R_+-4}.
\end{equation}
This centering fixes the Euclidean translation freedom.
For every prescribed $a<1/4$ and finite $Q\ge0$,
\[
 \int_{\{r\ge L\}}(1+r)^Q\,\mathrm d\nu'
       \le C_{a,Q}e^{-aL^2}\qquad(L\ge0).
\]
The constants depend
only on $\beta,\epsilon_{\rm sl}$ and the background, with additional
dependence on $a,Q$ in the tail estimate, and are independent of $b,R_+$
and $\epsilon_+$ satisfying the hypotheses.
\end{proposition}

\begin{proof}
\label{part8-full-nonlinear-construction-in-the-adjustable-weight-norm}
Take the parameter neighborhood of Theorem~\ref{thm:gauge-inverse}
at $\beta$ and $k=4$. We will restrict it to obtain the prescribed
reference closeness.
The density identity
\[
 \ell=-(F-F_b)+\frac12\log\det(G_b^{-1}G)
\]
and the uniform bounds on the fixed preliminary family give
$\|\ell\|_{C^8(N_{R_+})}\le C\epsilon_+$ for sufficiently small $\epsilon_+$,
as ensured by increasing $R_{+,0}$. The constant is uniform on this preliminary family.

Choose smooth radial cutoffs $0\le\chi_i\le1$, $i\in\{0,1,3\}$,
with uniformly bounded derivatives, such that
\[
 \chi_i=1\ \hbox{on }N_{R_+-i-1},\qquad
 \operatorname{supp}\chi_i\subset N_{R_+-i}.
\]
Extend $\chi_0\ell$ by zero and define
\begin{equation}\label{eq:gauge-completion}
 Z':=\int_{\bar M}e^{\chi_0\ell}\,\mathrm d\mu_b,
 \qquad
 \mathrm d\mu':=(Z')^{-1}e^{\chi_0\ell}\,\mathrm d\mu_b.
\end{equation}
Since $\mu_b$ is centered,
\[
 |\log Z'|\le C\epsilon_+,\qquad
 \vec m':=\int\vec x\,\mathrm d\mu',\qquad
 |\vec m'|\le C\epsilon_+ .
\]
Indeed $|e^{\chi_0\ell}-1|\le C\epsilon_+$ and the first Gaussian
moment of $\mu_b$ is uniform. On $N_{R_+-1}$,
$\mu'=(Z')^{-1}\mu$.

All $C^{5,1/2}$ norms of vector fields in this proof are taken on
$N_{R_+}$. Use the reference exponential map and put
$P_U:=\exp_{\bar g}(\chi_1U)$ for $\|U\|_{C^{5,1/2}}$ sufficiently small.
\phantomsection\label{idx:gauge-map}
Choose this neighborhood so that the displacement is less than $1/10$.
For every $U$ in this neighborhood, the small $C^1$ norm makes each
$\exp_{\bar g}(t\chi_1 U)$, $0\le t\le1$, a local diffeomorphism.
These maps are the identity outside a compact set, hence proper coverings
of the connected reference product. A point outside this compact set and
its image has exactly one preimage, so each covering has one sheet.
Thus they are global diffeomorphisms.

Since $\operatorname{supp}\chi_3\subset N_{R_+-3}$,
\[
 P_U(\operatorname{supp}\chi_3)\Subset N_{R_+-1}.
\]
This inclusion persists for small variations of $U$ and the radius;
hence the pullback is defined near the support and its product with
$\chi_3$ extends smoothly by zero.
Keep $\mu'$ fixed when varying $U$. On that support,
where $\chi_1=1$, Taylor expansion gives
\begin{equation}\label{eq:gauge-nonlinear-expansions}
 \begin{split}
 P_U^*G-G_b
   &=h_0+L_UG_b+Q_b(U,h_0),
       \qquad h_0:=G-G_b,\\
 \int\vec x\,P_U^*(\mathrm d\mu')
   &=\vec m'-\int\chi_1 U_E\,\mathrm d\mu_b+\mathcal R_g(U).
 \end{split}
\end{equation}
On $N_{R_+-3}$ the metric remainder has the decomposition
\[
 Q_b(U,h_0)
 =\bigl(P_U^*G_b-G_b-L_UG_b\bigr)
       +\bigl(P_U^*h_0-h_0\bigr).
\]
The first term is quadratic in $U$; the second has
$C^{4,1/2}$ norm at most $C\epsilon_+\|U\|_{C^{5,1/2}}$.
For the moment identity, change variables to obtain
$\int(P_U^{-1})_E\,\mathrm d\mu'$ and use
\[
 (P_U^{-1})_E=\vec x-(\chi_1U_E)\circ P_U^{-1}.
\]
Subtracting $\vec x-\chi_1U_E$ gives a quadratic inverse-map
error bounded by $C\|U\|_{C^{5,1/2}}^2$. Replacing
$\mu'$ by $\mu_b$ in the linear term costs
$C\epsilon_+\|U\|_{C^{5,1/2}}$, since
$\|\mathrm d\mu'/\mathrm d\mu_b-1\|_{L^\infty}\le C\epsilon_+$.
The product and composition estimates give
\[
 \begin{split}
 \|Q_b(U,h_0)\|_{C^{4,1/2}}
       +|\mathcal R_g(U)|
   &\le C\bigl(\epsilon_+\|U\|_{C^{5,1/2}}+\|U\|_{C^{5,1/2}}^2\bigr),\\
 \|Q_b(U,h_0)-Q_b(V,h_0)\|_{C^{4,1/2}}
       &+|\mathcal R_g(U)-\mathcal R_g(V)|\\
   &\le C\bigl(\epsilon_++\|U\|_{C^{5,1/2}}+\|V\|_{C^{5,1/2}}\bigr)
                                      \|U-V\|_{C^{5,1/2}}.
 \end{split}
\]
The tensor norms are on $N_{R_+-3}$. Uniform bounds for the product
geometry, compact-family coefficients, cutoff derivatives and Gaussian
moments make $C$ independent of $R_+$.
The raw $C^8$ bound controls the $C^{7,1/2}$ norm and supplies
the derivatives needed in the composition and difference estimates.

Keep $G_b$ in
$\mathcal P_b=-\frac12\delta_{\bar X}L_{(\cdot)}G_b$ and solve
\begin{equation}\label{eq:gauge-augmented-equation}
 \begin{split}
 \mathcal P_bU
   &=\frac12\delta_{\bar X}
          \{\chi_3[h_0+Q_b(U,h_0)]\},\\
 \int\chi_1 U_E\,\mathrm d\mu_b&=\vec m'+\mathcal R_g(U),
       \qquad \Pi_{\mathcal K_0}U=0 .
 \end{split}
\end{equation}
The first source annihilates every reference Killing field $Y$:

\[
 \int\left\langle\delta_{\bar X}
      \{\chi_3(h_0+Q_b)\},Y\right\rangle\mathrm d\bar\mu
 =\int\left\langle\chi_3(h_0+Q_b),
                         -\frac12L_Y\bar g\right\rangle
                                           \mathrm d\bar\mu=0 .
\]
The first-moment condition and
$\mathcal K_0$ projection use the reference normalizations of
Theorem~\ref{thm:gauge-inverse}.

For the complete source data, put
\[
 T_U:=\frac12\chi_3\bigl[h_0+Q_b(U,h_0)\bigr],
 \qquad \vec m_U:=\vec m'+\mathcal R_g(U).
\]
The raw bounds, the uniformly bounded cutoff derivatives and the
local nonlinear estimates give
\[
 \begin{split}
 \|T_U\|_{C^{4,1/2}}+|\vec m_U|
   &\le C\bigl(\epsilon_++\epsilon_+\|U\|_{C^{5,1/2}}+\|U\|_{C^{5,1/2}}^2\bigr),\\
 \|T_U-T_V\|_{C^{4,1/2}}+|\vec m_U-\vec m_V|
   &\le C\bigl(\epsilon_++\|U\|_{C^{5,1/2}}+\|V\|_{C^{5,1/2}}\bigr)
                                      \|U-V\|_{C^{5,1/2}}.
 \end{split}
\]
These tensor norms are global, with the cutoff tensors extended
by zero. The fixed terms $h_0$ and $\vec m'$ cancel in
the difference estimate. All constants so far are independent of
$R_+$ and $\epsilon_{\rm sl}$ on the preliminary neighborhood.

Let $\mathcal TU$ be the global solution of the augmented linear
system with source $\delta_{\bar X}T_U$ and moment $\vec m_U$.
Equip $H^2(\mu_{\beta})\cap C^{5,1/2}(N_{R_+})$ with the norm
\[
 \|U\|_*:=\|U\|_{H^2(\mu_{\beta})}+\|U\|_{C^{5,1/2}}.
\]
Since $\operatorname{supp}T_U\subset N_{R_+-3}$,
Theorem~\ref{thm:gauge-inverse} applies at the prescribed $\beta$ with
$S=R_+$ and $k=4$.
Equations~\eqref{eq:G-1-2} and~\eqref{eq:G-1-3} give, respectively,
\[
 \begin{split}
 \|\mathcal TU\|_{H^2(\mu_{\beta})}
   &\le C(\beta)\bigl(\|T_U\|_{C^1}+|\vec m_U|\bigr),\\
 \|\mathcal TU\|_{C^{5,1/2}}
   &\le C(\beta)(1+R_+)^Qe^{\beta R_+^2}
                  \bigl(\|T_U\|_{C^{4,1/2}}+|\vec m_U|\bigr).
 \end{split}
\]
For the first estimate,
$|\delta_{\bar X}T_U|\le C(|\bar\nabla T_U|+(1+r)|T_U|)$
and $\int(1+r)^2\,\mathrm d\mu_{\beta}<\infty$,
independently of the support radius.
The parameter $b$, cutoffs and normalization conditions stay fixed
during the contraction. Thus the same linear inverse applies to
$T_U-T_V$ and $\vec m_U-\vec m_V$.
Adding these bounds and applying the source estimates gives
\[
 \begin{split}
 \|\mathcal TU\|_*
   &\le C(\beta)(1+R_+)^Qe^{\beta R_+^2}
             \bigl(\epsilon_++\epsilon_+\|U\|_*+\|U\|_*^2\bigr),\\
 \|\mathcal TU-\mathcal TV\|_*
   &\le C(\beta)(1+R_+)^Qe^{\beta R_+^2}
       \bigl(\epsilon_++\|U\|_*+\|V\|_*\bigr)\|U-V\|_*.
 \end{split}
\]
Enlarge $Q\ge0$ and $C(\beta)\ge1$ to include the fixed source
constants and the finitely many low-order pullback and density estimates
below. Thus $C(\beta)$ and $Q$ depend only on $\beta$, the reference
product and these fixed orders.
Increase $Q_\beta$ to absorb the fixed constants, so that
\[
 C(\beta)(1+R)^Qe^{\beta R^2}\le\mathcal C_\beta(R)
       \qquad(R\ge1).
\]
This choice depends only on $\beta$ and the fixed background and
orders, and is made before $\epsilon_{\rm sl}$ is prescribed.
Consequently,
\[
 \begin{split}
 \|\mathcal TU\|_*&\le\mathcal C_\beta(R_+)
             \bigl(\epsilon_++\epsilon_+\|U\|_*+\|U\|_*^2\bigr),\\
 \|\mathcal TU-\mathcal TV\|_*
   &\le\mathcal C_\beta(R_+)
       \bigl(\epsilon_++\|U\|_*+\|V\|_*\bigr)\|U-V\|_*.
 \end{split}
\]
By~\eqref{eq:D-3}, increasing $R_{+,0}$ puts the closed ball
$\|U\|_*\le2\mathcal C_\beta(R_+)\epsilon_+$ inside the low-order
neighborhood where the preceding expansions hold. On this ball,
\[
 \begin{split}
 \|\mathcal TU\|_*
 &\le\mathcal C_\beta(R_+)\epsilon_+
       \bigl(1+6\mathcal C_\beta(R_+)^2\epsilon_+\bigr)
       \le2\mathcal C_\beta(R_+)\epsilon_+,\\
 \operatorname{Lip}(\mathcal T)
 &\le5\mathcal C_\beta(R_+)^2\epsilon_+<\frac12.
 \end{split}
\]
Thus $\mathcal T$ preserves the ball and is a strict contraction.
Its fixed point satisfies
\begin{equation}\label{eq:gauge-displacement-bound}
 \|U\|_{H^2(\mu_{\beta})}
 +\|U\|_{C^{5,1/2}(N_{R_+})}
       \le2\mathcal C_\beta(R_+)\epsilon_+ .
\end{equation}
The second inverse factor in~\eqref{eq:D-3} comes from substituting
the ball radius into the Lipschitz estimate.
For every fixed polynomial cost $(1+R_+)^Q$ in the pullback, density
and coefficient estimates, the construction ball satisfies
\begin{equation}\label{eq:gauge-prefactor-smallness}
 \begin{split}
 \epsilon_+&\le\frac1{100\mathcal C_\beta(R_+)^2},\qquad
 \|U\|_*\le\frac1{50\mathcal C_\beta(R_+)},\\
 (1+R_+)^Q\bigl(\epsilon_++\|U\|_*+\|U\|_*^2\bigr)
       &\le\frac{C(1+R_+)^Q}{\mathcal C_\beta(R_+)}
       =C(1+R_+)^{Q-Q_\beta}e^{-\beta(R_++1)^2}.
 \end{split}
\end{equation}
With $\beta$ and $Q_\beta$ fixed, the right-hand side of
\eqref{eq:gauge-prefactor-smallness} tends to zero as $R_+\to\infty$,
uniformly on the construction ball. Increasing $R_{+,0}$ therefore
gives the required coefficient bounds.

The bound in \eqref{eq:gauge-displacement-bound} gives displacement
less than $1/10$ throughout the isotopy. The smooth raw data and local
elliptic regularity of the gauge equation give smoothness of $U$
on $\operatorname{supp}\chi_1$; outside the nonlinear source the
linear equation gives the same conclusion.

On $N_{R_+-4}$, the cutoff $\chi_3$ equals one, so
the first equation in~\eqref{eq:gauge-augmented-equation} gives
$\delta_{\bar X}(P_U^*G-G_b)=0$ there. Since $G_b$ lies in the
reference slice, this also gives
$\delta_{\bar X}(P_U^*G-\bar g)=0$.
The first-moment equation gives the exact center in
\eqref{eq:gauge-centered-probability}.

Write the product density as
$\mathrm d\mu_b=e^{v_b(y)}e^{-|\vec x|^2/4}
 \mathrm dV_{g_N}\mathrm d\vec x$, absorbing its fixed constant
in $v_b$. On the retained core the logarithmic density of the
actual pullback relative to $\mu_b$ is exactly
\begin{equation}\label{eq:gauge-pullback-density}
 \ell\circ P+\log\operatorname{Jac}_{\mathrm dV_{g_N}
                  \mathrm d\vec x}P
       +v_b(P_N)-v_b(y)-\frac{|P_E|^2-|\vec x|^2}{4}.
\end{equation}
Since $P_E-\vec x=(\chi_1U)_E$, the $C^3$ norm of the last term is
bounded by $C(1+R_+)(\|U\|_{C^3}+\|U\|_{C^3}^2)$.
The Jacobian requires the $C^4$ norm of $U$. By
\eqref{eq:gauge-prefactor-smallness} the resulting error is bounded
by $C(1+R_+)^Q/\mathcal C_\beta(R_+)$. The metric pullback is
controlled in the same way.
Both pullback estimates hold near $\operatorname{supp}\chi_3$, since
$P(\operatorname{supp}\chi_3)\subset N_{R_+-1}$. Moreover,
\[
 F^+-F_b
   =\frac12\log\det(G^+/G_b)
               -\log(\mathrm d\mu^+/\mathrm d\mu_b).
\]
The logarithmic determinant is locally Lipschitz in $C^3$, uniformly
on the preliminary family. Combining these estimates gives
\[
 \begin{split}
 &\|G^+-G_b\|_{C^3(N_{R_+-4})}
 +\|F^+-F_b\|_{C^3(N_{R_+-4})}
 +\left\|\log\frac{\mathrm d\mu^+}{\mathrm d\mu_b}
       \right\|_{C^3(N_{R_+-4})}\\
 &\hspace{2cm}\le
 \frac{C(1+R_+)^Q}{\mathcal C_\beta(R_+)}
 \le\frac12\epsilon_{\rm sl},
 \end{split}
\]
where the last inequality holds for every $R_+\ge R_{+,0}$ after
increasing $R_{+,0}$. Use the same increase to meet the prescribed
displacement, metric and logarithmic-density coefficient bounds on
the entire construction ball. Next restrict the inverse neighborhood
to an invariant neighborhood $B_{\rm gauge}$ on which
\[
 \|G_b-\bar g\|_{C^3(\bar M)}
 +\|F_b-\bar f\|_{C^3(\bar M)}
 +\left\|\log\frac{\mathrm d\mu_b}{\mathrm d\bar\mu}
      \right\|_{C^3(\bar M)}
 \le\frac12\epsilon_{\rm sl}.
\]
The compact-factor differences vanish at $b=0$, so this restriction
and the required compact-family coefficient bounds are uniform in
$R_+$ and $\epsilon_+$.
The three triangle inequalities, using
$\ell^+=\log(\mathrm d\mu^+/\mathrm d\mu_b)
          +\log(\mathrm d\mu_b/\mathrm d\bar\mu)$,
give~\eqref{eq:D-4}.

On the retained core the extended measure differs from the
actual one only by $(Z')^{-1}$, which proves the last identity in
\eqref{eq:gauge-centered-probability}. Globally the same pullback
formula uses $\chi_0\ell-\log Z'$ in place of $\ell$.
The same polynomial estimate, including the term
$\vec x\cdot(P_E-\vec x)$, gives the global comparison
\[
 \left\|\log\frac{\nu'}{\bar\mu}\right\|_{L^\infty}
       \le C|b|+\frac{C(1+R_+)^Q}{\mathcal C_\beta(R_+)}
       \le C(|b|+\epsilon_{\rm sl}).
\]
Outside the displacement support the cutoff term $\chi_0\ell$
and its normalization are $O(\epsilon_+)$; outside the extension support
only the compact density and this scalar normalization remain.
The constant is independent of $\epsilon_{\rm sl},b,R_+$.
The chosen neighborhood and minimum radius bound the logarithmic
ratio uniformly; Gaussian integration gives the asserted tails for
every $a<1/4$.
\end{proof}

The construction also gives the global tensor
\begin{equation}\label{eq:gauge-patched-tensor}
 \mathfrak h_U
   :=\chi_3(P_U^*G-G_b)+(1-\chi_3)L_UG_b,
       \qquad \delta_{\bar X}\mathfrak h_U=0.
\end{equation}
The first cutoff product is extended by zero using the support margin
above. The divergence identity follows from
\eqref{eq:gauge-augmented-equation}, and on $N_{R_+-4}$ this tensor
equals $P_U^*G-G_b$.

For a Ricci-flow comparison chart at a current time
$0\le s\le\log2$, the input pair in
Proposition~\ref{prop:restart-gauge} is $(G,F)=(g_s,f_s)$.
Proposition~\ref{prop:actual-comparison} supplies its raw bound
whenever $A(s)\le1$ and $R_+\le e^{s/2}(\rho S-C_1)$.
Namely, after increasing its constants, put
\phantomsection\label{idx:gauge-radius}
\begin{equation}\label{eq:D-2}
 \epsilon_+=C|\mathrm De(b)|+CS^Q\{\boldsymbol{\epsilon}(S)^{1/2}+e^{-c_4S^2}\}
       +C(1+R_+)^Qe^{c_1R_+^2/2}A(s)^{1/2}.
\end{equation}
At the endpoint $s=\log2$, this uses $A(s)=A$ and
$R_+\le\sqrt2(\rho S-C_1)$.

To estimate the variation of the gauge, keep the raw pair, $b$, the cutoffs and
$\mu'$ fixed. Let $U$ range over the construction ball
\[
 \|U\|_*:=\|U\|_{H^2(\mu_{\beta})}
             +\|U\|_{C^{5,1/2}(N_{R_+})}
       \le2\mathcal C_\beta(R_+)\epsilon_+.
\]
For every such $U$, put $G^+=P_U^*G$ on a neighborhood of
$\operatorname{supp}\chi_3$ and
$\nu'=P_U^*\mu'$. Formula
\eqref{eq:gauge-patched-tensor} defines $\mathfrak h(U,G,b)$
throughout this ball; its divergence need not vanish.

\phantomsection\label{idx:weak-gauge-data}
Choose a smooth radial cutoff $\chi_2$ with fixed-width gaps
\[
 \operatorname{supp}\chi_3\Subset\{\chi_2=1\}
 \Subset\operatorname{supp}\chi_2\Subset\{\chi_1=1\}.
\]
For variations with the cutoffs fixed, define
\[
 A_UY:=(\mathrm DP_U)^{-1}\mathrm D_UP_U[Y].
\]
Choose a smooth radial cutoff $\zeta$, equal to one near
$\operatorname{supp}\chi_2$ and supported in the interior of
$\{\chi_1=1\}$, with derivatives bounded independently of $R_+$.
Set
\[
 B_U:=\Id+\zeta(A_U-\Id).
\]
The construction bounds give
\[
 \|B_U-\Id\|_{C^1}+\|B_U^{-1}-\Id\|_{C^1}
       \le C\|U\|_{C^2(N_{R_+})}.
\]
In particular, $B_U$ is invertible and agrees with $A_U$
near $\operatorname{supp}\chi_2$, whereas $A_U$ vanishes
where $\chi_1=0$. Define the first-moment component by
\[
 \mathcal M_U Z
   :=\int(A_UB_U^{-1}Z)_E\,\mathrm d\nu'.
\]

Set $H^1:=H^1(\bar\mu;T\bar M)$, $H^{-1}:=(H^1)^*$, and
\[
 H^{-1}_{\mathcal K}
 :=\{F\in H^{-1}:F(K)=0\text{ for all }K\in\mathcal K\},
 \qquad
 \mathcal Y:=H^{-1}_{\mathcal K}\times\mathcal K_0\times\mathbb R^d.
\]
All tensor norms use $\bar g$; the finite-dimensional components have their fixed
Euclidean norms. On the strong construction ball define
\[
 \mathcal F(U):=\left(
 -\frac12\delta_{\bar X}\mathfrak h(U,G,b),\quad
 \Pi_{\mathcal K_0}U,\quad
 -\int\vec x\,\mathrm d(P_U^*\mu')\right).
\]
Its first component is the divergence of a whole $L^2(\bar\mu)$ symmetric
tensor, and hence belongs to $H^{-1}_{\mathcal K}$.

\begin{lemma}[Augmented inverse in $H^1$]\label{lem:augmented-weak-inverse}
The reference operator
\[
 \mathscr PV
 :=\left(\bar{\mathcal P}V,\Pi_{\mathcal K_0}V,
                    \int V_E\,\mathrm d\bar\mu\right)
 :H^1\longrightarrow\mathcal Y
\]
is an isomorphism. There is a tolerance $\epsilon_{\rm w}>0$,
depending only on the reference geometry and fixed cutoff profiles,
with the following property. If
\begin{equation}\label{eq:weak-coefficient-smallness}
 \begin{split}
 &\|G^+-\bar g\|_{C^1(\operatorname{supp}\chi_3)}
       +\|G_b-\bar g\|_{C^1}
       +\|B_U-\Id\|_{C^1}+\|B_U^{-1}-\Id\|_{C^1}\\
 &\quad+\|U\|_{C^1(N_{R_+})}
       +\left\|\frac{\nu'}{\bar\mu}-1\right\|_\infty
       +\|1-\chi_1\|_{L^2(\bar\mu)}
       \le\epsilon_{\rm w},
 \end{split}
\end{equation}
then
\[
 \mathbb L_U Z:=\left(
 -\frac12\delta_{\bar X}
  [\chi_3L_ZG^++(1-\chi_3)L_{B_U^{-1}Z}G_b],\quad
 \Pi_{\mathcal K_0}(B_U^{-1}Z),\quad\mathcal M_UZ\right)
\]
is an isomorphism $H^1\to\mathcal Y$, with inverse norm independent
of $R_+$. Moreover,
\[
 \mathrm D\mathcal F(U)[Y]=\mathbb L_U(B_UY),
\]
and this differential extends to $H^1\to\mathcal Y$ with
\begin{equation}\label{eq:weak-gauge-derivative-close}
 \|\mathrm D\mathcal F(U)-\mathscr P\|_{H^1\to\mathcal Y}
 \le\epsilon,\qquad
 \epsilon\|\mathscr P^{-1}\|_{\mathcal Y\to H^1}
 \le\frac12.
\end{equation}
Here $\epsilon>0$ is fixed by the reference inverse.

Under the hypotheses of Proposition~\ref{prop:restart-gauge}, the
left side of \eqref{eq:weak-coefficient-smallness} is bounded,
throughout its strong construction ball, by
\begin{equation}\label{eq:weak-construction-coefficients}
 C_{\rm w}\bigl(|b|+\epsilon_{\rm sl}+e^{-a(R_+-C)^2}\bigr)
 \qquad(0<a<1/8).
\end{equation}
For fixed $\beta,a$ and background data, $C_{\rm w}$ and $C$ are independent of
$\epsilon_{\rm sl},R_+,b$ and $\epsilon_+$. Thus the tolerance is met by
first decreasing $\epsilon_{\rm sl}$ and the parameter neighborhood,
and then increasing the minimum radius.
\end{lemma}
\begin{proof}
Gaussian integration by parts and approximation by compactly supported
smooth fields give
\[
 \|rV\|_2^2
 =2d\|V\|_2^2+4\sum_i\int x_i\langle V,\bar\nabla_iV\rangle
                       \,\mathrm d\bar\mu,
 \qquad \|rV\|_2\le C\|V\|_{H^1}.
\]
Expanding the symmetric gradient and using
$\Ric_{\bar g}+\bar\nabla^2\bar f=\bar g/2$ gives
\[
 2\|\delta_{\bar X}^*V\|_2^2
 =\|\bar\nabla V\|_2^2+\|\delta_{\bar X}V\|_2^2
                          -\frac12\|V\|_2^2.
\]
The multiplier estimate makes the weighted divergence continuous
on $H^1$, so both identities extend by density.
For $V\perp\mathcal K$, Proposition~\ref{prop:reference-gauge-inverse}
supplies $W\perp\mathcal K$ with
$\bar{\mathcal P}W=V$ and $\|W\|_{H^2}\le C\|V\|_2$. Hence
\[
 \|V\|_2^2
 =(\delta_{\bar X}^*W,\delta_{\bar X}^*V)_{L^2(\bar\mu)}
 \le C\|V\|_2\|\delta_{\bar X}^*V\|_2.
\]
Consequently
\begin{equation}\label{eq:M-8}
 \|V\|_{H^1}\le C\|\operatorname{sym}\bar\nabla V\|_2
     \quad(V\perp\mathcal K),
 \qquad \|rV\|_2\le C\|V\|_{H^1}.
\end{equation}
Lax--Milgram applied to
$\int\langle\delta_{\bar X}^*V,\delta_{\bar X}^*Y\rangle
\,\mathrm d\bar\mu$ solves the metric equation on
$H^1\cap\mathcal K^\perp$. Adding a field in $\mathcal K_0$ and
a constant translation supplies the other two components uniquely: the
translation moment vanishes on $\mathcal K_0$. This proves the
reference isomorphism.

For $Z,Y\in H^1$, integration by parts pairs the tensor in the
divergence term with $\delta_{\bar X}^*Y$.
The identity
\[
 (L_ZG)_{ij}=G_{aj}\bar\nabla_iZ^a+G_{ia}\bar\nabla_jZ^a
                         +Z^a\bar\nabla_aG_{ij}
\]
shows that the difference of the metric components is bounded in
$H^{-1}$ by $C\epsilon_{\rm w}\|Z\|_{H^1}$.
The difference of the symmetry projections is at most
$C\|B_U^{-1}-\Id\|_{C^0}\|Z\|_2$.
Differentiating the reference exponential gives
\[
 \|A_UB_U^{-1}-\chi_1\Id\|_{C^0}
             \le C\|U\|_{C^1(N_{R_+})}.
\]
Thus Cauchy--Schwarz gives
\[
 \begin{split}
 \left|\mathcal M_UZ-\int Z_E\,\mathrm d\bar\mu\right|
 &\le C\left(\|U\|_{C^1(N_{R_+})}
       +\left\|\frac{\nu'}{\bar\mu}-1\right\|_\infty
       +\|1-\chi_1\|_2\right)\|Z\|_2\\
 &\le C\epsilon_{\rm w}\|Z\|_{H^1}.
 \end{split}
\]
Combining the components yields
$\|\mathbb L_U-\mathscr P\|\le C\epsilon_{\rm w}$.
A fixed decrease of $\epsilon_{\rm w}$ makes the Neumann series
converge, with a uniform inverse bound.

For frozen data, pullback differentiation gives
\[
 \mathrm D_U\mathfrak h[Y]
 =\chi_3L_{A_UY}G^++(1-\chi_3)L_YG_b.
\]
Also, integration by parts in the density variation gives
\[
 -\mathrm D_U\!\left[\int\vec x\,\mathrm d(P_U^*\mu')\right][Y]
 =\int(A_UY)_E\,\mathrm d\nu'.
\]
The varying field has compact support in $\operatorname{supp}\chi_1$,
which justifies the integration by parts for every $U$ in the ball.
Since $A_U=B_U$ near
$\operatorname{supp}\chi_3$, these formulas give
$\mathrm D\mathcal F(U)[Y]=\mathbb L_U(B_UY)$. Multiplication by
$B_U$ is a bounded $H^1$ isomorphism and
\[
 \mathrm D\mathcal F(U)-\mathscr P
 = (\mathbb L_U-\mathscr P)B_U
                      +\mathscr P(B_U-\Id).
\]
Its norm is at most $C\epsilon_{\rm w}$. Decreasing the fixed
tolerance once more proves \eqref{eq:weak-gauge-derivative-close}.

On the construction ball, the compact-family errors are $O(|b|)$. The pullback and multiplier
errors are bounded by the right side of
\eqref{eq:gauge-prefactor-smallness}, and hence by
$C\epsilon_{\rm sl}$ at the radius threshold of
Proposition~\ref{prop:restart-gauge}. The support margin
$P_U(\operatorname{supp}\chi_3)\Subset N_{R_+-1}$ holds throughout the
ball, so the metric bound holds on all of $\operatorname{supp}\chi_3$.
The same calculation for the
probability uses \eqref{eq:gauge-pullback-density} with
$\chi_0\ell-\log Z'$ in place of $\ell$. It includes the term
$\vec x\cdot(P_E-\vec x)$ and gives globally
\[
 \left\|\log\frac{\nu'}{\bar\mu}\right\|_\infty
          \le C(|b|+\epsilon_{\rm sl}).
\]
Finally $1-\chi_1$ is supported where
$r\ge R_+-2$, so Gaussian integration gives
$\|1-\chi_1\|_2\le C_a e^{-a(R_+-C)^2}$ for any fixed $a<1/8$.
This proves \eqref{eq:weak-construction-coefficients}, with constants
independent of $\epsilon_{\rm sl}$ and $R_+$.
\end{proof}

\subsection{Higher regularity at the endpoint}
\label{part8-high-endpoint-bounds-without-increasing-the-incoming-order}
\label{part8-interior-slice-regularity-in-sufficient-detail}

\begin{lemma}
\label{lem:restart-coordinate-regularity}
Let $Q_*$ and $c_3$ be the even endpoint order and Gaussian loss fixed in
Subsection~\ref{subsec:physical-common-setup}, with the terminal and caloric
estimates prescribed through order $Q_*+4$. Let
$u(t)=\tilde F(t)-|t|f(t)$ as in Lemma~\ref{thm:caloric}, and put
$u^+=u(-1/2)\circ(\iota_{\log2}\circ P)$. Under the Ricci-flow comparison
hypotheses and the gauge smallness condition of
Proposition~\ref{prop:restart-gauge}, the following estimates hold on $N_{R_+-5}$:
\begin{equation}\label{eq:D-5}
 |\bar\nabla^jG^+|\le C_j(1+r)^{Q_j}
                         \qquad(0\le j\le Q_*+3),
\end{equation}
and
\begin{equation}\label{eq:D-6}
 |\bar\nabla^ju^+|+|\bar\nabla^jF^+|
       \le C_j(1+r)^{Q_j}e^{8c_3 r^2}
                         \qquad(0\le j\le Q_*+3).
\end{equation}
The constants depend on $\rho,c_1,Q_*,\beta,\epsilon_{\rm sl}$
and the corresponding comparison bounds, but not on $b$ or $R_+$.
\end{lemma}
\begin{proof}
Put $h^+:=G^+-\bar g$. The exact slice condition is
\[
 \delta_{\bar X}h^+=0
 \quad\Longleftrightarrow\quad
 \bar\nabla^i h^+_{ij}=h^+_{ij}\bar\nabla^i\bar f.
\]
Substituting this identity in the Ricci equation gives an elliptic
system for $h^+$. Its principal part at the reference metric is
$-\frac12(\bar\Delta h^++\bar\nabla^2\Tr_{\bar g}h^+)$;
the remaining second-order coefficients are $O(|h^+|)$, and the
drift grows like $1+r$. The metric closeness in~\eqref{eq:D-4}
preserves ellipticity. On nested reference balls of radius
$c/(1+r)$, with $c>0$ fixed sufficiently small, the rescaled drift
and reference coefficients have uniform bounds. These balls lie in
$N_{R_+-4}$ when their centers lie in $N_{R_+-5}$.

Start with the $C^3$ bound in~\eqref{eq:D-4} and apply interior
elliptic estimates successively on these nested balls. At each step,
the preceding metric bounds control the connection difference and
convert the intrinsic Ricci derivative bounds from~\eqref{eq:A-4}
into the reference derivatives needed on the right-hand side.
One additional intrinsic derivative gives their local H\"older bound.
The differentiated equation then gains the next metric derivative;
its highest-order terms remain in the same elliptic principal part.
All coefficient products and rescaling costs are finite powers of
$1+r$. This proves~\eqref{eq:D-5}, with local
$C^{Q_*+3,1/2}$ bounds and constants independent of $R_+$.

For the scalar bounds, set
\[
 \iota_+=\iota_{\log2}\circ P,\qquad
 G^+=2\iota_+^*g(-1/2).
\]
The definitions $F(t)=|t|f(t)$ and $u=\tilde F-F$ give
\[
 F^+=f(-1/2)\circ \iota_+=2F(-1/2)\circ \iota_+,
 \qquad u^+=(\tilde F-F)(-1/2)\circ \iota_+.
\]
Proposition~\ref{prop:terminal-comparison} bounds the intrinsic
derivatives of $F,f$ through order $Q_*+4$ with Gaussian
loss $c_3$. Lemma~\ref{lem:caloric-local-derivatives} gives polynomial
intrinsic bounds through the same order for $\tilde F$.
For a scalar $v$, naturality and invariance of the connection under
constant metric scaling give
\[
 (\nabla^{G^+})^j(v\circ \iota_+)
       =\iota_+^*((\nabla^{g(-1/2)})^jv),\qquad
 |(\nabla^{G^+})^j(v\circ \iota_+)|_{G^+}
       =2^{-j/2}|(\nabla^{g(-1/2)})^jv|_{g(-1/2)}\circ \iota_+.
\]
For $p=(y,\vec x)$, the Euclidean radius of $\iota_+(p)$ in the
original chart is $|P^E(p)|/\sqrt2$, since the compact map does
not change the Euclidean coordinate and $\iota_{\log2}$ dilates
it by $1/\sqrt2$. Thus
\[
 \frac{|P^E(p)|}{\sqrt2}\le\frac{r+1/10}{\sqrt2},\qquad
 \frac{c_3|P^E(p)|^2}{2}\le c_3r^2+C,\qquad r=|\vec x|.
\]
The displacement margin and
$R_+\le\sqrt2(\rho S-C_1)$ put the relevant points and local
collars in the scalar-estimate region. Hence the pulled-back intrinsic scalar derivatives are bounded by
a polynomial times $e^{c_3r^2}$, and therefore by the right-hand side
of~\eqref{eq:D-6}.

Finally, write
$\mathcal B=\nabla^{G^+}-\bar\nabla=(G^+)^{-1}*\bar\nabla G^+$,
where $*$ denotes fixed tensor contractions. By~\eqref{eq:D-5},
the required derivatives of $\mathcal B$ have polynomial growth.
For a scalar $v$, the first conversion formulas are schematically
\[
 \begin{aligned}
 \bar\nabla^2v&=(\nabla^{G^+})^2v+\mathcal B*\mathrm dv,\\
 \bar\nabla^3v&=(\nabla^{G^+})^3v
       +\mathcal B*(\nabla^{G^+})^2v
       +(\bar\nabla\mathcal B+\mathcal B*\mathcal B)*\mathrm dv.
 \end{aligned}
\]
In order $j$, each further term is a product of connection
derivatives through order $j-2$ with one lower intrinsic derivative
of $v$. Consequently the conversion changes only the polynomial
factor. Applying it to $u^+$ and $F^+$ proves~\eqref{eq:D-6}.
For the local H\"older estimates at order $Q_*+3$, use the
metric H\"older bounds just proved and the intrinsic scalar derivative
of order $Q_*+4$ on an enlarged ball of radius $c/(1+r)$.
The remaining rescaling factors are again polynomial.
\end{proof}
\subsection{Smallness of the endpoint soliton defect}
\label{part8-uniform-oscillator-control-and-the-completed-block-defect}
We work on the completed Ricci-flow block $[-1,-1/2]$ with the
actual endpoint gauge of Proposition~\ref{prop:restart-gauge}.
Write
\[
 \iota_+=\iota_{\log2}\circ P,\qquad
 u^+=u(-1/2)\circ \iota_+,\qquad u=\tilde F-F.
\]
Thus $G^+=2\iota_+^*g(-1/2)$ and
$\mu^+=\iota_+^*\mu_{-1/2}$. The measures used here are restrictions
and pullbacks of the actual terminal probability. We assume $0\le A=A(\log2)\le1$.

Prescribe $0<p_3<1$ and choose the endpoint order so that
\[
 p_3<\frac{Q_*+1}{Q_*+2+n/2}.
\]
This order, together with the other prescribed comparison orders, is
fixed before the auxiliary time and metric accuracy. The derived
loss satisfies $16c_3<1/4$. Constants below are uniform in
$b,S,R_+$ and depend on the fixed comparison and gauge inputs and
on the stated interpolation exponents.

\begin{lemma}\label{lem:restart-hessian}
Fix
\[
 p_3<\theta<\frac{Q_*+1}{Q_*+2+n/2}.
\]
For the comparison and gauge inputs above, there are $C,Q$ and a
minimum radius $\underline R_+$ such that every endpoint satisfying
the preceding hypotheses with $R_+\ge\underline R_+$ obeys
\[
 |\nabla_{G^+}^2u^+|(z)
 \le C(1+r)^Q A^\theta
       \exp\!\left[\left(\frac\theta8
                      +8c_3(1-\theta)\right)r^2\right],
 \qquad z\in N_{R_+-7},\quad r=|\vec x(z)|.
\]
Consequently,
\begin{equation}\label{eq:restart-hessian-bound}
 \|\nabla_{G^+}^2u^+\|_{L^2(\mu^+);N_{R_+-7}}
       \le C A^\theta\le C A^{p_3}.
\end{equation}
The tensor norms use $G^+$. The constants depend on the fixed inputs
and $\theta$, independently of $b$ and $R_+$.
\end{lemma}
\begin{proof}
Put $V:=\mathrm du^+$. For $z\in N_{R_+-7}$ and $r=|\vec x(z)|$,
choose a reference coordinate ball of radius
$\varrho=c/(1+r)$, with fixed $c>0$ sufficiently small that
$B_{2\varrho}(z)\subset N_{R_+-5}$ and the doubled ball lies in a
uniform reference atlas. The squared Euclidean radius varies by a
bounded amount on this ball. Thus the reference Gaussian density
is comparable to $e^{-r^2/4}$ there. The metric and measure
comparisons~\eqref{eq:D-4}, followed by~\eqref{eq:G-2-1}, give
\[
 \begin{split}
 \|V\|_{L^2(B_{2\varrho}(z);\mathrm dV_{\bar g})}^2
 &\le Ce^{r^2/4}\int_{B_{2\varrho}(z)}
                         |\mathrm du^+|_{G^+}^2\,\mathrm d\mu^+\\
 &=\frac C2e^{r^2/4}\int_{\iota_+(B_{2\varrho}(z))}
             |\mathrm du(-1/2)|_{g(-1/2)}^2\,\mathrm d\mu_{-1/2}
 \le Ce^{r^2/4}A^2.
 \end{split}
\]
The factor $1/2$ comes from $G^+=2\iota_+^*g(-1/2)$.

Equation~\eqref{eq:D-6} bounds the reference derivatives of $V$
through order $Q_*+2$ on the doubled ball. Rescale its
coordinates to $B_2$, and denote the component functions by $V_r$.
Changing between component and covariant norms and rescaling costs
only powers of $1+r$. Increasing the fixed polynomial exponent gives
\[
 \|V_r\|_{L^2(B_2)}\le C(1+r)^Q e^{r^2/8}A,
 \qquad
 \|V_r\|_{C^{Q_*+2}(B_2)}
       \le C(1+r)^Q e^{8c_3r^2}.
\]
Ordinary interior interpolation on the fixed balls gives
\[
 \|V_r\|_{C^1(B_1)}
 \le C\|V_r\|_{L^2(B_2)}^\theta
       \|V_r\|_{C^{Q_*+2}(B_2)}^{1-\theta}.
\]
The strict upper bound on $\theta$ is the interpolation condition in
dimension $n$; the lower-order term is absorbed by
$\|V_r\|_{L^2(B_2)}\le C\|V_r\|_{C^{Q_*+2}(B_2)}$.
Substituting the preceding bounds and returning to reference derivatives
gives the asserted radial estimate for $|V|+|\bar\nabla V|$ at $z$. The connection difference is
uniformly bounded by~\eqref{eq:D-4}, and
\[
 \nabla_{G^+}^2u^+
 =\bar\nabla V-(\Gamma_{G^+}-\Gamma_{\bar g})V.
\]
Uniform metric equivalence proves the pointwise assertion.

Finally $\mu^+\le C\bar\mu$ on this region. Squaring the bound and
combining its exponential with the Gaussian density leaves
\[
 \exp\!\left[-(1-\theta)
                    \left(\frac14-16c_3\right)r^2\right].
\]
Since $16c_3<1/4$, this factor absorbs every fixed polynomial.
Integration over $N\times\mathbb R^d$ therefore bounds the squared
norm in~\eqref{eq:restart-hessian-bound} by $CA^{2\theta}$,
independently of $R_+$. Taking square roots and using
$0\le A\le1$ and $\theta>p_3$ finishes the proof.
\end{proof}

The following estimate is analogous to \cite[Lemma~6.11]{FLC}.

\begin{lemma}
\label{lem:restart-lichnerowicz}
For the completed block and the endpoint gauge above,
\begin{equation}\label{eq:restart-lichnerowicz-bound}
 \|\tilde T(-1/2)\|_{L^2(\mu_{-1/2});\iota_+(N_{R_+-9})}
       \le C A,
 \qquad \tilde T=|t|\Ric+\nabla^2\tilde F-g/2.
\end{equation}
Here the tensor norm uses the metric $g(-1/2)$. The
constant depends on the curvature and coordinate bounds
on the block and the fixed cutoff profile, and is uniform in the
varying radii and compact parameter.
\end{lemma}
\begin{proof}
Choose a time-independent smooth cutoff $\chi$ on $M$ which,
in the original chart $\iota$, is one on
$N_{(R_+-8)/\sqrt2}$ and supported in
$N_{(R_+-7)/\sqrt2}$. Extend it by zero on the closed manifold.
The cutoff depends only on the Euclidean radius in that chart.
The compact map $\psi_b(\log2)$ preserves this radius, the
endpoint dilation is $1/\sqrt2$, and $P$ has displacement less
than $1/10$. Consequently $\chi=1$ on
$\iota_+(N_{R_+-9})$. Moreover,
\[
 \frac{R_+-7}{\sqrt2}
       \le\rho S-C_1-\frac7{\sqrt2},
\]
so its support has a fixed collar inside the propagated Ricci-flow
cylinder throughout $[-1,-1/2]$. The transition width is fixed.
The coordinate bounds in the original chart therefore give
\[
 |\mathrm d\chi|_{g(t)}+|\Delta_{g(t)}\chi|\le C
       \qquad(-1\le t\le-1/2),\qquad \partial_t\chi=0.
\]

By \eqref{eq:G-2-3}, $\tilde T$ solves
$\partial_t\tilde T=\Delta_L\tilde T$. The evolution of the
norm under Ricci flow cancels the Ricci terms in $\Delta_L$,
giving
\[
 (\partial_t-\Delta_g)|\tilde T|_g^2
 =-2|\nabla\tilde T|_g^2
       +4\langle\Rm_g(\tilde T),\tilde T\rangle_g.
\]
Set
$\mathcal E_\chi(t)=\int\chi^2|\tilde T|_g^2\,\mathrm d\mu_t$.
Conjugate-heat duality gives its derivative by integrating
$(\partial_t-\Delta_g)(\chi^2|\tilde T|^2)$.
The curvature is bounded on the support. The mixed term is
bounded by
\[
 8\chi|\mathrm d\chi|\,|\tilde T|\,|\nabla\tilde T|
 \le\chi^2|\nabla\tilde T|^2
                       +C|\mathrm d\chi|^2|\tilde T|^2,
\]
and the remaining cutoff term is bounded by
$C|\tilde T|^2$ on its support. Thus, since $0\le\chi\le1$,
\[
 \frac{\mathrm d}{\mathrm dt}\big((t+1)\mathcal E_\chi(t)\big)
 \le (1+C(t+1))\int_M|\tilde T(t)|_g^2\,\mathrm d\mu_t .
\]
The initial term vanishes. Integrating from $-1$ to $-1/2$ and
using~\eqref{eq:G-2-1} gives $\mathcal E_\chi(-1/2)\le CA^2$,
which proves~\eqref{eq:restart-lichnerowicz-bound}.
\end{proof}

\begin{proposition}\label{prop:restart-defect}
Fix $p_3$ and $Q_*$ as above. Every completed block satisfying
the comparison hypotheses and Proposition~\ref{prop:restart-gauge}, with
$R_+$ above the thresholds of the preceding two lemmas, obeys
\begin{equation}\label{eq:D-9}
 \|\Phi(G^+,F^+)\|_{L^2(\mu^+);N_{R_+-9}}
       \le C A^{p_3}.
\end{equation}
The norm uses $G^+$. The constant depends on $p_3$ and the prescribed
comparison and gauge inputs, independently of the admissible radii.
\end{proposition}
\begin{proof}
Recall $T=|t|\Ric+\nabla^2F-g/2$ and
$T=\tilde T-\nabla^2u$. Constant metric rescaling does not change
the connection, so at the endpoint
\[
 \iota_+^*(\nabla_{g(-1/2)}^2u)=\nabla_{G^+}^2u^+,
 \qquad
 \iota_+^*T(-1/2)=-\frac12\Phi(G^+,F^+).
\]
At the endpoint, for any covariant two-tensor $Q$,
\[
 |\iota_+^*Q|_{G^+}=\frac12|Q|_{g(-1/2)}\circ \iota_+,
       \qquad \mu^+=\iota_+^*\mu_{-1/2}.
\]
Consequently the Lichnerowicz lemma controls
$2\iota_+^*\tilde T$ by $CA$ in the required actual
endpoint norm. The Hessian lemma controls
$2\nabla_{G^+}^2u^+$ by $CA^{p_3}$. Combining them on $N_{R_+-9}$ gives
$C(A+A^{p_3})\le C A^{p_3}$ for $0\le A\le1$, proving
\eqref{eq:D-9}.
\end{proof}
\subsection{Complete extension and stationary comparison}
\label{part8-stationary-extension-and-high-order-smallness}
We now construct from the actual endpoint pair $X^+$ a complete
pair $X^R$ to which Theorem~\ref{thm:logarithmic} applies.
The auxiliary probability from Proposition~\ref{prop:restart-gauge}
will control its center. Write the endpoint difference as
\[
 G^+=\bar g+h^+,\qquad F^+=\bar f+k^+,\qquad
 \ell^+=\log\frac{\mu^+}{\bar\mu}
       =\frac12\log\det(\Id+\bar g^{-1}h^+)-k^+.
\]
All these fields are defined on the actual endpoint chart. Fix
$\Lambda>1$ and a stationary core radius $R$ satisfying
\begin{equation}\label{eq:D-10}
 \Lambda R+10\le R_+.
\end{equation}
In particular, the entire cutoff collar below lies in the exact
slice and in the endpoint-defect domain.

\begin{lemma}\label{lem:restart-extension}
Given $\Lambda>1$, $0<c_2<1/8$ and $1<p_1<p$, there are
\[
 \bar\epsilon_{\rm sl}=\bar\epsilon_{\rm sl}(\Lambda,p,p_1)>0,
 \qquad R_*=R_*(\Lambda,c_2,p,p_1)>0
\]
with the following property. Assume the conclusions of
Propositions~\ref{prop:endpoint-gauge} and~\ref{prop:restart-defect},
with $0\le A\le1$, $\epsilon_{\rm sl}\le\bar\epsilon_{\rm sl}$,
$R\ge R_*$ and~\eqref{eq:D-10}. The estimate constant in the defect bound also depends on the
prescribed comparison and gauge bounds.

Choose a smooth radial cutoff $0\le\chi_R\le1$, equal to one on $N_R$
and zero outside $N_{\Lambda R}$, with
$|\bar\nabla^j\chi_R|\le C_{j,\Lambda}R^{-j}$, and set
\begin{equation}\label{eq:stationary-extension-definition}
 \begin{split}
 h^R&:=\chi_Rh^+,\qquad
 Z_R:=\int_{\bar M}e^{\chi_R\ell^+}\,\mathrm d\bar\mu,\\
 k^R&:=\frac12\log\det(\Id+\bar g^{-1}h^R)
                       -\chi_R\ell^++\log Z_R,\\
 X^R&:=(g^R,f^R)=(\bar g+h^R,\bar f+k^R),\qquad
 \mu^R:=Z_R^{-1}e^{\chi_R\ell^+}\bar\mu.
 \end{split}
\end{equation}
Products with $\chi_R$ are extended by zero. The metric $g^R$ is
smooth and complete, and $\mu^R$ is a probability measure. The pair
$X^R$ satisfies the hypotheses of Theorem~\ref{thm:logarithmic} for the
prescribed $p_1$.
On $N_R$,
\begin{equation}\label{eq:stationary-extension-core}
 g^R=G^+,\qquad f^R=F^++\log Z_R,\qquad
 \mu^R=Z_R^{-1}\mu^+.
\end{equation}
Define its defect, slice error and center by
\[
 D_R:=\|\Phi(X^R)\|_{L^2(\mu^R)},\qquad
 a_{{\rm sl},R}:=\|(\Id-\Pi_{\rm sl})h^R\|_{H^2(\bar\mu)},\qquad
 \vec m^{\,R}:=\int_{\bar M}\vec x\,\mathrm d\mu^R.
\]
Then
\begin{equation}\label{eq:stationary-extension-preliminary}
 D_R\le C\bigl(A^{p_3}+e^{-c_2R^2}\bigr),\qquad
 a_{{\rm sl},R}+|\vec m^{\,R}|\le Ce^{-c_2R^2}.
\end{equation}
The constants are uniform in $b,S,R_+,R$ after the prescribed
parameters have been fixed.
\end{lemma}
\begin{proof}
Let $\epsilon>0$ be the tolerance in
Theorem~\ref{thm:logarithmic} for the prescribed $p_1$; it depends only
on the background and $p,p_1$. Take $R_*\ge1$. By~\eqref{eq:D-4}, the sum of the
$C^2$ norms of $h^+,k^+$ and $\ell^+$ on $N_{R_+-4}$ is at most
$\epsilon_{\rm sl}$. Since $R\ge1$, multiplication by $\chi_R$
preserves this bound up to a factor depending only on $\Lambda$.
For sufficiently small $\epsilon_{\rm sl}$, the metric $\bar g+h^R$
is uniformly equivalent to $\bar g$. It equals $\bar g$ outside a
compact set and is therefore complete. Moreover,
$e^{-\epsilon_{\rm sl}}\le Z_R\le e^{\epsilon_{\rm sl}}$, so
$|\log Z_R|\le\epsilon_{\rm sl}$. The determinant formula gives
$\mu^R$ in~\eqref{eq:stationary-extension-definition}, its mass one,
and the core identities~\eqref{eq:stationary-extension-core}.
It also bounds the global $C^2$ norm of $k^R$ by
$C_\Lambda\epsilon_{\rm sl}$. Consequently
\[
 \|h^R\|_{C^2}
 +\sup_{\bar M}\left\{
   \frac{|k^R|}{1+|\vec x|^2}
   +\frac{|\bar\nabla k^R|}{1+|\vec x|}
   +|\bar\nabla^2k^R|\right\}
 \le C_\Lambda\epsilon_{\rm sl}.
\]
Choose $\bar\epsilon_{\rm sl}$ so that
$C_\Lambda\bar\epsilon_{\rm sl}\le\epsilon$. This verifies all
hypotheses of Theorem~\ref{thm:logarithmic}.

A constant shift of the potential does not change $\Phi$, so on
$N_R$ the squared defect integral is $Z_R^{-1}$ times the actual
one. Proposition~\ref{prop:restart-defect} bounds its square root
by $CA^{p_3}$, since $N_{\Lambda R}\subset N_{R_+-9}$.
Outside $N_{\Lambda R}$ the pair is
$(\bar g,\bar f+\log Z_R)$ and its defect vanishes.
On the transition annulus, the local Ricci and Hessian formulas
involve at most two derivatives of $h^R,k^R$; the only unbounded
background coefficient is $\mathrm d\bar f=O(1+r)$.
Consequently
\[
 |\Phi(X^R)|\le C(1+r)^Q,
 \qquad C^{-1}\bar\mu\le\mu^R\le C\bar\mu,
 \qquad r=|\vec x|.
\]
The Gaussian integral of its square on $r\ge R$ gives an
$L^2$ bound $C(1+R)^Q e^{-R^2/8}$.

Exact slicing holds throughout the cutoff collar, so
\[
 \delta_{\bar X}h^R
       =(\mathrm d\chi_R)\lrcorner h^+.
\]
The reference projection satisfies
\[
 \|(\Id-\Pi_{\rm sl})v\|_{H^2(\bar\mu)}
       \le C\|\delta_{\bar X}v\|_{H^1(\bar\mu)}.
\]
Indeed,
\[
 (\Id-\Pi_{\rm sl})v
       =\delta_{\bar X}^*\bar{\mathcal P}^{-1}\delta_{\bar X}v,
\]
where the inverse is taken on $\mathcal K^\perp$.
Apply~\eqref{eq:G-1-11} with $j=1$ and use
$\delta_{\bar X}^*=-\operatorname{sym}\bar\nabla$.
For $v=h^R$, the source is compatible by integration by parts,
since it is the divergence of a compactly supported symmetric tensor.
This source equals $(\mathrm d\chi_R)\lrcorner h^+$ and is supported
on the transition annulus. Differentiating it once uses only the already controlled
first derivatives of $h^+$ and two derivatives of the cutoff.
Its Gaussian $H^1$ norm is again bounded by
$C(1+R)^Q e^{-R^2/8}$.

Let $\nu'=P^*\mu'$ be the exactly centered probability from
\eqref{eq:gauge-centered-probability}. On $N_{R_+-4}$,
$\nu'=(Z')^{-1}\mu^+$, where $Z'$ normalizes the
auxiliary probability. Both $Z'$ and $Z_R$ are bounded above and
below independently of the radii. Hence on $N_R$,
$\mu^R=(Z'/Z_R)\nu'$, and exact centering gives
\[
 \vec m^{\,R}
 =-\frac{Z'}{Z_R}\int_{\bar M\setminus N_R}\vec x\,\mathrm d\nu'
       +\int_{\bar M\setminus N_R}\vec x\,\mathrm d\mu^R.
\]
The probability $\mu'$ has bounded logarithmic density relative to
$\mu_b$, and $P$ has displacement less than $1/10$.
Changing variables by $P$ therefore bounds the discarded first
moment of $\nu'$ by a Gaussian tail with every fixed exponent
strictly below $1/4$. The same statement for $\mu^R$ follows
from its bounded density ratio relative to $\bar\mu$.

For all three estimates choose a preliminary exponent
$c_2<c_+<1/8$. The preceding integrals, including the fixed
displacement losses, are bounded by
$C(1+R)^Q e^{-c_+R^2}$. The strict gap $c_+-c_2$ absorbs this
polynomial for all $R$ above a fixed threshold. Combining these
bounds with the actual core defect proves
\eqref{eq:stationary-extension-preliminary}.
\end{proof}

For the improvement step, assume
\begin{equation}\label{eq:D-13}
 A^{p_3}\le e^{-c_2R^2}.
\end{equation}
The preliminary estimates then imply
\begin{equation}\label{eq:D-12}
 D_R\le Ce^{-c_2R^2},\qquad
 a_{{\rm sl},R}+|\vec m^{\,R}|\le Ce^{-c_2R^2}.
\end{equation}

\begin{lemma}\label{lem:restart-entropy-localization}
\label{part8-entropy-tails-are-converted-at-the-energy-scale}
In the preceding setup,
\begin{equation}\label{eq:restart-entropy-localization}
 |\log Z_R|+
 |\mathcal W(2g(-1/2),f(-1/2))-\mathcal W(X^R)|
       \le Ce^{-3R^2/16}.
\end{equation}
Constants are uniform in the varying radii.
\end{lemma}
\begin{proof}
Set $W:=\mathcal W(2g(-1/2),f(-1/2))$ and
$\iota_+=\iota_{\log2}\circ P$. Since $P$ and $P^{-1}$ move Euclidean
radius by less than $1/10$ and $\psi_b$ maps $N$ onto itself,
\[
 \iota_+(N_R)\supset
 \iota\bigl(N_{(R-1/10)/\sqrt2}\bigr).
\]
All these regions lie in the Ricci-flow chart by~\eqref{eq:D-10}
and $R_+\le\sqrt2(\rho S-C_1)$. Increase the minimum radius
so that $(R-1/10)/\sqrt2\ge R_c$.
At the endpoint $t=-1/2$, the exponent in~\eqref{eq:B-6}
therefore becomes
\[
 \frac{(1-1/8)((R-1/10)/\sqrt2)^2}{4|t|}
       =\frac7{32}(R-1/10)^2.
\]
The strict gap $7/32-3/16$ absorbs the fixed linear collar loss.
Equation~\eqref{eq:B-6}, including its entropy weights, gives
\begin{equation}\label{eq:actual-endpoint-entropy-tail}
 \int_{M\setminus \iota_+(N_R)}
 \left(1+|f(t)|+|t||\mathrm R_{g(t)}|
                      +|t||\mathrm df(t)|_{g(t)}^2\right)\,\mathrm d\mu_t
       \le Ce^{-3R^2/16}.
\end{equation}

On $N_R$, $e^{\chi_R\ell^+}\bar\mu=\mu^+$, whose mass is
$\mu_t(\iota_+(N_R))$. On its complement the unnormalized extended
density is bounded by $C\bar\mu$. Hence
\[
 |Z_R-1|
 \le\mu_t(M\setminus \iota_+(N_R))
       +\int_{\bar M\setminus N_R}e^{\chi_R\ell^+}\,\mathrm d\bar\mu
 \le Ce^{-3R^2/16}.
\]
The second term is bounded at every fixed rate below $1/4$;
reserve a rate above $3/16$ to absorb its polynomial factor.
For large $R$, $Z_R\in[1/2,3/2]$, and the same bound holds
for $|\log Z_R|$ and $|Z_R^{-1}-1|$.

Temporarily remove the constant shift and put
$\widetilde f^R=f^R-\log Z_R$ and
$\mathrm d\widetilde\mu^R=e^{\chi_R\ell^+}\mathrm d\bar\mu$.
Consider the unnormalized entropy integral
\[
 \widetilde W_R
 =\int_{\bar M}
 \left(\mathrm R_{g^R}+|\mathrm d\widetilde f^R|_{g^R}^2
                         +\widetilde f^R-n\right)
                                  \,\mathrm d\widetilde\mu^R.
\]
After pullback by $\iota_+$, the entropy densities agree on $N_R$.
Equation~\eqref{eq:actual-endpoint-entropy-tail} bounds the entropy integral
outside $\iota_+(N_R)$.
On the cutoff annulus $N_{\Lambda R}\setminus N_R$, the bounded
$C^2$ input and the stated cutoff derivatives give
\[
 |\mathrm R_{g^R}|+|\mathrm d\widetilde f^R|_{g^R}^2
                         +|\widetilde f^R|+n\le C(1+r)^Q.
\]
Outside $N_{\Lambda R}$ this is the corresponding reference
integrand. Both regions have density at most $C\bar\mu$.
Their entropy integrals are therefore bounded by
$C(1+R)^Q e^{-R^2/4}\le Ce^{-3R^2/16}$.
Thus $|W-\widetilde W_R|\le Ce^{-3R^2/16}$, and
$|\widetilde W_R|\le C$.

Finally, the exact effect of the normalization is
\[
 \mathcal W(X^R)=Z_R^{-1}\widetilde W_R+\log Z_R.
\]
The bounds for $Z_R$ and $\widetilde W_R$ give
\eqref{eq:restart-entropy-localization}.
\end{proof}

\begin{proposition}
\label{prop:restart-stationary-comparison}
Given $1<p_1<p$, put $p_2=p_1/p$ and use the extension of
Lemma~\ref{lem:restart-extension} for this $p_1$.
There are a low-error threshold
and a minimum radius, depending on $p_1,\Lambda,c_2,p_3$ and the
prescribed comparison and gauge bounds, with the following property.
Under the extension hypotheses and~\eqref{eq:D-13}, the endpoint
entropy gap satisfies
\begin{equation}\label{eq:restart-entropy}
 E_+^{1/p_1}\le Ce^{-c_2R^2},\qquad
 E_+=\Theta_N-\mathcal W(2g(-1/2),f(-1/2))\ge0.
\end{equation}
Define $(b_+,a):=\Pi_Kh^R$, where $h^R=\chi_R(G^+-\bar g)$.
Choose the low-error threshold so that the closed parameter ball
of radius $\epsilon_{\rm sl}$ lies in the representative domain of
Proposition~\ref{prop:cutoffjets}. Then~\eqref{eq:D-4} gives
\[
 |(b_+,a)|=\|\Pi_Kh^R\|_{L^2(\bar\mu)}
       \le\|h^R\|_{L^2(\bar\mu)}
       \le\epsilon_{\rm sl}.
\]
Here $\Pi_K$ is orthogonal in the fixed metric $L^2$ parameter norm
of Definition~\ref{def:internal-parameters}, and
$\operatorname{supp}\chi_R\subset N_{\Lambda R}\subset N_{R_+-4}$
by~\eqref{eq:D-10}.
Thus $(b_+,a)$ lies in the fixed representative domain of
Proposition~\ref{prop:cutoffjets}, and the positive, mass-one, centered
representative is defined. Write
\[
 \hat X_{b_+,a}=(\bar g+\hat h(b_+,a),\bar f+\hat k(b_+,a))
\]
and put
\[
 H^R=(h^R,k^R)-(\hat h(b_+,a),\hat k(b_+,a)).
\]
Then
\begin{equation}\label{eq:D-14}
 \begin{split}
 |a|^2+|\mathrm De(b_+)|
       &\le Ce^{-c_2p_2R^2},\\
 \|H^R\|_{H^2(\bar\mu)}&\le Ce^{-c_2p_2R^2}.
 \end{split}
\end{equation}
The constants have the same dependencies as the thresholds and are
uniform in the compact parameter and the admissible radii.
\end{proposition}
\begin{proof}
Apply Theorem~\ref{thm:logarithmic} to $X^R$, using
Lemma~\ref{lem:restart-extension} and the defect, slice and center
bounds~\eqref{eq:D-12}. This gives
$|\mathcal W(X^R)-\Theta_N|^{1/p_1}\le Ce^{-c_2R^2}$.
Entropy monotonicity and~\eqref{eq:entropy-limit} give $E_+\ge0$.
Lemma~\ref{lem:restart-entropy-localization} and the triangle
inequality for the $1/p_1$ power imply
\[
 E_+^{1/p_1}
       \le Ce^{-c_2R^2}+Ce^{-3R^2/(16p_1)}\le Ce^{-c_2R^2}.
\]
Here $c_2<1/8$ and $p_1<p\le3/2$ give $3/(16p_1)>c_2$,
so the entropy-localization error after taking the $1/p_1$ power
is absorbed by the first term. This proves~\eqref{eq:restart-entropy}.
The metric--potential comparison estimate~\eqref{eq:L-6-2} gives the two
bounds in~\eqref{eq:D-14}, since
$D_R^{p_2}+a_{{\rm sl},R}+|\vec m^{\,R}|
\le Ce^{-c_2p_2R^2}$ because $0<p_2<1$.
\end{proof}
\subsection{Endpoint comparison at fixed order}
\label{subsec:restart-strong-seed}
We convert the weighted stationary comparison into an endpoint comparison
at order $100n$, retaining the actual potential constant.
Use the common analytical inputs above, with the local comparison
and gauge hypotheses, and the extension and stationary-comparison
hypotheses of Lemma~\ref{lem:restart-extension} and
Proposition~\ref{prop:restart-stationary-comparison}.

\begin{theorem}[Endpoint improvement]\label{thm:restart}
For the preceding inputs, there are an error threshold $\bar\epsilon>0$,
a minimum radius, and constants $C,Q$ with the following property.
Suppose an incoming comparison chart on $N_S$ has parameter
$b$ in the common comparison neighborhood of Section~\ref{sec:physical}
and in $B_{\rm gauge}(\beta,\epsilon_{\rm sl})$ from
Proposition~\ref{prop:restart-gauge}, with
$\boldsymbol{\epsilon}(S)\le\bar\epsilon$.
Let $0\le A=A(\log2)\le1$ and let $\epsilon_+$ be given by~\eqref{eq:D-2}
at the endpoint. Require $S,R_+,R$ to exceed the minimum radius and
assume
\[
 \mathcal C_\beta(R_+)^2\epsilon_+\le\frac1{100},\qquad
 R_+\le\sqrt2(\rho S-C_1),\qquad \Lambda R+10\le R_+,
\]
\[
 A^{p_3}\le e^{-c_2R^2},\qquad S'\ge1,\qquad S'+4\le R.
\]
Then the endpoint gives $(b_+,a)$ in the stationary representative domain,
with $|(b_+,a)|\le\epsilon_{\rm sl}$, and an endpoint comparison
chart $\iota_+=\iota_{\log2}\circ P$ on a neighborhood
of $\overline N_{S'}$, at order $100n$, satisfying
\begin{equation}\label{eq:D-18}
 \begin{gathered}
 E_+^{1/p_1}\le Ce^{-c_2R^2},\qquad
 |\mathrm De(b_+)|\le Ce^{-c_2p_2R^2},\\
 \boldsymbol{\epsilon}_+(S')
 \le C(1+S')^Q
       \exp\!\left(\frac{c_1(S')^2-c_2p_2R^2}{2}\right).
 \end{gathered}
\end{equation}
Here $E_+=\Theta_N-\mathcal W(2g(-1/2),f(-1/2))\ge0$,
and $\boldsymbol{\epsilon}_+$ is the norm~\eqref{eq:seed-error}
for the endpoint pair $(2g(-1/2),f(-1/2))$ in the chart $\iota_+$,
relative to $X_{b_+}$ and with its actual potential constant.
The thresholds and constants depend only on the prescribed inputs
and background, uniformly in $b$ and all admissible radii.
\end{theorem}
\begin{proof}
Proposition~\ref{prop:actual-comparison} supplies the raw endpoint bound,
and Proposition~\ref{prop:restart-gauge} gives the endpoint gauge.
Lemma~\ref{lem:restart-coordinate-regularity} supplies the high coordinate
bounds. Proposition~\ref{prop:restart-defect} and
Lemma~\ref{lem:restart-extension} construct the complete extension
$X^R$ with defect, slice error and center bounded by $Ce^{-c_2R^2}$.
Lemma~\ref{lem:restart-entropy-localization} and
Proposition~\ref{prop:restart-stationary-comparison} give the first
two estimates in~\eqref{eq:D-18}.
With $(b_+,a)$ and $\hat X_{b_+,a}$ as in that proposition, put
\[
 H^R=(h^R,k^R)-(\hat h(b_+,a),\hat k(b_+,a)).
\]
Equation~\eqref{eq:D-14} gives
\[
 \|H^R\|_{H^2(\bar\mu)}\le Ce^{-c_2p_2R^2},\qquad
 |a|\le Ce^{-c_2p_2R^2/2}.
\]
On $N_R$, the core identities~\eqref{eq:stationary-extension-core},
Lemma~\ref{lem:restart-coordinate-regularity},
Proposition~\ref{prop:cutoffjets} and the bound for $\log Z_R$ in
Lemma~\ref{lem:restart-entropy-localization} imply
\[
 |\bar\nabla^j H^R(y,\vec x)|
 \le C_j(1+r)^{Q_j}e^{8c_3r^2},\qquad
 r=|\vec x|,\quad 0\le j\le Q_*.
\]

At a point of $N_{S'}$, take a reference ball of radius $c/(1+r)$,
with $c>0$ fixed and small. Its doubled ball lies in $N_R$ by
$S'+4\le R$, and $N_R\subset N_{R_+-5}$ by $\Lambda R+10\le R_+$.
The Gaussian density varies by a bounded factor on this ball.
Let $V_r$ denote the components of $H^R$ rescaled to $B_2$.
The preceding estimates give
\[
 \|V_r\|_{L^2(B_2)}
 \le C(1+r)^Qe^{r^2/8}\|H^R\|_{L^2(\bar\mu)},\qquad
 \|V_r\|_{C^{Q_*}(B_2)}
 \le C(1+r)^Qe^{8c_3r^2}.
\]
Since $Q_*\ge1000n$, one has
$\frac12<(Q_*-100n)/(Q_*+n/2)$.
Ordinary local interpolation therefore yields
\[
 \|V_r\|_{C^{100n}(B_1)}
 \le C\|V_r\|_{L^2(B_2)}^{1/2}
           \|V_r\|_{C^{Q_*}(B_2)}^{1/2}.
\]
Here the inhomogeneous term is absorbed using
$\|V_r\|_{L^2(B_2)}\le C\|V_r\|_{C^{Q_*}(B_2)}$.
Returning to reference derivatives costs only a polynomial in $1+r$.
The Gaussian exponents satisfy
\[
 \frac{c_1}{2}-\left(\frac1{16}+4c_3\right)=12c_3>0.
\]
Taking the supremum over the centers controls the full integer norm
with uniform constants, giving
\[
 \begin{split}
 \|H^R\|_{C^{100n}(N_{S'})}
 &\le C(1+S')^Qe^{c_1(S')^2/2}\|H^R\|_{L^2(\bar\mu)}^{1/2}\\
 &\le C(1+S')^Qe^{(c_1(S')^2-c_2p_2R^2)/2}.
 \end{split}
\]

The actual endpoint difference on $N_{S'}$ is exactly
\[
 (G^+,F^+)-X_{b_+}
 =H^R+(\hat X_{b_+,a}-X_{b_+})-(0,\log Z_R).
\]
By Proposition~\ref{prop:cutoffjets} and
Lemma~\ref{lem:restart-entropy-localization}, respectively,
\[
 \|\hat X_{b_+,a}-X_{b_+}\|_{C^{100n}(N_{S'})}
 \le C(1+S')^Q|a|,\qquad
 |\log Z_R|\le Ce^{-3R^2/16}.
\]
Since $3/16>c_2p_2/2$, these terms and the preceding interpolation
bound prove the last estimate in~\eqref{eq:D-18} for the actual pair.

Finally, use the endpoint parabolic rescaling
\[
 g_+(t')=2g(t'/2),\qquad
 f_+(t')=f(t'/2),\qquad
 \iota_+=\iota_{\log2}\circ P.
\]
The constructed maps are defined on a neighborhood of $\overline N_{S'}$,
and the endpoint identities give
\[
 \iota_+^*
   (g_+(-1),f_+(-1))=(G^+,F^+)
       \quad\hbox{on }N_{S'}.
\]
The factor $2^{n/2}$ in the volume form cancels the change in the
probability normalization, so this rescaling preserves the normalized
probability and introduces no additive potential constant.
Thus $\iota_+$ provides order-$100n$ endpoint comparison data
with the normalized actual potential.
\end{proof}

The endpoint construction is equivariant under every fixed center-preserving
reference isometry $Q$.
\phantomsection\label{par:endpoint-equivariance}
The reference metric and probability, radial cutoffs and kernel decomposition
are invariant, and orthogonal projection satisfies
\begin{equation}\label{eq:endpoint-projection-equivariance}
 \Pi_K\!\left[\chi_R(Q^*G-\bar g)\right]
 =Q^*\Pi_K\!\left[\chi_R(G-\bar g)\right].
\end{equation}
The gauge and extension constructions and the stationary representatives
are therefore equivariant. Consequently
$\iota_+\mapsto\iota_+\circ Q$ and $b_+\mapsto Q\cdot b_+$
preserve the endpoint estimates, radii and constants in the invariant
parameter norm. Pullback by $Q$ preserves both probability masses and the vanishing
auxiliary first moment.

\section{Entropy summation and radius control}\label{sec:radii}

This section chooses compatible parameters and comparison radii, then proves
summability and quantitative bounds for the entropy, defect integral, and radii
uniformly on admissible finite initial segments. Section~\ref{sec:motion}
uses these bounds to confine the compact parameter and establish continuation.

On blocks of length $\log2$, propagation along the Ricci flow from the order-$100n$ initial comparison data
increases the available radius, while the current entropy drop limits
the auxiliary comparison gauge. Their minimum leaves an inner endpoint
radius for stationary comparison at order $100n$. In Gaussian-tail
variables this minimum becomes a delayed maximum, so earlier drops
can still determine later radii.

\subsection{Choice of parameters}
\label{subsec:parameter-selection}
\label{part8-a-noncircular-set-of-strict-inequalities}

Choose the parameters below in the stated order. They remain fixed throughout
the continuation argument in Section~\ref{sec:motion} and the qualitative
convergence argument in Subsection~\ref{part9-one-fixed-embedding-not-a-collection-of-unrelated-gauges}.
For quantitative convergence, Subsection~\ref{subsec:quantitative-fixed-gauge}
makes a separate choice for each prescribed derivative order and rate,
using the same modified flow and fixed embedding.
The incoming coordinate order is $100n$ and the recurrence power is
$3/2$. Table~\ref{tab:iteration-data} summarizes the construction.

\begin{enumerate}[label=\textup{\arabic*.},ref=\arabic*,leftmargin=*]
\item \emph{Entropy exponents.}\label{stage:entropy-exponents}
Let $p$ be the stationary exponent in Theorem~\ref{thm:logarithmic}.
Choose $1<p_1<p$ sufficiently close to $p$, set
\[
 p_2:=p_1/p,\qquad p_4:=(p_2+2-p_1)/2,
\]
and choose $p_3<p_2$ sufficiently close to $p_2$ that
\begin{equation}\label{eq:E-1}
 \begin{gathered}
 0<p_4<p_3<p_2<1<p_1<p,\\
 p_3(p_1+p_4)>2,\qquad p_2+p_1p_3>2.
 \end{gathered}
\end{equation}
These choices are possible since, as $p_1\uparrow p$ and $p_3\uparrow1$,
$p_4\to(3-p)/2<1$ and the last two left-hand sides tend to
$(p+3)/2>2$ and $1+p>2$. Moreover,
$p_4=1-p_1(p-1)/(2p)>(3-p)/2\ge3/4$.

\item \emph{Radius ratios.}\label{stage:geometric-retentions}
Choose $\rho<1$ and $\Lambda>1$ sufficiently close to one that
\begin{equation}\label{eq:E-2}
 0<\rho<1,\qquad \rho^2p_3>p_4,\qquad
 \frac32<\frac{2\rho^2p_3}{\Lambda^2},\qquad \Lambda^2p_3<1.
\end{equation}
Indeed, after choosing $\rho$ with $\rho^2p_3>p_4>3/4$, the interval
$1<\Lambda^2<\min\{4\rho^2p_3/3,1/p_3\}$ is nonempty.
All intermediate radius ratios are the functions of $\rho$ fixed in
Subsection~\ref{subsec:physical-common-setup}.

\item \emph{Gaussian rates.}\label{stage:gaussian-rates}
Choose $0<c_2<1/8$ sufficiently close to $1/8$, and set $c_1=1/4-c_2$, so that
\begin{equation}\label{eq:E-3}
 c_1p_3<c_2p_2,\qquad c_1\Lambda^2p_3<c_2.
\end{equation}
At $c_2=c_1=1/8$ these inequalities follow from $p_3<p_2$ and
$\Lambda^2p_3<1$. Set $c_3=(c_1-1/8)/32$, $c_6:=(c_2+1/8)/2$, and define the
entropy-radius coefficient
\begin{equation}\label{eq:E-4}
 \Lambda':=\frac12\left(\frac{\Lambda^2p_3}{c_2}+\frac1{c_1}\right),
 \qquad \frac{\Lambda^2p_3}{c_2}<\Lambda'<\frac1{c_1},
 \qquad c_2\Lambda'>p_4.
\end{equation}
The bounds follow from~\eqref{eq:E-3} and $\Lambda^2p_3>p_4$.
Finally, set
\begin{equation}\label{eq:E-5}
 c_5:=\frac14(c_2p_2/p_3-c_1)>0,
 \qquad c_5<\frac12(c_2p_2/p_3-c_1).
\end{equation}

\item \emph{Derivative orders.}\label{stage:baseline-regularity}
Choose an even integer $Q_*$ such that
\[
 Q_*\ge1000n,\qquad
 p_3<\frac{Q_*+1}{Q_*+2+n/2}.
\]
Prescribe intrinsic comparison estimates
through order $Q_*+4$, stationary representatives of
Proposition~\ref{prop:cutoffjets} through order $Q_*$, and
incoming reference-coordinate regularity of order $100n$.

\item \emph{Parameters for the Ricci flow estimates.}\label{stage:comparison-package}
Fix $(\rho,c_1,Q_*)$ and the preceding derivative orders,
including reference-coordinate order $100n-2$. Proposition~\ref{prop:terminal-comparison}
then supplies $\eta\in(0,1/3)$ and $\epsilon_{\rm met}>0$, chosen in that order.
The comparison estimates of Section~\ref{sec:physical} hold uniformly
on a common neighborhood of zero. Choose an invariant closed ball
\[
 \overline B_{\rm cmp},\qquad B_{\rm cmp}:=B_{r_{\rm cmp}}^{K_N}(0),
\]
inside that neighborhood, using the reference $L^2(\bar\mu)$ norm.
Fix the common input-error threshold, minimum radius $S_*$, propagation
collar $C_1$, Gaussian boundary-decay rate $c_4>0$, and transition
collar $C_2\ge1$. These depend on the preceding comparison inputs; in particular,
\[
 r_{\rm cmp}=r_{\rm cmp}(N,g_N,d,Y;\rho,c_1,Q_*),
\]
with $\eta,\epsilon_{\rm met}$ fixed, independently of the later gauge choices.

\item \emph{Gauge loss and family inverse.}\label{stage:gauge-loss}
Choose $\beta>0$ so that
\begin{equation}\label{eq:E-6}
 4\beta\rho^2<\min\{c_5/4,c_4/2\},
 \qquad \Lambda'(c_1/2+2\beta)<1/2,
\end{equation}
and
\begin{equation}\label{eq:E-7}
 \beta<1-8c_6,\qquad
 c_2<c_6<(1-\beta)/8.
\end{equation}
This is possible because $c_4,c_5>0$, $\Lambda' c_1<1$ and $c_6<1/8$.
Apply Theorem~\ref{thm:gauge-inverse} at this $\beta$ and
$k=4$, and fix its constants and $Q_\beta$ as in
Proposition~\ref{prop:restart-gauge}.

\item \emph{Localized pairings and the outer stationary neighborhood.}
\label{stage:parameter-pairings}\label{par:parameter-neighborhoods}
Choose an orthonormal basis $\{e_i\}$ of $K_N$ in the metric
$L^2(\bar\mu)$ norm and a nonnegative smooth radial function
\[
 \bar\varphi\in C_c^\infty(B_1(0)),\qquad
 \int_{\mathbb R^d}\bar\varphi\,\mathrm d\gamma_d=1.
\]
Extend $e_i$ constantly in the Euclidean factor and set
\[
 T_i:=\bar\varphi e_i,\qquad
 \mathcal O(G)_i:=\int_{N_1}\langle G-\bar g,T_i\rangle_{\bar g}
                     \,\mathrm d\bar\mu.
\]
Identify $K_N$ with its coordinates in this basis. By
Lemma~\ref{lem:compact-family} and the normalization of $\bar\varphi$,
\[
 \mathcal O(G_b)_i
 =\langle g_b-g_N,e_i\rangle_{L^2(\mu_N)}=b_i.
\]
The fixed tensors $T_i$ give
$|\mathcal O(G)-\mathcal O(G')|\le C\|G-G'\|_{C^0(N_1)}$.
Choose $r_{\rm out}>0$ so that the closure of
$B_{\rm out}:=B_{r_{\rm out}}^{K_N}(0)$ lies in the stationary-family
domain at the prescribed orders and
$\{|b|\le r_{\rm out},\ |a|\le r_{\rm out}\}$ lies in the
representative domain of Proposition~\ref{prop:cutoffjets}.
These balls are invariant under the reference group.

\item \emph{Stationary accuracy and iteration neighborhoods.}
\label{stage:iteration-neighborhoods}
Let $\epsilon_{\rm w}$ be the tolerance in
Lemma~\ref{lem:augmented-weak-inverse}, and take $C_{\rm w}$ and
$C$ in \eqref{eq:weak-construction-coefficients} with $a=1/16$.
Choose $0<\epsilon_{\rm sl}\le\bar\epsilon_{\rm sl}(\Lambda,p,p_1)$
from Lemma~\ref{lem:restart-extension}, with its closed parameter ball
in the representative domain, and require
\[
 \epsilon_{\rm sl}<\frac12r_{\rm out},\qquad
 3C_{\rm w}\epsilon_{\rm sl}\le\epsilon_{\rm w}.
\]
Apply Proposition~\ref{prop:restart-gauge} at $(\beta,\epsilon_{\rm sl})$
and choose $r_{\rm gauge}>0$ such that
\[
 \overline B_{r_{\rm gauge}}^{K_N}(0)
       \subset B_{\rm gauge}(\beta,\epsilon_{\rm sl}),\qquad
 3C_{\rm w}r_{\rm gauge}\le\epsilon_{\rm w}.
\]
Set
\[
 r_{\rm con}:=\frac12\min\{r_{\rm out},r_{\rm cmp},r_{\rm gauge}\},
 \qquad B_{\rm con}:=B_{r_{\rm con}}^{K_N}(0).
\]
Then
\[
 \overline B_{\rm con}\subset
 B_{\rm out}\cap B_{\rm cmp}\cap B_{\rm gauge}.
\]
The projection bound
in Proposition~\ref{prop:restart-stationary-comparison} now implies
$|(b_+,a)|\le\epsilon_{\rm sl}<r_{\rm out}/2$, and therefore $b_+\in B_{\rm out}$.

\item \emph{Spatial and entropy thresholds.}\label{stage:iteration-thresholds}
Increase the minimum radius from Proposition~\ref{prop:restart-gauge} so that
\[
 R_{+,0}>\max\{10,C\},\qquad
 3C_{\rm w}e^{-(R_{+,0}-C)^2/16}\le\epsilon_{\rm w}.
\]
Stages~\ref{stage:iteration-neighborhoods}--\ref{stage:iteration-thresholds}
ensure the smallness condition in Lemma~\ref{lem:augmented-weak-inverse}
throughout the ball used in Proposition~\ref{prop:restart-gauge}.
Choose $R_{\min}\ge R_{+,0}$ and $0<E_*<1$ with
\[
 \begin{gathered}
 \rho\sqrt{p_3}R_{\min}-C_1-C_2\ge R_{+,0},\qquad
 Ce^{-c_5p_3R_{\min}^2}<r_{\rm con}/8,\\
 \frac{p_3}{c_2}\log\frac1{\sqrt{E_*}}\ge R_{\min}^2,
 \qquad
 \sqrt{\Lambda'\log\frac1{\sqrt{E_*}}}-C_2\ge R_{+,0}.
 \end{gathered}
\]
Here $C$ in the exponential bound is the pairing constant from
Stage~\ref{stage:parameter-pairings}. Increase $R_{\min}$ and then
decrease $E_*$ and the input-error threshold as required below
before choosing the starting chart.
\end{enumerate}

The structural parameters and neighborhoods depend only on the
background and the preceding choices, independently of spatial radii,
iteration indices and the number of blocks.

\begin{table}[htbp]
\centering\small
\renewcommand{\arraystretch}{1.15}
\begin{tabular}{@{}p{0.07\linewidth}p{0.35\linewidth}p{0.51\linewidth}@{}}
\toprule
Stage & Main inputs & Main outputs \\
\midrule
1 & $p$ & $p_1,p_2,p_3,p_4$ \\
2 & $p_3$ and the preceding margins & $\rho,\Lambda$ \\
3 & Exponents and geometric margins & $c_2,c_1,c_3,c_5,c_6,\Lambda'$ \\
4 & $p_3,n$ & $Q_*$ and prescribed derivative orders \\
5 & $\rho,c_1,Q_*$ & $\eta,\epsilon_{\rm met},r_{\rm cmp},c_4$ and comparison thresholds \\
6 & $c_4,c_5$ and Gaussian margins & $\beta,Q_\beta$ and family inverse \\
7 & Stationary family and pairings supported in $N_1$ & $\bar\varphi,T_i,\mathcal O,r_{\rm out},B_{\rm out}$ \\
8 & Gauge and accuracy inputs & $\epsilon_{\rm sl},r_{\rm gauge},r_{\rm con},B_{\rm con}$ \\
9 & All preceding choices & $R_{+,0},R_{\min},E_*$ and input thresholds \\
\bottomrule
\end{tabular}
\caption{The parameter choices for continuation and qualitative convergence.
Subsection~\ref{subsec:quantitative-fixed-gauge} makes a separate construction
with exponents, radius ratios, derivative orders and Gaussian losses chosen
for each prescribed derivative order and convergence rate.}
\label{tab:iteration-data}
\end{table}

\subsection{Choice of radii and endpoint comparison}\label{part8-physical-input-and-causal-low-gauge-radius}
The radius must be chosen from the part of the flow already constructed.
The following proposition verifies the gauge smallness condition on every interval $[0,u]$ of a block and collects the endpoint consequences of Theorem~\ref{thm:restart}.
Stage~\ref{stage:iteration-neighborhoods} places the endpoint
parameter in $B_{\rm out}$; Section~\ref{sec:motion} proves its
confinement in $B_{\rm con}$.

For the initial chart of radius $S_0$, set $R_{-1}:=S_0/\sqrt{p_3}$.
The incoming comparison chart on block $j$ has radius
\begin{equation}\label{eq:E-8}
 S_j:=\sqrt{p_3}\,R_{j-1}.
\end{equation}
Write $\boldsymbol{\epsilon}_j:=\boldsymbol{\epsilon}(S_j)$ for its error
in~\eqref{eq:seed-error}. The incoming data satisfy
\begin{equation}\label{eq:E-9}
 \boldsymbol{\epsilon}_j+|\mathrm De(b_j)|
 \le e^{-c_5S_j^2}.
\end{equation}
Propagation along the Ricci flow starts from the given order-$100n$
comparison chart, before the low-order spatial gauge is applied.
Section~\ref{subsec:comparison-length} specifies how such charts are
inherited in the recursive construction.
The block is rescaled to $[-1,-1/2]$; the underlying closed flow
supplies the required interval back to $-2$.

\phantomsection\label{idx:prefix-drop}
For $s_j:=s_0+j\log2$, write $E_j:=E(s_j)$ and, for
$0\le u\le\log2$, define
\[
 A_j(u)^2:=E_j-E(s_j+u)
          =\mathcal W(s_j+u)-\mathcal W(s_j).
\]
Thus $A_j=A_j(\log2)$ and $E_j-E_{j+1}=A_j^2$.
All drops use the actual terminal potential and are nonnegative by
entropy monotonicity. We use the convention $\log(1/0)=+\infty$.
There are two restrictions on the comparison radius at time $u$.
Propagation along the Ricci flow makes $e^{u/2}(\rho S_j-C_1)$ available. The
current-drop term in the gauge smallness condition requires the entropy cap
$\sqrt{\Lambda'\log(1/A_j(u))}$, by~\eqref{eq:E-11}.
Both depend only on the constructed segment. Reserve the final fixed
transition collar after taking their minimum:
\begin{equation}\label{eq:E-10}
 L_j(u):=\min\left\{
 e^{u/2}(\rho S_j-C_1),\,
 \sqrt{\Lambda'\log\frac1{A_j(u)}}\right\}-C_2.
\end{equation}
The propagation branch increases and the
entropy branch decreases, so there is at most one switch. The resulting
radius is locally absolutely continuous: while its entropy branch is
active, $A_j(u)$ is bounded away from zero and is at most
$\sqrt{E_*}<1$, so that branch is an absolutely continuous function
of the accumulated drop. The branches agree at a switch.
Derivatives of moving cutoffs and the estimates below are
understood almost everywhere at a switch.

\begin{proposition}
\label{prop:causal-radii}
Fix $j\ge0$ and use the data of Subsection~\ref{subsec:parameter-selection}.
After increasing $R_{\min}$ and decreasing $E_*$ and the input-error
threshold in terms of these data, the following holds.
Suppose $b_j\in B_{\rm con}$,
$R_{j-1}\ge R_{\min}$, $0\le E_j\le E_*$, and the initial comparison chart satisfies
\eqref{eq:E-8}--\eqref{eq:E-9}. Then:
\begin{enumerate}[label=\textup{(\roman*)},leftmargin=*]
\item For every $u\in[0,\log2]$, the raw comparison on $N_{L_j(u)}$
      satisfies the low-gauge smallness condition~\eqref{eq:D-3}, uniformly in $u$. The resulting chart has the exact slice and centered
      auxiliary probability of Proposition~\ref{prop:restart-gauge}.
\item At a block endpoint retain
\begin{equation}\label{eq:E-12}
 R_j^2:=\min\left\{\frac32R_{j-1}^2,
          \frac{p_3}{c_2}\log\frac1{A_j}\right\}.
\end{equation}
      Then $R_j\ge R_{\min}$. With $\Delta_j:=e^{-c_2R_j^2}$,
      including $\Delta_{-1}=e^{-c_2R_{-1}^2}$, the one-step recurrence is
\begin{equation}\label{eq:E-14}
 \Delta_j=\max\{\Delta_{j-1}^{3/2},A_j^{p_3}\}.
\end{equation}
\item Put $R_+=L_j(\log2)$, $S_{j+1}=\sqrt{p_3} R_j$ and
      $L_{j+1}(0)=\rho S_{j+1}-C_1-C_2$, as determined by~\eqref{eq:E-10}
      at zero initial drop. These radii satisfy
\begin{equation}\label{eq:E-13}
 \Lambda R_j+10\le L_j(\log2)=R_+,
\end{equation}
      as well as $L_{j+1}(0)\le S_{j+1}-1$ and $S_{j+1}+4\le R_j$.
      If the next compact parameter belongs to $B_{\rm con}$ and its
      own gauge smallness condition holds, its low-order gauge change satisfies
\begin{equation}\label{eq:D-19}
 \operatorname{supp}(P_{\rm new}-\Id)
 \subset N_{L_{j+1}(0)-1}\Subset N_{S_{j+1}}
 \Subset N_{R_j}\Subset N_{R_+-4}.
\end{equation}
      The outer extension cutoff and every displaced cutoff support
      remain in the old exact-slice and auxiliary-measure agreement region.
\item The block endpoint gives $b_+\in B_{\rm out}$ and comparison
      data $\iota_+$ on a neighborhood of $\overline N_{S_{j+1}}$ at
      order $100n$, with the actual potential constant and
\[
 \boldsymbol{\epsilon}_+(S_{j+1})+|\mathrm De(b_+)|
       \le e^{-c_5S_{j+1}^2}.
\]
      The endpoint entropy gap satisfies
\begin{equation}\label{eq:E-15}
 E_{j+1}^{1/p_1}\le C\Delta_j.
\end{equation}
\end{enumerate}
The endpoint estimates also hold after any fixed center-preserving
reference isometry, by~\eqref{eq:endpoint-projection-equivariance}.
The constants and thresholds are independent of $j$ and are uniform
on every constructed part of this block, including a final partial block.
\end{proposition}
\begin{proof}
Apply Proposition~\ref{prop:actual-comparison} to the incoming chart.
On $N_{L_j(u)}$ it gives~\eqref{eq:D-2} with $L=L_j(u)$ and current
drop $A_j(u)$, before the low-order gauge change. To apply
Proposition~\ref{prop:restart-gauge}, we verify~\eqref{eq:D-3}.

From~\eqref{eq:E-10} and $C_2\ge1$,
\[
 L+1\le e^{u/2}(\rho S_j-C_1)\le\sqrt2\rho S_j.
\]
Thus $\mathcal C_\beta(L)^2$ is bounded by a fixed polynomial in $S_j$
times $e^{4\beta\rho^2S_j^2}$. The initial-error contribution is at most
\[
 C(1+S_j)^Q
 e^{-(c_5/2-4\beta\rho^2)S_j^2}.
\]
The reduced-gradient and Gaussian boundary-error contributions are bounded respectively by
\[
 C(1+S_j)^Q e^{-(c_5-4\beta\rho^2)S_j^2},
 \qquad
 C(1+S_j)^Q e^{-(c_4-4\beta\rho^2)S_j^2}.
\]
By~\eqref{eq:E-6}, all three exponents are positive. Increasing $R_{\min}$
makes each contribution at most $1/800$.

For $A=A_j(u)>0$, the entropy branch gives
$L+1\le\sqrt{\Lambda'\log(1/A)}$. The term involving $A$ in~\eqref{eq:D-2} is at most
$C(1+L)^Qe^{c_1(L+1)^2/2}A^{1/2}$. By~\eqref{eq:D-1},
$\mathcal C_\beta(L)^2=(1+L)^{2Q_\beta}e^{2\beta(L+1)^2}$.
Multiplying gives the bound
\begin{equation}\label{eq:E-11}
 C\bigl(1+\log(1/A)\bigr)^Q
 A^{1/2-\Lambda'(c_1/2+2\beta)}.
\end{equation}
Its power
is positive by~\eqref{eq:E-6}; since $A_j(u)\le\sqrt{E_*}$,
decreasing $E_*$ makes this contribution
at most $1/800$, uniformly in the segment. At zero drop it is exactly
zero. Thus the four contributions satisfy
\[
 \mathcal C_\beta(L)^2\epsilon_+\le1/200<1/100.
\]
Choose the same thresholds so that
\[
 \rho\sqrt{p_3}R_{\min}-C_1-C_2\ge R_{+,0},\qquad
 \sqrt{\Lambda'\log(1/\sqrt{E_*})}-C_2\ge R_{+,0}.
\]
Both branches in~\eqref{eq:E-10}, after reserving the collar $C_2$,
then give $L_j(u)\ge R_{+,0}>10$. Thus the coefficient accuracy holds
throughout every current segment, and
$N_4\subset N_{L_j(u)-4}$ lies in the exact reference slice.
In particular the fixed pairing core $N_1$ has the $N_3$ collar
used for the localized pairings in Proposition~\ref{prop:comparison-length}.
This proves (i).

The choice of $E_*$ ensures
$(p_3/c_2)\log(1/\sqrt{E_*})\ge R_{\min}^2$. Since $A_j^2\le E_j\le E_*$
and $3/2>1$, both entries of the minimum in~\eqref{eq:E-12} are at least
$R_{\min}^2$. Exponentiating that minimum gives~\eqref{eq:E-14}. This includes zero
drops under the convention $\log(1/0)=+\infty$. In particular
$A_j^{p_3}\le\Delta_j$, exactly the drop condition required by
Theorem~\ref{thm:restart}.

The propagation and entropy branches must both contain the proposed
stationary extension. For the propagation branch use
$R_j\le\sqrt{3/2} R_{j-1}$ and the strict gap
\[
 \sqrt2\rho\sqrt{p_3}-\Lambda\sqrt{3/2}>0.
\]
Multiplying this gap by $R_{j-1}\ge R_{\min}$ absorbs
$\sqrt2C_1+C_2+10$. Hence
\[
 \Lambda R_j+10
 \le\sqrt2(\rho S_j-C_1)-C_2.
\]
For a nonzero drop, the entropy branch uses
$R_j\le\sqrt{p_3/c_2}\sqrt{\log(1/A_j)}$ and
$\sqrt{\Lambda'}-\Lambda\sqrt{p_3/c_2}>0$. Increasing the lower
threshold for $\log(1/A_j)$ absorbs $C_2+10$ and gives
\[
 \Lambda R_j+10
 \le\sqrt{\Lambda'\log(1/A_j)}-C_2.
\]
For zero drop the second branch is infinite. Taking the minimum
proves~\eqref{eq:E-13}.

Since $p_3<1$ and $\rho<1$, a further fixed increase of $R_{\min}$
gives $S_{j+1}+4\le R_j$ and $L_{j+1}(0)\le S_{j+1}-1$,
and puts $L_{j+1}(0)$ above the gauge minimum radius. Indeed,
$L_{j+1}(0)=\rho\sqrt{p_3}R_j-C_1-C_2\ge R_{+,0}$
by $R_j\ge R_{\min}$ and the preceding choice of $R_{\min}$.
Conditional on the next parameter and its smallness condition, the low construction
has support in $N_{L_{j+1}(0)-1}$ and displacement less than $1/10$
throughout its isotopy. The three radial gaps in~\eqref{eq:D-19} are
at least $2$, $4$, and $6$. Moreover the outer extension cutoff is
supported in $N_{L_{j+1}(0)}$, with a collar of width at least one
inside $N_{S_{j+1}}$. The displacement bound, also valid for inverse
images, leaves positive margins for every cutoff support inside
$N_{R_+-4}$. This is precisely the old exact-slice region where
$\nu'=(Z')^{-1}\mu^+$.

Theorem~\ref{thm:restart} now gives the order-$100n$ endpoint chart
$\iota_+$ with its actual potential constant.
Its parameter belongs to $B_{\rm out}$ by
Stage~\ref{stage:iteration-neighborhoods}. For $S_{j+1}=\sqrt{p_3}R_j$,
\[
 \exp\!\left(\frac{c_1S_{j+1}^2-c_2p_2R_j^2}{2}\right)
 =\exp\!\left[-\frac12\left(\frac{c_2p_2}{p_3}-c_1\right)
 S_{j+1}^2\right]=e^{-2c_5S_{j+1}^2}.
\]
The choice of $c_5$ in
Stage~\ref{stage:gaussian-rates}, applied to~\eqref{eq:D-18}, gives
\[
 \boldsymbol{\epsilon}_+(S_{j+1})+|\mathrm De(b_+)|
 \le C(1+S_{j+1})^Q
 e^{-\frac12(c_2p_2/p_3-c_1)S_{j+1}^2}
       +Ce^{-(c_2p_2/p_3)S_{j+1}^2}
 \le e^{-c_5S_{j+1}^2}.
\]
The strict margin~\eqref{eq:E-5} absorbs both the polynomial and the
sum, uniformly for all sufficiently large $R_j$. This proves the
endpoint comparison assertion at order $100n$. The endpoint entropy estimate in
Theorem~\ref{thm:restart} gives~\eqref{eq:E-15}.
\end{proof}

\subsection{Summation on every finite initial segment}
\label{part9-current-prefix-length-and-delayed-summation}

The following estimate adapts the entropy-power summation of
\cite[Corollary~6.28]{FLC} to the inherited radius restriction.

\begin{lemma}[Summability of the delayed recurrence]\label{lem:delayed-sum}
For any compatible exponents fixed in
Stage~\ref{stage:entropy-exponents} and any $C_0>0$, there are
$\epsilon_*,\vartheta>0$ and $C$, depending only on these data,
with the following property. Suppose nonnegative finite sequences
satisfy, for $0\le j<J$,
\[
 E_j-E_{j+1}=A_j^2,\qquad
 \Delta_j=\max\{\Delta_{j-1}^{3/2},A_j^{p_3}\},\qquad
 E_{j+1}\le C_0\Delta_j^{p_1}.
\]
If $E_0+\Delta_{-1}\le\epsilon_*$, then
\begin{equation}\label{eq:M-14}
 \sum_{j=0}^{J-1}(A_j^{p_4}+\Delta_j^{p_4})
 \le C(E_0^\vartheta+\Delta_{-1}^{p_4}).
\end{equation}
The constants are independent of $J$.
\end{lemma}
\begin{proof}
Put $q=p_1p_3$ and $\vartheta=1-(2-p_3p_4)/q$. Stage~\ref{stage:entropy-exponents}
gives $q<2$ and $0<\vartheta<1$. Since $A_j^2\le E_0$, we may
choose $\epsilon_*$ so that $A_j\le1$ and
$\Delta_j\le\max\{\Delta_{-1},E_0^{p_3/2}\}\le1/2$ throughout.

At each index with $\Delta_j=A_j^{p_3}$, including ties, the endpoint
bound gives $E_j=E_{j+1}+A_j^2\le(C_0+1)A_j^q$. When $A_j>0$,
concavity and $q(\vartheta-1)+2=p_3p_4$ imply
\[
 E_j^\vartheta-E_{j+1}^\vartheta
 \ge\vartheta E_j^{\vartheta-1}A_j^2
 \ge c\Delta_j^{p_4}.
\]
A zero drop contributes zero. Since $E_j$ is nonincreasing,
the sum of $\Delta_j^{p_4}$ over these indices is at most $CE_0^\vartheta$.

On each intervening string of indices the radius is inherited.
Let $i$ be the index just before the string, taking $i=-1$ for an
initial string. Then $\Delta_{i+h}=\Delta_i^{(3/2)^h}$ throughout
the string. Since $(3/2)^h\ge1+h/2$ and $\Delta_i\le1/2$, its sum is
bounded by
\[
 \sum_{h\ge1}\Delta_i^{p_4(3/2)^h}
 \le\Delta_i^{p_4}\sum_{h\ge1}2^{-hp_4/2}
 \le C\Delta_i^{p_4}.
\]
Each preceding index occurs at most once, including when the final
string ends at $J-1$. Adding the initial contribution
$C\Delta_{-1}^{p_4}$ and using $A_j\le\Delta_j$ proves~\eqref{eq:M-14}.
\end{proof}

\subsection{Quantitative decay and radius bounds}
\label{subsec:quantitative-iteration}

To obtain polynomial decay, we add the preceding Gaussian tail to the entropy. We then use
the entropy-to-length calculation of \cite[Corollary~6.28]{FLC}
and \cite[Proposition~6.19]{LYZII}.

\begin{lemma}[Quantitative decay]\label{lem:quantitative-decay}
Fix $p_1>0$, $0<p_3<1$ and $C_0>0$ with $q:=p_1p_3\in(1,2)$.
For nonnegative finite sequences satisfying the three recurrences in
Lemma~\ref{lem:delayed-sum}, there is $\epsilon_*>0$, depending only
on $p_1,p_3,C_0$, such that $E_0+\Delta_{-1}\le\epsilon_*$ implies
\begin{equation}\label{eq:quantitative-entropy-decay}
 E_j\le F_j:=E_j+\frac12\Delta_{j-1}^{2/p_3}
 \le C(1+j)^{-q/(2-q)},\qquad 1\le j\le J.
\end{equation}
The constant depends only on $p_1,p_3,C_0$. Suppose also that a
nonnegative absolutely continuous
$E(s)$ on blocks $s_j=s_0+j\log2$ satisfies $E(s_j)=E_j$ and
$E'=-2\mathcal D^2$, allowing a final partial block. For every
$2-q<a<1$ and every finite terminal time $S$ in this interval,
\begin{equation}\label{eq:physical-length-rate}
 \int_s^S\mathcal D(u)^a\,\mathrm du
 \le C_a(1+s)^{-(a+q-2)/(2-q)}
\end{equation}
at sufficiently late modified times $s\le S$. The constants are
independent of the number of blocks and the final partial-block
duration; the waiting time may depend on $s_0$. After continuation,
the same bound holds with $S=\infty$.
\end{lemma}
\begin{proof}
Set $F_0=E_0+\frac12\Delta_{-1}^{2/p_3}$. Since $A_j^2\le E_0$,
all tails are bounded by
$\max\{\Delta_{-1},E_0^{p_3/2}\}$. Choose $\epsilon_*$ small enough
that $\Delta_{j-1}^{1/p_3}\le1/2$ for every $j$ and $F_0\le1$.
The maximum recurrence gives
$\Delta_j^{2/p_3}\le\Delta_{j-1}^{3/p_3}+A_j^2$. Hence
\begin{equation}\label{eq:augmented-entropy-drop}
 F_j-F_{j+1}
 \ge\frac12A_j^2+\frac12
       (\Delta_{j-1}^{2/p_3}-\Delta_{j-1}^{3/p_3})
 \ge\frac12A_j^2+\frac14\Delta_{j-1}^{2/p_3}.
\end{equation}
For $j\ge1$, the preceding endpoint estimate and $p_1<2/p_3$ give
$F_j\le(C_0+1/2)\Delta_{j-1}^{p_1}$. Together
with~\eqref{eq:augmented-entropy-drop} this gives
$F_j-F_{j+1}\ge cF_j^{2/q}$ for $1\le j<J$.
If $F_{j+1}>0$, integration yields
\[
 F_{j+1}^{1-2/q}-F_j^{1-2/q}
 =\left(\frac2q-1\right)\int_{F_{j+1}}^{F_j}x^{-2/q}\,\mathrm dx
 \ge c\left(\frac2q-1\right).
\]
Telescoping from $j=1$ and using the uniform bound on $F_1$ proves
\eqref{eq:quantitative-entropy-decay}. If some $F_j=0$, all subsequent
quantities vanish, and the conclusion is immediate.

For a block of duration $u\le\log2$, H\"older's inequality gives
\[
 \int_{s_j}^{s_j+u}\mathcal D(s)^a\,\mathrm ds
 \le u^{1-a/2}
   \left(\frac{E_j-E(s_j+u)}2\right)^{a/2}.
\]
In particular, a completed block contributes at most $C_aA_j^a$.
For $n\ge1$, another application of H\"older gives
\[
 \sum_{j=n}^{2n-1}A_j^a
 \le n^{1-a/2}\left(\sum_{j=n}^{2n-1}A_j^2\right)^{a/2}
 \le n^{1-a/2}E_n^{a/2}
 \le C_a n^{-(a+q-2)/(2-q)}.
\]
The same bound holds for a truncated group, including its last
partial block: there are at most $n$ terms, and the sum of their
drops is at most $E_n$. Summing the groups starting at
$n,2n,4n,\ldots$ is possible because $a+q-2>0$, with a constant
independent of the terminal time. For $s\in[s_n,s_{n+1}]$, start at
$s_n$. Once $s\ge2|s_0|+2$, $1+n$ and $1+s$ are comparable by fixed
factors. This proves~\eqref{eq:physical-length-rate}. Monotone
convergence permits $S\to\infty$ after continuation has been proved.
\end{proof}

\begin{lemma}\label{lem:quantitative-radii}
Consider a finite initial sequence of blocks from Proposition~\ref{prop:causal-radii}
with $E_0+\Delta_{-1}\le\epsilon_*$ as in Lemma~\ref{lem:quantitative-decay}.
With $q=p_1p_3$, one has
\begin{equation}\label{eq:quantitative-endpoint-radius}
 R_{j-1}^2\ge
 \frac{p_3q}{2c_2(2-q)}\log(1+j)-C.
\end{equation}
For every time $u$ reached in block $j$, including a final partial block,
\begin{equation}\label{eq:quantitative-causal-radius}
 L_j(u)^2\ge
 \frac{\rho^2p_3^2q}{2c_2(2-q)}\log(1+j)
                  -C\sqrt{\log(2+j)}.
\end{equation}
Fixed collars and reset displacements preserve the leading
coefficient. A constant remainder is available after strictly
decreasing that coefficient.
\end{lemma}
\begin{proof}
By~\eqref{eq:quantitative-entropy-decay},
\[
 \Delta_{j-1}^{2/p_3}\le2F_j
 \le C(1+j)^{-q/(2-q)}.
\]
Taking logarithms and using
$R_{j-1}^2=c_2^{-1}\log(1/\Delta_{j-1})$ proves
\eqref{eq:quantitative-endpoint-radius}. Also $A_j(u)^2\le E_j$ gives
\[
 \log\frac1{A_j(u)}
 \ge\frac{q}{2(2-q)}\log(1+j)-C,
\]
with the left side infinite at a zero drop. Since
$S_j=\sqrt{p_3}R_{j-1}$ and $e^{u/2}\ge1$, the propagation and
entropy branches in~\eqref{eq:E-10} have respective squared leading
coefficients at least
\[
 \frac{\rho^2p_3^2q}{2c_2(2-q)},\qquad
 \frac{\Lambda' q}{2(2-q)}.
\]
The second is larger, since
$c_2\Lambda'>\Lambda^2p_3>\rho^2p_3^2$. Thus both radii are bounded
below by the square root of the first coefficient times
$\sqrt{\log(1+j)}$, minus a fixed constant. Taking their minimum,
subtracting $C_2$, and squaring gives
\eqref{eq:quantitative-causal-radius}; increasing $C$ covers the
bounded initial range.

At a junction,
$L_{j+1}(0)=\rho\sqrt{p_3}R_j-C_1-C_2$ has the same bound.
The common collars~\eqref{eq:D-19} and displacement below $1/10$
subtract only another fixed constant from the radius, so the estimate
holds throughout the reset path. Such subtractions produce
$O(\sqrt{\log(2+j)})$ in the squared radius. A strictly smaller
logarithmic coefficient absorbs this remainder above a fixed
threshold.
\end{proof}

Section~\ref{sec:motion} combines these bounds with chart-motion
estimates to confine the compact parameter.

\section{Chart motion and continuation}\label{sec:motion}
Weak motion of the comparison charts controls finitely many localized
metric pairings. We first normalize the charts on one block, then
compare the gauges at one endpoint. Summing these estimates on finite
sequences of blocks confines the compact parameter and gives continuation.

\subsection{Comparison charts on one block}
\label{part9-exact-transport-of-the-constructed-raw-charts}
Fix the construction data of Subsection~\ref{subsec:parameter-selection},
an index $j$, and the interval $I_j=[s_j,s_{j+1}]$, where
$s_{j+1}=s_j+\log2$. Assume that the
incoming order-$100n$ comparison chart has parameter $b_j\in B_{\rm con}$
and satisfies the hypotheses of Proposition~\ref{prop:causal-radii}.
Use its fixed center and reference coordinates; $b_j$ is constant in
these raw coordinates throughout the propagation in Section~\ref{sec:physical}.
Put $u:=s-s_j$ and $L=L_j(u)$.

On the fixed manifold $M$ carrying the modified flow, \eqref{eq:modified-diffeomorphism} gives
\begin{equation}\label{eq:M-1}
 \partial_sg^{z_*}=2\Phi,\qquad
 \partial_sf^{z_*}=\Tr_{g^{z_*}}\Phi,\qquad
 \partial_s\mathrm d\mu^{z_*}=0,\qquad E'=-2\mathcal D^2.
\end{equation}
Use the parabolic normalization of this block.
Write $g_s,f_s,\mu_s,\Phi_s$ for the quantities in
Proposition~\ref{prop:actual-comparison} at elapsed time $u$. For
$\iota_s(y,\vec x)=\iota(\psi_{b_j}(u)y,e^{-u/2}\vec x)$, put
$J_s^{\rm raw}:=\Xi_s^{-1}\iota_s$. Differentiating gives
\begin{equation}\label{eq:M-2}
 \begin{gathered}
 W_s:=V_s+\nabla^{g_s}f_s,\qquad
 V_s:=((\mathrm D\psi_{b_j}(u))^{-1}\partial_u\psi_{b_j}(u),-\vec x/2),\\
 \dot g_s=2\Phi_s+L_{W_s}g_s,
 \qquad \partial_s\mu_s=L_{W_s}\mu_s.
 \end{gathered}
\end{equation}
The compact part of $V_s$ is the negative auxiliary gradient. The
Euclidean part cancels the model gradient $\vec x/2$; spatial scalar
normalizations do not affect either gradient.

\begin{proposition}[Weak chart motion on one block]\label{prop:coherent-maps}
There is a piecewise smooth path
\[
 R_s\in\operatorname{Isom}(N,g_N)\times O(d),
 \qquad R_{s_j}=\Id,
\]
such that, with $P_s^{(j)}=P_{U_s}$, the normalized charts
\[
 J_s^{(j)}
 =\bigl(J_s^{{\rm raw},(j)}\circ P_s^{(j)}\bigr)\circ R_s^{(j)}
\]
exist throughout $I_j$. Suppress the block superscript within this block.
Their auxiliary probabilities have mass one and zero Euclidean first
moment. Put $Y_s:=(\mathrm DR_s)^{-1}\dot R_s\in\mathcal K_0$,
$U_s^R=R_s^*U_s$ and $W_s^R=R_s^*W_s+Y_s$.
The global extension
\[
 Z_s^R=B_{U_s^R}\dot U_s^R+\chi_2P_{U_s^R}^*W_s^R
\]
agrees with $(\mathrm DJ_s)^{-1}\partial_sJ_s$ on $N_{L-4}$ and satisfies
\begin{equation}\label{eq:M-9}
 \Pi_{\mathcal K_0}(B_{U_s^R}^{-1}Z_s^R)=0,\qquad
 \|Z_s^R\|_{H^1(\bar\mu)}
 \le C\mathcal D(s)+Ce^{-c_2L^2}(1+|\dot L|).
\end{equation}
The constants depend only on the fixed construction data, uniformly
in $j$ and the duration of an initial subinterval of $I_j$.
The maps are smooth in space, piecewise smooth in time, and continuous
at radius switches, where differentiated identities hold almost everywhere. The reference background in these coordinates is
\[
 b(s)=R_s^{(j)}\cdot b_j\in B_{\rm con},\qquad b(s_j)=b_j.
\]
\end{proposition}

The map $P_{U_s}$ gives the spatial gauge, and composition with $R_s$
normalizes the chart by reference isometries. The parameter $b(s)$ transports the
incoming background; the stationary endpoint parameter is defined
in Subsection~\ref{subsec:inherited-endpoint-coordinates}.
We first estimate the metric, densities and velocity before this normalization.
\begin{lemma}\label{lem:motion-low-comparison}
At a fixed time $s\in I_j$, put
\[
 L=L_j(s-s_j),\qquad G=g_s,\qquad F=f_s,\qquad P=P_U,
 \qquad G^+=P^*G,\quad\mu^+=P^*\mu_s,\quad
 \nu'=P^*\mu'.
\]
Use the cutoffs of Proposition~\ref{prop:restart-gauge} and set
$\Omega_L=N_{L-5/2}$. Then
\[
 \operatorname{supp}\chi_3\subset\Omega_L,\qquad
 P(\Omega_L)\subset N_{L-1},\qquad
 \nu'=(Z')^{-1}\mu^+\quad\text{on }\Omega_L.
\]
For the block parameter $b=b_j$, one has
\[
  \|G^+-\bar g\|_{C^1(\Omega_L)}
 +\left\|\log\frac{\mathrm d\mu^+}{\mathrm d\bar\mu}
                   \right\|_{L^\infty(\Omega_L)}
 +
 \left\|\log\frac{\mathrm d\nu'}{\mathrm d\bar\mu}
                   \right\|_{L^\infty(\bar M)}
 \le C(|b|+\epsilon_{\rm sl}).
 \]
Consequently the metric and actual density on $\Omega_L$, and the
auxiliary probability on $\bar M$, are uniformly equivalent to their
reference counterparts. Before normalization by reference isometries,
\[
 \sum_{k=0}^2|\bar\nabla^kW_s|_{\bar g}\le C(1+r)
       \quad\text{on }N_L,\qquad r=|\vec x|.
\]
The constant depends only on the background, fixed comparison inputs,
$\beta,Q_\beta$ and cutoff profiles. It is independent of
$\epsilon_{\rm sl},L,b,j$ and the duration of the current segment.
The gauge-radius threshold $R_{+,0}$ may depend on
$\epsilon_{\rm sl}$.
\end{lemma}
\begin{proof}
Propositions~\ref{prop:causal-radii} and~\ref{prop:restart-gauge} give
\[
 \|G-G_b\|_{C^8(N_L)}+\|F-F_b\|_{C^8(N_L)}\le \epsilon_+,
 \qquad \mathcal C_\beta(L)^2\epsilon_+\le\frac1{100},
\]
and
\[
 \|U\|_{C^{5,1/2}(N_L)}
       \le2\mathcal C_\beta(L)\epsilon_+\le\frac1{50\mathcal C_\beta(L)}.
\]
The minimum radius already chosen in Proposition~\ref{prop:restart-gauge}
and~\eqref{eq:gauge-prefactor-smallness} give
\[
 \epsilon_+\le\frac1{100\mathcal C_\beta(L)^2}\le\epsilon_{\rm sl},
 \qquad
 \frac{C(1+L)^Q}{\mathcal C_\beta(L)}\le C\epsilon_{\rm sl}
\]
for the fixed polynomial costs below and all $L\ge R_{+,0}$.

The support of $\chi_3$ lies in $N_{L-3}$, with a collar of width
$1/2$ inside $\Omega_L$. The displacement bound $|P_E-\vec x|<1/10$
gives $P(\Omega_L)\subset N_{L-1}$. Thus $\chi_0(P(x))=1$ on
$\Omega_L$, and the definition
\eqref{eq:gauge-completion} proves the asserted measure identity
on $\Omega_L$. Ordinary pullback estimates, the bounds uniform in $b$
and $\|G_b-\bar g\|_{C^1}\le C|b|$ give
\[
 \|P^*G-\bar g\|_{C^1(\Omega_L)}
 \le C\bigl(|b|+\epsilon_++\|U\|_{C^2(N_L)}\bigr)
 \le C(|b|+\epsilon_{\rm sl}).
\]

The pullback formula~\eqref{eq:gauge-pullback-density} and the global
density estimate in Proposition~\ref{prop:restart-gauge} give
\[
 \left\|\log\frac{\mathrm d\nu'}{\mathrm d\bar\mu}\right\|_\infty
 \le C(|b|+\epsilon_{\rm sl}).
\]
Since $\mu^+=Z'\nu'$ on $\Omega_L$, it suffices to bound $Z'$. The $C^8$ comparison gives
$\|\log(\mathrm d\mu_s/\mathrm d\mu_b)\|_{C^0(N_L)}\le C\epsilon_+$,
and hence $|\log Z'|\le C\epsilon_+$ by~\eqref{eq:gauge-completion}.
Since $\epsilon_+\le\epsilon_{\rm sl}$, it follows that
$|\log(\mathrm d\mu^+/\mathrm d\bar\mu)|\le C(|b|+\epsilon_{\rm sl})$
on $\Omega_L$.

Finally, \eqref{eq:M-2} gives $W=V+G^{-1}\mathrm dF$.
The raw $C^8$ bound gives uniform bounds for $G^{-1}$ and its first
two reference derivatives on $N_L$. Since
$F_b=f_b(y)+r^2/4$, the same bound yields
\[
 |\mathrm dF|\le C(1+r),\qquad
 |\bar\nabla\mathrm dF|+|\bar\nabla^2\mathrm dF|\le C.
\]
The product rule thus bounds the first two spatial derivatives of
$G^{-1}\mathrm dF$ by $C(1+r)$. The compact component of $V$ has
uniformly bounded spatial derivatives by the compact comparison
construction; its Euclidean component is $-\vec x/2$.
This proves the velocity estimate.
\end{proof}

\subsection{Normalization and weak chart motion}
\label{part9-differentiate-the-patched-slice-before-taking-divergence}
\label{part9-reference-symmetries-normalize-the-total-generator}
\label{part9-the-reference-weak-estimate-and-its-high-local-upgrade}
Near a time $s_*$ on either smooth branch of $L$, use the fixed space
$H^2(\mu_\beta)\cap C^{5,1/2}(N_{L(s_*)})$ with the sum norm of
Proposition~\ref{prop:endpoint-gauge}. The cutoff supports and their
displaced arguments have a common collar. Translating the cutoff profiles
smoothly with $L$, the pullbacks in~\eqref{eq:gauge-nonlinear-expansions}
give a $C^1$ source map into $C^{4,1/2}$ and the moment space by the
raw $C^8$ bound and smooth time dependence. The augmented inverse varies
smoothly with its cutoff moment row by the inverse operator identity.
Since $\|\mathrm D_U\mathcal T_{s_*}\|<1/2$, the implicit function theorem
gives $C^1$ dependence of $U_s$ on each branch. At a switch the fixed
points agree by continuity and uniqueness.

For the unnormalized moving chart $J^{\rm raw}_s\circ P_{U_s}$,
the chain rule gives
\[
 Z^{\rm true}:=P_{U_s}^*W_s+(\mathrm DP_{U_s})^{-1}\partial_sP_{U_s}.
\]
Use the cutoffs and multipliers defined before
Lemma~\ref{lem:augmented-weak-inverse}. Their time derivatives satisfy
$\partial_s\chi_i=(\partial_L\chi_i)\dot L$, with spatial derivatives
bounded by $C|\dot L|$. On $N_{L-4}$, the displacement cutoff is
identically one on a spacetime neighborhood, so
$Z^{\rm true}=P_U^*W+A_U\dot U$. The global extension rule is
\begin{equation}\label{eq:M-4}
 Z:=B_U\dot U+\chi_2P_U^*W,\qquad
 v:=B_U^{-1}(\chi_2P_U^*W),\qquad \dot U=B_U^{-1}Z-v.
\end{equation}
Thus $Z=Z^{\rm true}$ on $N_{L-4}$, while~\eqref{eq:M-4}
defines $Z$ on all of $\bar M$.

The reference pairings, projections and orthogonal complements in
Lemma~\ref{lem:compact-family} are invariant under center-preserving
reference isometries, and $R^*X_b=X_{R\cdot b}$.
The reference exponential and radial cutoffs give
\[
 P_{R^*U}(x)
 =\exp_x\bigl(\chi_1(x)(\mathrm DR_x)^{-1}U(Rx)\bigr)
 =(R^{-1}\circ P_U\circ R)(x).
\]
Uniqueness of the constructed gauge therefore preserves equivariance
and the reference symmetry conditions. Differentiating this identity
along the reference flow generated by $Y\in\mathcal K_0$ gives
\[
 \mathrm D_UP_U[L_YU]=\mathrm DP_UY-Y\circ P_U,
 \qquad A_U(L_YU)=Y-P_U^*Y\quad\text{on }\{\chi_1=1\}.
\]

Now let $R_s$ be a smooth path in the reference group, with $R_{s_j}=\Id$ and $Y_s=(\mathrm DR_s)^{-1}\dot R_s$.
Write $Y=Y_s$ when suppressing time. The transformed maps satisfy
\[
 \begin{gathered}
 U^R=R^*U,\qquad P_{U^R}=R^{-1}\circ P_U\circ R,\\
 W^R=R^*W+Y,\qquad \dot U^R=R^*\dot U+L_YU^R,\\
 (J^{\rm raw}\circ R)\circ P_{U^R}=J^{\rm raw}\circ P_U\circ R.
 \end{gathered}
\]
The identity for $\dot U^R$ holds in $H^1(\bar\mu)$. For
$Y\in\mathcal K_0$, the product rule and the spatial bounds on reference
Killing fields give
\[
 |L_YU|+|\bar\nabla L_YU|
 \le C\|Y\|_{\mathcal K_0}(1+r)
                   \sum_{i=0}^2|\bar\nabla^iU|.
\]
Since $(1+r)^2\,\mathrm d\bar\mu/\mathrm d\mu_\beta$ is bounded,
\[
 \|L_YU\|_{H^1(\bar\mu)}
 \le C_\beta\|Y\|_{\mathcal K_0}\|U\|_{H^2(\mu_\beta)}.
\]
Approximation by smooth compactly supported fields and this bound justify
differentiation of the pullback in $H^1(\bar\mu)$, and hence of its
projection onto $\mathcal K_0$.

Reapplying~\eqref{eq:M-4} to the transformed data defines
\[
 Z^R=B_{U^R}\dot U^R+\chi_2P_{U^R}^*(R^*W+Y).
\]
On the common inner core, the preceding gauge derivative identity gives
\[
  Z^{R,\rm true}
 =P_{U^R}^*(R^*W+Y)+A_{U^R}(R^*\dot U+L_YU^R)
 =R^*Z^{\rm true}+Y.
 \]
The global extensions are defined by~\eqref{eq:M-4} in each coordinate system.

The transformed background is $R_s^*X_{b_j}=X_{b(s)}$, with
$b(s)=R_s\cdot b_j$. It stays in $B_{\rm con}$ by invariance.
Its variations satisfy
\[
 \dot G_b=L_{Y_s}G_b,\qquad
 h_b:=\partial_s\log\frac{\mathrm d\mu_b}{\mathrm d\bar\mu}
      =Y_s\!\left(\log\frac{\mathrm d\mu_b}{\mathrm d\bar\mu}\right).
\]
This uses $R_s^*\bar\mu=\bar\mu$. The norm of $Y_s$ is its fixed
finite-dimensional Gaussian norm; rotation fields have linear growth.
For the next lemma suppress the superscripts on the transformed fields,
so $W$ denotes $W^R$ and $Z$ denotes its cutoff extension $Z^R$.
With $P=P_U$, the patched tensor satisfies
\begin{equation}\label{eq:M-3}
 \mathfrak h=\mathfrak h(U,G,b),\qquad
 \delta_{\bar X}\mathfrak h=0.
\end{equation}

\begin{lemma}\label{lem:differentiated-gauge}
At every smooth time of the fixed block the extended generator
satisfies a compatible metric equation
\[
 -\frac12\delta_{\bar X}
     \bigl(\mathrm D_U\mathfrak h[B_U^{-1}Z]\bigr)
 =\delta_{\bar X}\bigl(\chi_3 P^*\Phi+\frac12T_{\rm tail}\bigr).
\]
The tensor $T_{\rm tail}$ consists of fixed-width annular terms and
an exterior term containing $U\dot G_b$. The first-moment equation is
\[
 \mathcal M_U Z=\vec m_{\rm tail}.
\]
Use the exponent $c_6=(c_2+1/8)/2$ fixed in
Subsection~\ref{subsec:parameter-selection}, so that
$c_2<c_6<(1-\beta)/8$. Then
\begin{equation}\label{eq:motion-tail-sources}
 \|T_{\rm tail}\|_{L^2(\bar\mu)}+|\vec m_{\rm tail}|
 \le C\bigl(1+\|Y_s\|_{\mathcal K_0}\bigr)
                  e^{-c_6L^2}(1+|\dot L|).
\end{equation}
The constant and required minimum radius depend only on the fixed
construction data, independently of $Y_s,L,\dot L$, the block index
and the duration of the initial subinterval.
\end{lemma}
\begin{proof}
First freeze the actual metric $G$, compact background $G_b$, and all
cutoffs, and vary only $U$. Put $P_\epsilon:=P_{U+\epsilon X}$.
On a neighborhood of $\operatorname{supp}\chi_3$, the support margin in
Proposition~\ref{prop:restart-gauge} gives
\[
 \left.\frac{\mathrm d}{\mathrm d\epsilon}\right|_{\epsilon=0}
       P_\epsilon^*G
   =L_{A_UX}G^+.
\]
Consequently
\[
 \mathrm D_U\mathfrak h[X]
   =\chi_3L_{A_UX}G^++(1-\chi_3)L_XG_b.
\]
Products with $\chi_3$ and their variations extend by zero.

The same support margin holds while the radius varies. The local
pullback differentiation formula for the actual metric is
\[
 \partial_s(P^*G)
 =P^*(\partial_sG)
        +L_{A_U\dot U+Y_L}(P^*G),
 \qquad
 Y_L:=(\mathrm DP)^{-1}\mathrm D_{\chi_1}P[\dot\chi_1].
\]
Here $Y_L$ is the velocity caused by varying the displacement cutoff;
it vanishes near $\operatorname{supp}\chi_3$, where $\chi_1=1$.
Thus the product rule for \eqref{eq:M-3} gives
\[
  \partial_s\mathfrak h
 ={}\dot\chi_3(G^+-G_b-L_UG_b)
       +\chi_3\bigl(\partial_sG^+-\dot G_b\bigr)
    +(1-\chi_3)\bigl(L_{\dot U}G_b+L_U\dot G_b\bigr).
 \]
Use $\dot G=2\Phi+L_WG$ and
$P^*L_WG=L_{P^*W}G^+$. The two terms containing $\dot U$ in
this full time derivative are precisely
\[
 \chi_3L_{A_U\dot U}G^++(1-\chi_3)L_{\dot U}G_b
       =\mathrm D_U\mathfrak h[\dot U].
\]
With the cutoffs and background still fixed, the two transport
variations are
\[
 \begin{split}
 \mathrm D_G\mathfrak h[L_WG]&=\chi_3L_{P^*W}G^+,\\
 \mathrm D_U\mathfrak h[v]
       &=\chi_3L_{A_Uv}G^++(1-\chi_3)L_vG_b.
 \end{split}
\]
Here $B_U=A_U$ and $\chi_2=1$ near
$\operatorname{supp}\chi_3$, so $A_Uv=P^*W$ there. Subtracting
therefore cancels the entire interior transport:
\begin{equation}\label{eq:M-5}
 \mathrm D_G\mathfrak h[L_WG]-\mathrm D_U\mathfrak h[v]
       =-(1-\chi_3)L_vG_b.
\end{equation}
Substitution of $\dot U=B_U^{-1}Z-v$ now gives the full tensor identity
\[
 \begin{aligned}
 \partial_s\mathfrak h
 ={}&\mathrm D_U\mathfrak h[B_U^{-1}Z]+2\chi_3P^*\Phi
          -(1-\chi_3)L_vG_b\\
 &+\dot\chi_3(G^+-G_b-L_UG_b)
       -\chi_3\dot G_b+(1-\chi_3)L_U\dot G_b.
 \end{aligned}
\]
The remaining transport term is supported where $1-\chi_3$ and $v$
overlap. Since $\delta_{\bar X}\mathfrak h=0$,
the reference divergence of the displayed identity vanishes.

Every compact background lies in the reference slice. Differentiating
that condition gives $\delta_{\bar X}\dot G_b=0$. Since
$\mathrm d\chi_3$ has only Euclidean components and $\dot G_b$ has
only compact indices,
\[
 \delta_{\bar X}(\chi_3\dot G_b)
 =\chi_3\delta_{\bar X}\dot G_b
          +(\mathrm d\chi_3)\lrcorner\dot G_b=0.
\]
Taking divergence of the tensor identity proves the metric equation with
\[
 T_{\rm tail}=-(1-\chi_3)L_vG_b
   +\dot\chi_3(G^+-G_b-L_UG_b)+(1-\chi_3)L_U\dot G_b.
\]
The first two terms have fixed-width annular support: $v$ is
supported in $\operatorname{supp}\chi_2$, and $\dot\chi_3$ is
supported in the transition of $\chi_3$. The last term is supported outside $\{\chi_3=1\}$; the
coordinate-field tail estimate below controls its exterior part.

Each source annihilates the reference Killing fields by
$\delta_{\bar X}^*X=0$. For the exterior tensor, compact cutoff
approximation and the tail estimate below justify integration by parts.

The function $h_b$ depends only on the compact variable. Set
$W^{\rm e}:=\chi_0W$ on $N_L$ and extend it by zero; it equals
$W$ on $N_{L-1}$. Lemma~\ref{lem:motion-low-comparison}, applied
to the original raw generator, and the fixed finite-dimensional
bounds for $Y_s$ give
\[
 \sum_{i=0}^2|\bar\nabla^iW|
 +\sum_{i=0}^2|\bar\nabla^iW^{\rm e}|
 \le C\bigl(1+\|Y_s\|_{\mathcal K_0}\bigr)(1+r)
                  \quad\text{on }N_L.
\]
Since
$\mu_s=e^\ell\mu_b$, the actual density identity in \eqref{eq:M-2}
gives
\[
 \dot\ell=\operatorname{div}_{\mu_b}W+W\ell-h_b,
 \qquad L_X\mu_b=(\operatorname{div}_{\mu_b}X)\mu_b.
\]
The extended probability is
$\mu'=(Z')^{-1}e^{\chi_0\ell}\mu_b$. Differentiating
its scalar and background factors separately yields
\[
 \begin{split}
 \partial_s\mu'
  &=(\dot\chi_0\ell+\chi_0\dot\ell+h_b
                         -\partial_s\log Z')\mu',\\
 L_{W^{\rm e}}\mu'
  &=\bigl(\operatorname{div}_{\mu_b}W^{\rm e}
       +\chi_0W^{\rm e}\ell+(W^{\rm e}\chi_0)\ell\bigr)
                                                \mu'.
 \end{split}
\]
Their difference is the exact global identity
\begin{equation}\label{eq:completed-density-motion}
 \partial_s\mu'-L_{W^{\rm e}}\mu'
       =(q_s-\partial_s\log Z')\mu',
\end{equation}
where
\[
 q_s:=\chi_0(\dot\ell-W^{\rm e}\ell)+h_b
       -\operatorname{div}_{\mu_b}W^{\rm e}
       +(\dot\chi_0-W^{\rm e}\chi_0)\ell.
\]
Terms involving $\ell$ are extended by zero where their cutoffs
vanish. Thus $q_s=0$ on $N_{L-1}$ and $q_s=h_b$ outside $N_L$.
On the intervening annulus the formula for $\dot\ell$ and the low
bounds give
\[
 |q_s|\le C\bigl(1+\|Y_s\|_{\mathcal K_0}\bigr)
                   (1+r)^Q(1+|\dot L|).
\]
The cutoff velocity $Y_L$ above is supported in the transition of
$\chi_1$ and is a bounded multiplier of $\dot\chi_1U$.
Pulling back \eqref{eq:completed-density-motion} gives
\[
 \partial_s\nu'
 =L_{A_U\dot U+Y_L+P^*W^{\rm e}}\nu'
       +(q_s\circ P-\partial_s\log Z')\nu'.
\]
On the support of $v$, $A_U=B_U$, so
$A_Uv=\chi_2P^*W$. Substituting \eqref{eq:M-4} yields
\[
 \partial_s\nu'
 =L_{A_UB_U^{-1}Z+e_W+Y_L}\nu'
       +(q_s\circ P-\partial_s\log Z')\nu',
 \qquad e_W:=P^*W^{\rm e}-\chi_2P^*W.
\]
The field $e_W$ is supported in a fixed-width outer annulus,
where $\chi_2$ differs from one.
For any field $X$ with the indicated integrability, integration by
parts gives
\begin{equation}\label{eq:M-6}
 \int\vec x\,L_X(\mathrm d\nu')
       =-\int X_E\,\mathrm d\nu'.
\end{equation}
Differentiating the exact zero first moment consequently gives
\[
  0={}-\mathcal M_UZ-\int(e_W+Y_L)_E\,\mathrm d\nu'
       +\int\vec x\,(q_s\circ P)\,\mathrm d\nu'
 -(\partial_s\log Z')\int\vec x\,\mathrm d\nu'.
 \]
The final integral vanishes by centering. We obtain the first-moment equation with
\[
 \vec m_{\rm tail}=-\int(e_W+Y_L)_E\,\mathrm d\nu'
       +\int\vec x\,(q_s\circ P)\,\mathrm d\nu'.
\]

At every required fixed spatial order, the identities for the
background variations give
\[
 |\bar\nabla^i\dot G_b|+|\bar\nabla^ih_b|
       \le C_i\|Y_s\|_{\mathcal K_0}.
\]
The compact dependence and linearity in $W,Y_s$ give polynomial
bounds for the annular terms, including $e_W$ and $q_s$, with factor
$(1+\|Y_s\|_{\mathcal K_0})(1+|\dot L|)$.
Their $L^2(\bar\mu)$ tails have
Gaussian rate $1/8$, and their probability moments have rate $1/4$,
before fixed polynomial and collar losses. The displacement is below
$1/10$, and Lemma~\ref{lem:motion-low-comparison}, preserved under
reference isometries, compares the auxiliary density with $\bar\mu$.

Choose $a$ such that $c_6<a<(1-\beta)/8$.
For the exterior tensor $(1-\chi_3)L_U\dot G_b$, factor out
$\|Y_s\|_{\mathcal K_0}$ and apply
\eqref{eq:gauge-coordinate-tail} to the remaining uniformly bounded
coefficient multipliers. It gives
\[
 \|(1-\chi_3)L_U\dot G_b\|_{L^2(\bar\mu)}
 \le C_a\|Y_s\|_{\mathcal K_0}e^{-a(L-C)^2}
                         \|U\|_{H^2(\mu_{\beta})}.
\]
The construction bounds the last norm uniformly. The same estimate
handles the cutoff velocity $Y_L$, with $|\dot L|$ in place of
$\|Y_s\|_{\mathcal K_0}$. The exterior density term $h_b$ has a
Gaussian first-moment tail with the factor $\|Y_s\|_{\mathcal K_0}$.
Combining these bounds yields
\[
 \|T_{\rm tail}\|_{L^2(\bar\mu)}+|\vec m_{\rm tail}|
 \le C_a\bigl(1+\|Y_s\|_{\mathcal K_0}\bigr)(1+|\dot L|)
                 (1+L)^Qe^{-a(L-C)^2}.
\]
Since $a>c_6$, this is bounded by the right-hand side of
\eqref{eq:motion-tail-sources} above a minimum radius independent of
$Y_s$ and $\dot L$.
\end{proof}

The metric equation leaves the reference Killing directions undetermined.
Translations are fixed by the center equation; the remaining symmetry
condition is imposed by a path of reference isometries.

\begin{lemma}\label{lem:reference-group-normalization}
On $I_j$ there is a reference-isometry path $R_s$ with $R_{s_j}=\Id$
for which $\Pi_{\mathcal K_0}(B_{U^R}^{-1}Z^R)=0$.
It is smooth on each smooth piece of $L$ and continues through its
switching corner. Its generator satisfies $\|Y_s\|_{\mathcal K_0}\le C$,
and the induced compact metric and logarithmic-density variations have
uniform bounds at every required fixed spatial order.
\end{lemma}
\begin{proof}
Return to the original raw generator $W_s$ in~\eqref{eq:M-2}.
Since $R$ preserves $\bar\mu$ and $\mathcal K_0$,
$\Pi_{\mathcal K_0}U^R=0$ and $\Pi_{\mathcal K_0}\dot U^R=0$.
For this proof write $U$ for $U^R$ and use its transformed multipliers;
the original raw field remains denoted by $W$.
Projecting the displayed formula for $Z^R$ gives
\[
 \Pi_{\mathcal K_0}(B_U^{-1}Z^R)
 =\Pi_{\mathcal K_0}(B_U^{-1}\chi_2P_U^*(R^*W))+\mathcal R_UY,
 \qquad
 \mathcal R_UY:=\Pi_{\mathcal K_0}(B_U^{-1}\chi_2P_U^*Y).
\]
Use the fixed Gaussian norm on $\mathcal K_0$. Its fields satisfy
$|Y|+|\bar\nabla Y|\le C(1+r)\|Y\|_{\mathcal K_0}$.
Since $\Pi_{\mathcal K_0}Y=Y$, write
\[
 (\mathcal R_U-\Id)Y
 =\Pi_{\mathcal K_0}\left[
    \chi_2(B_U^{-1}P_U^*Y-Y)+(\chi_2-1)Y\right].
\]
For the first term, use
$P_U^*Y=(\mathrm DP_U)^{-1}Y\circ P_U$.
The displacement, its first derivative and the $C^1$ bound for
$B_U^{-1}-\Id$ give
\[
 |B_U^{-1}P_U^*Y-Y|
       \le C\|U\|_{C^2(N_L)}(1+r)\|Y\|_{\mathcal K_0}.
\]
Finite Gaussian moments bound its projected norm by
$C\|U\|_{C^2(N_L)}\|Y\|_{\mathcal K_0}$.
For the second term, $1-\chi_2$ is supported where $r\ge L-C$;
the polynomial growth of $Y$ gives
\[
 \|(1-\chi_2)Y\|_2
 \le C(1+L)^Q e^{-(L-C)^2/8}\|Y\|_{\mathcal K_0}.
\]
Consequently
\[
 \|\mathcal R_U-\Id\|_{\mathcal K_0\to\mathcal K_0}
  \le C\|U\|_{C^2(N_L)}
       +C(1+L)^Q e^{-(L-C)^2/8}<\frac12.
\]
Here~\eqref{eq:gauge-prefactor-smallness} and the chosen minimum radius
give the last inequality. Hence $\|\mathcal R_U^{-1}\|\le2$.

Setting this projection to zero yields
\begin{equation}\label{eq:M-7}
 Y_s=(\mathrm DR_s)^{-1}\dot R_s
 =-\mathcal R_U^{-1}\Pi_{\mathcal K_0}
       \bigl(B_U^{-1}\chi_2P_U^*(R_s^*W_s)\bigr).
\end{equation}
The linear growth bound in Lemma~\ref{lem:motion-low-comparison}
is preserved by reference isometries and by the small displacement. The coefficients of $B_U^{-1}$
are bounded, and $\bar\mu$ has all polynomial moments. Hence the
Gaussian $L^2$ norm inside the projection is at most $C$, independently
of the radius. The inverse bound gives $\|Y\|_{\mathcal K_0}\le C$.

On each smooth piece of $L(s)$ the right-hand side depends smoothly
on $R$ through the unique gauge branch, so the ODE has a local solution.
It is uniformly bounded on the compact reference group and therefore
extends to the end of that piece. At a switch, $L$ and the gauge branch
are continuous. The right-hand side depends on $L$, not on $\dot L$,
so the solution continues across the corner and then along the next
piece. 
The formulas for $\dot G_b$ and $h_b$ above now give their uniform
spatial derivative bounds, since the family is smooth on the fixed
compact parameter neighborhood and $\|Y_s\|_{\mathcal K_0}\le C$.
\end{proof}

\begin{proof}[Proof of Proposition~\ref{prop:coherent-maps}]
Apply Lemma~\ref{lem:reference-group-normalization} and express all fields
in its normalized coordinates, suppressing the superscript $R$.
A tensor $T\in L^2(\bar\mu)$ defines a compatible functional by
\[
 \langle\delta_{\bar X}T,X\rangle
 =\int\langle T,\delta_{\bar X}^*X\rangle\,\mathrm d\bar\mu,
 \qquad \|\delta_{\bar X}T\|_{H^{-1}}\le C\|T\|_2.
\]
Lemma~\ref{lem:motion-low-comparison} compares the metric and actual
probability on the entire support of $\chi_3$, and hence
\[
 \|\chi_3P^*\Phi\|_{L^2(\bar\mu)}\le C\mathcal D(s).
\]
Lemma~\ref{lem:reference-group-normalization} gives the symmetry condition. Substitute the uniform bound for
$\|Y_s\|_{\mathcal K_0}$ into~\eqref{eq:motion-tail-sources} and apply
Lemma~\ref{lem:augmented-weak-inverse} to the compatible metric and
first-moment equations of Lemma~\ref{lem:differentiated-gauge}. Its coefficient
tolerance holds by Stages~\ref{stage:iteration-neighborhoods}--\ref{stage:iteration-thresholds}.
This gives
\[
 \|Z_s\|_{H^1}
 \le C\mathcal D(s)+Ce^{-c_6L^2}(1+|\dot L|)
 \le C\mathcal D(s)+Ce^{-c_2L^2}(1+|\dot L|),
\]
by $c_6>c_2$, proving~\eqref{eq:M-9}.

Proposition~\ref{prop:causal-radii}(i) gives the spatial gauge throughout the
block. Its local branches agree by uniqueness, and
Lemma~\ref{lem:reference-group-normalization} extends the normalization
across the radius switch.
\end{proof}

\subsection{Comparison at a block endpoint}
\label{part9-resets-on-the-actual-common-collar}
\label{subsec:inherited-endpoint-coordinates}

Complete the block $I_j$ by Proposition~\ref{prop:coherent-maps}. On its
full endpoint domain write $J_-,G_-,F_-$ for the limits from the left
at $s_{j+1}$, given by
\[
 J_-:=J_{s_{j+1}^-},\qquad
 G_-:=J_-^*g^{z_*}(s_{j+1}),\qquad
 F_-:=J_-^*f^{z_*}(s_{j+1}).
\]
Define the endpoint parameters directly in these coordinates by
\begin{equation}\label{eq:modified-endpoint-projection}
 (b_{j+1},a)
 :=\Pi_K\!\left[\chi_{R_j}(G_--\bar g)\right].
\end{equation}
The tensor in brackets is extended by zero. The radius inclusions
\eqref{eq:E-13}--\eqref{eq:D-19} give
\[
 \operatorname{supp}\chi_{R_j}\subset N_{\Lambda R_j}
       \Subset N_{L_j(\log2)-4},\qquad S_{j+1}+4\le R_j.
\]
Thus the projection uses the full endpoint chart before its restriction
to $N_{S_{j+1}}$.

Freeze $Q=R_{s_{j+1}}^{(j)}$. The endpoint equivariance
\eqref{eq:endpoint-projection-equivariance} shows that
\eqref{eq:modified-endpoint-projection} equals the projection in the
original block coordinates followed by the action of $Q$ on both
kernel components. Apply Theorem~\ref{thm:restart} and
Proposition~\ref{prop:causal-radii} to the fixed incoming background
and pull their endpoint estimates back by this frozen $Q$.
This gives $b_{j+1}\in B_{\rm out}$, the entropy estimate
\eqref{eq:E-15}, and
\begin{equation}\label{eq:modified-endpoint-estimates}
 \|G_--G_{b_{j+1}}\|_{C^{100n}(N_{S_{j+1}})}
 +\|F_--F_{b_{j+1}}\|_{C^{100n}(N_{S_{j+1}})}
 +|\mathrm De(b_{j+1})|\le e^{-c_5S_{j+1}^2}.
\end{equation}
The transported incoming background is $Q\cdot b_j$, and
$b_{j+1}$ is the projection in~\eqref{eq:modified-endpoint-projection}.

Define the candidate incoming chart by restriction:
\begin{equation}\label{eq:inherited-incoming-chart}
 J_{s_{j+1}}^{{\rm raw},(j+1)}
       :=J_-\big|_{N_{S_{j+1}}},\qquad
 \iota_{j+1}:=\Xi_{s_{j+1}}\circ
                 J_-\big|_{N_{S_{j+1}}}.
\end{equation}
The first map takes values in $M$; the second
in the Ricci-flow slice at $t=-e^{-s_{j+1}}$. After the endpoint
parabolic rescaling, the latter has the order-$100n$ comparison just
obtained, including the actual potential constant. Both the chart and
$b_{j+1}$ are defined before testing $b_{j+1}\in B_{\rm con}$.

If this membership holds, keep $b_{j+1}$ fixed in the next raw
comparison and construct its spatial gauge. With
$P_{\rm new}=P_{j+1}(0)$ and the next reference-group equation
initialized at the identity, its initial moving chart is
\begin{equation}\label{eq:inherited-initial-gauge}
 J_{s_{j+1}}^{(j+1)}
 =J_{s_{j+1}}^{{\rm raw},(j+1)}\circ P_{\rm new},\qquad
 R_{s_{j+1}}^{(j+1)}=\Id.
\end{equation}
Thus the junction map is $P_{\rm new}$, with $b_{j+1}$ fixed by
\eqref{eq:modified-endpoint-projection}. The same transition identity
holds after composition on the left with $\Xi_{s_{j+1}}$.

Write $\mu^+$ for the actual probability of
$(G_-,F_-)$, and $\nu'$ for its centered
auxiliary probability. Pullback by $Q$ preserves their agreement up to $(Z')^{-1}$ and the
zero auxiliary first moment. The next proposition controls the annular
divergence source and first moment created by the new cutoffs.

\begin{proposition}[Comparison at one junction]\label{prop:gauge-junction}
Use the outgoing radii of Proposition~\ref{prop:causal-radii}, with
$S'=S_{j+1}$ and $L=L_{j+1}(0)$, and its common-collar
inclusions~\eqref{eq:E-13}--\eqref{eq:D-19}. Suppose in addition that
$b_{j+1}\in B_{\rm con}$, $R_{+,0}\le L\le S'-1$, and the new data satisfy
\[
 \|G_--G_{b_{j+1}}\|_{C^8(N_L)}
 +\|F_--F_{b_{j+1}}\|_{C^8(N_L)}\le \epsilon_+,
 \qquad \mathcal C_\beta(L)^2\epsilon_+\le\frac1{100}.
\]
Then Proposition~\ref{prop:restart-gauge} constructs
\[
 P_{\rm new}=\exp_{\bar g}(\chi_{1,\rm new}U_{\rm reset}).
\]
For $P_v=\exp_{\bar g}(v\chi_{1,\rm new}U_{\rm reset})$ and
$Z_v^{\rm reset}=(\mathrm DP_v)^{-1}\partial_vP_v$, one has
\begin{equation}\label{eq:M-12}
 \|U_{\rm reset}\|_{H^1(\bar\mu)}\le Ce^{-c_2L^2},
 \qquad
 \int_0^1\|Z_v^{\rm reset}\|_{H^1(\bar\mu)}\,\mathrm dv
       \le Ce^{-c_2L^2}.
\end{equation}
Throughout $0\le v\le1$, the displaced cutoff supports and their
inverse images remain inside the old exact slice and auxiliary-measure
agreement region. The constant depends only on the fixed construction
data, independently of the block index and the radii.
\end{proposition}
\begin{proof}
Apply Proposition~\ref{prop:restart-gauge}, using
$B_{\rm con}\subset B_{\rm gauge}$. In this proof all cutoffs refer
to radius $L$; write $P=P_{\rm new}$ and $U=U_{\rm reset}$.
Denote the new extended probability by $\mu'_{\rm new}$;
its pullback by $P$ is centered. The actual probability remains $\mu^+$,
and the old centered auxiliary probability remains $\nu'$.
The radial cutoffs satisfy
\[
 \chi_i=1\ \text{on }N_{L-i-1},\qquad
 \operatorname{supp}\chi_i\subset N_{L-i},\qquad i=0,1,3.
\]
Their derivatives have the bounds of Proposition~\ref{prop:restart-gauge}.
The collars~\eqref{eq:D-19} and displacement below $1/10$ place all
supports, including $N_L$, and their images and inverse images inside
the old exact slice and auxiliary-measure agreement region throughout
the isotopy $P_v=\exp_{\bar g}(v\chi_1U)$.

Set $h:=G_--G_{b_{j+1}}$ and define the unpulled auxiliary
probability by
\[
 \ell=\log\frac{\mu^+}{\mu_{b_{j+1}}},\qquad
 \mu'_{\rm new}
 =(Z'_{\rm new})^{-1}e^{\chi_0\ell}\mu_{b_{j+1}},
 \qquad
 Z'_{\rm new}=\int e^{\chi_0\ell}\,\mathrm d\mu_{b_{j+1}}.
\]
Freeze the data on $\operatorname{supp}\chi_0$, the parameter $b_{j+1}$,
and the cutoffs. For a variable field $V$, put
\[
 \mathfrak h(V)
 =\chi_3(P_V^*G_--G_{b_{j+1}})
       +(1-\chi_3)L_VG_{b_{j+1}}.
\]
The first term extends by zero by the common-collar inclusion above.
The augmented reset equation is
\[
 \mathcal F(V)=
 \left(-\frac12\delta_{\bar X}\mathfrak h(V),
       \Pi_{\mathcal K_0}V,
       -\int\vec x\,\mathrm d(P_V^*\mu'_{\rm new})\right),
 \qquad \mathcal F(U)=0.
\]
Writing $\vec m:=\int\vec x\,\mathrm d\mu'_{\rm new}$ gives
\[
 \mathcal F(0)=
 \left(-\frac12\delta_{\bar X}(\chi_3h),
       0,-\vec m\right).
\]

Every compact background lies in the same fixed reference slice.
Consequently, on the whole new cutoff support,
\[
 \delta_{\bar X}h
 =\delta_{\bar X}(G_--\bar g)
       -\delta_{\bar X}(G_{b_{j+1}}-\bar g)=0,
\]
and hence
\begin{equation}\label{eq:M-11}
 \delta_{\bar X}(\chi_3h)
       =(\mathrm d\chi_3)\lrcorner h.
\end{equation}
This one-form is supported in $L-4\le r\le L-3$ and is compatible:
for every reference Killing field $X$,
\[
 \bigl\langle\delta_{\bar X}(\chi_3h),X\bigr\rangle_{L^2(\bar\mu)}
 =\bigl\langle\chi_3h,
                         \delta_{\bar X}^*X\bigr\rangle_{L^2(\bar\mu)}=0.
\]
We estimate the annular one-form in $L^2(\bar\mu)$ and use
$L^2\hookrightarrow H^{-1}$.

For the first-moment equation, Proposition~\ref{prop:restart-gauge} gives
$\nu'=(Z')^{-1}\mu^+$ on the old agreement region.
On $N_{L-1}$ one has
$\mu'_{\rm new}=a_g\nu'$, where
$a_g=Z'/Z'_{\rm new}$. Since the old probability is centered,
\begin{equation}\label{eq:restart-annular-moment}
 \vec m
 =\int_{\bar M\setminus N_{L-1}}\vec x\,
                                  \mathrm d\mu'_{\rm new}
       -a_g\int_{\bar M\setminus N_{L-1}}\vec x\,
                                  \mathrm d\nu'.
\end{equation}
The low logarithmic bounds and the old bounded displacement give
uniform Gaussian tails for both extended probabilities. Mass one and
$\nu'(N_{L-1})\ge1/2$ for large $L$ imply
\[
 0<a_g=\frac{\mu'_{\rm new}(N_{L-1})}
                 {\nu'(N_{L-1})}\le2.
\]
Thus $|\vec m|$ is bounded by the first exterior moment of
$\mu'_{\rm new}$ plus twice that of $\nu'$.

Use the fixed margin $c_2<c_6<1/8$ from
Subsection~\ref{subsec:parameter-selection}. The polynomial bounds on
$h$ and the fixed cutoff derivatives give
\[
 \bigl\|\delta_{\bar X}(\chi_3h)\bigr\|_{L^2(\bar\mu)}
 \le C(1+L)^Q e^{-(L-4)^2/8}.
\]
The probability first-moment tails have rate $1/4$. Since $c_6<1/8$,
the fixed collar losses and polynomial factors give
\[
 \bigl\|\delta_{\bar X}(\chi_3h)\bigr\|_{L^2(\bar\mu)}
       +|\vec m|
 \le Ce^{-c_6L^2},\qquad
 \|\mathcal F(0)\|_{\mathcal Y}
 \le Ce^{-c_2L^2}.
\]
Here $\mathcal Y=H^{-1}_{\mathcal K}\times\mathcal K_0\times\mathbb R^d$
is the target space of Lemma~\ref{lem:augmented-weak-inverse}, and
all constants are uniform in the construction neighborhood and $L$.

The norm bound for the ball in Proposition~\ref{prop:restart-gauge}
places the entire segment $vU$, $0\le v\le1$, in the construction ball
of Lemma~\ref{lem:augmented-weak-inverse}. The final assertion of that
lemma gives, with the inputs frozen as above,
\[
 \|\mathrm D\mathcal F(vU)
                   -\mathscr P\|_{H^1\to\mathcal Y}
 \le\epsilon,\qquad
 \epsilon\|\mathscr P^{-1}\|_{\mathcal Y\to H^1}\le\frac12.
\]
Here $\mathscr P$ is the fixed reference augmented operator
from Lemma~\ref{lem:augmented-weak-inverse}. The fundamental theorem of
calculus therefore gives
\[
 \mathscr PU
 =-\mathcal F(0)-\int_0^1
 \bigl(\mathrm D\mathcal F(vU)
                       -\mathscr P\bigr)
                            U\,\mathrm dv.
\]
Apply its reference inverse and absorb the last term to obtain
\[
 \|U\|_{H^1(\bar\mu)}
 \le 2\|\mathscr P^{-1}\|_{\mathcal Y\to H^1}
                    \|\mathcal F(0)\|_{\mathcal Y}
 \le Ce^{-c_2L^2}.
\]
The fixed cutoff derivative bounds give the displacement estimate
\[
 \|\chi_1U\|_{H^1(\bar\mu)}
 \le Ce^{-c_2L^2}.
\]
The construction bounds make $P_v=\exp_{\bar g}(v\chi_1U)$
diffeomorphisms on the common collar, with uniformly bounded first
two derivatives and inverses. Differentiating the reference exponential
expresses
\[
 Z_v^{\rm reset}=(\mathrm DP_v)^{-1}\partial_vP_v
\]
as a multiplier of $\chi_1U$ whose coefficients
and first derivatives are uniformly bounded. The product rule yields
\[
 \|Z_v^{\rm reset}\|_{H^1(\bar\mu)}
 \le C\|\chi_1U\|_{H^1(\bar\mu)}
 \le Ce^{-c_2L^2},\qquad 0\le v\le1.
\]
Integration proves~\eqref{eq:M-12}. The inherited map identities show that
$J_-\circ P_v$ joins the actual endpoint charts.
\end{proof}

\subsection{Weak motion on finite sequences of blocks}
\label{subsec:comparison-length}\label{subsec:weak-metric-observations}
Start with an order-$100n$ incoming chart and $b_0=0$ satisfying
the initialization conditions below. Given the incoming chart and
$b_j\in B_{\rm con}$, propagate the raw comparison with this fixed
parameter and construct its spatial gauges using
Proposition~\ref{prop:causal-radii}. Proposition~\ref{prop:coherent-maps}
then supplies $R_s^{(j)}$, starting at the identity, and the normalized
charts $J_s^{(j)}$ on $I_j$. Put
\[
 G(s):=(J_s^{(j)})^*g^{z_*}(s),\qquad s_j<s<s_{j+1}.
\]
Project the normalized endpoint metric by~\eqref{eq:modified-endpoint-projection}
and define the candidate incoming chart by~\eqref{eq:inherited-incoming-chart}.
Its order-$100n$ error and residual satisfy~\eqref{eq:E-9} with $j+1$
in place of $j$, and $R_j\ge R_{\min}$, $E_{j+1}\le E_0\le E_*$.
If $b_{j+1}\in B_{\rm con}$, Proposition~\ref{prop:causal-radii}
supplies the new gauge smallness condition, lower radius and collars,
with current drop zero at the initial time. Construct $P_{\rm new}$,
use~\eqref{eq:inherited-initial-gauge}, and continue. Otherwise stop
at this candidate endpoint. Proposition~\ref{prop:gauge-junction}
controls every junction actually constructed by this rule.

At a constructed junction, with $P_j(0)$ the new spatial gauge,
\[
 J_{s_j^+}=J_{s_j^-}\circ P_j(0),\qquad
 G(s_j^+)=P_j(0)^*G(s_j^-).
\]
The two values are taken before and after the new gauge change. Compose $J_{s_j^-}$ with
the gauge isotopy $P_{j,v}$ to obtain the comparison path, and denote
its generator by $Z_{j,v}^{\rm reset}$. On block interiors write
$J_s$ for $J_s^{(j)}$ and $Z_s$ for the normalized extension from
Proposition~\ref{prop:coherent-maps}.

Consider any finite initial sequence of this process, allowing a final
partial block and including only junctions already constructed.
All $H^1$ norms below use the reference metric and probability $\bar\mu$.

\begin{proposition}
\label{prop:comparison-length}
For such a finite sequence and $0\le u\le\log2$ in its $j$th block,
\begin{equation}\label{eq:M-13}
 \int_0^u\bigl(\mathcal D(s_j+v)+\mathcal D(s_j+v)^{p_2}
              +\|Z_{s_j+v}\|_{H^1}\bigr)\;\mathrm dv
 \le C\bigl(A_j(u)^{p_4}+\Delta_{j-1}^{p_4}\bigr).
\end{equation}
Every constructed incoming junction satisfies
\[
 \int_0^1\|Z_{j,v}^{\rm reset}\|_{H^1}\,\mathrm dv
       \le C\Delta_{j-1}^{p_4},\qquad j\ge1.
\]
The first spatial gauge is included in the initial pairing error.

Consequently, for every such segment ending at $S=s_J+u$,
\begin{equation}\label{eq:finite-prefix-motion}
 \begin{split}
 &\int_{s_0}^{S}\bigl(\mathcal D(s)+\mathcal D(s)^{p_2}+\|Z_s\|_{H^1}\bigr)
              \;\mathrm ds\\
 &\quad+\sum_{\text{junctions }j\text{ in the segment}}
       \int_0^1\|Z_{j,v}^{\rm reset}\|_{H^1}\;\mathrm dv
 \le\omega(E_0+\Delta_{-1}).
 \end{split}
\end{equation}
At every point of the concatenated chart path, including the endpoint before the next gauge change, its pulled-back metric $G$ satisfies
\begin{equation}\label{eq:M-15}
 |\mathcal O(G)-\mathcal O(G(s_0))|
       \le\omega(E_0+\Delta_{-1}).
\end{equation}
The constants and modulus depend only on the fixed construction data,
uniformly in the number of blocks, the final partial-block duration
and the endpoint of the finite sequence.
\end{proposition}
\begin{proof}
Since $E'=-2\mathcal D^2$, for
$\alpha\in\{1,p_2\}\subset(p_4,2)$ and $u>0$,
\[
 \int_0^u\mathcal D(s_j+v)^\alpha\,\mathrm dv
 \le u^{1-\alpha/2}\left(\frac{A_j(u)^2}{2}\right)^{\alpha/2}
 \le C A_j(u)^{p_4}.
\]
This controls the defect term in~\eqref{eq:M-9}.

Let $L_{\min}=\min_{0\le v\le u}L_j(v)$. The increasing propagation
branch is bounded below by its initial value and the decreasing entropy
branch by its current value. Thus, for a fixed collar constant $C$,
\[
 L_{\min}\ge\min\left\{
 \rho\sqrt{p_3} R_{j-1}-C,
 \sqrt{\Lambda'\log(1/A_j(u))}-C\right\}.
\]
The second entry is infinite at zero drop. The inequalities
$\rho^2p_3>p_4$ and $c_2\Lambda'>p_4$ imply
\[
 \begin{split}
 e^{-c_2(\rho\sqrt{p_3} R_{j-1}-C)^2}
       &\le C e^{-c_2p_4 R_{j-1}^2}=C\Delta_{j-1}^{p_4},\\
 e^{-c_2(\sqrt{\Lambda'\log(1/A_j(u))}-C)^2}
       &\le C A_j(u)^{p_4}.
 \end{split}
\]
The strict coefficient gaps absorb the linear collar terms. Hence
\[
 \int_0^u e^{-c_2L_j(v)^2}\,\mathrm dv
 \le (\log2)e^{-c_2L_{\min}^2}
 \le C(\Delta_{j-1}^{p_4}+A_j(u)^{p_4}).
\]

Before the possible switch in~\eqref{eq:E-10}, the active radius
increases; afterward it decreases. On either monotone part the
change-of-variable formula for absolutely continuous functions gives
the integral of $e^{-c_2r^2}$ over the corresponding interval of radii.
Each interval is contained in $[L_{\min},\infty)$, so
\[
 \int_0^u e^{-c_2L_j(v)^2}|\dot L_j(v)|\,\mathrm dv
 \le2\int_{L_{\min}}^\infty e^{-c_2r^2}\,\mathrm dr
 \le C e^{-c_2L_{\min}^2}.
\]
The last inequality uses the fixed lower radius threshold. Continuity
of $L_j$ excludes a jump at the switch, so this contribution has the
same bound $C(\Delta_{j-1}^{p_4}+A_j(u)^{p_4})$.

For $j\ge1$, the incoming reset radius is
$L_j(0)=\rho\sqrt{p_3} R_{j-1}-C_1-C_2$.
Estimate~\eqref{eq:M-12} and the first of the preceding tail bounds
give the asserted junction integral. The smooth-block estimates and
\eqref{eq:M-9} prove~\eqref{eq:M-13}.

For the completed blocks, Lemma~\ref{lem:delayed-sum} sums both
$A_j^{p_4}$ and $\Delta_j^{p_4}$. The shifted tails contribute also
$\Delta_{-1}^{p_4}$. A final partial block costs at most
$CA_J(u)^{p_4}\le CE_0^{p_4/2}$, together with its preceding tail.
Hence the total is bounded by
\[
 C\bigl(E_0^\vartheta+\Delta_{-1}^{p_4}+E_0^{p_4/2}\bigr),
\]
which is a uniform vanishing modulus of $E_0+\Delta_{-1}$.
Summing the constructed junctions proves~\eqref{eq:finite-prefix-motion}.

The tensors $T_i$ defining $\mathcal O$ are supported in $N_1$.
Since $L_j(u)\ge R_{+,0}>10$, the collar $N_3$ lies in
$N_{L_j(u)-4}$. The extension $Z_s$ agrees there with the true
chart generator, so $\partial_sG=2\Phi+L_{Z_s}G$ on the pairing support.
Differentiating the fixed pairings gives
\[
 \frac{\mathrm d}{\mathrm ds}\mathcal O(G(s))_i
 =\int_{N_1}\langle 2\Phi+L_{Z_s}G,T_i\rangle_{\bar g}
                                      \,\mathrm d\bar\mu.
\]
The uniform $C^1$ metric bound and the metric and actual-measure
comparisons of Lemma~\ref{lem:motion-low-comparison} imply
$|L_{Z_s}G|_{\bar g}\le C(|Z_s|+|\bar\nabla Z_s|)$ and
$\|\Phi\|_{L^2(N_1,\bar\mu)}\le C\mathcal D(s)$. Hence
\[
 \left|\frac{\mathrm d}{\mathrm ds}\mathcal O(G(s))\right|
       \le C\bigl(\mathcal D(s)+\|Z_s\|_{H^1}\bigr).
\]

At a constructed junction, displacement below $1/10$ gives
$P_v^{\pm1}(N_3)\subset N_4$. The common collar and the uniform
$C^2$ bounds for $P_v$ and its inverse give a uniform $C^1$ bound
for $P_v^*G$ on $N_1$. Apply the same pairing calculation to
$\partial_v(P_v^*G)=L_{Z_v^{\rm reset}}(P_v^*G)$ and integrate:
\[
 |\mathcal O(P_1^*G)-\mathcal O(G)|
       \le C\int_0^1\|Z_v^{\rm reset}\|_{H^1}\,\mathrm dv.
\]
This estimate holds on every initial part of the junction path.
Summation over the constructed blocks and junctions, using
\eqref{eq:finite-prefix-motion}, proves~\eqref{eq:M-15}, also for a sequence
ending before the next gauge change.
\end{proof}

\subsection{Parameter confinement and continuation}
\label{part9-nested-neighborhoods-and-a-genuine-no-exit-proof}

We use the finite-sequence pairing estimate to confine the compact parameter,
as in \cite[Corollary~6.18]{LYZII}.

Choose $S_0$ sufficiently large for all fixed thresholds and collars,
with $R_{-1}=S_0/\sqrt{p_3}\ge R_{\min}$ and
$\Delta_{-1}=e^{-c_2S_0^2/p_3}$ sufficiently small. Keep this radius
fixed. A sufficiently late chart from the selected tangent sequence
then has arbitrarily small $C^{100n}$ pair error on a neighborhood of
$\overline N_{S_0}$. Set $b_0=0$. Since $\mathrm De(0)=0$, this
error satisfies~\eqref{eq:E-9}; the entropy limit~\eqref{eq:entropy-limit}
likewise makes $E_0\le E_*$ arbitrarily small. Choose the chart late enough
also to supply the rescaled interval $[-2,0)$.

Apply the first spatial gauge. At this fixed radius the gauge of the exact
model is the identity, and the constructed gauge depends continuously on
the data. Thus $|\mathcal O(G(s_0))|$ tends to zero with the selected-chart
error. Choose the radius and then the chart so that
\begin{equation}\label{eq:observation-initial-smallness}
 |\mathcal O(G(s_0))|+\omega(E_0+\Delta_{-1})<r_{\rm con}/8.
\end{equation}
These choices use only the fixed construction data, independently of
the number of blocks.

\begin{proposition}\label{prop:continuation}
For this initialization, the recursive construction continues for
every $j\ge0$, and
\[
 |b_j|<r_{\rm con}/2\qquad(j\ge0).
\]
\end{proposition}
\begin{proof}
Suppose $k\ge1$ is the first candidate index with
$|b_k|\ge r_{\rm con}/2$. The preceding parameters lie in the
half-radius ball, so blocks $0,\ldots,k-1$ and the junctions preceding
them have been constructed. Their last endpoint defines
$b_k\in B_{\rm out}$. Apply~\eqref{eq:M-15} to the finite sequence
ending at $G(s_k^-)$ before the next gauge change.

By Stage~\ref{stage:parameter-pairings}, $\mathcal O(G_{b_k})=b_k$.
In the inherited coordinates, \eqref{eq:modified-endpoint-estimates} gives
\[
 |\mathcal O(G(s_k^-))-b_k|
 \le C\|G(s_k^-)-G_{b_k}\|_{C^0(N_1)}
 \le Ce^{-c_5S_k^2}<r_{\rm con}/8,
\]
since $S_k^2=p_3R_{k-1}^2\ge p_3R_{\min}^2$ by
Proposition~\ref{prop:causal-radii} and Stage~\ref{stage:iteration-thresholds}.
By~\eqref{eq:M-15} and~\eqref{eq:observation-initial-smallness},
$|\mathcal O(G(s_k^-))|<r_{\rm con}/8$.
Thus $|b_k|<r_{\rm con}/4$, a contradiction.
Every candidate lies in the half-radius ball, so the next gauge and
block exist.
\end{proof}

\phantomsection\label{par:continued-radii}
The minimum retained radii tend to infinity. Indeed, the uniform
finite-sequence sum~\eqref{eq:M-14}, now applied for all numbers of
blocks, gives $A_j,\Delta_j\to0$ and hence
$R_j^2=c_2^{-1}\log(1/\Delta_j)\to\infty$.
For $0\le u\le\log2$, the radius formula gives
\[
 L_j(u)\ge\min\left\{\rho\sqrt{p_3}R_{j-1}-C,
                  \sqrt{\Lambda'\log(1/A_j)}-C\right\}\longrightarrow\infty.
\]
The second term is infinite when $A_j=0$. The junction radii
$L_{j+1}(0)=\rho\sqrt{p_3}R_j-C_1-C_2$ also diverge, and their
isotopies have displacement below $1/10$. Thus this growth is uniform
throughout the chart paths. Finally, Proposition~\ref{prop:comparison-length}
and monotone convergence give
\begin{equation}\label{eq:length-tail}
 I_{p_2}(s)
 :=\int_s^\infty\mathcal D(u)^{p_2}\,\mathrm du<\infty,
 \qquad I_{p_2}(s)\longrightarrow0.
\end{equation}

\section{Intrinsic estimates and convergence in a fixed gauge}\label{sec:uniqueness}

The radius and length estimates of
Lemmas~\ref{lem:quantitative-decay}--\ref{lem:quantitative-radii},
containment of entire potential sublevels in the comparison charts,
and intrinsic interpolation allow us to integrate the modified-flow
equations in one fixed embedding. We use the fixed-subset and
reconstruction arguments of \cite[Section~6.7]{LYZII}, following
\cite[Section~7]{FLC}.

\subsection{Intrinsic estimates on potential sublevels}\label{subsec:intrinsic-sublevels}
Write $(g_s,f_s)=(g^{z_*}(s),f^{z_*}(s))$ for the modified flow after
the continuation in Section~\ref{sec:motion}. Derivatives in this
subsection are intrinsic to $g_s$, and $p_2$ belongs to the construction
under consideration.

\begin{lemma}\label{lem:quantitative-sublevel-coverage}
For every finite $B$, there are $R_B$ and $s_B$ such that
\begin{equation}\label{eq:whole-sublevel-coverage}
 \{f_s\le B\}\subset J_s(N_{R_B}),\qquad s\ge s_B,
\end{equation}
with a fixed surrounding neighborhood in the controlled comparison
region. For quantitative coverage put $q:=p_1p_3$ and prescribe
\[
 0<c'<\frac{\rho^2p_3^2q}{10c_2(2-q)}.
\]
There are a low-order accuracy threshold $\bar\epsilon>0$ and
$0<\vartheta<1$ such that, if the retained charts and junction paths
satisfy~\eqref{eq:D-4} with right side at most $\bar\epsilon$, then
\begin{equation}\label{eq:quantitative-sublevel-coverage}
 \{f_s\le c'\log(1+s)\}
 \subset J_s\bigl(N_{\vartheta(L_j(u)-C)}\bigr),
 \qquad s=s_j+u,\quad 0\le u\le\log2,
\end{equation}
for all sufficiently large $s$. Here $C$ is a fixed collar width.
Every point of this sublevel has a controlled ambient neighborhood
in the region of the intrinsic estimates in Section~\ref{sec:physical}.
The moving charts on each block and both endpoint charts at a junction
satisfy the slice equation on their retained cores. Throughout the
junction isotopy the sublevel and its surrounding collar remain in the
old controlled region with the stated metric comparison.

The constants depend only on the fixed construction data and the
displayed sublevel parameters. Waiting times may also depend on the
initialization. The additional accuracy is imposed before the minimum
radii and final thresholds of the quantitative construction are fixed.
\end{lemma}
\begin{proof}
Let $z_H(s)$ be a terminal $H_n$-center on the modified manifold and
put $p_s:=J_s(y_0,0)$. Lemma~\ref{lem:terminal-tails} gives
\[
 f_s(x)\ge\frac15d_{g_s}(z_H(s),x)^2-C.
\]
Metric and actual-potential comparison give $J_s(N_1)$ probability
at least $C^{-1}$ and diameter at most $C$. The second-moment bound yields
\[
 H_n\ge C^{-1}\bigl(d_{g_s}(z_H(s),p_s)-C\bigr)_+^2,
 \qquad d_{g_s}(z_H(s),p_s)\le C.
\]
Thus an entire sublevel of height $B$ lies in the ambient ball about
$p_s$ of radius $\sqrt{5B+C}+C$. If $J_s^*g_s\ge a\bar g$,
every path from $p_s$ leaving $J_s(N_R)$ has length at least
$\sqrt aR$, by Lemma~\ref{lem:coverage}. Taking $R_B$ larger than
this ball and its fixed collars proves~\eqref{eq:whole-sublevel-coverage}
once the minimum radii have grown, as on
p.~\pageref{par:continued-radii}.

For the second assertion choose $a,\vartheta<1$ close enough to one that
\[
 5c'<a\vartheta^2
              \frac{\rho^2p_3^2q}{2c_2(2-q)},
\]
and choose $\bar\epsilon$ to ensure
$a\bar g\le J_s^*g_s\le a^{-1}\bar g$ on the exact slice and along
the junction paths. Lemma~\ref{lem:quantitative-radii}, with the fixed
collars subtracted, gives
\[
 (L_j(u)-C)^2\ge
 \frac{\rho^2p_3^2q}{2c_2(2-q)}\log(1+s)
             -C\sqrt{\log(2+s)}-C.
\]
Here the fixed block length makes $1+j$ and $1+s$ comparable after
a waiting time. The strict coefficient gap implies
\[
 \sqrt a\,\vartheta(L_j(u)-C)
   >\sqrt{5c'\log(1+s)+C}+C+2.
\]
The preceding center and first-exit argument proves
\eqref{eq:quantitative-sublevel-coverage}.

The gap to $N_{L_j(u)-4}$ tends to infinity, leaving ambient balls
of fixed radius around the sublevel. In the ungauged Ricci-flow
coordinates their radial coordinates are at most
$e^{-u/2}(\vartheta(L_j(u)-C)+1)\le\vartheta\rho S_j+C$.
Thus the fixed surrounding cylinders lie in $N_{\rho_2S_j}$, where
the earlier interval used for curvature estimates and the later interval
used for terminal-density estimates in Section~\ref{sec:physical} apply.
At junctions the same radius bound,
the collars~\eqref{eq:D-19}, and displacement below $1/10$ preserve
these strict margins.
\end{proof}

\begin{lemma}\label{lem:sublevel-velocity}\label{lem:weighted-intrinsic-velocity}
For every finite $B$, the construction of Sections~\ref{sec:radii}--\ref{sec:motion}
with the parameters fixed in Subsection~\ref{subsec:parameter-selection} gives a time $s_B$
such that, for every integer $k\ge0$,
\begin{equation}\label{eq:sublevel-velocity}
 \sum_{\ell=0}^k\bigl(|\nabla^\ell\Phi|
       +|\nabla^\ell\partial_sg_s|+|\nabla^\ell\partial_sf_s|\bigr)
 \le C_{B,k}\mathcal D(s)^{p_2}
 \quad\text{on }\{f_s\le B\},\qquad s\ge s_B.
\end{equation}
The constant depends on $B,k$ and the fixed data; the time $s_B$ is
independent of $k$.

For expanding sublevels fix $k$ and $c'$ as in
Lemma~\ref{lem:quantitative-sublevel-coverage}. There are a finite rough
order and a positive Gaussian loss, depending on $k,p_2$ and the
radius ratios, such that a construction supplying
\eqref{eq:terminal-arbitrary-rough} at these inputs and the accuracy
in that lemma satisfies
\begin{equation}\label{eq:weighted-intrinsic-velocity}
 \sum_{\ell=0}^k\bigl(|\nabla^\ell\Phi|
       +|\nabla^\ell\partial_sg_s|+|\nabla^\ell\partial_sf_s|\bigr)
 \le C_k e^{f_s/2}\mathcal D(s)^{p_2}
 \quad\text{on }\{f_s\le c'\log(1+s)\}
\end{equation}
for all sufficiently large $s$. Its constant depends on the background,
$k,p_2,c'$ and the fixed comparison bounds and margins, but
not on the point or block. The waiting time may also depend on the
flow and selected tangent chart. The rough order and loss precede
the choices of $\eta,\epsilon_{\rm met}$ and the associated thresholds;
the incoming order remains $100n$.
\end{lemma}
\begin{proof}
On the collared regions just obtained, metric comparison, local volume
bounds and~\eqref{eq:A-4} give uniformly controlled intrinsic coordinate
balls. On each fixed larger core,~\eqref{eq:fixed-core-rough} gives
all curvature and potential derivatives on the same cylinders. Its earlier
interval for curvature estimates and later interval for terminal-density
estimates are available also at a block endpoint. Enlarging the derivative
order changes only
the constants, not the core or its waiting time. Since the actual
density is $(4\pi)^{-n/2}e^{-f_s}$ and $f_s$ is bounded on this core,
the ordinary local $L^2$ norm of $\Phi$ is at most $C_B\mathcal D(s)$.
Local interpolation with sufficiently high rough order therefore gives
\eqref{eq:sublevel-velocity}, as in the differentiated estimate of
\cite[Lemma~6.21]{LYZII}.

For the expanding assertion write $x=J_s(y,\vec x)$ and $r=|\vec x|$.
The low comparison of the actual pair gives, on its surrounding collar,
\begin{equation}\label{eq:weighted-velocity-low-data}
 f_s(x)=f_N+r^2/4+O(1),\qquad |\nabla f_s|\le C(1+r).
\end{equation}
Choose $l$ with $(l-k)/(l+n/2)>p_2$, and take the Gaussian loss
in~\eqref{eq:terminal-arbitrary-rough} through order $l+2$ small
enough that the exponent $\delta$ below lies in $(0,1/8)$. Parabolic
normalization and the bounded gauge displacement give, on the controlled neighborhoods,
\begin{equation}\label{eq:weighted-velocity-rough}
 \sum_{i=0}^{l}|\nabla^i\Phi|
       \le C(1+r)^Qe^{\delta r^2},
\end{equation}
where the two additional derivative orders account for $\nabla^2f_s$ in $\Phi$.

On balls $B=B_{g_s}(x,2a/(1+r))$, with fixed small $a>0$, the local
geometry is uniform after rescaling; their doubles lie in the
controlled neighborhood. The potential oscillation is bounded
by~\eqref{eq:weighted-velocity-low-data}. Consequently
\begin{equation}\label{eq:local-defect-unweighting}
 \|\Phi\|_{L^2(B;\mathrm dV_{g_s})}^2
 =(4\pi)^{n/2}\int_B e^{f_s}|\Phi|^2\,\mathrm d\mu^{z_*}
 \le C e^{f_s(x)}\mathcal D(s)^2.
\end{equation}
For $\mathcal D(s)\le1$, interpolate on the rescaled ball as in
\cite[Lemma~6.8, \textup{(6.20)}]{LYZII}. The choice of $l$ allows the
$L^2$ exponent $p_2$ for every derivative through order $k$.
Rescaling contributes powers of $1+r$, so
\eqref{eq:weighted-velocity-rough} and~\eqref{eq:local-defect-unweighting} give
\[
 \sum_{i=0}^k|\nabla^i\Phi|(x)
 \le C(1+r)^Q e^{(1-p_2)\delta r^2+p_2 f_s(x)/2}\mathcal D(s)^{p_2}.
\]
After division by $e^{f_s(x)/2}$, the exponential factor is
$e^{-(1-p_2)(1/8-\delta)r^2+O(1)}$ by~\eqref{eq:weighted-velocity-low-data}.
It absorbs the polynomial uniformly in $r$. For $\mathcal D(s)\ge1$,
\eqref{eq:weighted-velocity-rough} gives the bound directly since
$(1+r)^Q e^{\delta r^2}\le C e^{f_s(x)/2}$. Finally,
$\partial_sg_s=2\Phi$ and $\partial_sf_s=\Tr_{g_s}\Phi$ prove the
two velocity estimates.
\end{proof}

\subsection{Convergence in a fixed embedding}
\label{part9-one-fixed-embedding-not-a-collection-of-unrelated-gauges}

The preceding intrinsic estimates permit the fixed-subset arguments
of \cite[Lemma~6.21 and Propositions~6.22--6.23]{LYZII}.

\begin{proposition}\label{prop:fixed-limit}\label{prop:embedding}
There are a connected open set $U\subset M$ and one smooth,
time-independent open embedding $\Psi:\bar M\to M$, with image $U$,
such that
\begin{equation}\label{eq:M-20}
 \Psi^*g_s\longrightarrow\bar g,\qquad
 \Psi^*f_s\longrightarrow\bar f
 \quad\text{in }C^\infty_{\rm loc}(\bar M),\qquad
 \Psi^*\mathrm d\mu^{z_*}=\mathrm d\bar\mu.
\end{equation}
The limiting metric on $U$ is complete and its potential is proper.
For every finite $B$, the entire sublevel $\{f_s\le B\}$ eventually
lies in a fixed compact subset of $U$.
\end{proposition}
\begin{proof}
By~\eqref{eq:M-1}, the modified pair satisfies
$\partial_sg_s=2\Phi$, $\partial_sf_s=\Tr_{g_s}\Phi$, and its
normalized probability measure $\mu^{z_*}$ is time-independent, smooth
and positive.
The global Gaussian estimate gives $f_s\ge-C$, while the scalar lower
bound of the original closed flow gives $\mathrm R_{g_s}\ge-Ce^{-s}$.
Lemma~\ref{lem:sublevel-velocity} supplies all intrinsic derivative
estimates on each entire bounded potential sublevel after one waiting
time, and~\eqref{eq:length-tail} supplies the integrable factor
$\mathcal D^{p_2}$.

Apply the potential barriers and fixed-domain convergence arguments
of \cite[Lemma~6.21 and Proposition~6.22]{LYZII}. The selected tangent
sequence supplies bounded-potential points at arbitrarily late times,
so the resulting open set $U$ is nonempty. It carries a smooth limiting
pair $(g_\infty,f_\infty)$, complete on each component, with proper
potential. Every entire bounded sublevel of $f_s$ eventually lies
in a fixed compact subset of $U$.

Compose the selected measured tangent embeddings with
$\Xi_{\sigma_i}^{-1}$. Their pulled-back pairs converge to
$(\bar g,\bar f)$ on exhausting collared compact domains. For each
fixed compact $K\subset\bar M$, the pulled-back potentials are uniformly
bounded on $K$. Sublevel confinement therefore puts its images in a
fixed compact subset of $U$ for all sufficiently large $i$.
By the argument of \cite[Proposition~6.23]{LYZII}, a subsequential
limit of these embeddings identifies the complete reference shrinker
with one component of $U$ and pulls back the time-independent density
to $\bar\mu$. This component has mass one, and positivity of
$\mu^{z_*}$ excludes any other component. This gives $\Psi$
and~\eqref{eq:M-20}.
\end{proof}

\subsection{Quantitative convergence on expanding sublevels}
\label{subsec:quantitative-fixed-gauge}

The radius estimates supply the lower bound required in
\cite[Lemma~6.24 and the proof of Theorem~6.20]{LYZII}.

\begin{proposition}\label{prop:quantitative-fixed-gauge}
The embedding $\Psi$ in Proposition~\ref{prop:fixed-limit} has the
following property. For every integer $k\ge0$ and every
\[
 0<\gamma<\frac{p-1}{2-p},\qquad 0<c<2\gamma,
\]
there are $C$ and $\bar s$ such that, for every integer $j\ge\bar s$
and every real $s\ge j$,
\begin{equation}\label{eq:quantitative-fixed-marking}
 [\Psi^*g_s-\bar g]_k+[\Psi^*f_s-\bar f]_k
 \le C e^{\bar f/2}(1+s)^{-\gamma}
 \quad\text{on }\Omega_j:=\{\bar f\le c\log j\}.
\end{equation}
Here $[T]_k=\sum_{\ell=0}^k|\bar\nabla^\ell T|_{\bar g}$.
The constant $C$ depends on the fixed background data and on
$p,k,\gamma,c$; the waiting time $\bar s$ may also depend on the
flow, the singular point and the initial normalization of the
modified flow. The same embedding works for all these parameters.
\end{proposition}
\begin{proof}
Choose $\gamma<\gamma_1<(p-1)/(2-p)$ and a compatible auxiliary
construction with $q=p_1p_3$ such that
\begin{equation}\label{eq:quantitative-parameter-margins}
 \frac{p_2+q-2}{2-q}>\gamma_1,\qquad
 \frac{\rho^2p_3^2q}{10c_2(2-q)}>2\gamma_1.
\end{equation}
Indeed, as $p_1\uparrow p$, $p_3\uparrow1$, $\rho\uparrow1$ and
$c_2\uparrow1/8$, the two left sides tend to
\[
 \frac{p-1}{2-p},\qquad
 \frac{4p}{5(2-p)}>\frac{2(p-1)}{2-p}
       \quad(1<p\le3/2).
\]
The strict inequalities in Stages~\ref{stage:entropy-exponents}--\ref{stage:gaussian-rates}
are compatible with these limits. Choose
\[
 2\gamma_1<c'
       <\frac{\rho^2p_3^2q}{10c_2(2-q)}.
\]
Prescribe the finite rough order and Gaussian loss required for
Lemma~\ref{lem:weighted-intrinsic-velocity} at order $k$ before choosing
the comparison time, metric accuracy and final thresholds, in the order
of Subsection~\ref{subsec:parameter-selection}. When fixing $\epsilon_{\rm sl}$,
also impose the accuracy of Lemma~\ref{lem:quantitative-sublevel-coverage}.
The input order remains
$100n$. Initialize along the selected tangent sequence and apply
Proposition~\ref{prop:continuation}. This auxiliary construction supplies
intrinsic estimates for the same modified flow in the fixed embedding
$\Psi$ of Proposition~\ref{prop:fixed-limit}.

Lemmas~\ref{lem:quantitative-decay}--\ref{lem:quantitative-radii} and
\ref{lem:quantitative-sublevel-coverage} place the entire sublevel
$\{f_s\le c'\log(1+s)\}$ in the comparison chart, with the collars
required for
\eqref{eq:weighted-intrinsic-velocity}, and
\begin{equation}\label{eq:quantitative-length-input}
 I_{p_2}(s)=\int_s^\infty\mathcal D(u)^{p_2}\,\mathrm du
       \le C(1+s)^{-\gamma_1}.
\end{equation}
The comparison between block indices and modified time holds beyond
twice the auxiliary starting time, so this starting time affects only
the waiting threshold.

Put $G=\Psi^*g_s$, $F=\Psi^*f_s$, and $h=G-\bar g$.
The density identity~\eqref{eq:M-20} gives, as in \cite[(6.78)]{LYZII},
\[
 F-\bar f=\frac12\log\det(\bar g^{-1}G).
\]
Fix $j$. By smooth convergence on the compact set $\Omega_j$, one has
$\sup_{\Omega_j}[h]_k<1/2$ on a terminal interval in $[j,\infty)$.
There the eigenvalues of $\bar g^{-1}G$ lie in $[1/2,3/2]$.
Differentiating the determinant identity gives
$[F-\bar f]_k\le C_{n,k}[h]_k$, and its zeroth-order bound gives
$F\le\bar f+C_n\le c\log j+C_n<c'\log(1+s)$ for large $j$.

The intrinsic estimate~\eqref{eq:weighted-intrinsic-velocity} thus
applies on $\Psi(\Omega_j)$. In converting derivatives to $\bar\nabla$,
the connection difference $\nabla^G-\bar\nabla$ is a contraction of
$G^{-1}$ with $\bar\nabla h$, so the metric bound controls its derivatives
through order $k-1$ when $k\ge1$. Consequently
$[\partial_sG]_k=2[\Psi^*\Phi]_k\le C_k e^{\bar f/2}\mathcal D(s)^{p_2}$.
Integrating to infinity and using the determinant identity yields
\begin{equation}\label{eq:fixed-weighted-length}
 [h]_k+[F-\bar f]_k
       \le C_k e^{\bar f/2}I_{p_2}(s)
       \le C_k e^{\bar f/2}(1+s)^{-\gamma_1}.
\end{equation}
On $\Omega_j$ the right-hand side is at most $C_kj^{c/2-\gamma_1}<1/4$
for large $j$. If the maximal terminal interval had left endpoint
greater than $j$, continuity of $\sup_{\Omega_j}[h]_k$ would give
the value $1/2$ there, contradicting this bound. Hence
\eqref{eq:fixed-weighted-length} holds for all $s\ge j$.
Since $\gamma_1>\gamma$, this proves~\eqref{eq:quantitative-fixed-marking}.
\end{proof}

\subsection{Time-preserving reconstruction}
\label{part9-complete-negative-time-measured-limit-and-calibration}

An ODE gives time-independent spatial identifications on each
rescaled interval.

\begin{proposition}[Uniqueness of tangent flow]\label{prop:time-preserving-reconstruction}
Let $\Psi$ be the embedding obtained above. For every sequence
$\tau_i\downarrow0$, put $\sigma_i:=-\log\tau_i$ and
$F_i:=\Xi_{\sigma_i}\circ\Psi$. On every compact product region and
every compact interval $[a_-,a_+]\Subset(0,\infty)$,
\[
 F_i^*(\tau_i^{-1}g(-\tau_i a))\longrightarrow ag_N+g_E,
 \qquad F_i^*f(-\tau_i a)\longrightarrow
 f_N+\frac{|\vec x|^2}{4a}
\]
smoothly in space and time. The pulled-back terminal probabilities
converge to
$\mu_N\otimes(4\pi a)^{-d/2}e^{-|\vec x|^2/(4a)}\,\mathrm d\vec x$.
These time-independent spatial maps exhaust every bounded-distance
region of the rescaled flows and identify their entire negative-time
measured limit with the negative-time restriction of $\mathcal C_N^d$, preserving its
terminal probability and time vector. The induced identification of spacetime completions
identifies every measured tangent flow with $\mathcal C_N^d$.
\end{proposition}
\begin{proof}
Enlarge the interval to contain $1$, and set
$G(s)=\Psi^*g_s$, $F(s)=\Psi^*f_s$, and
$s_i(a):=\sigma_i-\log a$. Solve
\begin{equation}\label{eq:M-18}
 \frac{\mathrm d}{\mathrm da}\varphi_i(a)
 =-a^{-1}\nabla^{G(s_i(a))}F(s_i(a))\circ\varphi_i(a),
 \qquad\varphi_i(1)=\Id.
\end{equation}
The modified-flow equations upgrade compact smooth convergence of
$G,F$ to space-time convergence. The vector fields in~\eqref{eq:M-18}
therefore converge smoothly to $-a^{-1}\bar\nabla\bar f$, whose flow is
\[
 \varphi_\infty(a)(y,\vec x)=(y,a^{-1/2}\vec x).
\]
The reference trajectories of a fixed core lie in a larger compact
product region with a positive collar. Smooth dependence gives existence
throughout the interval and $\varphi_i\to\varphi_\infty$ smoothly in
space and time.

The gradient transformation and~\eqref{eq:modified-diffeomorphism} give
\[
 \mathrm D(\Xi_s\Psi)\nabla^{G(s)}F(s)
 =\nabla^{\tilde g(s)}\tilde f(s)\circ\Xi_s\Psi,
\]
and hence
\[
 \begin{split}
 \frac{\mathrm d}{\mathrm da}
       \bigl(\Xi_{s_i(a)}\Psi\varphi_i(a)\bigr)
 &=a^{-1}\nabla^{\tilde g(s_i(a))}\tilde f(s_i(a))
                    \circ\Xi_{s_i(a)}\Psi\varphi_i(a)\\
 &\quad-a^{-1}\mathrm D(\Xi_{s_i(a)}\Psi)
       \nabla^{G(s_i(a))}F(s_i(a))\circ\varphi_i(a)=0.
 \end{split}
\]
Its value at $a=1$ is $F_i$, so
$\Xi_{s_i(a)}\Psi\varphi_i(a)=F_i$ exactly. Since
$e^{-s_i(a)}=\tau_i a$, this proves
\begin{equation}\label{eq:M-19}
 \begin{split}
 F_i^*(\tau_i^{-1}g(-\tau_i a))
   &=a\varphi_i(a)^*G(s_i(a))\longrightarrow ag_N+g_E,\\
 F_i^*f(-\tau_i a)
   &=F(s_i(a))\circ\varphi_i(a)
       \longrightarrow f_N+\frac{|\vec x|^2}{4a}.
 \end{split}
\end{equation}
For $\mu_i(a)=\mu_{-\tau_i a}$, the corresponding exact identity is
\begin{equation}\label{eq:reconstruction-probability}
  F_i^*\mathrm d\mu_i(a)
 =(4\pi a)^{-n/2}e^{-F(s_i(a))\circ\varphi_i(a)}
          \,\mathrm dV_{a\varphi_i(a)^*G(s_i(a))}
 =\varphi_i(a)^*\mathrm d\bar\mu.
 \end{equation}
This gives the stated Gaussian limit. Since $F_i$ is independent of
$a$, the maps $(x,-a)\mapsto(F_i(x),-a)$ preserve the time vector.

Put $g_i(a):=\tau_i^{-1}g(-\tau_i a)$ and $p_i:=F_i(y_0,0)$.
On every fixed larger core, $F_i^*g_i(a)$ is uniformly close to
$ag_N+g_E$. The first-exit argument of Lemma~\ref{lem:coverage},
whose radial lower bound is uniform for $a\in[a_-,a_+]$, gives
\[
 B_{g_i(a)}(p_i,B)\subset F_i(N_R)
\]
for each fixed $B$, sufficiently large fixed $R$, and all large $i$.
A larger collar keeps minimizing paths between inner points in the
controlled region. Thus the comparisons give ambient distances and
exhaust the bounded-distance regions.

Mass one and~\eqref{eq:reconstruction-probability} give
\[
 \mu_i(a)(M\setminus F_i(N_R))
 =1-\int_{N_R}\varphi_i(a)^*\mathrm d\bar\mu.
\]
The limiting integrals are at least
$\bar\mu(N_{R/\sqrt{a_+}})$, uniformly in $a$, and tend to one as
$R\to\infty$. This proves probability tightness on the whole time
interval. Positive mass on a fixed core and the second-moment bound
also keep $p_i$ a uniformly bounded distance from the normalized
$H_n$-centers.

The reference slabs are complete with bounded curvature. The coverage,
tightness and time-vector identity allow the identification argument
in the proof of \cite[Corollary~6.25]{LYZII} to identify the regular
spacetime charts, their conjugate kernels and the completion.
Exhausting compact negative-time intervals identifies the negative-time measured flow;
passing to its spacetime completion identifies $\mathcal C_N^d$, with its distinguished
terminal probability.
Since $\tau_i\downarrow0$ was arbitrary, this holds for every
measured tangent flow.
\end{proof}

\begin{proof}[Completion of the proof of Theorem~\ref{thm:main}]
The entropy limit~\eqref{eq:entropy-limit},
Proposition~\ref{prop:continuation} and~\eqref{eq:length-tail} supply
the hypotheses of Proposition~\ref{prop:fixed-limit}, which gives
\eqref{eq:main-convergence}. Proposition~\ref{prop:quantitative-fixed-gauge}
gives all derivative orders and admissible rates in this embedding.
Proposition~\ref{prop:time-preserving-reconstruction} identifies the
measured tangent flow, including its spacetime completion.
\end{proof}

\subsection{Finite quotients and ancient solutions}\label{sec:finite-quotients}
Let $\Gamma\le\operatorname{Isom}(N,g_N)\times O(d)$ be finite and act
freely on $\bar M=N\times\R^d$. Write
$\pi:\bar M\to\bar M_\Gamma:=\bar M/\Gamma$ and normalize the quotient by
\begin{equation}\label{eq:quotient-normalization}
 \pi^*\bar g_\Gamma=\bar g,\qquad
 \pi^*\bar f_\Gamma=\bar f-\log|\Gamma|,\qquad
 \bar\mu_\Gamma=\pi_*\bar\mu,\qquad
 \Theta_\Gamma:=\Theta_N-\log|\Gamma|.
\end{equation}
The Euclidean action may be nontrivial, so the quotient need not be a
product. If that action is trivial, it is $(N/\Gamma)\times\R^d$.
The following results use the quotient and ancient-flow frameworks of
\cite[Theorems~8.15--8.16 and~8.18--8.19]{FLC} and the generalized quotient formulation
of \cite[Theorem~4.16]{FLQ}.

\begin{theorem}[Finite quotients]\label{thm:quotient-uniqueness}
Under the flow hypotheses of Theorem~\ref{thm:main}, suppose that one
measured tangent flow at $z_*$ is the canonical shrinking flow associated
with $(\bar M_\Gamma,\bar g_\Gamma,\bar f_\Gamma)$. Every measured tangent
flow at $z_*$ is isometric to this flow, preserving terminal probability
and the time vector. There is a smooth, time-independent open embedding
$\Psi_\Gamma:\bar M_\Gamma\to M$ such that
$\Psi_\Gamma^*\mu^{z_*}=\bar\mu_\Gamma$ and the pulled-back modified pair
converges smoothly on compact sets to $(\bar g_\Gamma,\bar f_\Gamma)$.
With the product exponent $p$ of Theorem~\ref{thm:logarithmic}, this
embedding satisfies the following estimate for every integer $k\ge0$,
$0<\gamma<(p-1)/(2-p)$ and $0<c<2\gamma$. There are $C,\bar s<\infty$
such that, for every integer $j\ge\bar s$ and every $s\ge j$,
\begin{equation}\label{eq:quotient-quantitative}
 [\Psi_\Gamma^*g^{z_*}(s)-\bar g_\Gamma]_k
 +[\Psi_\Gamma^*f^{z_*}(s)-\bar f_\Gamma]_k
 \le Ce^{\bar f_\Gamma/2}(1+s)^{-\gamma}
 \quad\text{on }\{\bar f_\Gamma\le c\log j\}.
\end{equation}
Here $[T]_k:=\sum_{\ell=0}^k
|\nabla_{\bar g_\Gamma}^{\ell}T|_{\bar g_\Gamma}$.
The constants have the dependence stated in Theorem~\ref{thm:main},
with additional dependence on the fixed quotient.
\end{theorem}

\begin{proof}
The normalized lift of a probability-normalized pair
$X_\Gamma=(g_\Gamma,f_\Gamma)$ is
$\widetilde X:=(\pi^*g_\Gamma,\pi^*f_\Gamma+\log|\Gamma|)$, with
$\mathrm d\widetilde\mu=|\Gamma|^{-1}\pi^*\mathrm d\mu_\Gamma$.
Covering integration gives
\[
 \mathcal W(\widetilde X)-\Theta_N
   =\mathcal W(X_\Gamma)-\Theta_\Gamma,\qquad
 \|\Phi(\widetilde X)\|_{L^2(\widetilde\mu)}
   =\|\Phi(X_\Gamma)\|_{L^2(\mu_\Gamma)}.
\]
The slice projections commute with $\Gamma$, and
Proposition~\ref{prop:kernel} gives
\[
 K^\Gamma=\bigoplus_{j\ge0}
       \bigl(\mathcal H_j\otimes\mathscr E_{j/2}(N)\bigr)^\Gamma.
\]
Thus invariants are taken after forming the tensor products; the radial
direction $(|\vec x|^2-2d)g_N$ remains invariant. The unique normal graph
and the finite-order reduction in Subsection~2.2 are equivariant.
The scalar corrections restoring mass and first moment may be averaged
over $\Gamma$, since the constraints are linear and their prescribed
values are invariant.
The reduced functions and parameter norms are invariant. Their full
gradients at invariant parameters therefore belong to the invariant
subspaces, so restricting the gradient inequalities does not change their
norms or the admissible exponent $p$. The stationary estimates of
Theorem~\ref{thm:logarithmic} consequently descend.

Write $\rho_E:\Gamma\to O(d)$ for the Euclidean action. The first moment
of an invariant probability belongs to
$V:=(\R^d)^{\rho_E(\Gamma)}$. Restrict translations to $V$ and reference
Killing fields to $\mathcal K_0^\Gamma$. The gauge equations, including
their normalization rows, commute with $\Gamma$; uniqueness restricts
the inverses of Appendix~\ref{sec:inverse} and
Lemma~\ref{lem:augmented-weak-inverse} to invariant data. Radial cutoffs,
dilations and reference exponential maps preserve equivariance, as do
the forward comparisons by uniqueness. The reference-group equation
of Lemma~\ref{lem:reference-group-normalization} starts at $\Id$ in the
identity component of the centralizer of $\Gamma$. The endpoint
projection commutes with $\Gamma$, and the inherited maps in
Subsection~\ref{subsec:inherited-endpoint-coordinates} remain equivariant.
Every chart change therefore descends.

Use $N_R/\Gamma$ as comparison domains. The local calculations in
Sections~\ref{sec:physical}--\ref{sec:motion} lift to their product covers;
the exterior probability estimates remain on the ambient flow.
On a block starting at $t=-1$, the normalized lift adds $|t|\log|\Gamma|$
to $F=|t|f$ and the constant $\log|\Gamma|$ to its forward caloric
extension $\tilde F$ by~\eqref{eq:C-4}. Thus the spatial derivatives
of $u=\tilde F-F$ lift unchanged. We apply the terminal Harnack inequality
and use the mass-one identity for the actual potential downstairs.
The normalized covering measure accounts for multiplicity, and the
collars, displacement margins, and radius recurrence are unchanged.
For $b\in K_N^\Gamma$ the invariant pairings still give
$\mathcal O(G_b)=b$. Proposition~\ref{prop:continuation} therefore confines
the compact parameter and gives continuation. The decay and radius estimates of
Subsection~6.4 then apply. Subsections~8.1--8.4 give the fixed embedding,
the exact probability identity, all orders and rates in that embedding,
and time-preserving reconstruction, as in
\cite[Theorems~8.15--8.16]{FLC}.
Only the comparison domains require covers. The constant potential
shift changes the radial offset and the multiplicative constants in
the weighted estimate, but not its exponents.
\end{proof}

\begin{theorem}[Ancient solutions]\label{thm:ancient-uniqueness}
Let $(M^n,g(t))$, $-\infty<t\le0$, be a complete ancient Ricci flow
with bounded curvature on every compact time interval. Suppose that its
pointed Nash entropy $\mathcal N_{x,t}(\tau)$ is at least $-Y$ at every regular
base point and every positive backward scale, in the sense of
\cite[Definitions~2.7 and~2.20]{FLSt}. Fix $(x_*,0)$ and its conjugate
heat kernel, and suppose that one measured tangent flow at infinity
based there is the canonical shrinking flow of
$(\bar M_\Gamma,\bar g_\Gamma,\bar f_\Gamma)$.

Use $s=-\log(-t)$ and define $(g^{x_*}(s),f^{x_*}(s))$ and its fixed
probability $\mu^{x_*}$ by
\eqref{eq:modified-diffeomorphism}--\eqref{eq:modified-gauge}, with
the kernel based at $(x_*,0)$ and $\Xi_0=\Id$, corresponding to $t=-1$.
There is one smooth, time-independent open embedding
$\Psi_\Gamma:\bar M_\Gamma\to M$ such that
\begin{equation}\label{eq:ancient-convergence}
 \Psi_\Gamma^*g^{x_*}(s)\longrightarrow\bar g_\Gamma,\qquad
 \Psi_\Gamma^*f^{x_*}(s)\longrightarrow\bar f_\Gamma
 \quad\text{in }C^\infty_{\rm loc},\qquad
 \Psi_\Gamma^*\mu^{x_*}=\bar\mu_\Gamma
 \quad(s\longrightarrow-\infty).
\end{equation}
For every integer $k\ge0$, $0<\gamma<(p-1)/(2-p)$ and $0<c<2\gamma$,
there are $C,\bar s<\infty$ such that, for every integer $j\ge\bar s$
and every $s\le-j$,
\begin{equation}\label{eq:ancient-quantitative}
 [\Psi_\Gamma^*g^{x_*}(s)-\bar g_\Gamma]_k
 +[\Psi_\Gamma^*f^{x_*}(s)-\bar f_\Gamma]_k
 \le Ce^{\bar f_\Gamma/2}(1+|s|)^{-\gamma}
 \quad\text{on }\{\bar f_\Gamma\le c\log j\}.
\end{equation}
The embedding is independent of $k,\gamma,c$. The constant $C$ has
the dependence in Theorem~\ref{thm:quotient-uniqueness}; the threshold
may also depend on the ancient flow, the base point, and the finite-time
normalization of $\Xi$. Every measured tangent flow at infinity based
at $(x_*,0)$ is the canonical quotient flow, preserving its probability
and time vector.
\end{theorem}

\begin{proof}
Write $(g_s,f_s):=(g^{x_*}(s),f^{x_*}(s))$ and use $\mathcal W(s)$
and $\mathcal D(s)$ for its entropy and defect. We use the ancient modified-flow formulation and argument of
\cite[(7.14) and Theorems~7.10--7.11, 8.18--8.19]{FLC}, with
the following changes to the preceding proof.

The complete-flow estimates recalled in \cite[Section~11, before
Theorem~11.2]{FLSt} give the Gaussian and entropy tails in
\eqref{eq:terminal-tail}. Smooth measured convergence follows from
\cite[Theorem~11.3 and Example~11.5]{FLSt}. The scalar curvature is nonnegative by~\cite[Corollary~2.5]{Chen09}.

For~\eqref{eq:C-4}, add $C+n(t+1)/2$ to both $F$ and $\tilde F$ so
that the initial datum is nonnegative and the forward equation is
homogeneous. Evolve increasing smooth truncations of $F(-1)+C$ by the
minimal heat kernel. Take the truncations to have slopes in $[0,1]$,
to equal the identity up to $R$, and to be constant above $2R$.
Each datum is constant outside a compact set by the Gaussian lower
bound for $f$. The shifted terminal supersolution
bounds the evolutions above, so their limit, after removing the shift,
solves~\eqref{eq:C-4}. For each truncation the gradient and Hessian
evolution equations bound the forward derivatives of orders one and two
on the slab. In the calculation of Lemma~\ref{thm:caloric} the spatial
cutoff errors involve $|\mathrm du|^2$ and $|\tilde T|\,|\mathrm du|$
for this truncation. Their annular integrals vanish by
\eqref{eq:terminal-tail}, and $\int_{-1}^t\!\int |T|^2\,
\mathrm d\mu_r\,\mathrm dr$ is finite by entropy monotonicity.
The initial gradients converge in $L^2(\mu_{-1})$ by dominated convergence.
For the difference $v$ of two truncated evolutions the caloric Bochner
identity gives
\[
 \int|\mathrm dv(t)|^2\,\mathrm d\mu_t
 +2\int_{-1}^t\!\int|\nabla^2v|^2\,\mathrm d\mu_r\,\mathrm dr
 =\int|\mathrm dv(-1)|^2\,\mathrm d\mu_{-1}.
\]
Hence the gradients and Hessians converge in the norms of
\eqref{eq:G-2-1}, and the initial terms of the truncated identities tend to zero.
This gives~\eqref{eq:G-2-1} with the same constants on every block.

The selected tangent sequence and the entropy tails identify the
entropy limit with $\Theta_\Gamma$. Hence
\begin{equation}\label{eq:ancient-entropy}
 E(s):=\mathcal W(s)-\Theta_\Gamma
 =2\int_{-\infty}^s\mathcal D(u)^2\,\mathrm du\ge0,
 \qquad E'(s)=2\mathcal D(s)^2.
\end{equation}
Fix $S<-1$ with $E(S)$ small. Starting at selected times
$s_i\to-\infty$, choose initial radii tending slowly to infinity
so that~\eqref{eq:E-9} holds with $b_0=0$, and both the initial
pairing and $\Delta_{-1}$ tend to zero. Construct blocks forward
toward $S$, stopping at a possible parameter exit. The estimate
$A_j(u)^2\le E(s_j+u)\le E(S)$ supplies the small-drop condition
of Proposition~\ref{prop:causal-radii}.
Lemma~\ref{lem:restart-entropy-localization} and the proof of
Proposition~\ref{prop:restart-stationary-comparison} bound the absolute
entropy difference by
$Ce^{-c_2p_1R_j^2}+Ce^{-3R_j^2/16}\le C\Delta_j^{p_1}$.
Thus the same charts, collars and order-$100n$ restart give
\begin{equation}\label{eq:ancient-recurrence}
 E_{j+1}-E_j=A_j^2,\qquad
 \Delta_j=\max\{\Delta_{j-1}^{3/2},A_j^{p_3}\},\qquad
 E_{j+1}\le C\Delta_j^{p_1}.
\end{equation}

Put $q=p_1p_3$ and use $\vartheta=1-(2-p_3p_4)/q$ from
Lemma~\ref{lem:delayed-sum}. At indices with
$\Delta_j=A_j^{p_3}>0$, concavity gives
\[
 E_{j+1}^\vartheta-E_j^\vartheta
 \ge\vartheta E_{j+1}^{\vartheta-1}A_j^2
 \ge c\Delta_j^{p_4}.
\]
Sum the intervening inherited-radius strings as in that lemma to obtain
\begin{equation}\label{eq:ancient-finite-sum}
 \sum_{j<J}(A_j^{p_4}+\Delta_j^{p_4})
 \le C\bigl(E_J^\vartheta+\Delta_{-1}^{p_4}\bigr).
\end{equation}
A final partial block contributes at most $CE(S)^{p_4/2}$.
The proof of Proposition~\ref{prop:comparison-length} now bounds the
total change of the pairings by $\omega(E(S)+\Delta_{-1})$. Choose
$S$ and then the initial chart so that this bound and the initial
pairing together are less than $r_{\rm con}/8$. The endpoint comparison
\eqref{eq:modified-endpoint-estimates} gives the other
$r_{\rm con}/8$ bound in the proof of Proposition~\ref{prop:continuation}.
Thus every such construction reaches $S$.

Let $J$ count the complete blocks before $S$. For decay, set $x_j=\Delta_j^{2/p_3}$, $\delta=x_{-1}$ and
$F_j=E_j-x_{j-1}/2$. The recurrence gives
$x_{j-1}\le\max\{\delta,E_j\}$ and, for sufficiently small data,
\[
 F_{j+1}-F_j
 \ge\frac12A_j^2+\frac14x_{j-1}
 \ge cE_{j+1}^{2/q}.
\]
Indeed $x_j\le x_{j-1}^{3/2}+A_j^2$ and
$E_{j+1}^{2/q}\le Cx_j$. For $G_j=F_j+\delta/2$ we have
$E_j/2\le G_j\le E_j+\delta/2$ and
$G_{j+1}-G_j\ge c(G_{j+1}-\delta/2)_+^{2/q}$.
Apply the inverse-power calculation of
Lemma~\ref{lem:quantitative-decay} with reversed indices until the reversed sequence falls below
$\delta$; it remains below $\delta$ thereafter. This gives
\begin{equation}\label{eq:ancient-finite-decay}
 E_j\le C\bigl((1+J-j)^{-q/(2-q)}+\delta\bigr).
\end{equation}
For fixed $s<S-\log2$, take the first block endpoint at or after $s$.
Then $E(s)\le E_j$ and $J-j=(S-s)/\log2+O(1)$.
Sending the selected starting time to $-\infty$ makes $\delta\to0$,
so
\begin{equation}\label{eq:ancient-decay}
 E(s)\le C(1+|s|)^{-q/(2-q)}\qquad(s\le2S).
\end{equation}

The dyadic H\"older calculation of Lemma~\ref{lem:quantitative-decay},
now on intervals extending toward $-\infty$, yields
\begin{equation}\label{eq:ancient-length}
 \int_{-\infty}^s\mathcal D(u)^a\,\mathrm du
 \le C_a(1+|s|)^{-(a+q-2)/(2-q)},\qquad 2-q<a<1.
\end{equation}
On any compact interval in $(-\infty,2S]$, choose the starting chart so early
that $\delta$ is smaller than the right side of~\eqref{eq:ancient-decay}.
Using $x_{j-1}\le\max\{\delta,E_j\}$ and
$A_j(u)^2\le E(s_j+u)$, the two-branch calculation of
Lemma~\ref{lem:quantitative-radii} gives, throughout blocks and resets,
\begin{equation}\label{eq:ancient-radii}
 L_j(u)^2\ge
 \frac{\rho^2p_3^2q}{2c_2(2-q)}\log(1+|s_j|)
                    -C\sqrt{\log(2+|s_j|)}-C.
\end{equation}
The constants are independent of the starting time. Thus
Lemmas~\ref{lem:quantitative-sublevel-coverage}--\ref{lem:weighted-intrinsic-velocity}
give entire-sublevel coverage and the weighted velocity estimates with
$\log(1+s)$ replaced by $\log(1+|s|)$.

Apply the proof of Proposition~\ref{prop:fixed-limit}, integrating
absolute velocities from $-\infty$ instead of toward $+\infty$.
Equation~\eqref{eq:ancient-length} supplies the integrable factor,
and $\mathrm R\ge0$ supplies the scalar lower bound. At each fixed time
the Gaussian lower bound places every closed sublevel of $f_s$ in a
bounded ball of the complete metric $g_s$. These sublevels are compact
by Hopf--Rinow. The potential barriers therefore give a complete limit
on an open subset of $M$;
the selected measured tangent identifies it with the quotient model.
The fixed probability and mass-one argument give
\eqref{eq:ancient-convergence}, fixing $\Psi_\Gamma$ once and for all.

For the quantitative conclusion, choose $\gamma_1>\gamma$ and the
auxiliary parameters as in the proof of
Proposition~\ref{prop:quantitative-fixed-gauge}. Then
\eqref{eq:ancient-length} with $a=p_2$ decays faster than
$(1+|s|)^{-\gamma_1}$, and~\eqref{eq:ancient-radii} supplies the
same strict sublevel-coverage margin. On
$\Omega_j=\{\bar f_\Gamma\le c\log j\}$, that proof, integrating
from $-\infty$ in the already fixed embedding, gives
\[
 [\Psi_\Gamma^*g_s-\bar g_\Gamma]_k+
 [\Psi_\Gamma^*f_s-\bar f_\Gamma]_k
 \le C_ke^{\bar f_\Gamma/2}(1+|s|)^{-\gamma_1}.
\]
Initially this holds sufficiently far in the past. Its supremum is
at most $C_kj^{c/2-\gamma_1}<1/4$ for $s\le-j$ and large $j$,
so the same first-failure argument extends it to $s=-j$.
Weakening the exponent proves~\eqref{eq:ancient-quantitative}.

Finally apply the proof of
Proposition~\ref{prop:time-preserving-reconstruction} with
$\tau_i\to\infty$ and $\sigma_i=-\log\tau_i\to-\infty$.
The reconstruction ODE, probability identity and time-vector
cancellation are unchanged. The exhaustion and tightness arguments,
in the complete-flow framework of \cite[Theorem~11.3]{FLSt}, identify
every measured tangent flow at infinity with the canonical quotient
flow.
\end{proof}

\appendix
\refstepcounter{section}
\section*{Appendix \thesection. The gauge inverse and coordinate tails}
\label{sec:inverse}\label{part6-the-compatible-source-inverse}
\addcontentsline{toc}{section}{Appendix \thesection. The gauge inverse and coordinate tails}

The gauge equation corrects the metric divergence and fixes the
reference symmetries. We prove reference and family inverse estimates
for its global solution, then the stronger Gaussian estimates used
to cut it off in
Section~\ref{part8-low-slice-first-with-an-explicit-numerical-hypothesis}.

We use the reference operator $\bar{\mathcal P}$ and fixed symmetry
spaces $\mathcal K,\mathcal K_0$ from
Section~\ref{part8-low-slice-first-with-an-explicit-numerical-hypothesis},
identifying vector fields and one-forms using $\bar g$.
All orthogonality conditions below use $\bar\mu$ and these fixed spaces.
A source is reference-compatible if it annihilates $\mathcal K$ in the
$L^2(\bar\mu)$ pairing.

For a compactly supported symmetric tensor $T$, the source
$F=\delta_{\bar X}T$ satisfies
\[
 (F,X)_{L^2(\bar\mu)}=(T,\delta_{\bar X}^*X)_{L^2(\bar\mu)}=0
 \qquad(X\in\mathcal K).
\]
In particular, a cutoff is always taken inside the divergence:
\[
 \delta_{\bar X}(\chi T)
    =\chi\delta_{\bar X}T+(\mathrm d\chi)\lrcorner T.
\]

\begin{proposition}\label{prop:reference-gauge-inverse}
The reference-normalized solution of
$\bar{\mathcal P}V=F$, $V\perp\mathcal K$, satisfies
\begin{equation}\label{eq:G-1-11}
 \|V\|_{H^{j+2}(\bar\mu)}\le C_j\|F\|_{H^j(\bar\mu)},
 \qquad j=0,1,2,\ldots,
\end{equation}
for every reference-compatible source of the indicated regularity.
The constant $C_j$ depends only on $j$ and the reference product.
\end{proposition}
\begin{proof}
The spectral proof of \cite[Proposition~3.10]{LZ26} uses the divergence
form of the source only to establish orthogonality to the kernel.
It therefore gives~\eqref{eq:G-1-11} for every reference-compatible source.
\end{proof}

\begin{proposition}
\label{prop:family-reference-inverse}
Let $G_b=g_b+g_E$ and $\mu_b=\mu_b^N\otimes\gamma_d$ be the
fixed smooth compact family. There are a parameter neighborhood,
$S_{\rm inv}>0$, and $C>0$ such that the following holds for every
$b$ in that neighborhood and every $S\ge S_{\rm inv}$.
Let $0\le\chi_{1,S}\le1$ be a smooth radial cutoff equal to one on
$N_{S-2}$ and zero outside $N_{S-1}$. For every
$F\in L^2(\bar\mu)$ satisfying
$(F,X)_{L^2(\bar\mu)}=0$ for all $X\in\mathcal K$, and every
$\vec m\in\R^d$, there is a unique $V\in H^2(\bar\mu)$ solving
\begin{equation}\label{eq:family-reference-system}
 \mathcal P_bV:=-\frac12\delta_{\bar X}L_VG_b=F,
 \qquad \Pi_{\mathcal K_0}V=0,
 \qquad \int\chi_{1,S} V_E\,\mathrm d\mu_b=\vec m.
\end{equation}
It satisfies
\begin{equation}\label{eq:family-reference-bound}
 \|V\|_{H^2(\bar\mu)}
       \le C\bigl(\|F\|_{L^2(\bar\mu)}+|\vec m|\bigr).
\end{equation}
The constants are uniform for $S\ge S_{\rm inv}$. The parameter
neighborhood is chosen using $g_b\to g_N$ in $C^2(N)$ and
$\mathrm d\mu_b^N/\mathrm d\mu_N\to1$ in $C^0(N)$ as $b\to0$,
independently of $S$.
\end{proposition}
\begin{proof}
Lower the unknown with the family metric:
$\omega=G_b(V,\cdot)$, and put
$C_b:=\nabla^{G_b}-\bar\nabla$. The covariant derivative formula
gives
\begin{equation}\label{eq:G-1-4}
 \begin{split}
 L_VG_b&=2\operatorname{sym}\bar\nabla\omega-2C_b\omega,\\
 \mathcal P_bV&=\bar{\mathcal P}\omega+B_b\omega,
 \qquad B_b\omega:=\delta_{\bar X}(C_b\omega).
 \end{split}
\end{equation}
Here $(C_b\omega)_{ij}=(C_b)^a{}_{ij}\omega_a$ is symmetric.
All indices of $C_b$ are compact and its coefficients are independent
of $\vec x$. Thus the Euclidean drift contraction vanishes, and
$B_b$ contains at most one compact derivative, giving
\[
 \|B_b\omega\|_{L^2(\bar\mu)}
       \le C|b|\|\omega\|_{H^2(\bar\mu)}.
\]
Integration by parts against $\delta_{\bar X}^*X=0$ makes
$B_b\omega$ reference-compatible. Cutoff approximation extends this
identity to $H^2$ fields, since the errors vanish by Gaussian tails.

Let $L^2_{\mathcal K^\perp}$ denote this fixed compatible target.
Using $\bar g$ to identify one-forms and vector fields, the reference
augmented operator of Lemma~\ref{lem:augmented-weak-inverse} is
\[
 \mathscr P\omega
 =\left(\bar{\mathcal P}\omega,
        \Pi_{\mathcal K_0}(\bar g^{-1}\omega),
        \int(\bar g^{-1}\omega)_E\,\mathrm d\bar\mu\right).
\]
By Proposition~\ref{prop:reference-gauge-inverse}, it is an isomorphism
$H^2(\bar\mu)\longrightarrow
L^2_{\mathcal K^\perp}\times\mathcal K_0\times\R^d$.
The first component is inverted on $\mathcal K^\perp$; the remaining
components determine its Killing coefficients. The first-moment
component is the identity on constant translations and vanishes on
$\mathcal K_0$.

For the family metric and cutoff first moment, use
\[
 \mathscr P_{b,S}\omega
 :=\left(\bar{\mathcal P}\omega+B_b\omega,
        \Pi_{\mathcal K_0}(G_b^{-1}\omega),
        \int\chi_{1,S}(G_b^{-1}\omega)_E\,\mathrm d\mu_b\right).
\]
Raising indices changes the $\mathcal K_0$ component by $O(|b|)$.
The Euclidean component is unchanged, while
$\|\mathrm d\mu_b^N/\mathrm d\mu_N-1\|_{C^0}\le C|b|$;
thus the full first-moment component also changes by $O(|b|)$.
For each fixed $0<a<1/8$, Cauchy--Schwarz gives
\[
 \left|\int(1-\chi_{1,S})(G_b^{-1}\omega)_E\,\mathrm d\mu_b\right|
 \le C\|\omega\|_{L^2(\bar\mu)}\bar\mu(\bar M\setminus N_{S-2})^{1/2}
 \le C_a e^{-a(S-2)^2}\|\omega\|_{H^2(\bar\mu)}.
\]
The strict exponent absorbs the Gaussian-tail polynomial. Hence
\[
 \|\mathscr P_{b,S}-\mathscr P\|
       \le C|b|+C_a e^{-a(S-2)^2}.
\]
Choose $S_{\rm inv}$ large and then $|b|$ small so that a Neumann series
inverts $\mathscr P_{b,S}$ uniformly. Apply it to $(F,0,\vec m)$ and
raise indices with $G_b^{-1}$ to obtain the assertion.
\end{proof}

\begin{theorem}[Gauge inverse and tail estimates]\label{thm:gauge-inverse}
Fix $0<\beta<1$, $k\ge3$, and uniformly
bounded cutoff derivatives through the required finite orders, and let
$\mu_\beta$ be the Gaussian measure of Definition~\ref{def:gauge-reference-data}.
There are $b_0>0$, $S_{\rm inv}$, $C=C(\beta,k)$ and $Q=Q(\beta,k)$,
depending also on the fixed smooth compact family and cutoff bounds,
such that the following holds
for $|b|\le b_0$ and $S\ge S_{\rm inv}$.
For every symmetric $T\in C^{k,1/2}$ supported in $N_{S-3}$ and
$\vec m\in\R^d$, the system
\begin{equation}\label{eq:G-1-1}
 \mathcal P_bV=\delta_{\bar X}T,\qquad
 V\perp_{\bar\mu}\mathcal K_0,\qquad
 \int\chi_{1,S} V_E\,\mathrm d\mu_b=\vec m
\end{equation}
has a unique solution in $H^2(\bar\mu)$. It satisfies
\begin{equation}\label{eq:G-1-2}
 \|V\|_{H^2(\mu_\beta)}
 \le C(\beta)\bigl(\|\delta_{\bar X}T\|_{L^2(\mu_\beta)}+|\vec m|\bigr)
 \le C(\beta)\bigl(\|T\|_{C^1}+|\vec m|\bigr)
\end{equation}
and
\begin{equation}\label{eq:G-1-3}
 \|V\|_{H^2(\bar\mu)}+\|V\|_{C^{k+1,1/2}(N_S)}
 \le C(\beta,k)(1+S)^Q e^{\beta S^2}
                      \bigl(\|T\|_{C^{k,1/2}}+|\vec m|\bigr).
\end{equation}
The parameter neighborhood and constants are independent of $S$.
For a multiplier $B$ with
$|B|+|\bar\nabla B|\le C_B(1+r)^{Q_B}$, where $r=|\vec x|$, and every
$0<a<(1-\beta)/8$,
\begin{equation}\label{eq:gauge-coordinate-tail}
 \|BV\|_{H^1(\bar\mu;\,r\ge L)}
 \le C(\beta,a,B)e^{-aL^2}\|V\|_{H^2(\mu_\beta)},
 \qquad L\ge0.
\end{equation}
\end{theorem}

\begin{proof}
Write $\|\cdot\|_\beta$ for the $L^2(\mu_\beta)$ norm.
We first invert $\bar{\mathcal P}$ in $H^2(\mu_\beta)$, then perturb
it by $B_b$ in this space.

For $0<\beta<1$ and reference-compatible $F$, let
$\bar{\mathcal P}\omega=F$ with $\omega\perp_{\bar\mu}\mathcal K$.
We first prove
\[
 \|\omega\|_{H^2(\mu_\beta)}
       \le C(\beta)\|F\|_{L^2(\mu_\beta)}.
\]
On functions and reference one-forms,
\[
 -\Delta_{\bar f}=-\Delta_N-\Delta_{\vec x}
                         +\frac12\vec x\cdot\partial_{\vec x},
\]
where $\Delta_N$ is the connection Laplacian on the indicated bundle.
Integration by parts in the Euclidean variables gives
\[
 \|rv\|_\beta\le C(\beta)\|v\|_{H^1(\mu_\beta)},\qquad
 \|r\bar\nabla v\|_\beta+\|r^2v\|_\beta
       \le C(\beta)\|v\|_{H^2(\mu_\beta)}.
\]
Thus the form domain $\{v:v,\bar\nabla v,rv\in L^2(\mu_\beta)\}$
is $H^1(\mu_\beta)$, with equivalent norms.
Multiplication by a constant times $e^{-\beta r^2/8}$ is unitary
from $L^2(\mu_\beta)$ to $L^2(\mathrm d\mu_N\,\mathrm d\vec x)$ and
conjugates $-\Delta_{\bar f}$ to
\begin{equation}\label{eq:G-1-6}
 -\Delta_N-\Delta_{\vec x}
       +\frac{1-\beta}{2}
                   (\vec x\cdot\partial_{\vec x}+d/2)
       +\frac{\beta(2-\beta)}{16}r^2-\frac d4.
\end{equation}
Unmarked norms in the next display use the unweighted product measure.
For real compactly supported fields,
\[
 \bigl(e^{-\beta r^2/8}(-\Delta_{\bar f})(e^{\beta r^2/8}v),v\bigr)
   =\|\bar\nabla v\|^2
      +\frac{\beta(2-\beta)}{16}\|rv\|^2-\frac d4\|v\|^2.
\]
The same identity holds for the formal adjoint. The coefficients
satisfy the bounds of \cite[Lemma~4.1]{LYZII} with weight
$(1+r^2)^{1/2}$. Bounded product geometry gives uniformly controlled
rescaled charts on balls of radius comparable to $(1+r)^{-1}$ and
a cover with bounded overlap. The truncated multiplier argument and
local elliptic estimates in that lemma therefore give, for this operator
and its adjoint, the common minimal and maximal domain
\[
 \{v:v,\bar\nabla^2v,r\bar\nabla v,r^2v\in L^2\},
\]
with the corresponding graph estimate and compact inclusion into
$L^2$.

The Gaussian multiplication argument of \cite[Lemma~4.5]{LYZII},
applied with $1-\kappa=\beta$, identifies this domain with
$H^2(\mu_\beta)$: use the preceding multiplication estimates and
the product rule for the conjugation and its inverse. The form on
$H^1(\mu_\beta)$ is bounded and coercive after a sufficiently large
shift. Its inclusion into $L^2(\mu_\beta)$ is compact by local Rellich
compactness and the uniform tail bound
$\|v\|_{L^2(\mu_\beta;\,r>R)}\le R^{-1}\|rv\|_\beta$.
Consequently every fixed shift is Fredholm of index zero, both from
$H^1(\mu_\beta)$ to its dual and from $H^2(\mu_\beta)$ to
$L^2(\mu_\beta)$, and
\begin{equation}\label{eq:G-1-8}
 \|v\|_{H^2(\mu_\beta)}
 \le C(\beta,\lambda)
       (\|(-\Delta_{\bar f}-\lambda)v\|_\beta+\|v\|_\beta).
\end{equation}
Since $H^2(\mu_\beta)\subset H^2(\bar\mu)$, the kernel of
$-\Delta_{\bar f}-\lambda$ is its finite-dimensional reference kernel:
compact-eigenfield/Hermite expansion gives the converse inclusion
in every Gaussian space. Its cokernel is represented by
\begin{equation}\label{eq:G-1-9}
 R_\beta X,\qquad
 R_\beta:=\frac{\mathrm d\bar\mu}{\mathrm d\mu_\beta}
        =\beta^{-d/2}e^{-(1-\beta)r^2/4},
\end{equation}
where $X$ runs through that reference kernel. These independent
functionals annihilate the range and exhaust the cokernel by
index zero. Thus the reference-normalized reduced resolvents are
bounded $L^2(\mu_\beta)\to H^2(\mu_\beta)$.
Write $(-\Delta_{\bar f}-\lambda)^\#$ for this reduced resolvent on compatible
sources, normalized by orthogonality to its kernel in $L^2(\bar\mu)$.

For $F\in L^2(\mu_\beta;T^*\bar M)$, integration by parts gives
\[
 \int(-\delta_{\bar X}F)\varphi\,\mathrm d\mu_\beta
 =\int\left\langle F,\mathrm d\varphi+
       \frac{1-\beta}{2}\varphi\sum_i x_i\,\mathrm dx_i
       \right\rangle\mathrm d\mu_\beta.
\]
The Gaussian multiplication estimate makes this a bounded functional
on $H^1(\mu_\beta)$. It annihilates $R_\beta$, since
$\mathrm dR_\beta=-(1-\beta)R_\beta\sum_i x_i\,\mathrm dx_i/2$.
The form Fredholm inverse therefore gives the unique
reference-mean-zero solution
\[
 -\Delta_{\bar f}u=-\delta_{\bar X}F,\qquad
 \|u\|_{H^1(\mu_\beta)}\le C(\beta)\|F\|_\beta.
\]
The form kernel consists of constants, by the reference $H^1$
energy identity and cutoff approximation; its cokernel is
$\R R_\beta$ by index zero. Hence $\Pi_{\rm ex}F=\mathrm du$ is bounded
in $L^2(\mu_\beta)$. Reference weak uniqueness identifies $\Pi_{\rm ex}$
with the usual orthogonal exact-form projection, and hence
$\delta_{\bar X}(\Id-\Pi_{\rm ex})F=0$.

The reference shrinker identities are
\begin{equation}\label{eq:G-1-10}
 2\bar{\mathcal P}=-\Delta_{\bar f}-\mathrm d\delta_{\bar X}-\frac12,
 \qquad -\Delta_{\bar f}\mathrm du
       =-\mathrm d\Delta_{\bar f}u-\frac12\mathrm du.
\end{equation}
Thus $\bar{\mathcal P}=-\Delta_{\bar f}$ on exact forms and
$\bar{\mathcal P}=(-\Delta_{\bar f}-1/2)/2$ on reference co-closed forms.
For $F$ compatible with $\mathcal K$, define
\[
 \omega=(-\Delta_{\bar f})^\#\Pi_{\rm ex}F
       +2(-\Delta_{\bar f}-\frac12)^\#(\Id-\Pi_{\rm ex})F.
\]
On one-forms, the kernel of $-\Delta_{\bar f}$ consists of Euclidean
translations. In the co-closed subspace, its $1/2$ eigenspace consists of the
nontranslation reference Killing forms, since
$(\bar{\mathcal P}v,v)_{\bar\mu}
=\|\delta_{\bar X}^*v\|_{\bar\mu}^2$.
The remaining $1/2$ eigenforms are exact and orthogonal to
$(\Id-\Pi_{\rm ex})F$. Each projected source is therefore orthogonal to the kernel of its
shifted operator. Reference spectral invariance
preserves the two subspaces, so~\eqref{eq:G-1-10} proves the equation
and weighted estimate. The reduced normalizations give
$\omega\perp\mathcal K$.

Use the augmented maps of
Proposition~\ref{prop:family-reference-inverse} in this fixed
weighted space. The proved inverse and finite Killing conditions make
$\mathscr P$ an isomorphism on $H^2(\mu_\beta)$.
Equation~\eqref{eq:G-1-4} has only compact
indices in $C_b$, and consequently
\[
 \|B_b\omega\|_\beta
       \le C|b|\|\omega\|_{H^1(\mu_\beta)}.
\]
Weighted approximation in the symmetric-divergence identity gives
reference compatibility. The reference Killing conditions are
continuous in this space and their perturbations have norm
$O_\beta(|b|)$. The cutoff first-moment error is bounded by
\[
 C(\beta)|b|+
 C(\beta)\left(\int_{r\ge S-2}R_\beta^2\,
                                      \mathrm d\mu_\beta\right)^{1/2},
\]
whose second term tends to zero as $S\to\infty$, uniformly on
the weighted unit ball.
The same Neumann argument, first choosing $S_{\rm inv}$ and then a
neighborhood independent of $S\ge S_{\rm inv}$, gives
\[
 \|V\|_{H^2(\mu_\beta)}
       \le C(\beta)(\|F\|_\beta+|\vec m|).
\]
Reference uniqueness identifies this weighted solution with that
of Proposition~\ref{prop:family-reference-inverse}.

Since the divergence has just one linear Euclidean coefficient,
\[
 \|\delta_{\bar X}T\|_{L^2(\mu_\beta)}
       \le C(\beta)\|T\|_{C^1},
\]
independently of the support radius. This proves \eqref{eq:G-1-2}.

On concentric reference balls of radii $c/(1+r)$ and $2c/(1+r)$,
the rescaled system is uniformly elliptic with bounded coefficients
through the required orders. These bounds are uniform in $b$ and $r$:
$G_b$ has fixed compact-factor bounds, and its drift is bounded after
rescaling. On the larger ball, the Gaussian density is comparable to
$e^{-\beta r^2/4}$, so
\[
 \|V\|_{L^2(B_{2c/(1+r)};\,\mathrm dV_{\bar g})}
 \le C(\beta)e^{\beta r^2/8}\|V\|_{L^2(\mu_\beta)}.
\]
Interior Schauder estimates therefore give
\[
 \|V\|_{C^{k+1,1/2}(B_{c/(1+r)})}
 \le C(1+r)^Q
 \left(e^{\beta r^2/8}\|V\|_{L^2(\mu_{\beta})}
                          +\|T\|_{C^{k,1/2}}\right).
\]
Here the source uses $k-1$ derivatives of $\delta_{\bar X}T$,
hence only $k$ derivatives of $T$. Raising indices with $G_b$ uses
uniformly bounded coefficients on $N$. Together with the
$H^2(\mu_\beta)$ bound, this gives the global estimate
\[
 \sum_{j=0}^{k+1}|\bar\nabla^jV|(y,\vec x)
 \le C(1+r)^Q e^{\beta r^2/8}
              \bigl(\|T\|_{C^{k,1/2}}+|\vec m|\bigr).
\]
The full H\"older norm on
$N_S$ follows by using these balls for nearby pairs and the supremum
bound for pairs separated by at least $c/(1+S)$. This costs at most
a fixed power of $1+S$.
The strict inequality $\beta/8<\beta$ absorbs these polynomial and
fixed-radius costs and proves \eqref{eq:G-1-3}.

Finally,
$\mathrm d\bar\mu/\mathrm d\mu_\beta
 =\beta^{-d/2}e^{-(1-\beta)r^2/4}$.
Apply this identity to $BV$ and
$\bar\nabla(BV)$, absorbing their fixed polynomial factors into
the strict gap $a<(1-\beta)/8$. Taking square roots proves
\eqref{eq:gauge-coordinate-tail}.
\end{proof}

\clearpage
\phantomsection
\section*{Index of notation}
\addcontentsline{toc}{section}{Index of notation}
\label{sec:notation-index}
\begingroup
\small\linespread{1.0}\selectfont
\setlength{\columnsep}{1.8em}
\setlength{\multicolsep}{0.7\baselineskip}
\newcommand{\notationentry}[2]{\par\noindent\hangindent=1em\hangafter=1 #1, #2\par}
\begin{multicols}{2}
\notationentry{$A_j$, $A_j(u)$}{\pageref{idx:prefix-drop}}
\notationentry{$p,p_1,p_2,p_3,p_4$}{\pageref{thm:logarithmic}, \pageref{stage:entropy-exponents}}
\notationentry{$b,b_+$}{\pageref{eq:compactfamily}, \pageref{thm:restart}}
\notationentry{$B_{\rm cmp},r_{\rm cmp}$}{\pageref{stage:comparison-package}}
\notationentry{$B_{\rm out},r_{\rm out}$}{\pageref{stage:parameter-pairings}}
\notationentry{$B_{\rm gauge},B_{\rm con},r_{\rm gauge},r_{\rm con}$}{\pageref{stage:iteration-neighborhoods}}
\notationentry{$\bar\varphi,T_i,\mathcal O$}{\pageref{stage:parameter-pairings}}
\notationentry{$B_b$}{\pageref{eq:G-1-4}}
\notationentry{$A_U,B_U$}{\pageref{idx:weak-gauge-data}}
\notationentry{$\mathcal M_U$}{\pageref{idx:weak-gauge-data}}
\notationentry{$\mathcal C_\beta(R)$}{\pageref{eq:D-1}}
\notationentry{$\mathcal C_N^d$}{\pageref{idx:product-flow}}
\notationentry{$\mathcal D$}{\pageref{eq:physical-defect}}
\notationentry{$e(b)$}{\pageref{eq:compactfamily}}
\notationentry{$E(s),E_j,E_+$}{\pageref{eq:modified-flow}, \pageref{idx:prefix-drop}, \pageref{thm:restart}}
\notationentry{$E(s)$ (ancient limit)}{\pageref{eq:ancient-entropy}}
\notationentry{$F,\tilde F,u$}{\pageref{idx:caloric-data}}
\notationentry{$F_b,\hat f_b$}{\pageref{idx:product-background}, \pageref{eq:compact-physical-forward}}
\notationentry{$f(t),f_\tau$}{\pageref{eq:terminal-probability}, \pageref{lem:terminal-tails}}
\notationentry{$F_2$}{\pageref{eq:reduced-energy}}
\notationentry{$\bar f,f_N,f_E$}{\pageref{eq:model-potential}, \pageref{idx:model-potentials}}
\notationentry{$g^{z_*},f^{z_*}$}{\pageref{eq:modified-gauge}}
\notationentry{$g^{x_*},f^{x_*},\mu^{x_*}$}{\pageref{thm:ancient-uniqueness}}
\notationentry{$g_\infty,f_\infty$}{\pageref{prop:fixed-limit}}
\notationentry{$F_i,\varphi_i$}{\pageref{prop:time-preserving-reconstruction}}
\notationentry{$\bar g,g_N,g_E$}{\pageref{eq:model}}
\notationentry{$\boldsymbol{\epsilon}(S),\boldsymbol{\epsilon}_+(S')$}{\pageref{eq:seed-error}, \pageref{thm:restart}}
\notationentry{$\boldsymbol{\epsilon}_j$}{\pageref{eq:E-9}}
\notationentry{$[T]_k$}{\pageref{eq:main-quantitative}}
\notationentry{$\gamma,c$}{\pageref{eq:main-quantitative}}
\notationentry{$H^2$}{\pageref{def:weighted-sobolev}}
\notationentry{$H^2(\mu)$}{\pageref{def:weighted-sobolev}}
\notationentry{$H^2(\mu_\beta)$}{\pageref{idx:adjustable-weights}}
\notationentry{$(h,k)$}{\pageref{eq:chart}}
\notationentry{$\mathcal H_j$}{\pageref{prop:kernel}}
\notationentry{$\iota_s,\iota_+$}{\pageref{eq:C-15}, \pageref{thm:restart}}
\notationentry{$J_s$}{\pageref{prop:coherent-maps}}
\notationentry{$Y_s$}{\pageref{prop:coherent-maps}}
\notationentry{$K,K_j$}{\pageref{prop:kernel}}
\notationentry{$K_N,K_+$}{\pageref{prop:kernel}}
\notationentry{$N_R$}{\pageref{idx:product-regions}}
\notationentry{$\mathcal K,\mathcal K_0$}{\pageref{idx:killing-spaces}}
\notationentry{$\mathcal L_{\bar f},\mathcal L_N$}{\pageref{prop:kernel}}
\notationentry{$\mathcal L_f$}{\pageref{idx:weighted-operators}}
\notationentry{$\mathbb L_U,\mathscr P$}{\pageref{lem:augmented-weak-inverse}}
\notationentry{$\mathscr P_{b,S}$}{\pageref{prop:family-reference-inverse}}
\notationentry{$c_1,c_2,c_3,c_5,c_6$}{\pageref{stage:gaussian-rates}, \pageref{eq:E-5}}
\notationentry{$c_4$}{\pageref{prop:actual-comparison}, \pageref{stage:comparison-package}}
\notationentry{$Q_*$}{\pageref{stage:baseline-regularity}}
\notationentry{$Q_\beta$}{\pageref{eq:D-1}}
\notationentry{$\beta$}{\pageref{eq:D-1}, \pageref{eq:E-6}}
\notationentry{$\epsilon_+$}{\pageref{prop:restart-gauge}}
\notationentry{$\epsilon_{\rm met},\epsilon_{\rm sl}$}{\pageref{subsec:physical-common-setup}, \pageref{prop:restart-gauge}}
\notationentry{$\rho,\rho_1,\rho_2,\rho_3$}{\pageref{subsec:physical-common-setup}}
\notationentry{$\Gamma,\rho_E$}{\pageref{sec:finite-quotients}}
\notationentry{$\Lambda$}{\pageref{stage:geometric-retentions}}
\notationentry{$c'$}{\pageref{lem:quantitative-sublevel-coverage}}
\notationentry{$\Lambda'$}{\pageref{eq:E-4}}
\notationentry{$\chi_0,\chi_1,\chi_2,\chi_3$}{\pageref{prop:restart-gauge}, \pageref{idx:weak-gauge-data}, \pageref{eq:M-3}}
\notationentry{$L_j(u)$}{\pageref{eq:E-10}}
\notationentry{$\ell$}{\pageref{eq:density-ratio}}
\notationentry{$\vec m$}{\pageref{thm:logarithmic}, \pageref{eq:family-reference-system}}
\notationentry{$\boldsymbol{\mu}(g,\tau)$}{\pageref{eq:nash-hypothesis}}
\notationentry{$\mathcal N_{x,t}(\tau)$}{\pageref{thm:ancient-uniqueness}}
\notationentry{$\bar{\mathcal P},\mathcal P_b$}{\pageref{idx:killing-spaces}, \pageref{eq:gauge-augmented-equation}}
\notationentry{$P_U$}{\pageref{idx:gauge-map}}
\notationentry{$R_{\min}$}{\pageref{stage:iteration-thresholds}}
\notationentry{$R_j$}{\pageref{eq:E-12}}
\notationentry{$r_s$}{\pageref{eq:C-16}}
\notationentry{$R_s$}{\pageref{prop:coherent-maps}}
\notationentry{$R_+$}{\pageref{idx:gauge-radius}}
\notationentry{$s,t$}{\pageref{idx:modified-time}}
\notationentry{$S_j$ (comparison radii, including $S_0$)}{\pageref{eq:E-8}}
\notationentry{$S_*$}{\pageref{subsec:physical-common-setup}}
\notationentry{$S_{\rm inv}$}{\pageref{prop:family-reference-inverse}, \pageref{thm:gauge-inverse}}
\notationentry{$T,\tilde T$}{\pageref{idx:caloric-data}}
\notationentry{$w(g,f)$}{\pageref{eq:physical-defect}}
\notationentry{$\mathcal W(g,f)$}{\pageref{eq:physical-defect}}
\notationentry{$W_s$}{\pageref{eq:M-2}}
\notationentry{$a\in K_+$}{\pageref{thm:jets}}
\notationentry{$X_b=(G_b,F_b)$}{\pageref{idx:product-background}}
\notationentry{$\hat X_{b,a}$}{\pageref{prop:cutoffjets}}
\notationentry{$Z_b(s)$}{\pageref{prop:compact-forward}}
\notationentry{$Z_s$}{\pageref{prop:coherent-maps}, \pageref{eq:M-4}}
\notationentry{$z_*;z_H$}{\pageref{eq:terminal-probability}, \pageref{lem:terminal-tails}}
\notationentry{$\vec x,x_i$}{\pageref{idx:euclidean-coordinates}}
\notationentry{$\Delta_f$}{\pageref{idx:weighted-operators}}
\notationentry{$\delta_{\bar X},\delta_{\bar X}^*$}{\pageref{idx:weighted-divergence}}
\notationentry{$\delta_g,\delta_X,\delta_X^*$}{\pageref{idx:weighted-divergence}}
\notationentry{$\Theta_N$}{\pageref{eq:model}}
\notationentry{$\mu_X$}{\pageref{eq:weighted-conventions}}
\notationentry{$\bar\mu,\mu_N,\gamma_d$}{\pageref{eq:model-potential}}
\notationentry{$\mu_t,\mu^{z_*}$}{\pageref{eq:terminal-probability}}
\notationentry{$\mu_\beta$}{\pageref{idx:adjustable-weights}}
\notationentry{$\mu',\nu',Z',\vec m'$}{\pageref{prop:restart-gauge}}
\notationentry{$\Xi_s$}{\pageref{eq:modified-diffeomorphism}}
\notationentry{$\Pi_{\mathrm{sl}}$}{\pageref{idx:slice-projections}}
\notationentry{$\Phi$}{\pageref{eq:physical-defect}}
\notationentry{$\Phi_b$}{\pageref{eq:compactfamily}, \pageref{idx:product-background}}
\notationentry{$\Psi$}{\pageref{thm:main}, \pageref{prop:embedding}}
\notationentry{$\Psi_\Gamma$}{\pageref{thm:quotient-uniqueness}, \pageref{thm:ancient-uniqueness}}
\notationentry{$\omega$}{\pageref{idx:vanishing-modulus}}
\notationentry{$\Delta_j$}{\pageref{eq:E-14}}

\end{multicols}
\endgroup
\clearpage
\begingroup\linespread{0.99}\small
\phantomsection

\endgroup

\medskip
Yu Li, Institute of Geometry and Physics, University of Science and Technology of China, No.~96
Jinzhai Road, Hefei, Anhui Province, 230026, China; Hefei National Laboratory, No.~5099 West
Wangjiang Road, Hefei, Anhui Province, 230088, China; Email: \href{mailto:yuli21@ustc.edu.cn}{yuli21@ustc.edu.cn}.
\end{document}